\documentclass{article}

\usepackage{macro}

\AtEveryBibitem{%
  \clearname{editor}%
  \clearname{editora}%
  \clearname{editorb}%
  \clearname{editorc}%
}

\hypersetup{
    colorlinks=true,
    linkcolor=purple,
    citecolor=forestgreen,
}

\title{Optimal Estimation in Orthogonally Invariant Generalized Linear Models: Spectral Initialization and Approximate Message Passing}

\author{
    Yihan Zhang\thanks{School of Mathematics, University of Bristol. Email: \href{mailto:yihan.zhang@bristol.ac.uk}{\texttt{yihan.zhang@bristol.ac.uk}}.}
    \and
    Hong Chang Ji\thanks{Department of Mathematics, Sungkyunkwan University. Email: \href{mailto:hcji@skku.edu}{\texttt{hcji@skku.edu}}.}
    \and
    Ramji Venkataramanan\thanks{Department of Engineering, University of Cambridge. Email: \href{mailto:rv285@cam.ac.uk}{\texttt{rv285@cam.ac.uk}}.}
    \and
    Marco Mondelli\thanks{Institute of Science and Technology Austria. Email: \href{mailto:marco.mondelli@ist.ac.at}{\texttt{marco.mondelli@ist.ac.at}}.}
}

\begin{document}
\maketitle

\pagenumbering{roman}

\begin{abstract}
We consider the problem of parameter estimation from a  generalized linear model with a random design matrix that is orthogonally invariant in law. 
Such a model allows the design have an arbitrary distribution of singular values and only assumes that its singular vectors are generic. It is a vast generalization of the i.i.d.\ Gaussian design  typically considered in the theoretical literature, and is motivated by the fact that real data often have a complex correlation structure so that  methods relying on i.i.d.~assumptions can be highly suboptimal.
Building on the paradigm of spectrally-initialized iterative optimization, this paper proposes optimal spectral estimators and combines them with an approximate message passing (AMP) algorithm, establishing rigorous performance guarantees for these two algorithmic steps. Both the spectral initialization and the subsequent AMP meet existing conjectures on the fundamental limits to estimation---the former on the optimal sample complexity for efficient weak recovery, and the latter on the optimal errors.  
Numerical experiments 
suggest the effectiveness of our methods and accuracy of our theory beyond orthogonally invariant data. 
\end{abstract}

\tableofcontents
\newpage

\pagenumbering{arabic}

\newcounter{asmpctr} 
\setcounter{asmpctr}{\value{enumi}}

\section{Introduction}
\label{sec:intro}

We study the canonical problem of parameter estimation from a generalized linear model (GLM) \cite{Nelder_Wedderburn}. Consider a dataset $ \cD = (X,y) \in \bbR^{n\times d} \times \bbR^n $ where the design matrix $ X = \matrix{x_1 & \cdots & x_n}^\top $ consists of $n$ covariate vectors each of dimension $d$, and  $y$ is a vector of response variables each given by 
\begin{align}
y_i &= q(\inprod{x_i}{\beta_*}, \varepsilon_i), \qquad 
1\le i\le n.
\notag 
\end{align}
Here, $q:\mathbb R^2\to\mathbb R$ is a link function, $\eps_i$ a noise variable,  $\beta_*\in\bbR^d$ an unknown vector of regression coefficients, and $\inprod{\cdot\, }{\cdot}$ denotes the Euclidean inner product. The task is to estimate $ \beta_* $ from $\cD$. The nonlinear link function $q$ generalizes linear regression and encompasses the most commonly used statistical models in engineering and the natural sciences, including
logistic regression \cite{Tolles_Meurer}, phase retrieval \cite{shechtman2015phase,fannjiang2020numerics} 
($y_i = |\inprod{x_i}{\beta_*}| + \eps_i$), and 1-bit compressive sensing \cite{boufounos20081} 
($y_i = \sgn(\inprod{x_i}{\beta_*} + \eps_i)$).  

In the classical regime where the sample size $n$ is much larger than the parameter dimension $d$, the regression coefficients in $\beta_*$ can be consistently estimated and standard inference procedures, such as confidence intervals and likelihood-based tests, provide accurate uncertainty quantification. However, in modern high-dimensional settings where $n$ and $d$ are both large and comparable ($n/d\to\delta\in(0,\infty)$), these classical results are no longer valid \cite{donoho2016high,sur2019modern}. This has led to a flurry of work aimed at establishing both statistical limits and efficient estimators in this high-dimensional regime, largely focusing on the setting where the design matrix $X$ is unstructured, i.e., its entries are i.i.d.\ Gaussian \cite{Barbier_etal,bayati2011lasso,donoho2016high, thrampoulidis2018precise, sur2019modern}. 

Estimation in the high-dimensional regime can be broken down into two steps: \emph{(i)} producing a pilot estimate $ \wh{\beta} \in \bbR^d $ that has a nontrivial correlation with $\beta_*$, and \emph{(ii)} refining the initialization $ \wh{\beta} $ via an iterative procedure.
More precisely, in step \emph{(i)} we aim to produce an estimate with a nonvanishing asymptotic overlap with $\beta_*$:
\begin{align}
    \liminf_{d\to\infty} \frac{|\langle\wh{\beta}, \beta_*\rangle|}{\|\wh{\beta}\|_2 \|\beta_*\|_2} > 0 . \label{eqn:overlap}
\end{align}
This criterion is often called `weak recovery' in the literature \cite{Mondelli_Montanari,Maillard_etal_threshold,Ma_etal}. Spectral methods \cite{ChenChiFanMa_spec_book} provide a simple and popular approach for obtaining a $\wh{\beta}$ that satisfies \Cref{eqn:overlap}.  A spectral estimator for a GLM is given by the principal eigenvector of the following data-dependent matrix:
\begin{equation}\label{eqn:specintro}
D = X^\top \diag(\cT(y_1), \cdots, \cT(y_n)) X=\sum_{i=1}^n \cT(y_i) x_i x_i^\top\in \bbR^{d\times d},    
\end{equation}
where $ \cT \colon \bbR \to \bbR $ is a user-defined preprocessing function. The performance of spectral estimators obtained via the matrix in \Cref{eqn:specintro} has been studied in significant detail \cite{Netrapalli_alt_min,candes2015phase,Chen_Candes_trim,Lu_Li,Mondelli_Montanari,Luo_Alghamdi_Lu}; specifically, for an i.i.d.\ Gaussian design and prior, it has been shown in \cite{Mondelli_Montanari} that a suitable choice of $\cT$ achieves weak recovery up to the information-theoretic limit for a class of link functions $q$, including phase retrieval. 
Moving on to step \emph{(ii)}, a variety of iterative procedures have been proposed, ranging from optimizing a suitable loss function (e.g., penalized log-likelihood) via gradient descent  or alternating minimization \cite{candes2015phase, Netrapalli_alt_min, yi2014alternating, chandrasekher2023sharp} to Kaczmarz methods \cite{wei2015solving}
and higher-order methods \cite{celentano2025state}.

In the Bayesian setting, where $\beta_*$ is generated according to a prior distribution independent of $X$, the posterior mean $\expt{\beta_* \mid X, y }$ minimizes the mean squared error. Although it is usually computationally infeasible, the asymptotic Bayes risk $ \lim_{d\to\infty} \frac{1}{d} \expt{\normtwo{\beta_* - \expt{\beta_* \mid X, y }}^2}$ has been precisely characterized  for i.i.d.\ Gaussian designs, under suitable model assumptions \cite{Barbier_etal}.  For this setting, it has been shown  \cite{Mondelli_Venkataramanan} that the spectral estimator in  \Cref{eqn:specintro}  can be refined via a Generalized Approximate Message Passing (GAMP) algorithm \cite{Rangan}. In certain instances, this combined procedure achieves the Bayes risk \cite{Barbier_etal}, and more broadly, GAMP estimators have been shown to be optimal among general first order methods (GFOMs) \cite{Celentano_Montanari_Wu,Montanari_Wu}. 

\begin{figure}[t]
    \centering
    \begin{subfigure}{0.48\textwidth}
        \includegraphics[width=\linewidth]{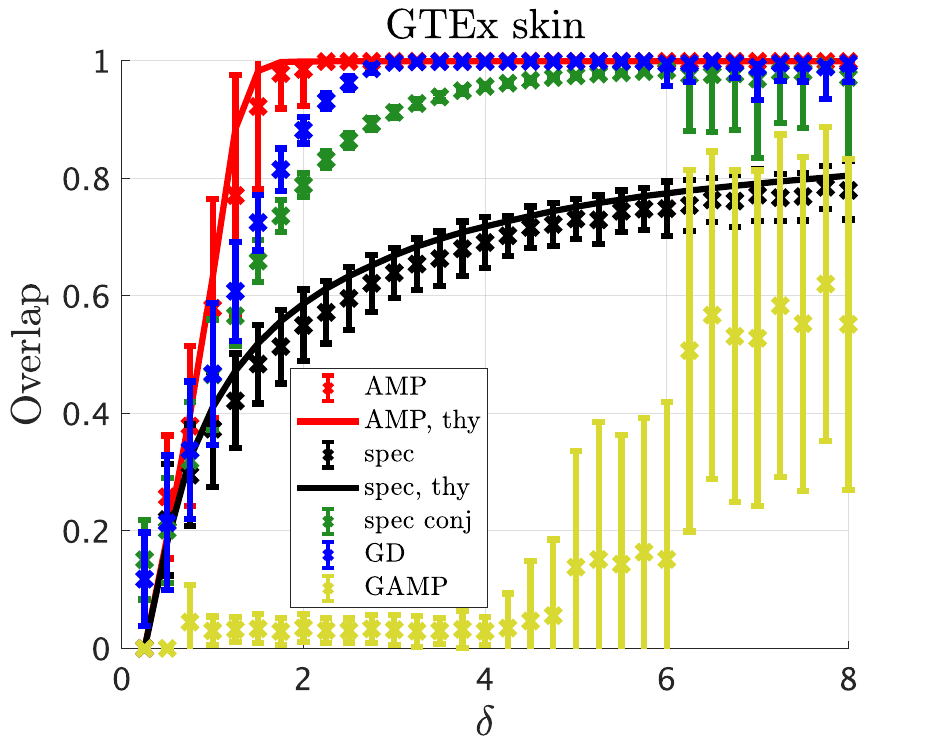}
    \end{subfigure}
    \hspace{1em}
    \begin{subfigure}{0.48\textwidth}
        \includegraphics[width=\linewidth]{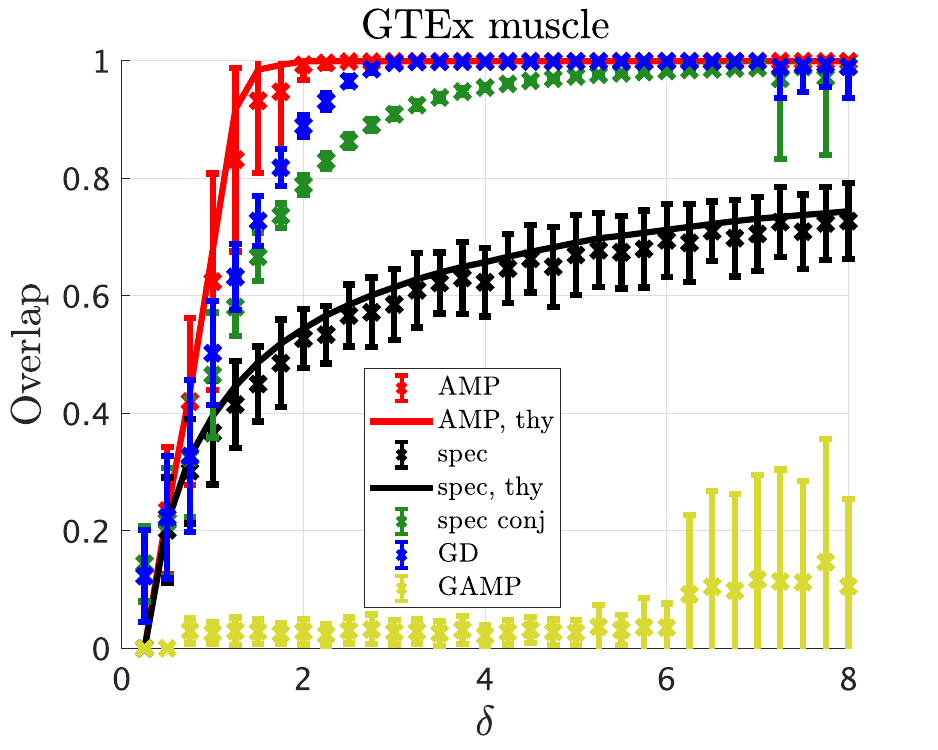}
    \end{subfigure}
    \caption{Performance of $5$ algorithms  on phase retrieval whose design matrix is given by one of two real datasets from GTEx, see \Cref{sec:experiments} for details. 
    The GAMP algorithm \cite{Mondelli_Venkataramanan} (in yellow) with provable optimality on Gaussian data is no longer as effective. 
    The proposed Bayes-GVAMP (`AMP' in red) with spectral initialization (`spec' in black) shows dominant performance which is well predicted by our theory under the orthogonal invariance assumption. 
    Gradient descent (`GD' in blue) and GAMP are both initialized with the spectral estimator from \cite[Equation (20)]{Luo_Alghamdi_Lu} (`spec conj' in green). 
    They achieve the same weak recovery threshold as spectrally initialized Bayes-GVAMP but worse overlaps.}
    \label{fig:GTEx}
\end{figure}

The results above provide a clear picture of how to perform optimal estimation in GLMs when the design matrix $X$ has  i.i.d.\ Gaussian entries. Some of these asymptotic results, e.g.\ those for spectral methods and certain convex penalized estimators, can be extended to the broader setting of Gaussian designs with arbitrary covariance \cite{Zhang_Ji_Venkataramanan_Mondelli, zhao2022asymptotic,CelentanoMontanariWei}.
However, in many applications the data have a complex structure and the correlations across the entries of $X$ (including dependence across rows) are not adequately captured by a Gaussian design. 
As a consequence, methods relying on Gaussian assumptions on the design  may be highly suboptimal.
This is highlighted in \Cref{fig:GTEx}: the algorithms that are provably optimal when $X$ has i.i.d.\ Gaussian entries \cite{Barbier_etal,Mondelli_Venkataramanan} fall short of the performance promised by the theory when $X$ is taken from computational genomics \cite{GTEx} datasets.

A common way to capture the structure in the data and retain mathematical tractability is to model the design matrix $X$ as random with an \emph{orthogonally invariant} law. This corresponds to the singular vectors of $X$ being generic (or, more formally, Haar distributed over the orthogonal group) while its  spectrum can be arbitrary, which allows for row-column dependencies as well as a range of other structures in the design. 
The universality of this design choice has been both observed empirically 
\cite{oymak2014case,ccakmak2019memory,Ma_etal} and proved formally \cite{Dudeja_Bakhshizadeh,Dudeja_Lu_Sen,wang2024universality,Dudeja_Sen_Lu,Wang_Zhong_Fan} for a broad class of well-known ensembles, including those obtained from discrete cosine/sine transforms and Hadamard-Walsh transforms; see \cite[Section 1.3]{Dudeja_Lu_Sen}. 

Due to their flexibility, orthogonally invariant designs have been recently investigated in the statistical literature, e.g., in the context of both linear regression \cite{Li_Sur}
and principal component analysis \cite{fan2020approximate,barbier2023fundamental,Dudeja_Liu_Ma}. 
For estimation in a GLM with orthogonally invariant design, variants of AMP algorithms have been proposed \cite{Schniter_Rangan_Fletcher,VKM}, and their fixed points have been shown to match the Bayes risk conjectured using tools from statistical physics \cite{kabashima2008inference, takahashi2022macroscopic}. However, these AMP algorithms require an initialization that  satisfies \Cref{eqn:overlap} and is independent of the design matrix. Assuming the availability of such an initialization is unrealistic for several popular choices of nonlinear link functions $q$ (e.g., phase retrieval). Spectral estimators are a natural choice for initialization, but existing results are either restrictive in terms of the spectrum of the design matrix (unit singular values in \cite{Dudeja_Bakhshizadeh_Ma_Maleki}, correlated Gaussian design in \cite{Zhang_Ji_Venkataramanan_Mondelli}) or they remain at a conjectural level \cite{Maillard_etal_construction}. In short, obtaining rigorous performance guarantees on initialization methods has proved hard, let alone the combination of the initialization with iterative methods, such as AMP.


\paragraph{Contributions.} In this paper, we address several open questions on optimal estimation in orthogonally invariant GLMs: we provide optimal spectral estimators, show how they can be combined with an AMP algorithm, and establish rigorous asymptotic guarantees for both. Both the initialization and the iterative refinement via AMP meet existing conjectures on  fundamental limits to estimation---the former on the optimal sample complexity for efficient weak recovery \cite{Maillard_etal_threshold}, and the latter on the optimal errors \cite{kabashima2008inference, takahashi2022macroscopic}. 
\begin{figure}[t]
    \centering
    \begin{subfigure}{0.48\textwidth}
        \includegraphics[width=\linewidth]{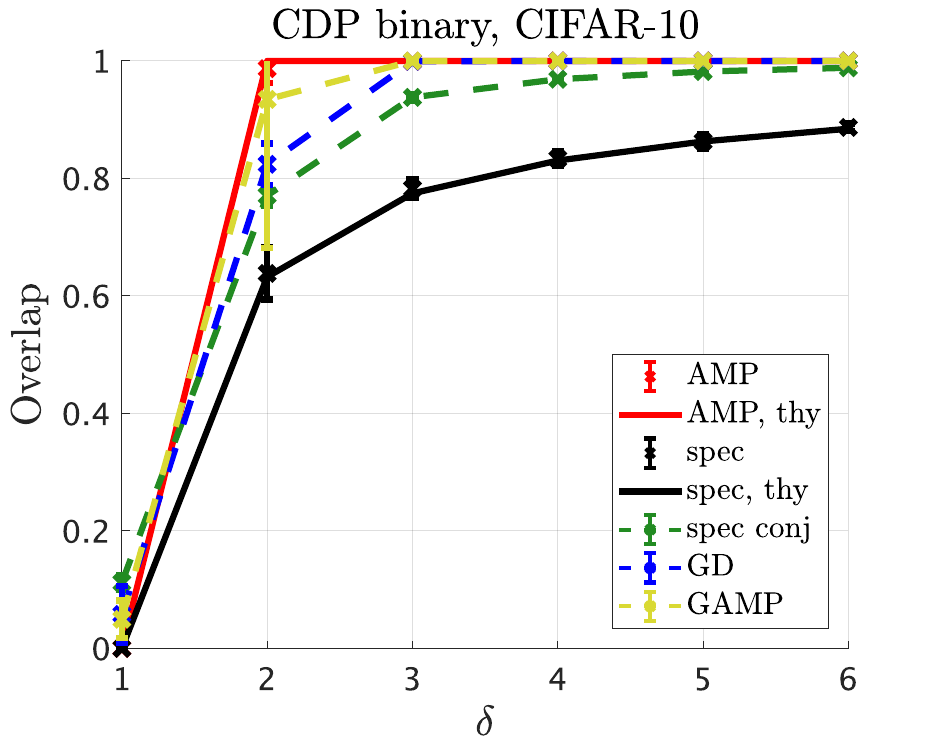}
    \end{subfigure}
    \hspace{1em}
    \begin{subfigure}{0.48\textwidth}
        \includegraphics[width=\linewidth]{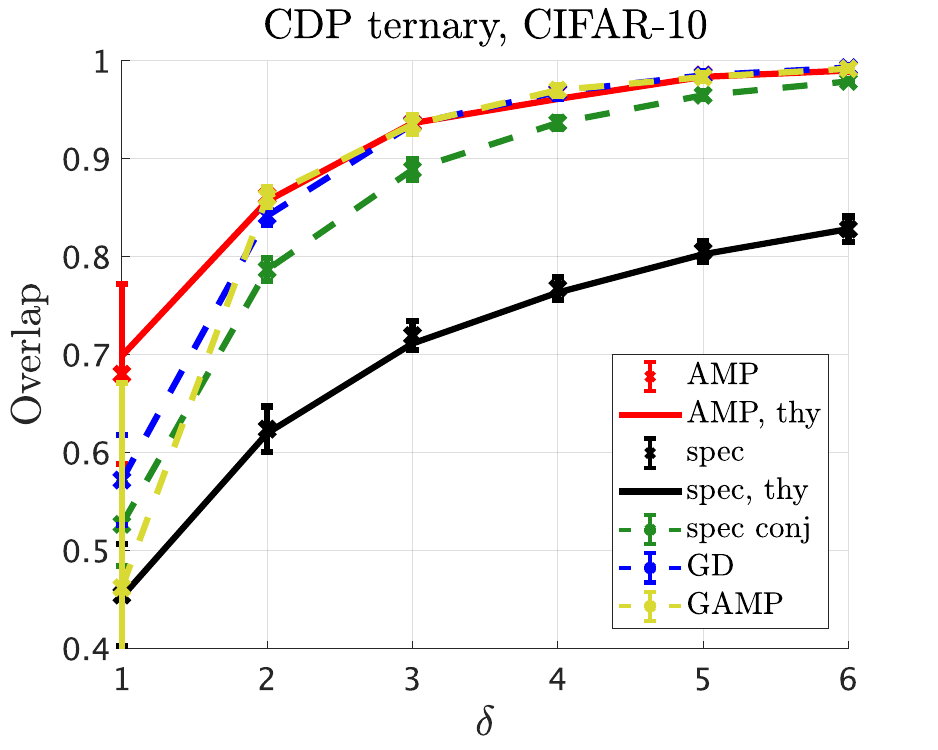}
    \end{subfigure}
    \caption{Performance of $5$ algorithms on phase retrieval whose design matrix is given by either binary or ternary coded diffraction patterns (CDP), see \Cref{sec:experiments} for details. The regression coefficients are obtained by applying standard preprocessing to a $32\times32$ image of a truck from the CIFAR-10 dataset.
    Our theoretical guarantees (`thy', solid curves) remain accurate when neither the design matrix is orthogonally invariant nor the prior is i.i.d., a strong sign of universality.}
    \label{fig:CDP}
\end{figure}
Our contributions are summarized as follows.
\begin{enumerate}
    \item We establish the asymptotic overlap of spectral estimators for a class of preprocessing functions (\Cref{thm:spec}). More precisely, we unveil a criticality condition \Cref{eqn:thr} under which
    the largest eigenvalue of $D$ in \Cref{eqn:specintro} escapes the bulk of the spectrum and the spectral estimator given by the top eigenvector of $D$ achieves weak recovery in the sense of \Cref{eqn:overlap}, a phenomenon similar to the celebrated BBP phase transition \cite{BBAP}. 

    \item Leveraging the characterization of \Cref{thm:spec}, we identify a preprocessing function that achieves weak recovery with the smallest possible sample complexity (\Cref{thm:opt_thr}). Notably, the critical sample complexity above which spectral estimators achieve nonvanishing asymptotic overlap coincides with the one conjectured in \cite{Maillard_etal_threshold} to be optimal among all efficient estimators. 
    
    \item We design a general class of AMP iterations which, following the terminology of \cite{Cademartori_Rush}, we call Generalized Vector AMP (GVAMP). GVAMP iteratively refines the spectral initialization by applying nonlinearities both to the singular values of the design matrix and to the entries of the iterates. We prove that the dynamics (and, hence, the performance) of this class of algorithms can be tracked via a state evolution recursion, a succinct set of low-dimensional deterministic equations (\Cref{thm:spec_GVAMP}). This resolves the issue, common in related work, of requiring an initialization that achieves weak recovery, while being independent of the design matrix. 

    \item Finally, via a suitable specialization of the generic GVAMP from \Cref{thm:spec_GVAMP}, we obtain an algorithm---spectrally initialized Bayes-GVAMP---whose fixed points match the Bayes risk (\Cref{thm:replica}), as conjectured in \cite{kabashima2008inference,takahashi2022macroscopic}. 
    This offers an end-to-end algorithmic solution to the problem of generalized linear regression, with a rigorous performance characterization and conjectured optimality among a large class of efficient algorithms. 
\end{enumerate}
Due to the aforementioned universality of orthogonally invariant designs, we expect our findings to hold for design matrices common in applications. This is showcased in \Cref{fig:GTEx,fig:CDP}: the empirical results match our theoretical predictions even when $X$ is not generated from an orthogonally invariant law, but taken from computational imaging \cite{candes_CDP} or quantitative genomics \cite{GTEx} datasets; our proposed algorithms consistently outperform the approaches stemming from i.i.d.\ assumptions on the design. 

\section{Related work}
\label{sec:related_work}

\paragraph{Spectral methods.}
A line of work has shown the effectiveness and, in some cases, optimality of spectral methods for weak recovery in generalized linear models with i.i.d.\ Gaussian designs \cite{Lu_Li,Mondelli_Montanari,Luo_Alghamdi_Lu,Dudeja_Bakhshizadeh}, correlated Gaussian designs \cite{Zhang_Ji_Venkataramanan_Mondelli} and partial Haar designs \cite{Ma_etal,Dudeja_Bakhshizadeh_Ma_Maleki}. 
Moreover, for i.i.d.\ Gaussian designs, spectral methods have been extended to mixed regression \cite{mixed-zmv-arxiv}, multi-index models \cite{kovavcevic2025spectral,defilippis2025optimal} and higher-order tensors \cite{Damian_Lee_Bruna,Damian_etal}. The 
optimal combination of linear and spectral estimators is obtained in \cite{mondelli2021optimalcombination}. 
The present paper further advances this body of literature by considering orthogonally invariant designs: we propose and analyze spectral estimators that are effective down to the conjectured computationally optimal weak recovery threshold $ \delta_* $ \cite{Maillard_etal_threshold}, i.e., the smallest $\delta$ above which weak recovery is achievable by efficient algorithms. 
The conjectured computational optimality of spectral methods is underpinned by a well-known heuristic connection to linearized AMP \cite{Maillard_etal_construction}. 
Specifically, considerations based on statistical physics postulate that Bayes-GVAMP is optimal among all polynomial-time algorithms. 
Spectral estimators, when viewed as power iteration, arise as the linearization of Bayes-GVAMP around its trivial fixed point. 
Then, $ \delta_* $ can be heuristically determined as the smallest $ \delta $ above which the trivial fixed point is repulsive. 
Our results regarding spectral estimators build upon and make rigorous this connection. 

\paragraph{Approximate Message Passing.}
Our generic formalism of GVAMP takes inspiration from a number of other iterative algorithms for Bayesian linear or generalized linear models including 
Adaptive Thouless--Anderson--Palmer (AdaTAP) \cite{Opper_Winther_ada1,Opper_Winther_ada2}, 
Expectation Propagation (EP) \cite{Minka,Seeger,Ma_etal,Ma_Xu_Maleki}, 
Expectation-Consistent approximation (EC) \cite{Opper_Winther_EC,Kabashima_Vehkapera}, 
S-AMP \cite{Cakmak_Winther_Fleury,Cakmak_Winther_Fleury_nonlin}, RI-GAMP \cite{VKM},
Orthogonal AMP \cite{Ma_Ping,Dudeja_Bakhshizadeh,Dudeja_Lu_Sen,Dudeja_Sen_Lu,Dudeja_Liu_Ma},
Vector AMP \cite{Rangan_Schniter_Fletcher,Schniter_Rangan_Fletcher,Takeuchi}. 
It is argued in \cite{Maillard_etal_high} via a non-rigorous Plefka expansion approach that the derivations of AdaTAP, EC and VAMP are equivalent. 
All the above algorithms operate on certain specific resolvent-type matrices of $X$ (e.g., the linear MMSE estimate) and apply the MMSE denoisers to the iterates. 
This feature is fundamentally due to the ways in which those algorithms are derived. 
Our generic formalism of GVAMP does not have MMSE estimates baked in and allows for arbitrary trace-free spectral transformations of $X$ and divergence-free entry-wise nonlinearities of the iterates. 
This level of generality is expected to be useful beyond the Bayes-optimal setting. 
VAMP-type algorithms with MMSE denoisers are also studied for more sophisticated models than GLMs \cite{pandit2020inference,Pandit_etal,xu2023approximate}. 
It is expected that our generic GVAMP formalism and the associated state evolution result can be extended to those settings. 
AMPs of the style of \cite{Bayati_Montanari,Bolthausen_AMP} developed for an i.i.d.\ design find applications to diverse problems as both an algorithmic primitive and a proof technique \cite{Mousavi_Maleki_Baraniuk,Montanari_Venkataramanan,Celentano_Fan_Mei,Celentano,Han}. The 
GVAMP proposed in this paper, when specialized to an i.i.d.\ Gaussian design, differs from the AMP of \cite{Bayati_Montanari,Bolthausen_AMP}: it is empirically observed that GVAMP has the same fixed points as AMP (left panel of \Cref{fig:gauss}) but converges faster. 

\paragraph{Bayes risk.}
Under an i.i.d.\ Gaussian design, the Bayes risk of generalized linear models is established in \cite{Barbier_etal} by the adaptive interpolation method \cite{Barbier_Macris}. 
To our knowledge, this method has not been successfully applied to orthogonally invariant designs. 
However, the Bayes risk 
is recently proved in \cite{Li_Fan_Sen_Wu} for linear models with orthogonally invariant design under an additional technical condition that restricts the support of the limiting spectral distribution of $X$ to be sufficiently narrow.
The conjectured Bayes risk of generalized linear models with orthogonally invariant design in \cite{kabashima2008inference,takahashi2022macroscopic,Maillard_etal_threshold} has not been rigorously justified in its full generality.
\cite{reeves2017additivity,gabrie2018entropy} extend such conjectures regarding Bayes risk to more complex architectures such as multilayer networks with orthogonally invariant weight matrices. 

\section{Preliminaries}
\label{sec:prelim}

\subsection{Notation and definitions}
\label{sec:notation}

All vectors are column vectors. 
For a symmetric matrix $ M \in \bbR^{d\times d} $, we denote by $ \lambda_1(M) \ge \cdots \ge \lambda_d(M) $ its (real) eigenvalues in non-increasing order and by $ v_1(M), \cdots, v_d(M) \in \bbS^{d-1} $ the corresponding eigenvectors (of unit $\ell_2$ norm). 
For a matrix $ V = \matrix{v_1 & \cdots & v_k} \in \bbR^{n\times k} $, we use $ \Pi_V \in \bbR^{n\times n} $ to denote the orthogonal projection onto $ \spn\{(v_i)_{i\in[k]}\} $, and  let  $ \Pi_{V^\perp} = I_n - \Pi_V $.

All $\lim$, $\liminf$, $\limsup$ hold almost surely, unless otherwise specified. 
We use $ \delta_x(\cdot) = \delta(\cdot - x) $ to denote the point mass at $x\in\bbC$. 
The orthogonal group in dimension $n$ is defined as $$ \bbO(n) \coloneqq \brace{O\in\bbR^{n\times n} : OO^\top = O^\top O = I_n}.$$ 
The unique bi-invariant Haar measure on it (normalized to be a probability measure) is denoted by $ \haar(\bbO(n)) $. 
The empirical spectral distribution of a matrix $M$ is defined as the empirical distribution of its eigenvalues if $M\in\bbR^{d\times d}$ is symmetric: 
\begin{align}
    \frac{1}{d} \sum_{i = 1}^d \delta_{\lambda_i(M)} , \notag 
\end{align}
or of its singular values if $M\in\bbR^{n\times d}$ is rectangular: 
\begin{align}
    \frac{1}{n\wedge d} \sum_{i = 1}^{n\wedge d} \delta_{\sigma_i(M)} . \notag 
\end{align}
Generically, we use $ p_\sfA, p_{\sfA,\sfB}, p_{\lr{\sfA \mid \sfB}} $ to denote the marginal, joint and conditional densities of random variables $\sfA, \sfB$.  

\begin{definition}[Pseudo-Lipschitz functions]
\label{def:pseudo_lip}
For $ p\ge1 $, a function $ f \colon \bbR^k \to \bbR $ is called a pseudo-Lipschitz function of order $p$ if there exists $ C>0 $ such that for all $ x,y\in\bbR^k $, 
\begin{align}
    \abs{ f(x) - f(y) } &\le C \normtwo{x - y} \paren{1 + \normtwo{x}^{p-1} + \normtwo{y}^{p-1}} . \notag 
\end{align}
Denote by $ \PL^p_{k\to1} $ the set of all such functions and $ \PL_{k\to1} = \bigcup_{p\ge1} \PL^p_{k\to1} $. 
\end{definition}

\begin{definition}[Polynomial growth]
\label{def:poly_grow}
For $p\ge1$, a function $f\colon\bbR^k\to\bbR^\ell$ is said to have polynomial growth of order $p$ if there exists $C>0$ such that 
\begin{align}
    \normtwo{f(x)} &\le C \paren{1 + \normtwo{x}^p} , \label{eqn:poly_growth}
\end{align}
for all $ x\in\bbR^k $. 
Denote by $ \PG^p_{k\to\ell} $ the set of all such functions and $ \PG_{k\to\ell} = \bigcup_{p\ge1} \PG^p_{k\to\ell} $. 
\end{definition}
Note that any pseudo-Lipschitz function has polynomial growth, i.e., $ \PL^p_{k\to1} \subset \PG^p_{k\to1} $. 

\begin{definition}[Wasserstein convergence, {\cite[Chapter 6]{Villani_book}}]
\label{def:wasserstein} 
For $p\ge1$, we say that the empirical distribution of the rows of a (possibly random) matrix $ \matrix{ v_1 & \cdots & v_k } = \matrix{v_{i,j}}_{(i,j)\in[n]\times[p]} \in \bbR^{n\times k} $ converges in Wasserstein-$p$ to a random vector $ \matrix{ \sfV_1 & \cdots & \sfV_k } \in \bbR^{1\times k} $ with finite $ \expt{\normtwo{\matrix{\sfV_1 & \cdots & \sfV_k}}^p} $, denoted by 
\begin{align}
\matrix{ v_1 & \cdots & v_k } &\xrightarrow{W_p} \matrix{ \sfV_1 & \cdots & \sfV_k } , \label{eqn:wass} 
\end{align}
if for any $ f \in \PL^p_{k\to1} $, 
it holds almost surely that 
\begin{align}
\lim_{n\to\infty} \frac{1}{n} \sum_{i=1}^n f(v_{i,1}, \cdots, v_{i,k}) &= \expt{ f(\sfV_1, \cdots, \sfV_k) } . \label{eqn:wass_f} 
\end{align}
We write $ \matrix{ v_1 & \cdots & v_k } \wto \matrix{ \sfV_1 & \cdots & \sfV_k } $ to mean $ W_p $ convergence for every fixed $p\ge1$, where the right-hand side has finite mixed moments of all orders. 
\end{definition}

To verify \Cref{eqn:wass}, it suffices to verify \Cref{eqn:wass_f} for either (a) any $ f \in \PG^p_{k\to1} $ continuous, or (b) $ f(v_1, \cdots, v_k) = \normtwo{\matrix{v_1 & \cdots & v_k}}^p $ and any bounded Lipschitz $ f $; 
see \cite[Appendix E]{fan2020approximate}. 

\subsection{Model}
\label{sec:model}

Let  $ X \in \bbR^{n\times d} $ be a design matrix, $ \beta_* \in \bbR^d $  a vector of regression coefficients,and $ \eps \in \bbR^n $  a vector of noise variables. 
Consider a generalized linear model 
\begin{align}
y &= q(X \beta_*, \eps) . 
\label{eqn:model} 
\end{align}
In the above display, the function $ q \colon \bbR^2 \to \bbR $ is applied component-wise to its vector inputs, i.e., $ q(z, \eps) = \matrix{ q(z_1, \eps_1) & \cdots & q(z_n, \eps_n) }^\top $, where $ z \coloneqq  X \beta_* $. 
Equivalently, 
\begin{align}
\lr{y \mid X \beta_*} &\sim Q(\cdot \mid  X \beta_*) , \label{eqn:model_Q}
\end{align}
for a conditional density $ Q\paren{\cdot \mid \cdot} $ applied component-wise, i.e., $ \prob{ y \mid z } = \prod_{i = 1}^n Q\paren{y_i \mid z_i} $. 

We assume that $ \beta_*, X, \eps $ are mutually independent each subject to the following distributional assumption. 
\begin{enumerate}[label=(A\arabic*)]
\setcounter{enumi}{\value{asmpctr}}

    \item\label[asmp]{asmp:signal} $ \beta_* \sim P_{\sfB_*}^{\ot d} $ where $ P_{\sfB_*} $ is a distribution on $\bbR$ with mean $0$ and variance $ \rho \in (0,\infty) $. 

    \item\label[asmp]{asmp:design} $ X $ is bi-orthogonally invariant in law, meaning $ X \eqqlaw \wt{O} X \wt{Q}^\top $ for any deterministic $ \wt{O} \in \bbO(n), \wt{Q} \in \bbO(d) $. 
    Equivalently, $X$ is identified with its singular value decomposition $ X = O \Lambda Q^\top $ where $ (O, Q) \sim \haar(\bbO(n)) \ot \haar(\bbO(d)) $ is independent of the rectangular diagonal matrix $ \Lambda  \in \bbR^{n \times d}$. 
    Letting $ \lambda \in \bbR^{n\wedge d} $ be the vector containing the diagonal elements of $ \Lambda $, we further assume that almost surely, the empirical distribution of $ \lambda $ converges in Wasserstein (of every order) to the law $\mu_{\sfLambda}$ of a random variable $ \sfLambda $ that is compactly supported on $ (0,\infty) $ and has positive variance. 
    Moreover, $ \max_{1\le i \le n\wedge d} \lambda_{i}, \max_{1\le i \le n\wedge d} \lambda_{i}^{-1} $ are uniformly bounded over $n,d$.
    Almost surely, for every $\varsigma > 0$, there exists $N(\varsigma)$ such that for every $n\wedge d\ge N(\varsigma)$, 
    \begin{align}
        & \brace{\lambda_i}_{1\le i\le n\wedge d} \subset \supp\mu_{\sfLambda} + (-\varsigma, \varsigma) . \label{eqn:strong} 
    \end{align}
    The following two limits hold
    \begin{align}
        & \lim_{z\nearrow\inf\supp\mu_{\sfLambda}} \expt{\frac{\sfLambda}{\sfLambda - z}} = 
        \lim_{z\searrow\sup\supp\mu_{\sfLambda}} \expt{\frac{\sfLambda}{z - \sfLambda}} = \infty . \label{eqn:hard} 
    \end{align}

    \item\label[asmp]{asmp:noise} $ \eps \in \bbR^n $ satisfies $ \eps \wto \sfE $ for a real-valued random variable $ \sfE $ with finite variance.  
    Denote by $ P_{\sfE} $ the law of $\sfE$. 

    \item\label[asmp]{asmp:proportional} $ n,d\to\infty $ with $ n/d\to\delta\in(0,\infty) $. All limits, unless otherwise specified, are taken with respect to this scaling. 

\setcounter{asmpctr}{\value{enumi}}
\end{enumerate}

We make a few comments on \Cref{asmp:design}. 
The condition \Cref{eqn:strong} therein implies that the empirical spectral distribution of $X$ converges \emph{strongly} to $ \mu_{\sfLambda} $ in the sense that the top and bottom eigenvalues of $X$ converge respectively to the upper and lower edges of the support of $ \sfLambda $, i.e., 
\begin{align}
    & 0 < \lim_{n\wedge d \to \infty} \min_{1\le i \le n\wedge d} \lambda_{i} = \inf\supp \sfLambda
    < \sup\supp \sfLambda = \lim_{n\wedge d \to \infty} \max_{1\le i \le n\wedge d} \lambda_{i} < \infty
    . \notag 
\end{align}
This in particular excludes outlying singular values to the left or right of the spectrum of $X$. 
Otherwise, such outliers may result in outliers in the spectrum of $D$ constructed as \Cref{eqn:specintro} for spectral estimators \cite[Theorem 2.5]{Belinschi_etal},  causing confusion with the informative outlier arising from $\beta_*$. 

The limits in \Cref{eqn:hard} require $ \sfLambda $ to have hard edges on the leftmost and rightmost boundaries of its support. 
Similar assumptions are common in the literature of high-dimensional statistics; see e.g.\ \cite[Assumption 6 `Thickness of the bulk edge']{Donoho_Gavish_Romanov} and \cite[(A6), (A7)]{Zhang_Ji_Venkataramanan_Mondelli}. 
We believe that this is a purely technical assumption and can be removed via a perturbation argument \cite[Appendix C]{Zhang_Ji_Venkataramanan_Mondelli}. 
Indeed, simulations in \Cref{sec:experiments} numerically validate our theoretical results for $ \sfLambda $ proportional to $ \Beta(3,1) $ and $ \Beta(1,2) $ for which the first and second limit in \Cref{eqn:hard} is finite, respectively. 
Similarly, since $ \inf \supp \Beta(3,1) = \inf\supp \Beta(1,2) = 0 $, the assumption that the left edge of the support of $\sfLambda$ is bounded away from $0$ is technical and may be removed. 
On the other hand, the bounded support assumption on $\sfLambda$ is likely needed, otherwise the spectrum of $D$ in \Cref{eqn:specintro} will be unbounded \cite[Theorem 5.9]{Belinschi_Bercovici_Ho}. 

The bi-orthogonal invariance assumption on the law of $X$ captures situations where the left and right singular vectors of $X$ are generic and independent of each other and of the singular values. 
This unlocks analytic tractability of the model under consideration. 
However, numerical experiments in \Cref{sec:phase_real} suggest that our theoretical predictions remain accurate for certain design ensembles that are not orthogonally invariant or that are even deterministic. 
Identifying suitable universality classes of $X$ and justifying within them the validity of our theory is left for future work. 
Progress in this direction includes \cite{Dudeja_Lu_Sen,Dudeja_Sen_Lu,Wang_Zhong_Fan}, although they are not directly applicable to our proposed algorithms. 

Assuming that $ q, P_{\sfB_*}, \law(\sfLambda), P_{\sfE} $ are known (see \Cref{sec:discuss} for an additional discussion on this requirement), our interest is in constructing estimators of $ \beta_*\in\bbR^d $ from the data $ (X,y) $, and  
the quality of an estimator $ \wh{\beta} \equiv \wh{\beta}(X,y) \in \bbR^d $ is measured by its asymptotic overlap with $ \beta_* $ defined in \Cref{eqn:overlap}. 

We introduce some notation used frequently in the rest of the paper. 
Given the random variable $ \sfLambda $ in \Cref{asmp:design}, we also define $\sfLambda_n, \sfLambda_d$ by 
\begin{align}
&&
    \law(\sfLambda_n) &\coloneqq \begin{cases}
        \law(\sfLambda) , & \textnormal{if } \delta \le 1 \\
        \delta^{-1} \law(\sfLambda) + (1 - \delta^{-1}) \delta_0 , & \textnormal{if } \delta > 1
    \end{cases} , &
    \law(\sfLambda_d) &\coloneqq \begin{cases}
        \delta \law(\sfLambda) + (1 - \delta) \delta_0 , & \textnormal{if } \delta \le 1 \\
        \law(\sfLambda) , & \textnormal{if } \delta > 1
    \end{cases} . & 
& \label{eqn:sfLambda_nd}
\end{align}
By definition, $ \sfLambda_n, \sfLambda_d $ are distributed according to the limiting spectral distribution of $ X $, accounting for $ n,d $ singular values, respectively. 
Also, the limiting spectral distributions of $ XX^\top, X^\top X $ are given respectively by the laws of $ \sfLambda_n^2, \sfLambda_d^2 $. 
Let $ (\ol{\kappa}_{2k})_{k\ge0} \subset \bbR $ and $ (\ol{m}_{2k})_{k\ge0} \subset\bbR_{\ge0} $ be the rectangular free cumulants and moments, respectively, of $ \sfLambda_n^2 $; see \Cref{sec:cumulant} for the relevant background. 
Let 
\begin{align}
&& 
    \sigma^2 &\coloneqq \ol{\kappa}_2 \rho , &
& \label{eqn:sigma2} 
\end{align}
where $ \rho $ is the variance of $ P_{\sfB_*} $ from \Cref{asmp:signal}. 
Note that $ \sigma^2 > 0 $ since \Cref{asmp:design} implies $ \ol{\kappa}_2 = \expt{\sfLambda_n^2} > 0 $. 
Let
\begin{align}
&& 
    (\sfB_*, \sfZ, \sfE) &\sim P_{\sfB_*} \ot \cN(0, \sigma^2) \ot P_{\sfE} , & 
    \sfY &= q(\sfZ, \sfE) . &
& \label{eqn:BZEY} 
\end{align}
We also define $ \ol{\sfZ} \coloneqq \sigma^{-1} \sfZ \sim \cN(0,1) $. 
It is easy to verify that $ \beta_* \wto \sfB_* $, $ \matrix{ z & \eps & y} \wto \matrix{\sfZ & \sfE & \sfY} $. 

\section{Main results}

This section presents our main formal results. 
\Cref{sec:spec} contains a precise characterization of the spectral properties of the matrix $D$ in \Cref{eqn:specintro} (in particular the overlap \eqref{eqn:overlap} of its principal eigenvector) and the optimal choice of preprocessing function. 
\Cref{sec:spec_GVAMP} introduces our generic GVAMP formalism and characterizes its dynamics by a state evolution recursion when initialized with spectral estimators from \Cref{sec:spec}. 
\Cref{sec:BVAMP} describes GVAMP with the optimal denoisers, leading to the central algorithmic deliverable of this paper---the spectrally initialized Bayes-GVAMP---with exact asymptotic guarantees. 
\Cref{sec:replica} provides evidence on its computational optimality by connecting its performance to the conjectured Bayes risk. 

\subsection{Optimal spectral estimators}
\label{sec:spec}

Fix any preprocessing function $ \cT \colon \bbR \to \bbR $. 
Let
\begin{subequations}
\label{eqn:spec_est}
\begin{align}
&& 
    T &\coloneqq \diag(\cT(y)) , & 
    D &\coloneqq X^\top T X , &
& \label{eqn:D} 
\end{align}
where $ \cT $ acts component-wise on $y$. 
The spectral estimator outputs 
\begin{align}
\wh{\beta}_{\textnormal{s}} &\coloneqq v_1(D) . \label{eqn:v1} 
\end{align}
\end{subequations}

We will consider preprocessing functions $\cT$ subject to the following condition. 
\begin{enumerate}[label=(A\arabic*)]
\setcounter{enumi}{\value{asmpctr}}

    \item\label[asmp]{asmp:preprocess} The preprocessing function $ \cT \colon \bbR \to \bbR $ takes the form
    \begin{align}
    &&
        \cT(y) &= \frac{\ol{g}(y)}{\ol{g}(y) + \gamma} \, , & 
        \textnormal{where }
        \gamma &\coloneqq \paren{\frac{\ol{\kappa}_4}{\ol{\kappa}_2^2} + \delta} \expt{\ol{g}(\sfY) \ol{\sfZ}^2} , &
    & \label{eqn:T} 
    \end{align}
    for some function $ \ol{g}$ bounded from above and satisfying
    \begin{align}
    &&
        \inf_{y\in\supp\sfY} \ol{g}(y) &> - \gamma , & 
        \expt{ \ol{g}(\sfY) } &= 0 , & 
        \expt{\ol{g}(\sfY) \ol{\sfZ}^2} &> 0 . & 
    & 
    \label{eqn:gbar} 
    \end{align}
    Furthermore, 
    \begin{align}
        \lim_{z\searrow\sup\supp\mu_{\sfT}} \expt{\frac{\sfT}{z - \sfT}} &= \infty , \label{eqn:T_edge} 
    \end{align}
    where $ \sfT = \cT(\sfY) $. 
    
\setcounter{asmpctr}{\value{enumi}}
\end{enumerate}

The first inequality in \Cref{eqn:gbar} ensures that $ \cT(y) $ in \Cref{eqn:T} is well-defined for any $y\in\supp\sfY$ in that the denominator $ \ol{g}(y) + \gamma $ is positive. 
Note that the positivity assumption on $ \expt{\ol{g}(\sfY) \ol{\sfZ}^2} $ implies $ \gamma > 0 $. 
This is because by the moment-cumulant relation (see \Cref{sec:cumulant}), 
\begin{align}
    \frac{\ol{\kappa}_4}{\ol{\kappa}_2^2} + \delta &= \frac{\ol{m}_4 - \ol{m}_2^2}{\ol{m}_2^2} > 0 , \label{eqn:k4-by-k2+delta} 
\end{align}
where the inequality follows by recalling the definition of $ \ol{m}_{2k} $ and applying Cauchy--Schwarz $ \expt{\sfLambda_n}^2 \le \expt{\sfLambda_n^2} $. 
The positivity of $ \expt{\ol{g}(\sfY) \ol{\sfZ}^2} $ is without loss of generality. 
If $ \expt{\ol{g}(\sfY) \ol{\sfZ}^2} < 0 $, one can instead consider $ -\ol{g}(y) $. 
Then, under the mapping $ \ol{g}(y) \mapsto - \ol{g}(y) $, $ \gamma \mapsto - \gamma $ and the function $ \cT $ remains unchanged. 
We exclude the case where $ \expt{\ol{g}(\sfY) \ol{\sfZ}^2} = 0 $ since that induces the trivial function $ \cT(y) = 1 $ for every $y$. 
The condition \Cref{eqn:T_edge} on the upper edge of $\supp\mu_{\sfT}$ is akin to \Cref{eqn:hard} imposed on $\supp\mu_{\sfLambda}$ and can potentially be removed. 

\Cref{asmp:preprocess} is mild in the sense that there exists a preprocessing function subject to it that attains the conjecturally optimal weak recovery threshold; see \Cref{thm:opt_thr}. 
The function class defined by \Cref{asmp:preprocess} is largely motivated by our proof technique based on linearized GVAMP; see \Cref{sec:heuristics}. 
Though higher asymptotic overlap can be achieved by considering a larger function class (see \Cref{sec:experiments,sec:discuss}), those preprocessing functions do not lead to a weak recovery threshold lower than that in \Cref{thm:opt_thr} and their advantage in overlap is superseded by subsequent refinement of Bayes-GVAMP (\Cref{thm:SE_Bayes_GVAMP}). 
On the other hand, spectral estimators with preprocessing functions satisfying \Cref{asmp:preprocess} are more convenient to refinement by GVAMP since they arise from linearization of GVAMP. 

\Cref{thm:spec} below characterizes the limit of the top two eigenvalues of $D$ and the asymptotic overlap of the spectral estimator, provided a criticality condition. 
To state these characterizations, we need a few definitions. 
Recall $ \sfT = \cT(\sfY) $ and write $ \sfS = \sfLambda_n^2 $ for brevity. 
Define the function $ \psi \colon (\sup\supp\sfT, \infty) \to \bbR $ as
\begin{align}
    \psi(a) &= a \expt{ \frac{\sfS}{\omega(a) - \sfS} } \expt{ \frac{1}{\omega(a) - \sfS} }^{-1} 
    = \omega(a) \expt{ \frac{\sfT}{a - \sfT} } \expt{ \frac{1}{a - \sfT} }^{-1} , \label{eqn:def_psi} 
\end{align}
where $ \omega(a) $ is implicitly defined as the unique solution $ \omega \in \bbR \setminus (\inf\supp\sfLambda^2, \sup\supp\sfLambda^2) $ to 
\begin{align}
    \expt{ \frac{\sfT}{a - \sfT} } &= \expt{ \frac{\sfS}{\omega - \sfS} } , \label{eqn:omega} 
\end{align}
or equivalently (by adding $1$ to both sides of \Cref{eqn:omega}), 
\begin{align}
    \expt{ \frac{a}{a - \sfT} } &= \expt{ \frac{\omega}{\omega - \sfS} } . \label{eqn:omega_equiv} 
\end{align}
Define 
\begin{align}
    \lambda^\circ &= \psi(a^\circ) , 
    \label{eqn:def_lambda2}
\end{align}
where $ a^\circ \in (\sup\supp\sfT, \infty) $ is the largest critical point of $ \psi $, i.e., the largest solution in $(\sup\supp\sfT, \infty)$ to  $ \psi'(a^\circ) = 0 $, or more explicitly, 
\begin{align}
    \frac{\expt{ \frac{\sfT}{a^\circ - \sfT} } \expt{\frac{1}{a^\circ - \sfT}}}{\expt{ \frac{\sfT}{(a^\circ - \sfT)^2} }} 
    + \frac{\expt{ \frac{\sfS}{\omega(a^\circ) - \sfS} } \expt{\frac{1}{\omega(a^\circ) - \sfS}}}{\expt{ \frac{\sfS}{(\omega(a^\circ) - \sfS)^2} }} 
    &= 1 . \label{eqn:a0_explicit}
\end{align}
Finally, let 
\begin{align}
    \eta &\coloneqq \delta \frac{
        \paren{ \frac{\ol{\kappa}_4}{\ol{\kappa}_2^2} + \delta }^2 \expt{\ol{g}(\sfY) \ol{\sfZ}^2}^2
        - \paren{ \frac{\ol{\kappa}_4}{\ol{\kappa}_2^2} + \delta } \expt{\ol{g}(\sfY)^2}
    }{
        \paren{ \frac{\ol{\kappa}_4}{\ol{\kappa}_2^2} + \delta }^3 \expt{\ol{g}(\sfY) \ol{\sfZ}^2}^2
        + \paren{ \frac{\ol{\kappa}_4}{\ol{\kappa}_2^2} + \delta }^2 \expt{\ol{g}(\sfY)^2 \ol{\sfZ}^2} 
        + \paren{ \frac{\ol{\kappa}_6}{\ol{\kappa}_2^3} - 2 \frac{\ol{\kappa}_4^2}{\ol{\kappa}_2^4} - \delta \paren{\frac{\ol{\kappa}_4}{\ol{\kappa}_2^2} + \delta} } \expt{\ol{g}(\sfY)^2}
    } . \label{eqn:eta} 
\end{align}

\begin{theorem}[Spectral estimator]
\label{thm:spec}
Consider the spectral estimator defined through \Cref{eqn:spec_est} using a preprocessing function subject to \Cref{asmp:preprocess}. 
Then we have almost surely, 
\begin{align}
    \lim_{d\to\infty} \lambda_2(D) &= \lambda^\circ , \label{eqn:lambda2_D} 
\end{align}
where $\lambda^\circ$ is defined in \Cref{eqn:def_lambda2}. 
Furthermore, if 
\begin{align}
    \paren{ \frac{\ol{\kappa}_4}{\ol{\kappa}_2^2} + \delta } \expt{\ol{g}(\sfY) \ol{\sfZ}^2}^2 &> \expt{\ol{g}(\sfY)^2} , \label{eqn:thr} 
\end{align}
then we have the following: 
\begin{align}
&&
    \lim_{d\to\infty} \lambda_1(D) &= \ol{\kappa}_2 > \lambda^\circ , &
    \lim_{d\to\infty} \frac{\abs{ \inprod{\beta_*}{v_1(D)} }^2}{\normtwo{\beta_*}^2} 
    &= \eta \in (0,1),&
& \label{eqn:lambda1_eta} 
\end{align}
where $\eta$ is defined in \Cref{eqn:eta}.
\end{theorem}
\Cref{thm:spec} follows from a sequence of results in \Cref{sec:lin_VAMP_SE} whose proofs span \Cref{sec:pf_spec_step1,app:right_edge,app:pf_lem:align}. 
See \Cref{sec:heu_lin,sec:heu_D} for a detailed overview of the proof. 


\begin{remark}[Universality]
The value of the asymptotic overlap $ \eta $ in \Cref{eqn:eta} enjoys a surprising universality property with respect to the limiting spectral distribution of $X$. 
Specifically, the only dependence on $ \sfLambda_n $ in the expression \Cref{eqn:eta} is through $ \ol{\kappa}_2, \ol{\kappa}_4, \ol{\kappa}_6 $ (or equivalently by the moment-cumulant relation, the first three moments of $ \sfLambda_n^2 $). 
This indicates that as far as spectral estimators subject to \Cref{asmp:preprocess} are concerned, the details of the spectral law of the design have minimal impact on the performance of the estimators. 
\end{remark}

The following result characterizes, among the family of preprocessing functions given by \Cref{asmp:preprocess}, the optimal spectral threshold for weak recovery, which is the minimum $\delta$ for which there exists a preprocessing function $\cT$ identified with $\ol{g}$ that satisfies \Cref{eqn:thr}. 

\begin{theorem}[Optimal spectral threshold and preprocessing]
\label{thm:opt_thr}
\,
\begin{enumerate}
    \item \label{itm:ach} Assume that 
    \begin{align}
            \paren{ \frac{\ol{\kappa}_4}{\ol{\kappa}_2^2} + \delta } \int \, \frac{\expt{Q(y\,|\,\sfZ) \paren{\ol{\sfZ}^2 - 1}}^2}{\expt{Q(y\,|\,\sfZ)}} \diff y &> 1 . \label{eqn:opt_thr} 
        \end{align}
    Then, by taking 
    \begin{align}
        \ol{g}(y) &= \frac{\expt{\paren{\ol{\sfZ}^2 - 1} Q(y \,|\, \sfZ)}}{\expt{Q(y \,|\, \sfZ)}}
        = \expt{\ol{\sfZ}^2 - 1 \mid \sfY = y} , 
        \label{eqn:olg_opt}
    \end{align}
    the corresponding preprocessing function defined through \Cref{asmp:preprocess} achieves weak recovery, in the sense that \Cref{eqn:thr} holds. 
    Furthermore, 
    $\ol{g}$ in \Cref{eqn:olg_opt} satisfies \Cref{eqn:gbar}
    provided: 
    \begin{align}
        \var{\expt{\ol{\sfZ}^2 \mid \sfY}} &= \expt{\expt{\ol{\sfZ}^2 - 1 \mid \sfY}^2} > 0 . \label{eqn:nontrivial}
    \end{align}

    \item \label{itm:conv} Conversely, if there is a preprocessing function subject to \Cref{asmp:preprocess} achieving weak recovery (in the sense that \Cref{eqn:thr} holds), then \Cref{eqn:opt_thr} must also hold.
\end{enumerate}
\end{theorem}

\Cref{thm:opt_thr} is proved in \Cref{app:pf_thm:opt_thr}. 
\Cref{itm:conv} follows from a Cauchy--Schwarz argument applied to the criticality condition \Cref{eqn:thr}. 
The optimal choice \Cref{eqn:olg_opt} of $ \ol{g} $ in \Cref{itm:ach} can be read off from the equality case of the Cauchy--Schwarz argument. 

\begin{remark}[Conjectured computationally optimal weak recovery threshold]
Using the moment-cumulant relation (see \Cref{sec:cumulant})
\begin{align}
\frac{\ol{\kappa}_4}{\ol{\kappa}_2^2} + \delta &= \delta \frac{\ol{m}_4'}{\ol{m}_2'^2} - 1 , \notag
\end{align}
we can write the condition \Cref{eqn:opt_thr} equivalently as 
\begin{align}
\delta &> \frac{\ol{m}_2'^2}{\ol{m}_4'} \brack{ 1 + \paren{ \int \frac{\expt{\paren{\ol{\sfZ}^2 - 1} Q(y \,|\, \sfZ)}^2}{ \expt{Q(y \,|\, \sfZ)} } \diff y }^{-1} } . \label{eqn:thr_mmt}
\end{align}
The above condition coincides with \cite[(11)]{Maillard_etal_threshold} which is conjectured to be the computationally optimal weak recovery threshold, i.e., no other polynomial-time algorithm is believed to achieve positive asymptotic overlap if \Cref{eqn:thr_mmt} is violated. 
We observe that the threshold given by \Cref{eqn:thr_mmt} (or equivalently \Cref{eqn:opt_thr}) exhibits remarkable universality with respect to the limiting spectral distribution of $X$ in that it only depends on the law of $ \sfLambda_d^2 $ through its first two moments (or equivalently its first two free cumulants). 
This implies that the threshold \Cref{eqn:thr_mmt} remains the same for all designs with matching first two moments of their limiting spectral distributions. 
\end{remark}

\begin{remark}[Effect of non-Gaussian design on spectral threshold]
Recall from \cite[Remark 5.1]{fan2020approximate} that when $X_{i,j}\iid\cN(0,1/d)$, $ \sfLambda_n^2 $ is distributed according to the Marchenko--Pastur law whose rectangular free cumulants are $ \ol{\kappa}_2 = 1 $ and $ \ol{\kappa}_{2k} = 0 $ for $k\ge2$. 
Using this in \Cref{eqn:opt_thr} recovers the optimal spectral threshold identified in \cite[Equation (20)]{Mondelli_Montanari} for i.i.d.\ Gaussian design. 
We therefore see clearly that the effect of non-Gaussian design on the optimal spectral threshold is through the cumulant ratio $ \ol{\kappa}_4 / \ol{\kappa}_2^2 $ on the left of \Cref{eqn:opt_thr}. 
\end{remark}

\subsection{Spectrally initialized GVAMP and its state evolution}
\label{sec:spec_GVAMP}

This section introduces spectrally initialized GVAMP and characterizes its limiting dynamics using a deterministic state evolution recursion. 
A generic GVAMP iteration updates its iterates according to the following rules: for every $t\ge0$, 
\begin{subequations}
\label{eqn:GVAMP}    
\begin{align}
&&
    r^{t+1} &= \Phi_t(X^\top X) \wt{r}^t + \wt{\Phi}_t(X)^\top \wt{p}^t , & 
    \wt{r}^{t+1} &= f_{t+1}(r^{t+1}; \beta_*, \Theta) , & 
& \label{eqn:GVAMP1} \\ 
&&
    p^{t+1} &= \Psi_t(XX^\top) \wt{p}^t + \wt{\Psi}_t(X) \wt{r}^t , & 
    \wt{p}^{t+1} &= g_{t+1}(p^{t+1}; z, \eps, \Xi) , &
& \label{eqn:GVAMP2}
\end{align}
\end{subequations}
with initialization $ \wt{r}^0 \in \bbR^d , \wt{p}^0 \in \bbR^n $. 
The matrices $ \Theta \in \bbR^{d\times d'} , \Xi \in \bbR^{n\times n'} $ are side information subject to \Cref{asmp:side}. 
\begin{enumerate}[label=(A\arabic*)]
\setcounter{enumi}{\value{asmpctr}}

    \item \label[asmp]{asmp:side} The matrices $ \Theta \in \bbR^{d\times d'} , \Xi \in \bbR^{n\times n'} $ are independent of each other and of $O,Q$ from \Cref{asmp:design}. 
    They satisfy the following Wasserstein limits: 
    \begin{align}
    &&
        \matrix{\beta_* & \Theta} &\wto \matrix{\sfB_* & \sfTheta^\top} , & \matrix{z & \eps & \Xi} &\wto \matrix{\sfZ & \sfE & \sfXi^\top} , & 
    & \notag 
    \end{align}
    as $d,n\to\infty$ and $d',n'$ are fixed.
    
\setcounter{asmpctr}{\value{enumi}}
\end{enumerate}
The functions $ \Phi_t, \Psi_t \colon \bbR \to \bbR $ apply to each eigenvalue of their symmetric input matrices separately, retaining the eigenvectors; 
and $ \wt{\Phi}_t, \wt{\Psi}_t \colon \bbR \to \bbR $ apply to each singular value of their rectangular input matrices separately, retaining the singular vectors. 
For instance, for a symmetric matrix $ A $, $ \Phi_t(A) $ is understood as $ \Phi_t(A) = \sum_i \Phi_t(\lambda_i(A)) v_i(A) v_i(A)^\top $. 
The functions $ f_{t+1} \colon \bbR^{2+d'}\to\bbR , g_{t+1} \colon \bbR^{3+n'} \to \bbR $ are applied row-wise to their inputs, e.g., $f_{t+1}(r^{t+1}; \beta_*, \Theta) = \matrix{f_{t+1}(r^{t+1}_1; \beta_{*,1}, \Theta_{1,:}) & \cdots & f_{t+1}(r^{t+1}_d; \beta_{*,d}, \Theta_{d,:})}^\top$. 
See \Cref{sec:heu_GVAMP} for motivation behind the formalism of \Cref{eqn:GVAMP}. 

The initializers $ \wt{r}^0, \wt{p}^0 $ are obtained from the family of spectral estimators described in \Cref{sec:spec}. 
Specifically, consider a matrix $D$ as in \Cref{eqn:D} constructed from a preprocessing function of the form \Cref{eqn:T} identified with some $ \ol{g} $. 
In particular, let $ \gamma $ be defined through \Cref{eqn:T} using this $ \ol{g} $. 
We will derive a state evolution result for the GVAMP iteration \Cref{eqn:GVAMP} initialized with (functions of) the leading eigenvector $ v_1(D) $. 
We consider initializers $ \wt{r}_0 = f_0(r^0) $ and $ \wt{p}^0 = g_0(p^0; y) $, where 
\begin{align}
&&
    r^0 &= s \, c_r \, \sqrt{d} \, v_1(D) , & 
    p^0 &= c_p \, \gamma (\gamma I_n + \ol{G})^{-1} X \, \sqrt{d} \, v_1(D) , & 
& \label{eqn:spec_GVAMP_init} 
\end{align}
with 
\begin{align}
&& 
    \ol{G} &= \diag(\ol{g}(y)) , & 
    s &= \sgn\paren{\inprod{v_1(D)}{\beta_*}} , & 
    c_r &= \frac{\sqrt{\rho + w_2}}{w_2} , & 
    c_p &= \frac{\sqrt{\rho + w_2}}{w_1 w_3} . & 
& \label{eqn:cr_cp} 
\end{align}
The constants $ w_1,w_2,w_3 $ in \Cref{eqn:cr_cp} along with $w_4$ are given below: 
\begin{align}
    w_1 &\coloneqq \rho \frac{
        \paren{ \frac{\ol{\kappa}_4}{\ol{\kappa}_2} + \delta \ol{\kappa}_2 } \expt{\ol{g}(\sfY)^2 \ol{\sfZ}^2} + \paren{\frac{\ol{\kappa}_6}{\ol{\kappa}_2^2} - \frac{\ol{\kappa}_4^2}{\ol{\kappa}_2^3} + \delta \frac{\ol{\kappa}_4}{\ol{\kappa}_2}} \expt{\ol{g}(\sfY) \ol{\sfZ}^2}^2
    }{
        \paren{ \frac{\ol{\kappa}_4}{\ol{\kappa}_2^2} + \delta }^2 \expt{\ol{g}(\sfY) \ol{\sfZ}^2}^2
        - \paren{ \frac{\ol{\kappa}_4}{\ol{\kappa}_2^2} + \delta } \expt{\ol{g}(\sfY)^2}
    } , \label{eqn:tau2/nu2_FP} \\
    w_2 &\coloneqq \frac{\rho}{\delta} 
    \frac{
        \paren{ \frac{\ol{\kappa}_4}{\ol{\kappa}_2^2} + \delta }^2 \expt{\ol{g}(\sfY)^2 \ol{\sfZ}^2} 
        + \frac{\ol{\kappa}_4}{\ol{\kappa}_2^2}
        \paren{ \frac{\ol{\kappa}_4}{\ol{\kappa}_2^2} + \delta }^2 \expt{\ol{g}(\sfY) \ol{\sfZ}^2}^2
        + \paren{ \frac{\ol{\kappa}_6}{\ol{\kappa}_2^3} - 2 \frac{\ol{\kappa}_4^2}{\ol{\kappa}_2^4} } \expt{\ol{g}(\sfY)^2} 
    }{
        \paren{ \frac{\ol{\kappa}_4}{\ol{\kappa}_2^2} + \delta }^2 \expt{\ol{g}(\sfY) \ol{\sfZ}^2}^2
        - \paren{ \frac{\ol{\kappa}_4}{\ol{\kappa}_2^2} + \delta } \expt{\ol{g}(\sfY)^2}
    } , \label{eqn:omega2/chi2_FP} \\
    w_3 &= 1 + \frac{\ol{\kappa}_4}{\delta \ol{\kappa}_2^2} , \qquad 
    w_4 = \delta \expt{\ol{g}(\sfY) \ol{\sfZ}^2} . 
    \label{eqn:w34}
\end{align}
Note from \Cref{eqn:T} that $ w_3 w_4 = \gamma $. 
Moreover, the denominators of $w_1,w_2$ are positive provided \Cref{eqn:thr} holds. 


\paragraph{State Evolution.} The dynamics of \Cref{eqn:GVAMP} are accurately predicted by a deterministic \emph{state evolution} recursion in the high-dimensional limit. Indeed, \Cref{thm:spec_GVAMP} below shows that  for $t \ge 1$, the empirical distribution of $r^t$ converges to the law of $b_t \sfB_* + \sfK_t$, where $b_t$ is a deterministic scalar  and $\sfK_t$ is a zero-mean Gaussian independent of $\sfB_*$. Similarly, the empirical distribution of $p^t$ converges to the law of $a_t \sfZ + \sfJ_t$, for a scalar $a_t$ and a zero-mean Gaussian $\sfJ_t$ independent of $\sfZ$. 
We recursively define these state evolution parameters as follows. 

Let 
\begin{align}
&& 
    a_0 &= \frac{1}{w_1} , & 
    b_0 &= \frac{1}{w_2} . & 
& \label{eqn:a0_b0}
\end{align}
Define 
\begin{align}
&&
    \sfR_0 &= b_0 \sfB_* + \sfK_0 , & 
    \sfP_0 &= a_0 \sfZ + \sfJ_0 , & 
    \wt{\sfR}_0 &= f_0(\sfR_0) , & 
    \wt{\sfP}_0 &= g_0(\sfP_0; \sfY) , & 
& \label{eqn:R0_P0} 
\end{align}
where 
\begin{align}
&&
    (\sfB_*, \sfK_0) &\sim P_{\sfB_*} \ot \cN(0, b_0) , & 
    (\sfZ, \sfJ_0) &\sim \cN(0,\sigma^2) \ot \cN(0, a_0) , & 
    \sfY &\sim Q\paren{\cdot \mid \sfZ} . & 
& \label{eqn:J0_K0} 
\end{align}
Furthermore, let 
\begin{align}
&&
    \wt{a}_0 &= \sigma^{-2} \expt{\wt{\sfP}_0 \sfZ} , &
    \wt{b}_0 &= \rho^{-1} \expt{\wt{\sfR}_0 \sfB_*} , &
& 
\label{eqn:wta0_wtb0}
\end{align}
and $ \wt{\sfP}_{-1} = \gamma^{-1} \ol{g}(\sfY) \sfP_0, \wt{a}_{-1} = \gamma^{-1} a_0 \expt{\ol{g}(\sfY) \ol{\sfZ}^2} $. 

Let $ (\sfJ_{r})_{0\le r\le t} $ and $ (\sfK_{r})_{0\le r\le t} $ be two centered Gaussian processes that are independent of each other and of $(\wt{\sfR}_0, \wt{\sfP}_0, \sfB_*, \sfZ, \sfE, \sfLambda_n, \sfLambda_d)$, with covariance structures specified below.
The law of $ \sfJ_0,\sfK_0 $ is given in \Cref{eqn:J0_K0}. 
For any $ t\ge0 $, the correlation between $ \sfJ_{t+1} $ (resp.\ $ \sfK_{t+1} $) and $ \sfJ_0 $ (resp.\ $ \sfK_0 $) is given by: 
\begin{subequations}
\label{eqn:SE_GVAMP}
\begin{align}
    \expt{\sfJ_{t+1} \sfJ_{0}}
    &= \ol{\kappa}_2^{-1} \expt{ \Psi_t(\sfLambda_n^2) \sfLambda_n^2 } \expt{\wt{\sfP}_t \wt{\sfP}_{-1}}
    - 2 \expt{ \Psi_t(\sfLambda_n^2) \sfLambda_n^2 } \wt{a}_t \wt{a}_{-1} \rho
    + \ol{\kappa}_2^{-1} \expt{ \Psi_t(\sfLambda_n^2) \sfLambda_n^4 } \wt{a}_t \wt{a}_{-1} \rho \notag \\
    &\quad + \expt{ (\ol{\kappa}_2^{-1} \sfLambda_n^2 - 1) \wt{\Psi}_t(\sfLambda_n) \sfLambda_n } \wt{a}_{-1} \wt{b}_t \rho 
    - a_{t+1} a_{0} \sigma^2 , \label{eqn:J_cov0} \\
    \expt{ \sfK_{t+1} \sfK_{0} }
    &= \ol{\kappa}_2^{-1} \expt{ \wt{\Phi}_t(\sfLambda_d) \sfLambda_d } \expt{\wt{\sfP}_t \wt{\sfP}_{-1}}
    + \expt{ (\ol{\kappa}_2^{-1} \sfLambda_d^2 - 1) \wt{\Phi}_t(\sfLambda_d) \sfLambda_d } \wt{a}_t \wt{a}_{-1} \rho \notag \\
    &\quad + \ol{\kappa}_2^{-1} \expt{ \Phi_t(\sfLambda_d^2) \sfLambda_d^2 } \wt{a}_{-1} \wt{b}_t \rho 
    - b_{t+1} b_{0} \rho . \label{eqn:K_cov0}
\end{align}
For every $0\le r,s\le t$, 
\begin{align}
    \expt{\sfJ_{r+1} \sfJ_{s+1}}
    &= \expt{ \Psi_r(\sfLambda_n^2) \Psi_s(\sfLambda_n^2) } \paren{ \expt{\wt{\sfP}_r \wt{\sfP}_s} - \wt{a}_r \wt{a}_s \sigma^2 }
    + \expt{ \Psi_r(\sfLambda_n^2) \Psi_s(\sfLambda_n^2) \sfLambda_n^2 } \wt{a}_r \wt{a}_s \rho \notag \\
    &\quad + \expt{ \wt{\Psi}_r(\sfLambda_n) \wt{\Psi}_s(\sfLambda_n) } \expt{ \wt{\sfR}_r \wt{\sfR}_s } 
    - a_{r+1} a_{s+1} \sigma^2 \notag \\
    &\quad + \expt{ \Psi_r(\sfLambda_n^2) \wt{\Psi}_s(\sfLambda_n) \sfLambda_n } \wt{a}_r \wt{b}_s \rho
    + \expt{ \Psi_s(\sfLambda_n^2) \wt{\Psi}_r(\sfLambda_n) \sfLambda_n } \wt{a}_s \wt{b}_r \rho , \label{eqn:J_cov} \\
    \expt{ \sfK_{r+1} \sfK_{s+1} }
    &= \expt{ \Phi_r(\sfLambda_d^2) \Phi_s(\sfLambda_d^2) } \expt{\wt{\sfR}_r \wt{\sfR}_s} 
    - b_{r+1} b_{s+1} \rho \notag \\
    &\quad + \expt{ \wt{\Phi}_r(\sfLambda_d) \wt{\Phi}_s(\sfLambda_d) } \paren{ \expt{\wt{\sfP}_r \wt{\sfP}_s} - \wt{a}_r \wt{a}_s \sigma^2 }
    + \expt{ \wt{\Phi}_r(\sfLambda_d) \wt{\Phi}_s(\sfLambda_d) \sfLambda_d^2 } \wt{a}_r \wt{a}_s \rho \notag \\
    &\quad + \expt{ \Phi_r(\sfLambda_d^2) \wt{\Phi}_s(\sfLambda_d) \sfLambda_d } \wt{a}_s \wt{b}_r \rho
    + \expt{ \Phi_s(\sfLambda_d^2) \wt{\Phi}_r(\sfLambda_d) \sfLambda_d } \wt{a}_r \wt{b}_s \rho . \label{eqn:K_cov} 
\end{align}
In \Cref{eqn:J_cov0,eqn:K_cov0,eqn:J_cov,eqn:K_cov}, the  scalars $a_{r+1}, b_{r+1}$ are defined as
\begin{align}
&&
    \wt{a}_r &= \frac{1}{\sigma^2} \expt{ \wt{\sfP}_r \sfZ } , &
    \wt{b}_r &= \frac{1}{\rho} \expt{ \wt{\sfR}_r \sfB_* } , &
& \label{eqn:ab_tilde} \\
&&
    a_{r+1} &= \frac{1}{\ol{\kappa}_2} \paren{ \expt{ \Psi_r(\sfLambda_n^2) \sfLambda_n^2 } \wt{a}_r + \expt{ \wt{\Psi}_r(\sfLambda_n) \sfLambda_n } \wt{b}_r } , & 
    b_{r+1} &= \expt{ \wt{\Phi}_r(\sfLambda_d) \sfLambda_d } \wt{a}_r , & 
& \label{eqn:ab}
\end{align}
and the random variables $(\wt{\sfP}_{r+1})_{0\le r\le t}, (\wt{\sfR}_{r+1})_{0\le r\le t}$ are defined as
\begin{align}
&&
    \wt{\sfP}_{r+1} &= g_{r+1}(\sfP_{r+1}; \sfZ, \sfE, \sfXi) , & 
    \wt{\sfR}_{r+1} &= f_{r+1}(\sfR_{r+1}; \sfB_*, \sfTheta), & 
& \label{eqn:wtP_wtR} 
\end{align}
where
\begin{align}
&& 
    \sfP_{r+1} &= a_{r+1} \sfZ + \sfJ_{r+1} , & 
    \sfR_{r+1} &= b_{r+1} \sfB_* + \sfK_{r+1} . &
& \label{eqn:PR} 
\end{align}
\end{subequations}
This completely specifies the laws of the Gaussian processes $ (\sfJ_t)_{t\ge0},(\sfK_t)_{t\ge0} $. 



We are now ready to present the state evolution result for spectrally initialized GVAMP which asserts that the random variables defined above identify the Wasserstein limits of the iterates of \Cref{eqn:GVAMP}. 
For this result to be valid, the matrix denoisers $ \Phi_t, \Psi_t $ and vector denoisers $ f_{t+1}, g_{t+1} $ are required to satisfy the conditions below. 

\begin{enumerate}[label=(A\arabic*)]
\setcounter{enumi}{\value{asmpctr}}

    \item \label[asmp]{asmp:tr_free} $ \Phi_t, \Psi_t $ are \emph{trace-free} for every $ t\ge0 $, that is, $ \expt{\Phi_t(\sfLambda_d^2)} = \expt{\Psi_t(\sfLambda_n^2)} = 0 $. 
    Moreover, $ \Phi_t, \wt{\Phi}_t \in \PG_{1\to1} $ are continuous with probability $1$ under the law of $ \sfLambda_d $; 
    $ \Psi_t, \wt{\Psi}_t \in \PG_{1\to1} $ are continuous 
    with probability $1$ under $ \sfLambda_n $. 

    \item \label[asmp]{asmp:div_free} $ f_{t+1}, g_{t+1} $ are \emph{divergence-free} with respect to their first inputs for every $ t\ge0 $, that is, 
    \begin{align}
    && 
         \expt{ f_{t+1}'(\sfR_{t+1}; \sfB_*, \sfTheta) } &= 0 , & 
        \expt{ g_{t+1}'(\sfP_{t+1}; \sfZ, \sfE, \sfXi) } &= 0 , & 
    & \notag 
    \end{align} 
    where $ f_{t+1}', g_{t+1}' $ denote the partial derivatives of $ f_{t+1}, g_{t+1} $ with respect to the first arguments. 

\setcounter{asmpctr}{\value{enumi}}
\end{enumerate}
Recalling \Cref{asmp:design} and using \Cref{prop:wto5}, \Cref{asmp:tr_free} 
implies that almost surely, 
\begin{align}
    \lim_{d\to\infty} \frac{1}{d} \tr(\Phi_t(X^\top X)) &= \lim_{n\to\infty} \frac{1}{n} \tr(\Psi_t(XX^\top)) = 0 . \notag 
\end{align}
Finally, we isolate two technical assumptions regarding the regularity of matrix and vector denoisers. 
\begin{enumerate}[label=(A\arabic*)]
\setcounter{enumi}{\value{asmpctr}}

\item \label[asmp]{asmp:remove_PhiPsi} For every $t\ge0$, 
\begin{align}
&&
    & \limsup_{d\to\infty} \normtwo{\Phi_t(\Lambda^\top \Lambda)} , & 
    & \limsup_{d\to\infty} \normtwo{\wt{\Phi}_t(\Lambda)} , &
    & \limsup_{n\to\infty} \normtwo{\Psi_t(\Lambda \Lambda^\top)} ,  &
    & \limsup_{n\to\infty} \normtwo{\wt{\Psi}_t(\Lambda)} & 
& \notag 
\end{align}
are all finite almost surely. 

\item \label[asmp]{asmp:remove} 
For every $t\ge0$, $ f_{t+1}, g_{t+1} $ are Lipschitz in all arguments. 
Furthermore, the partial derivatives of $ r \mapsto f_{t+1}(r;b,\theta) $ and $ p \mapsto g_{t+1}(p;z,e,\xi) $ are continuous with probability $1$ with respect to the law of $ (\sfR_{t+1}, \sfB_*, \sfTheta) $ and $ (\sfP_{t+1}, \sfZ, \sfE, \sfXi) $, respectively. 

\setcounter{asmpctr}{\value{enumi}}
\end{enumerate}

\begin{theorem}[State evolution of spectrally initialized GVAMP]
\label{thm:spec_GVAMP}
Consider the GVAMP iteration \Cref{eqn:GVAMP} initialized with $ \wt{r}^0, \wt{p}^0 $ given in \Cref{eqn:spec_GVAMP_init}. 
If \Cref{asmp:preprocess,asmp:side,asmp:tr_free,asmp:div_free,asmp:remove_PhiPsi,asmp:remove} and \Cref{eqn:thr} hold, then for any fixed $t\ge0$, 
\begin{align}
\matrix{
    r^0 & \cdots & r^{t} & 
    \wt{r}^0 & \cdots & \wt{r}^{t} & 
    \beta_* & \Theta
} 
&\xrightarrow{W_2} \matrix{
    \sfR_0 & \cdots & \sfR_{t} & 
    \wt{\sfR}_0 & \cdots & \wt{\sfR}_{t} & 
    \sfB_* & \sfTheta^\top
} , \label{eqn:SE_G_r_spec} \\
\matrix{
    p^0 & \cdots & p^{t} & 
    \wt{p}^0 & \cdots & \wt{p}^{t} & 
    z & \eps & \Xi
} 
&\xrightarrow{W_2} \matrix{
    \sfP_0 & \cdots & \sfP_{t} & 
    \wt{\sfP}_0 & \cdots & \wt{\sfP}_{t} & 
    \sfZ & \sfE & \sfXi^\top
} . \label{eqn:SE_G_p_spec}
\end{align}
\end{theorem}

The full proof of \Cref{thm:spec_GVAMP} is deferred to \Cref{app:pf_SE_GVAMP}. 
See \Cref{sec:heu_GVAMP} for a brief overview of the proof. 

\begin{remark}[Choice of initializer]
\label{rmk:init}
The initializer $ r^0 $ in \Cref{eqn:spec_GVAMP_init} is given by the spectral estimator up to a rescaling factor $c_r$, while $ p^0 $ is motivated by the heuristics in \Cref{sec:heu_init} based on linearized GVAMP. 
The rescaling factor $ c_r $ in \Cref{eqn:spec_GVAMP_init} ensures that $ r^0 \xrightarrow{W_2} \sfR_0 = b_0 \sfB_* + \sqrt{b_0} \sfG $ as $d\to\infty$, with $ \sfG \sim \cN(0,1) $ independent of $\sfB_*$; see \Cref{eqn:R0_P0}. That is, 
 the signal strength and the noise variance are both equal to $b_0$. 
The rescaling factor $c_p$ in $p^0$ has a similar effect. 
These special properties are immaterial for a generic GVAMP algorithm with spectral initialization. 
However, they are crucial for obtaining a succinct state evolution for the spectrally initialized Bayes-GVAMP in the next section. 
\end{remark}

\begin{remark}[Sign calibration]
The spectral initializer $r^0$ in \Cref{eqn:spec_GVAMP_init} contains a sign factor $s$ defined in \Cref{eqn:cr_cp}, which ensures that the sign of the eigenvector is such that $ \inprod{r^0}{\beta_*} \ge 0 $. 
However, $s$ is not computable by the statistician since it involves the unknown $\beta_*$. 
In practice, without the sign calibration, one can run the state evolution twice with $s=1$ and $s=-1$ and one of the trials is guaranteed to track the dynamics of spectrally initialized GVAMP accurately. 
Drawing inspiration from \cite[Appendix H]{Chen_Liu_Ma}, we believe that $s$ is statistically identifiable via e.g.\ moment methods under additional structures (such as asymmetry) of the link $q$ and prior $P_{\sfB_*}$, and non-identifiable otherwise (see also \cite[Remark 3.6]{amp-tutorial}). 
We do not address this issue in more detail but remark that for the purpose of understanding the overlap \Cref{eqn:overlap} which is defined with an absolute value, the sign calibration is inconsequential and the theoretical prediction of asymptotic overlap given by our state evolution result is accurate regardless of the value of $s\in\{-1,1\}$. 
\end{remark}

\subsection{Spectrally initialized Bayes-GVAMP and its state evolution}
\label{sec:BVAMP}

If the link function $q$ and the laws of  $\sfB_*, \sfLambda, \sfE$ in \Cref{asmp:signal,asmp:design,asmp:noise} are known, 
one can choose the matrix and vector denoisers defining GVAMP to maximize the asymptotic overlap between $r^t$ and $\beta_*$. 
We now present this special case of GVAMP, called Bayes-GVAMP, which uses $q, P_{\sfB_*}, \law(\sfLambda), P_{\sfE}$ to define the denoisers used in each iteration.

Define two functions $ g_{x1} \colon \bbR^2 \to \bbR $ and $ g_{z1} \colon \bbR^3 \to \bbR $ as 
\begin{align}
&&
    g_{x1}(r_0; \gamma) &= \frac{ \int x P_{\sfB_*}(x) \exp\paren{ - \frac{\gamma}{2} x^2 + r_0 x } \diff x }{ \int P_{\sfB_*}(x) \exp\paren{ - \frac{\gamma}{2} x^2 + r_0 x } \diff x }
    = \expt{ \sfB_* \mid \sfB_* + \frac{1}{\sqrt{\gamma}} \ol{\sfJ} = \frac{r_0}{\gamma} } , &
& \label{eqn:gx1} \\
&&
    g_{z1}(p_0, y_0; \tau) &= \frac{ \int z Q(y_0\,|\,z) \exp\paren{ -\frac{\tau}{2} \paren{z - \frac{p_0}{\tau}}^2 } \diff z }{ \int Q(y_0\,|\,z) \exp\paren{ -\frac{\tau}{2} \paren{z - \frac{p_0}{\tau}}^2 } \diff z }
    = \expt{ \sfZ \mid \sfP = p_0 , \sfY = y_0 } , & 
& \label{eqn:gz1}
\end{align}
for any $ r_0, \gamma, p_0, y_0, \tau \in \bbR $. 
The expectations in \Cref{eqn:gx1,eqn:gz1} are taken over random variables whose laws are specified as follows. 
In \Cref{eqn:gx1}, $ (\sfB_*, \ol{\sfJ}) \sim P_{\sfB_*} \ot \cN\paren{0, 1} $. 
In \Cref{eqn:gz1}, $ \sfZ = \frac{1}{\tau} \sfP + \frac{1}{\sqrt{\tau}} \ol{\sfK} $, where $ (\sfP, \ol{\sfK}) \sim \cN(0,(\tau - \sigma^{-2}) \tau \sigma^2) \ot \cN(0,1) $. 
Furthermore, $ \sfY \sim Q\paren{ \cdot \mid \sfZ } $. 
We use $ g_{x1}', g_{z1}' $ to denote the derivatives of $ g_{x1}, g_{z1} $ with respect to their first arguments. 
Elementary calculations show that 
\begin{align}
&& 
    g_{x1}'(r_0; \gamma) &= \var{ \sfB_* \mid \sfB_* + \frac{1}{\sqrt{\gamma}} \ol{\sfJ} = \frac{r_0}{\gamma} } , &
    g_{z1}'(p_0, y_0; \tau) &= \var{ \sfZ \mid \sfP = p_0, \sfY = y_0 } . &
& \label{eqn:g'} 
\end{align}
We also define functions $ g_{x2} \colon \bbR^d \times \bbR^n \times \bbR^2 \to \bbR^d $ and $ g_{z2} \colon \bbR^d \times \bbR^n \times \bbR^2 \to \bbR^n $ as
\begin{subequations}
\label{eqn:gx2_gz2}
\begin{align}
    g_{x2}(r,p;\gamma,\tau) &\coloneqq \frac{1}{\gamma} r + X^\top \paren{\frac{\gamma}{\tau} I_n + XX^\top}^{-1} \paren{ \frac{1}{\tau} p - \frac{1}{\gamma} X r } , \\
    g_{z2}(r,p;\gamma,\tau) &\coloneqq X g_{x2}(r,p,\gamma,\tau) , 
\end{align}
for any $ r\in\bbR^d , p \in \bbR^n $ and $ \gamma, \tau \in \bbR $. 
\end{subequations}

The Bayes-GVAMP algorithm maintains a collection of iterates that are updated for every $ t\ge1 $ as follows:
\begin{subequations}
\label{eqn:Bayes_GVAMP}
\begin{align}
&& 
    \wh{x}_1^t &= g_{x1}(r^t; \gamma_1^t) \in \bbR^d , & 
    \wh{z}_1^t &= g_{z1}(p^t, y; \tau_1^t) \in \bbR^n , & 
& \\
&& 
    v_1^t &= \expt{g_{x1}'(\sfR_t; \gamma_1^t)} \in \bbR , &
    c_1^t &= \expt{g_{z1}'(\sfP_t, \sfY; \tau_1^t)} \in \bbR , & 
& \label{eqn:B_v1_c1}\\
&&
    \wt{r}^t &= \frac{1}{v_1^t} \wh{x}_1^t - r^t \in \bbR^d , &
    \wt{p}^t &= \frac{1}{c_1^t} \wh{z}_1^t - p^t \in \bbR^n , & 
& \\
&&
    \gamma_2^t &= \frac{1}{v_1^t} - \gamma_1^t \in \bbR , &
    \tau_2^t &= \frac{1}{c_1^t} - \tau_1^t \in \bbR , & 
& \label{eqn:B_gamma2_tau2} \\
&&
    \wh{x}_2^t &= g_{x2}(\wt{r}^t, \wt{p}^t; \gamma_2^t, \tau_2^t) \in \bbR^d & 
    \wh{z}_2^t &= g_{z2}(\wt{r}^t, \wt{p}^t; \gamma_2^t, \tau_2^t) \in \bbR^n & 
& \\
&&
    v_2^t &= \expt{ \frac{1}{\tau_2^t \sfLambda_d^2 + \gamma_2^t} } \in \bbR , &
    c_2^t &= \expt{ \frac{\sfLambda_n^2}{\tau_2^t \sfLambda_n^2 + \gamma_2^t} } \in \bbR , & 
& \label{eqn:B_v2_c2} \\
&& 
    r^{t+1} &= \frac{1}{v_2^t} \wh{x}_2^t - \wt{r}^t \in\bbR^d , &
    p^{t+1} &= \frac{1}{c_2^t} \wh{z}_2^t - \wt{p}^t \in \bbR^n , & 
& \\
&&
    \gamma_1^{t+1} &= \frac{1}{v_2^t} - \gamma_2^t \in\bbR , &
    \tau_1^{t+1} &= \frac{1}{c_2^t} - \tau_2^t \in\bbR , & 
& \label{eqn:B_gamma1_tau1}
\end{align}
\end{subequations}
where 
\begin{align}
&&
    \sfR_{t} &= \gamma_1^{t} \sfB_* + \sqrt{\gamma_1^{t}} \ol{\sfJ}_{t} , & 
    \sfP_{t} &= \paren{\tau_1^{t} - \frac{1}{\sigma^2}} \sfZ + \sqrt{\tau_1^{t} - \frac{1}{\sigma^2}} \ol{\sfK}_{t} , & 
& \label{eqn:RP_BGVAMP} 
\end{align}
$ \ol{\sfJ}_{t} \sim \cN(0,1) $ is independent of $ \sfB_* $, and $ \ol{\sfK}_{t} \sim \cN(0,1) $ is independent of $ \sfZ, \sfE $. 
The above iteration is initialized with $ r^0, p^0 $ given in \Cref{eqn:spec_GVAMP_init} and 
\begin{align}
&&
    \gamma_1^0 &= \frac{1}{w_2}, & 
    \tau_1^0 &= \frac{1}{w_1} + \frac{1}{\sigma^2} . & 
& \label{eqn:RP0_BGVAMP}
\end{align}

The Bayes-GVAMP iteration \Cref{eqn:Bayes_GVAMP} can be shown to be a special case of the generic GVAMP iteration in \Cref{eqn:GVAMP}. 
Using the definition of $ g_{x2}, g_{z2} $ and applying elementary manipulations to \Cref{eqn:Bayes_GVAMP} bring us to 
\begin{subequations}
\label{eqn:BGVAMP}
\begin{align}
    r^{t+1} &= \brack{ \paren{ \frac{1}{v_2^t \gamma_2^t} - 1 } I_d - \frac{1}{v_2^t \gamma_2^t} X^\top \paren{\frac{\gamma_2^t}{\tau_2^t} I_n + XX^\top}^{-1} X } \paren{ \frac{1}{v_1^t} \wh{x}_1^t - r^t } \notag \\
    &\quad + \frac{1}{v_2^t \tau_2^t} X^\top \paren{\frac{\gamma_2^t}{\tau_2^t} I_n + XX^\top}^{-1} \paren{ \frac{1}{c_1^t} \wh{z}_1^t - p^t } , \\
    p^{t+1} &= \brack{ \frac{1}{c_2^t \tau_2^t} X X^\top \paren{\frac{\gamma_2^t}{\tau_2^t} I_n + XX^\top}^{-1} - I_n } \paren{ \frac{1}{c_1^t} \wh{z}_1^t - p^t } \notag \\
    &\quad + \frac{1}{c_2^t \gamma_2^t} X \brack{ I_d - X^\top \paren{\frac{\gamma_2^t}{\tau_2^t} I_n + XX^\top}^{-1} X } \paren{ \frac{1}{v_1^t} \wh{x}_1^t - r^t } . 
\end{align}
This corresponds to taking the denoisers in \Cref{eqn:GVAMP} to be 
\begin{align}
&& 
    \Phi_t(x) &= \paren{ \frac{1}{v_2^t \gamma_2^t} - 1 } - \frac{1}{v_2^t \gamma_2^t} \cdot \frac{x}{\frac{\gamma_2^t}{\tau_2^t} + x} , & 
    \wt{\Phi}_t(x) &= \frac{1}{v_2^t \tau_2^t} \cdot \frac{x}{\frac{\gamma_2^t}{\tau_2^t} + x^2} , & 
& \notag \\
&&
    &= \frac{1}{v_2^t} \cdot \frac{1}{\gamma_2^t + \tau_2^t x} - 1 , & 
    &= \frac{1}{v_2^t} \cdot \frac{x}{\gamma_2^t + \tau_2^t x^2} , &
& \label{eqn:B_Phi} \\
&&
    \Psi_t(x) &= \frac{1}{c_2^t \tau_2^t} \cdot \frac{x}{\frac{\gamma_2^t}{\tau_2^t} + x} - 1 , & 
    \wt{\Psi}_t(x) &= \frac{1}{c_2^t \gamma_2^t} \cdot x \paren{ 1 - \frac{x^2}{\frac{\gamma_2^t}{\tau_2^t} + x^2} } , & 
& \notag \\
&& 
    &= \frac{1}{c_2^t} \cdot \frac{x}{\gamma_2^t + \tau_2^t x} - 1 , & 
    &= \frac{1}{c_2^t} \cdot \frac{x}{\gamma_2^t + \tau_2^t x^2} , & 
& \label{eqn:B_Psi} \\
&& 
    f_t(r) &= \frac{1}{v_1^t} g_{x1}(r; \gamma_1^t) - r , & 
    g_t(p; y) &= \frac{1}{c_1^t} g_{z1}(p, y; \tau_1^t) - p . & 
& \label{eqn:B_fg}
\end{align}
\end{subequations}
The matrix denoisers $ \Phi_t, \Psi_t $ and the vector denoisers $ f_t,g_t $ above indeed satisfy the trace-free and divergence-free conditions required by \Cref{asmp:tr_free,asmp:div_free}, as shown in \Cref{lem:BGVAMP_free}. 

\begin{theorem}[State evolution of spectrally initialized Bayes-GVAMP]
\label{thm:SE_Bayes_GVAMP}
Consider the Bayes-GVAMP algorithm in \Cref{eqn:Bayes_GVAMP} initialized with \Cref{eqn:spec_GVAMP_init}. 
Then for any fixed $ t\ge0 $, 
\begin{align}
&&
    \matrix{r^{t} & \beta_*} &\xrightarrow{W_2} \matrix{\sfR_{t} & \sfB_*} , & 
    \matrix{p^{t} & z & \eps} &\xrightarrow{W_2} \matrix{\sfP_{t} & \sfZ & \sfE} , &
& \notag 
\end{align}
where the random variables on the right are defined in \Cref{eqn:RP_BGVAMP}. 
\end{theorem}

\begin{figure}[t]
    \centering
    \begin{subfigure}{0.32\textwidth}
        \includegraphics[width=\linewidth]{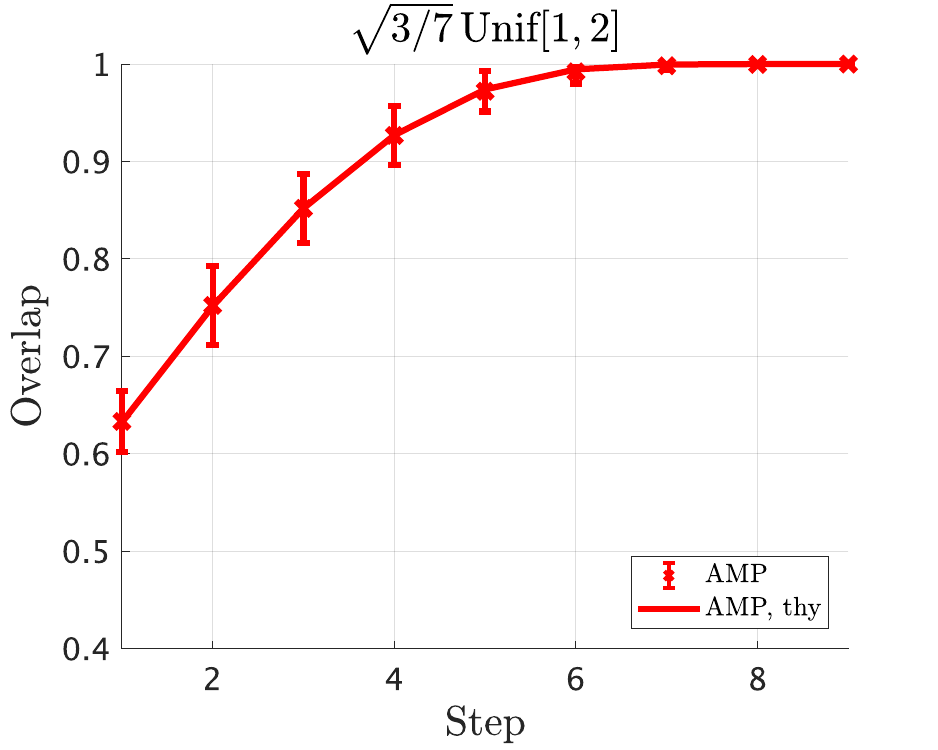}
    \end{subfigure}
    \begin{subfigure}{0.32\textwidth}
        \includegraphics[width=\linewidth]{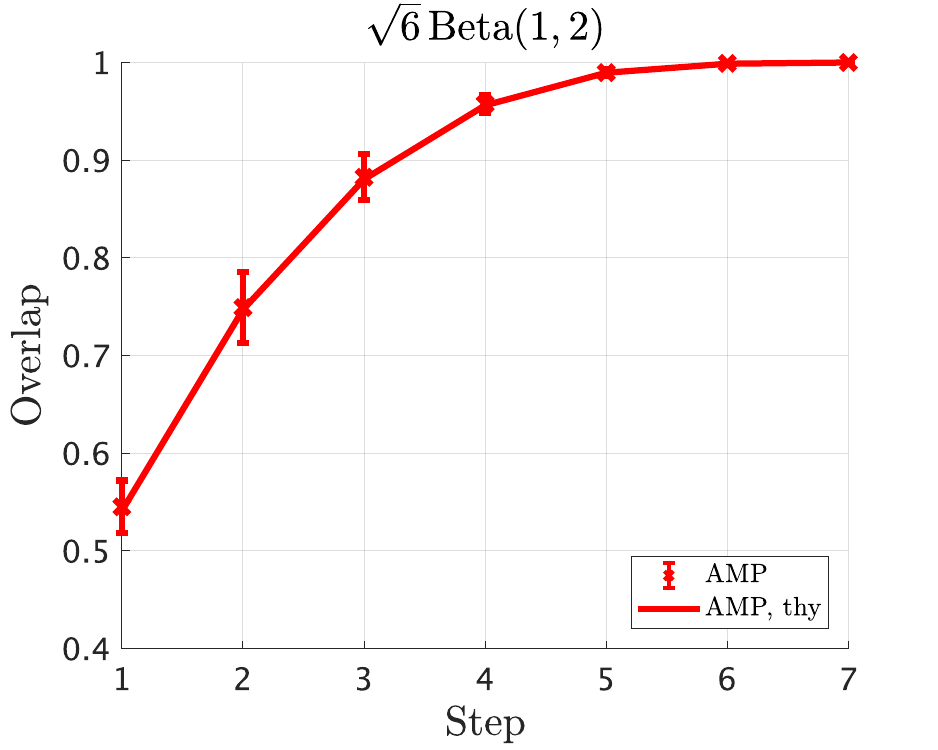}
    \end{subfigure}
    \begin{subfigure}{0.32\textwidth}
        \includegraphics[width=\linewidth]{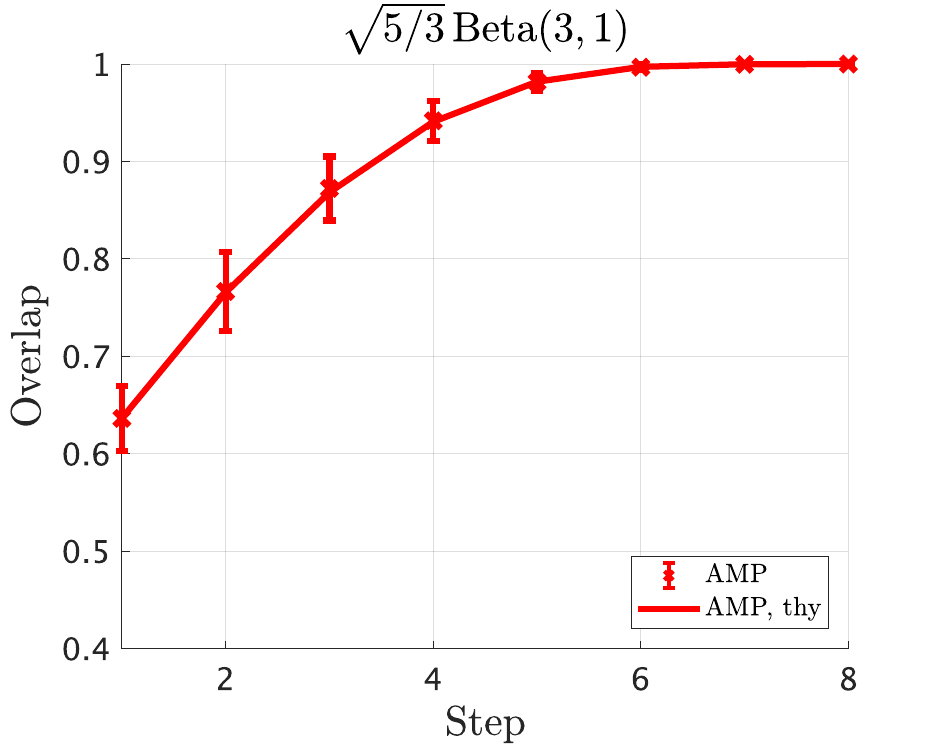}
    \end{subfigure}
    \caption{The trajectory of spectrally initialized Bayes-GVAMP in \Cref{eqn:Bayes_GVAMP} on phase retrieval is accurately tracked by its state evolution. Both trajectories are plotted for $3$ spectral distributions in \Cref{eqn:eg_Lambda} at $\delta = 2$. }
    \label{fig:amp_step}
\end{figure}

\Cref{thm:SE_Bayes_GVAMP} is proved in \Cref{app:pf_thm:SE_Bayes_GVAMP} by inducting on $t$ to simplify the state evolution recursion \Cref{eqn:SE_GVAMP} to \Cref{eqn:Bayes_GVAMP} under the configuration of \Cref{eqn:BGVAMP}. 
\Cref{fig:amp_step} numerically validates that the trajectory of overlaps achieved by spectrally initialized Bayes-GVAMP \Cref{eqn:Bayes_GVAMP} is accurately tracked by its state evolution. 
It is observed that Bayes-GVAMP typically converges rapidly (within $10$ steps) even at  moderate values of $\delta$. 

\subsection{Conjectured Bayes risk}
\label{sec:replica}

In this section, we recall the Bayes risk derived in \cite{Maillard_etal_threshold} using the non-rigorous replica method and show that Bayes-GVAMP has the same fixed points. 
Define two functions $ \cZ_0 \colon \bbR^{2} \to \bbR_{\ge0} $ and $ \cZ_\out \colon \bbR^2 \times \bbR_{\ge0} \to \bbR_{\ge0} $ as
\begin{align}
&& 
    \cZ_0(b,a) &= \int P_{\sfB_*}(x) \exp\paren{ -\frac{1}{2} a x^2 + b x } \diff x , & 
    \cZ_\out(y; \xi, v) &= \expt{ Q\paren{y\mid \sqrt{v} \, \sfOmega + \xi} } , & 
& \label{eqn:Z0_Zout} 
\end{align}
where the expectation is over $ \sfOmega \sim \cN(0,1) $. 
Denote by $ \partial_b \log \cZ_0(b, \cdot) $ and $ \partial_\xi \log \cZ_\out(\cdot; \xi, \cdot) $ the derivatives of the maps $ b \mapsto \log \cZ_0(b, \cdot) $ and $ \xi \mapsto \log \cZ_\out(\cdot; \xi, \cdot) $, respectively. 

Consider the following system of equations in the variables $q_x, q_y, \wh{q}_x, \wh{q}_y,  \gamma_x, \gamma_y$: 
\begin{subequations}
\label{eqn:replica}
\begin{align}
    & q_x = \expt{ \cZ_0\paren{\sqrt{\wh{q}_x}\,\sfXi, \wh{q}_x} \paren{ \partial_b \log \cZ_0\paren{\sqrt{\wh{q}_x}\,\sfXi, \wh{q}_x} }^2 } , \label{eqn:replica1} \\
    & q_z = \frac{1}{\wh{Q}_z + \wh{q}_z} \paren{ \frac{\wh{q}_z}{\wh{Q}_z} + \expt{ \int \cZ_\out\paren{ y; \sqrt{\frac{\wh{q}_z}{\wh{Q}_z (\wh{Q}_z + \wh{q}_z)}}\,\sfXi, \frac{1}{\wh{Q}_z + \wh{q}_z} } \paren{ \partial_\xi \log \cZ_\out\paren{ y; \sqrt{\frac{\wh{q}_z}{\wh{Q}_z (\wh{Q}_z + \wh{q}_z)}}\,\sfXi, \frac{1}{\wh{Q}_z + \wh{q}_z} } }^2 \diff y } } , \label{eqn:replica2} \\
    & \wh{q}_x = \frac{q_x}{\rho(\rho - q_x)} - \gamma_x , \label{eqn:replica3} \\
    & \wh{q}_z = \frac{q_z}{Q_z (Q_z - q_z)} - \gamma_z , \label{eqn:replica4} \\
    & \rho - q_x = \expt{ \frac{1}{\rho^{-1} + \gamma_x + \gamma_z \sfLambda_d^2} } , \label{eqn:replica5} \\
    & Q_z - q_z = \expt{ \frac{\sfLambda_n^2}{\rho^{-1} + \gamma_x + \gamma_z \sfLambda_n^2} } , \label{eqn:replica6}
\end{align}
\end{subequations}
where the expectations in the first two equations are over $ \sfXi \sim \cN(0,1) $ and the constants $ Q_z, \wt{Q}_z $ are defined as
\begin{align}
&& 
    Q_z &= \frac{\rho \expt{ \sfLambda_d^2 }}{\delta} , &
    \wh{Q}_z &= \frac{1}{Q_z} . & 
& \label{eqn:Qz} 
\end{align}
We also recall that $ \rho $ is the variance of $ P_{\sfB_*} $ (\Cref{asmp:signal}).
The system of equations \Cref{eqn:replica} is referred to as the replica saddle point equations in the literature \cite{Maillard_etal_threshold}. 

Define $ \wh{\beta}_{\textnormal{B}} \equiv \wh{\beta}_{\textnormal{B}}(X,y) \in \bbR^d $ (with subscript $\textnormal{B}$ for `Bayesian') as the posterior mean of $ \beta_* $: 
\begin{align}
    \wh{\beta}_{\textnormal{B}} &= \expt{\beta_* \mid X,y} = \frac{\int_{\bbR^d} \beta P_{\sfB_*}^{\ot d}(\beta) Q^{\ot d}\paren{y \mid X \beta} \diff \beta}{\int_{\bbR^d} P_{\sfB_*}^{\ot d}(\beta) Q^{\ot d}\paren{y \mid X \beta} \diff \beta} . \label{eqn:beta_B} 
\end{align}
It is well-known that $ \wh{\beta}_{\textnormal{B}} $ minimizes the mean squared error for estimating $ \beta_* $ given $(X,y)$ from the generalized linear model \Cref{eqn:model_Q}. 

\begin{conjecture}[{\cite[Conjecture 2.1]{Maillard_etal_threshold}}]
\label{conj:replica}
There exists a tuple 
\begin{align}
    & (\wh{q}_x^*, \wh{q}_z^*, \gamma_x^*, \gamma_z^*, q_x^*, q_z^*) 
    \in [0, \infty) \times [0, \infty) \times [0, \infty) \times [0, \infty) \times [0, \rho] \times [0, Q_z] \notag 
\end{align}
solving \Cref{eqn:replica} such that 
\begin{align}
    \lim_{d\to\infty} \frac{1}{d} \expt{\normtwo{\beta_* - \wh{\beta}_{\textnormal{B}}}^2} &= \rho - q_x^* , \notag 
\end{align}
where the expectation is over $ \beta_*, X, y $. 
\end{conjecture}

When the solution to \Cref{eqn:replica} is not unique, $ q_x^* $ is obtained by maximizing a certain free entropy functional. 
We do not define this functional formally due to its complexity, but instead refer the reader to \cite[Conjecture 2.1]{Maillard_etal_threshold}. 
In what follows, we verify that the set of solutions to \Cref{eqn:replica} is in one-to-one correspondence with the set of fixed points of the state evolution recursion for Bayes-GVAMP in \Cref{eqn:Bayes_GVAMP}. 

\begin{theorem}
\label{thm:replica}
The set of fixed points
\begin{align}
    & (\gamma_1, \tau_1, \gamma_2, \tau_2, v_1, v_2, c_1, c_2) \in 
    [0, \infty) \times \left[\frac{1}{\sigma^2}, \infty\right) \times \left[\frac{1}{\rho}, \infty\right) \times [0, \infty) \times [0,\rho] \times [0,\rho] \times [0,\sigma^2] \times [0,\sigma^2] \label{eqn:SE_FP}
\end{align}
of the recursion for $ (\gamma_1^t, \tau_1^t, \gamma_2^t, \tau_2^t, v_1^t, v_2^t, c_1^t, c_2^t) $ in \Cref{eqn:Bayes_GVAMP} is the same as the set of tuples
\begin{align}
    & \paren{\wh{q}_x, \wh{q}_z + \frac{1}{\sigma^2}, \gamma_x + \frac{1}{\rho}, \gamma_z, \rho - q_x, \rho - q_x, \sigma^2 - q_z, \sigma^2 - q_z} 
    \label{eqn:change_var} 
\end{align}
such that
\begin{align}
    & (\wh{q}_x, \wh{q}_z, \gamma_x, \gamma_z, q_x, q_z) 
    \in [0, \infty) \times [0, \infty) \times [0, \infty) \times [0, \infty) \times [0, \rho] \times [0, Q_z] \notag 
\end{align}
are the solutions to the replica saddle point equations \Cref{eqn:replica}. 
\end{theorem}

\Cref{thm:replica} is proved in \Cref{app:pf_thm:replica} by verifying the equivalence between \Cref{eqn:Bayes_GVAMP,eqn:replica} under the change of variables in \Cref{eqn:change_var}. 

\begin{figure}[t]
    \centering
    \begin{subfigure}{0.32\textwidth}
        \includegraphics[width=\linewidth]{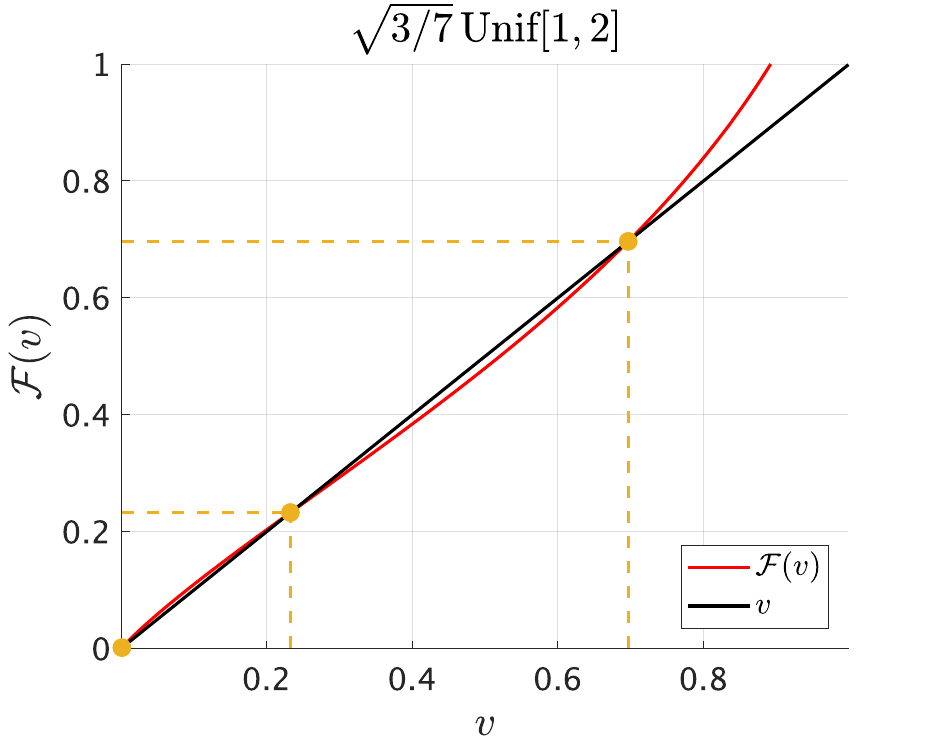}
    \end{subfigure}
    \begin{subfigure}{0.32\textwidth}
        \includegraphics[width=\linewidth]{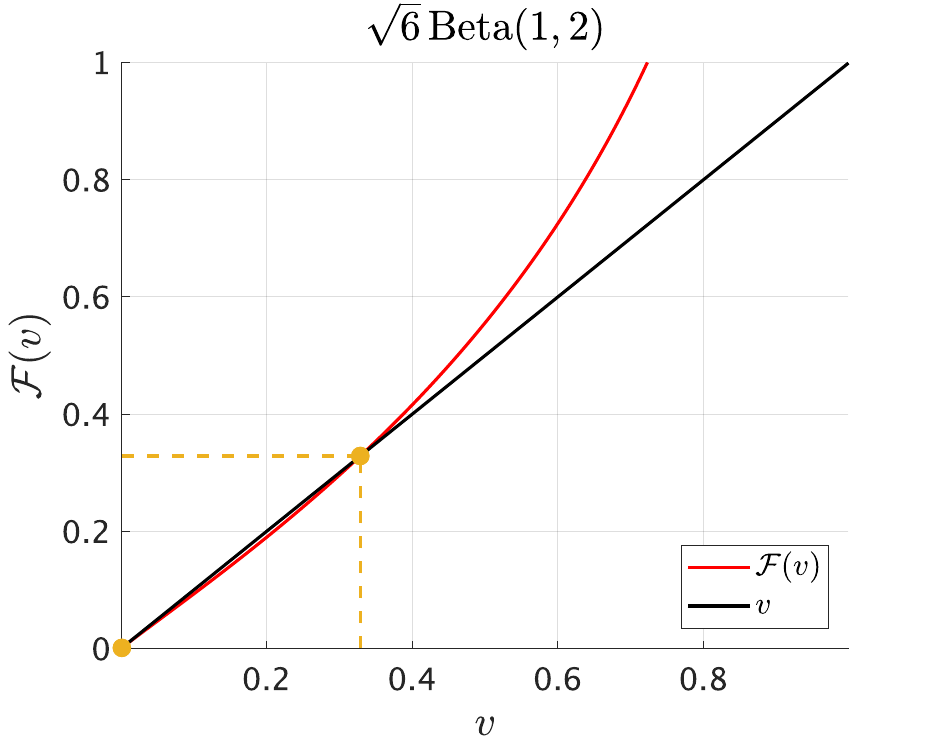}
    \end{subfigure}
    \begin{subfigure}{0.32\textwidth}
        \includegraphics[width=\linewidth]{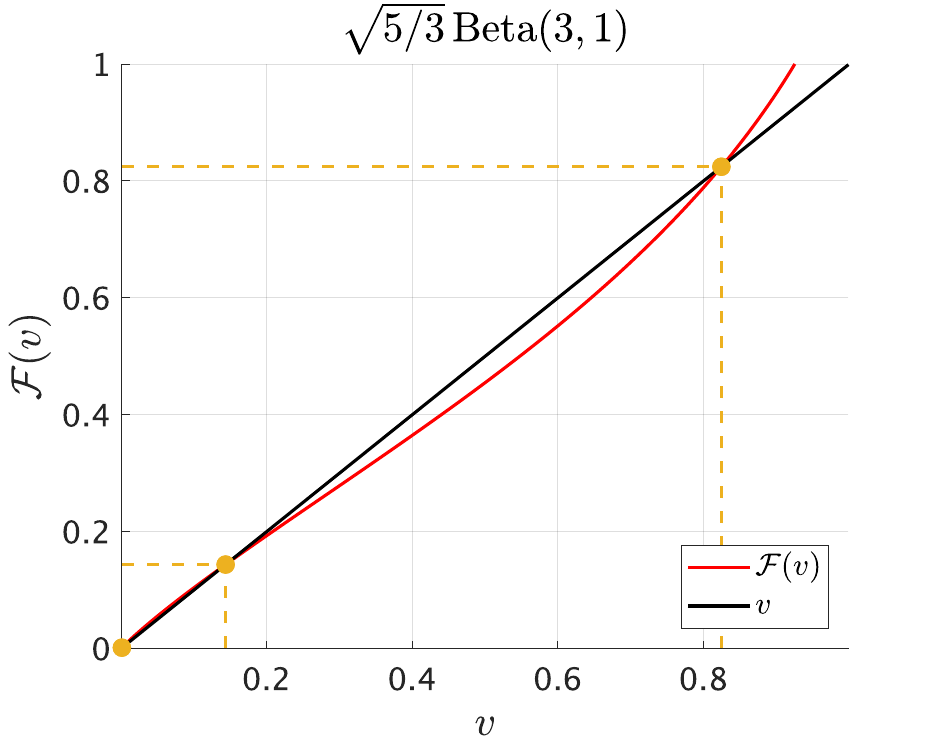}
    \end{subfigure}
    \caption{Plots of the function $ \cF $ and its fixed points, i.e.\ solution $v\in[0,\rho]$ to $ v = \cF(v) $, for $3$ spectral distributions \Cref{eqn:eg_Lambda}.  
    From the leftmost to the rightmost panel, $ \delta $ is taken to be $ 1.44 $, $ 1 $, $ 1.34 $, respectively.}
    \label{fig:replica}
\end{figure}

The tuple of fixed points in \Cref{thm:replica} may not be unique. 
Due to the equivalence between \Cref{eqn:SE_FP,eqn:change_var}, let us only examine \Cref{eqn:SE_FP}. 
For instance, consider the (noiseless) phase retrieval model $ q(z,\eps) = \abs{z} $ with Gaussian prior $ P_{\sfB_*} = \cN(0,1) $. 
Specializing \Cref{conj:replica} to this setting, the asymptotic Bayes risk is conjectured to be $ v \in [0,1] $ which solves the equation $ \cF(v) = v $ for a certain function $\cF$; 
see \Cref{sec:replica_eg} for the precise definition of $\cF$ and a derivation of this result. 
\Cref{fig:replica} shows the non-uniqueness of the solution $v$ under $3$ different spectral distributions of $ \sfLambda $ (see \Cref{sec:phase_synthetic} for their definitions) and certain values of $\delta$. 

\begin{figure}[t]
    \centering
    \begin{subfigure}{0.32\textwidth}
        \includegraphics[width=\linewidth]{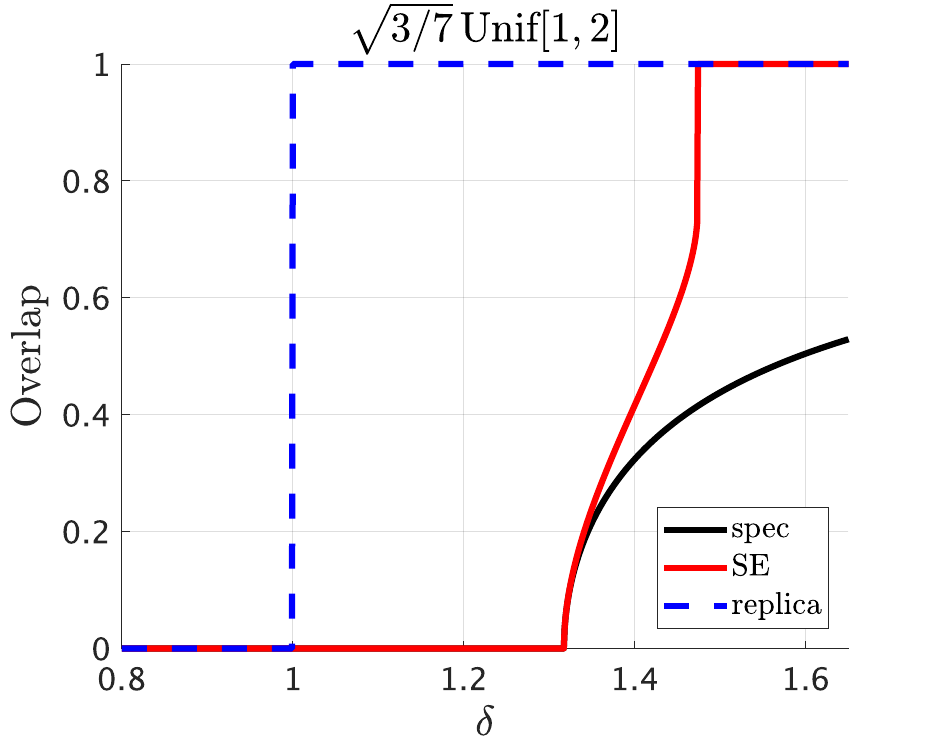}
    \end{subfigure}
    \begin{subfigure}{0.32\textwidth}
        \includegraphics[width=\linewidth]{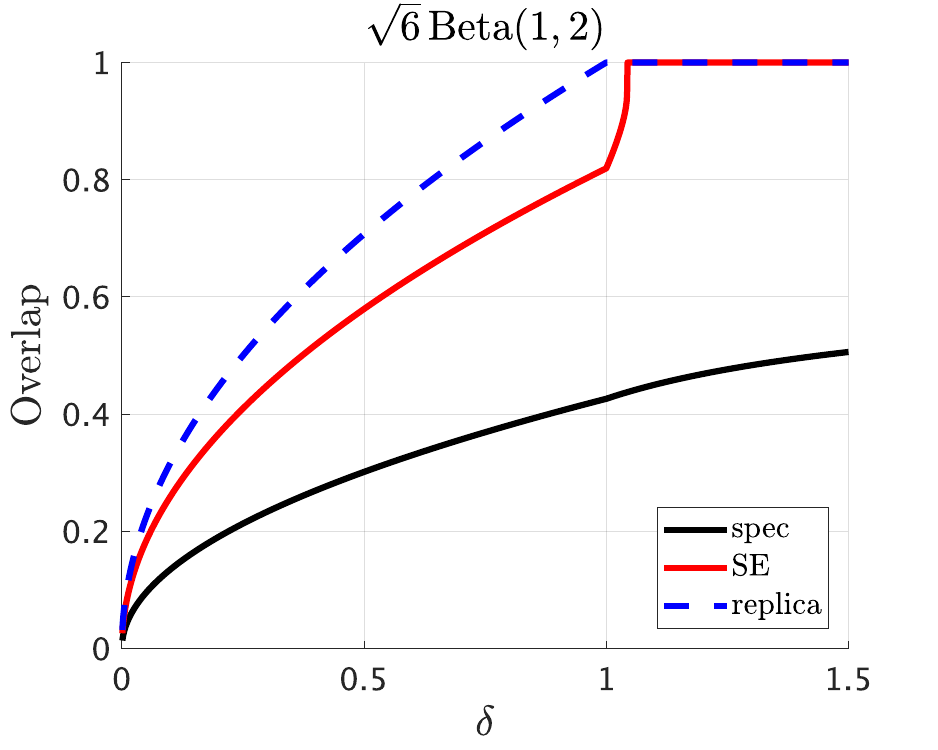}
    \end{subfigure}
    \begin{subfigure}{0.32\textwidth}
        \includegraphics[width=\linewidth]{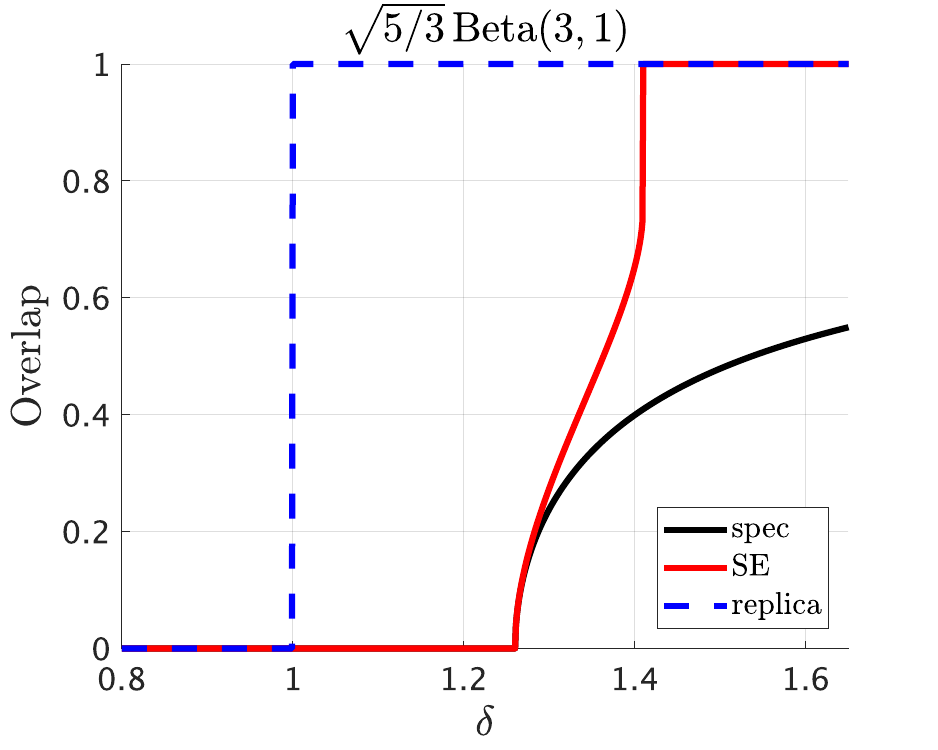}
    \end{subfigure}
    \caption{Asymptotic overlap \Cref{eqn:eta} achieved by the spectral estimator in \Cref{thm:opt_thr} (`spec' in black), overlap achieved by the state evolution of spectrally initialized Bayes-GVAMP in \Cref{eqn:Bayes_GVAMP} run till convergence (`SE' in red), and the conjectured Bayes risk \Cref{eqn:replica} expressed in terms of overlap (`replica' in blue). All curves are plotted for $3$ spectral distributions given in \Cref{eqn:eg_Lambda}. }
    \label{fig:gap}
\end{figure}

Under the same setting above, we also plot in \Cref{fig:gap} the overlap achieved by (the state evolution of) spectrally initialized Bayes-GVAMP at convergence, the conjectured Bayes-optimal overlap $ \sqrt{1 - v/\rho} $ (which can be seen by applying the tower property of conditional expectation to \Cref{conj:replica}), and for comparison, the asymptotic overlap of spectral estimator given by \Cref{thm:opt_thr}. 
It is observed that spectrally initialized Bayes-GVAMP in \Cref{eqn:Bayes_GVAMP} achieves the conjectured Bayes-optimal overlap when the solution $v$ is unique, and otherwise achieves an overlap corresponding to a worse solution $v$. 

For an i.i.d.\ Gaussian design which is unstructured, it has been well observed \cite{Barbier_etal} that multiple state evolution fixed points and the consequent performance gap between the optimal AMP algorithm \cite{Rangan,Mondelli_Venkataramanan} and the Bayes-optimal posterior mean \Cref{eqn:beta_B} may emerge in the presence of structured \emph{prior}. 
Here, \Cref{fig:replica,fig:gap} are exclusively under the unstructured Gaussian prior. 
Our finding therefore complements the above observation by providing another source for multiple fixed points and performance gaps, which is a structured \emph{design}. 
We note that such performance gaps can be drastic. 
Indeed, in both the leftmost and rightmost panels of \Cref{fig:gap}, there exists a range of $\delta$ in which the posterior mean achieves perfect recovery whereas Bayes-GVAMP is ineffective.

\section{Heuristics and proof overview}
\label{sec:heuristics}

In \Cref{sec:heu_lin}, we draw heuristic connections between linearized GVAMP and spectral estimators, thereby motivating spectral estimators of the form \Cref{eqn:spec_est} subject to \Cref{asmp:preprocess}. 
In \Cref{sec:heu_D}, we highlight the key difficulties in the spectral analysis of the random matrix $D$ and provide a proof overview of \Cref{thm:spec}. 
\Cref{sec:heu_GVAMP} motivates our generic GVAMP formalism. 
In \Cref{sec:heu_init}, we offer additional intuition on the design of the spectral initialization in  \Cref{eqn:spec_GVAMP_init} for GVAMP. 

\subsection{Simulating spectral estimators using linearized GVAMP}
\label{sec:heu_lin}

Let us consider the following GVAMP algorithm with linear denoising functions: 
\begin{align}
&&
    u^{t} &= \paren{ \ol{\kappa}_2^{-1} X X^\top - I_n } g(u^{t-1}) , & 
    v^{t} &= \ol{\kappa}_2^{-1} X^\top g(u^{t-1}) , & 
& \label{eqn:recur}
\end{align}
where 
\begin{align}
g(u) &= \paren{ \frac{G}{\tr(G)/n} - I_n } u, \notag 
\end{align}
and $ G = \diag(\wh{g}(y)) $ is a diagonal matrix obtained by applying an arbitrary function $ \wh{g}\colon \bbR \to \bbR $ to $y$ component-wise. 
As shown in \Cref{app:pf_prop:SE}, \Cref{eqn:recur} is a special case of the GVAMP algorithm in \Cref{eqn:GVAMP} and also of the RI-GAMP algorithm in \cite{VKM}. 
By \cite[Equation (16)]{Maillard_etal_construction}, with a specific choice of $ G $, the iteration \Cref{eqn:recur} arises as the linearization of the Bayes-GVAMP algorithm \cite{Schniter_Rangan_Fletcher} (see also \Cref{eqn:Bayes_GVAMP}) around its trivial fixed point. 
In the context of generalized linear models with partial Haar design, a similar iteration is derived in \cite[Equation (21a)]{Ma_etal} from expectation propagation and is used to heuristically analyze spectral estimators. 

Note that the iteration of $ u^t $ above operates on a random matrix with vanishing trace: 
\begin{align}
\lim_{n\to\infty} \frac{1}{n} \tr\paren{ \ol{\kappa}_2^{-1} X X^\top - I_n } &= 0 . \notag 
\end{align}
Furthermore, the denoising function $g$ has zero divergence since 
\begin{align}
\tr\paren{ \frac{G}{\tr(G)/n} - I_n } &= 0 . \notag 
\end{align}
For large $t$, we expect that $ u^t / \normtwo{u^t} $ converges in the following sense: 
\begin{align}
\lim_{t\to\infty} \lim_{n\to\infty} \frac{\abs{\inprod{u^t}{u^{t-1}}}}{\normtwo{u^t} \normtwo{u^{t-1}}} &= 1 . \label{eqn:heu_conv1} 
\end{align}
In other words, assuming that the limit
\begin{align}
    \gamma &= \lim_{t\to\infty} \lim_{d\to\infty} \frac{\inprod{u^t}{u^{t-1}}}{\normtwo{u^{t-1}}^2} \label{eqn:gamma_heu}
\end{align}
exists, then we expect that two adjacent iterates are collinear for large $t$, i.e., 
\begin{align}
    \lim_{t\to\infty} \lim_{n \to\infty} \frac{\normtwo{u^t - \gamma u^{t-1}}}{\normtwo{u^{t-1}}} &= 0 . \label{eqn:heu_conv2}
\end{align}
Note that $ \gamma \equiv \gamma(\wh{g}) $ depends on the function $\wh{g}$. 
We caution that $\gamma = 1$ need not hold, in which case \Cref{eqn:heu_conv1} or \Cref{eqn:heu_conv2} does not imply $ \lim_{t\to\infty} \lim_{n\to\infty} n^{-1} \normtwo{u^t - u^{t-1}}^2 = 0 $. 
See the discussion following \Cref{lem:SE_stay} for more details. 
The fixed point of the iteration \Cref{eqn:recur} can then be written as: 
\begin{align}
\gamma \wh{u} &= \paren{ \ol{\kappa}_2^{-1} X X^\top - I_n } \paren{ \frac{G}{\tr(G)/n} - I_n } \wh{u} \label{eqn:whu} 
\end{align}
for some unit vector $ \wh{u} \in \bbS^{n-1} $. 
This equation is closely related to an eigenequation of $D$ of the form \Cref{eqn:D} (with respect to some $\cT$ to be specified). 
Indeed, rearranging terms, we get 
\begin{align}
\brack{ (\gamma - 1) I_n + \frac{G}{\tr(G)/n} } \wh{u} &= \paren{ \ol{\kappa}_2^{-1} X X^\top } \paren{ \frac{G}{\tr(G)/n} - I_n } \wh{u} . \notag 
\end{align}
Left-multiplying both sides by 
\begin{align}
    X^\top \paren{ \frac{G}{\tr(G)/n} - I_n } \brack{ (\gamma - 1) I_n + \frac{G}{\tr(G)/n} }^{-1} , \notag 
\end{align}
we get
\begin{align}
X^\top \paren{ \frac{G}{\tr(G)/n} - I_n } \wh{u}
&= \ol{\kappa}_2^{-1} X^\top \paren{ \frac{G}{\tr(G)/n} - I_n } \brack{ (\gamma - 1) I_n + \frac{G}{\tr(G)/n} }^{-1} X X^\top \paren{ \frac{G}{\tr(G)/n} - I_n } \wh{u} . \notag 
\end{align}
We recognize that this is an eigenequation of the matrix 
\begin{align}
X^\top \paren{ \frac{G}{\tr(G)/n} - I_n } \brack{ (\gamma - 1) I_n + \frac{G}{\tr(G)/n} }^{-1} X \label{eqn:DD} 
\end{align}
corresponding to the eigenvector (up to rescaling)
\begin{align}
X^\top \paren{ \frac{G}{\tr(G)/n} - I_n } \wh{u} . \label{eqn:eigvec_heu} 
\end{align}
The matrix in \Cref{eqn:DD} is asymptotically equivalent to $D$ in \Cref{eqn:D} with $ \cT $ taking the form of \Cref{eqn:T} and $ \ol{g} $ therein given by 
\begin{align}
    \ol{g}(y) &= \frac{\wh{g}(y)}{\expt{\wh{g}(\sfY)}} - 1 , \notag 
\end{align}
which clearly satisfies the condition $ \expt{\ol{g}(\sfY)} = 0 $ required by \Cref{eqn:gbar}. 
We remark that the value of $ \gamma $ in \Cref{eqn:gamma_heu} can be derived using state evolution (see \Cref{sec:lin_VAMP_SE}, or more specifically \Cref{lem:technical}), leading to the explicit expression in \Cref{eqn:T}. 
This motivates the spectral estimator \Cref{eqn:spec_est} with respect to preprocessing functions subject to \Cref{asmp:preprocess}. 

\subsection{Analysis of $D$ and proof overview of \Cref{thm:spec}}
\label{sec:heu_D}

Our analysis of the matrix $D$ in \Cref{eqn:D} deviates from the standard random matrix theory approach. 
Recalling \Cref{asmp:design}, we can write this matrix as $ D = Q \Lambda^\top O^\top T O \Lambda Q^\top $ which formally resembles a spiked multiplicative model \cite{Ding_Ji}. 
There are, however, three important distinctions. 
In a standard spiked multiplicative model, it is typically assumed that \emph{(a)} $\Lambda,T$ are independent of the Haar matrices $O,Q$, \emph{(b)} $T$ is positive semidefinite which brings the convenient property that the spectra of $ Q \Lambda^\top O^\top T O \Lambda Q^\top $ and $ T^{1/2} O \Lambda \Lambda^\top O^\top T^{1/2} $ differ only at $0$, and \emph{(c)} $\Lambda,T$ contain spectral outliers at fixed locations. 
The last assumption can be equivalently understood as requiring the low rank spikes on $\Lambda,T$ to be aligned with the eigenvectors corresponding to the given outlying eigenvalues. 
In our setting, none of \emph{(a)}, \emph{(b)}, \emph{(c)} hold: $T$ involves the Haar randomness from $O,Q$ and can take negative values (see \Cref{thm:opt_thr}); the spike enters $T$ through a $1$-dimensional projection $X\beta_*$ along the unknown parameter vector $\beta_*$. 
Existing results on spiked multiplicative models are therefore inapplicable. 
Instead, we exploit and make rigorous the connections in \Cref{sec:heu_lin} between spectral estimators and linearized GVAMP by combining tools from AMP, random matrix theory and free probability. 

The precise characterization of asymptotic spectral properties of $D$ in \Cref{thm:spec} is accomplished in three steps. 
First, motivated by the linearization heuristics in \Cref{sec:heuristics}, we consider an instance \Cref{eqn:recur} of GVAMP in \Cref{eqn:GVAMP}, derive its state evolution and study its fixed points. 
These results are summarized in \Cref{lem:SE_stay} (proved in \Cref{sec:pf_spec_step1}). 
Second, using tools from random matrix theory and free probability, we determine the limit of the second (non-outlying) eigenvalue of $D$; see \Cref{lem:edge} (proved in \Cref{app:right_edge}). 
Taking the above two ingredients as input, the crux of the argument lies in \Cref{lem:align} (proved in \Cref{app:pf_lem:align}) which shows that the linearized GVAMP iterate is asymptotically aligned with $ v_1(D) $ of interest and the corresponding outlying eigenvalue $ \lambda_1(D) $ can be accurately located. 
Given this, the asymptotic overlap achieved by the spectral estimator is computed by invoking the state evolution result in the first step. 

We discuss in more detail the proof strategy for characterizing the right edge of the bulk of $D$ (\Cref{lem:edge}) as it is of independent interest in  random matrix theory and free probability. 
This part is in turn divided into three steps. 
\Cref{sec:edge1} sets the stage for the study of a `null'  counterpart of $D$ for which some assumptions in \Cref{sec:model,sec:spec} can be relaxed. 
Informally, $X$ and $T$ in the null matrix are independent. 
\Cref{sec:edge2} shows that the top $k$ eigenvalues (for any fixed $k\ge1$ relative to $n,d$) of the null matrix all converge to the supremum of the support of its limiting spectral distribution, a type of result known as strong convergence \cite{vanHandel}. 
By eigenvalue interlacing (\Cref{lem:conditioning}), this immediately implies that the second largest eigenvalue of $D$ has the same limit. 
A key component in this step in fact establishes a more general statement regarding the strong asymptotic freeness of Haar orthogonal and deterministic matrices, which implies strong convergence of matrix products involving such matrices. 
This is done by a reduction using functional calculus to known strong asymptotic freeness of GOEs and deterministic matrices \cite[Theorem 4.2]{Fan_Sun_Wang}. 
The limiting spectral distribution is expressed in terms of free multiplicative convolution, a characterization of whose right edge does not appear to be available. 
In \Cref{sec:edge3}, 
using tools from complex analysis and free probability, we develop such an explicit characterization for the right edge of the free convolution of more general infinite-dimensional operators in terms of the largest critical point of a certain functional. 
This functional can be loosely understood as the inverse Stieltjes transform of the free convolution. 
Finally, \Cref{sec:edge4} specializes this explicit characterization to the setting of \Cref{sec:model,sec:spec} under which further simplifications can be made, eventually leading to the characterization in \Cref{lem:edge}. 

\subsection{The generic GVAMP formalism and proof overview of \Cref{thm:spec_GVAMP}}
\label{sec:heu_GVAMP}

Our generic GVAMP formalism in  \Cref{eqn:GVAMP} draws inspiration from multiple sources (see \Cref{sec:related_work}) with the closest one in spirit being \cite{Dudeja_Sen_Lu}. 
In settings where the data matrix $ X \in \bbR^{d\times d} $ is symmetric, an abstract form of the VAMP iteration $ r^{t+1} = \Phi_t(X) f_t(r^t) $ is introduced in \cite{Dudeja_Sen_Lu}, where $ \Phi_t $ is trace-free and $f_t$ is divergence-free in the same senses as \Cref{asmp:tr_free,asmp:div_free} respectively. 
These two conditions are the key to a tractable state evolution without loss of optimality.  
In settings such as GLMs where $X\in\bbR^{n\times d}$ is typically rectangular, our GVAMP formalism can be understood as operating on an augmented block matrix $Y_t\in\bbR^{(n+d)\times(n+d)}$: 
\begin{align}
    \matrix{r^{t+1} \\ p^{t+1}}
    &= \underbrace{\matrix{\Phi_t(X^\top X) & \wt{\Phi}_t(X)^\top \\ \wt{\Psi}_t(X) & \Psi_t(XX^\top)}}_{Y_t} \matrix{f_t(r^t) \\ g_t(p^t)} , \label{eqn:rp_block}
\end{align}
while still imposing divergence-free conditions on the vector denoisers, but only imposing trace-free condition on the diagonal blocks of $Y_t$. 
Note that \Cref{eqn:rp_block} cannot be reduced to the symmetric version studied in \cite{Dudeja_Sen_Lu} for two reasons. 
First, $Y_t$ is not invariant under conjugation by $ \bbO(n+d) $, but is instead block-invariant under conjugation by $ \bbO(n) \times \bbO(d) $. 
Second, $Y_t$ is not symmetric since $ \wt{\Phi}_t, \wt{\Psi}_t $ may differ. 
This flexibility is particularly useful in applications to GLMs where the optimal choices of $ \wt{\Phi}_t, \wt{\Psi}_t $ are indeed different. 
Informally, this is because the statistical model \Cref{eqn:model} only reveals $y$ but not $ \beta_*$, which eventually breaks the symmetry between $n$- and $d$-dimensional iterates in \Cref{eqn:rp_block}. 
The proof of the state evolution result (\Cref{thm:spec_GVAMP}) for GVAMP follows the iterative Haar conditioning strategy of \cite{fan2020approximate} with careful treatment of the distinguished vectors $ \beta_* $ and $y$ (or more generally $z,\eps$) which $ f_t $ and $ g_t $ can take as inputs (though omitted from the notation in \Cref{eqn:rp_block}). 
The formalism \Cref{eqn:rp_block} can be easily seen equivalent to \Cref{eqn:GVAMP} and is general enough to include the Bayes-GVAMP \Cref{eqn:Bayes_GVAMP} first derived heuristically in \cite{Schniter_Rangan_Fletcher}. 
The general theory developed here allows us to rigorously track the dynamics of GVAMP with spectral initialization. 

\subsection{Design of spectral initialization}
\label{sec:heu_init}

To properly initialize GVAMP with spectral estimators, the $d$- and $n$-dimensional initializers $r^0,p^0$ need to respect a certain algebraic relation to ensure the correctness of state evolution. 
Consider again the linearized GVAMP algorithm in \Cref{eqn:recur} with some $ \wh{g} $. 
Let $ \ol{G} = \frac{G}{\tr(G)/n} - I_n $. 
Denote the fixed point of this iteration by the same notation with time indices removed. 
Let $ \wh{u}, \gamma $ be as in \Cref{eqn:gamma_heu,eqn:whu}. 
According to \Cref{eqn:eigvec_heu}, $ v_1(D) $ is aligned with $ \wh{v} = \ol{\kappa}_2^{-1} X^\top \ol{G} \wh{u} $ which is in turn aligned with $ v^t $ for large $d,t$. 
Then by the update rule of \Cref{eqn:recur}, 
\begin{align}
    \gamma \wh{u} &= (\ol{\kappa}_2^{-1} XX^\top - I) \ol{G} \wh{u}
    = X \wh{v} - \ol{G} \wh{u} . \notag 
\end{align}
Solving for $ \gamma \wh{u} $ gives
\begin{align}
    \gamma \wh{u} &= \gamma \paren{\gamma I_n + \ol{G}}^{-1} X \wh{v} , \notag 
\end{align}
thereby motivating the choice of $ c_p^{-1} p^0 $ in \Cref{eqn:spec_GVAMP_init}. 

\section{Examples and numerical experiments}
\label{sec:experiments}

We present numerical simulations of our proposed algorithms (as well as some existing methods from the literature) and empirical validate the corresponding theory on the phase retrieval and Poisson regression models, with both synthetic and real data. 
Some implementation details are deferred to \Cref{sec:details}. 

\begin{figure}[t]
    \centering
    \begin{subfigure}{0.32\textwidth}
        \includegraphics[width=\linewidth]{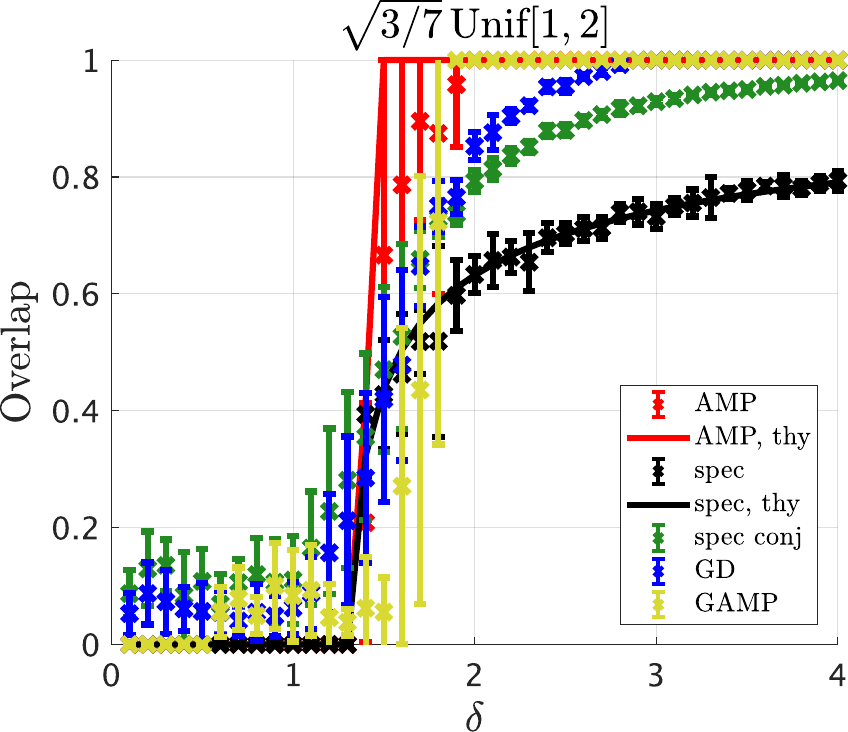}
    \end{subfigure}
    \begin{subfigure}{0.32\textwidth}
        \includegraphics[width=\linewidth]{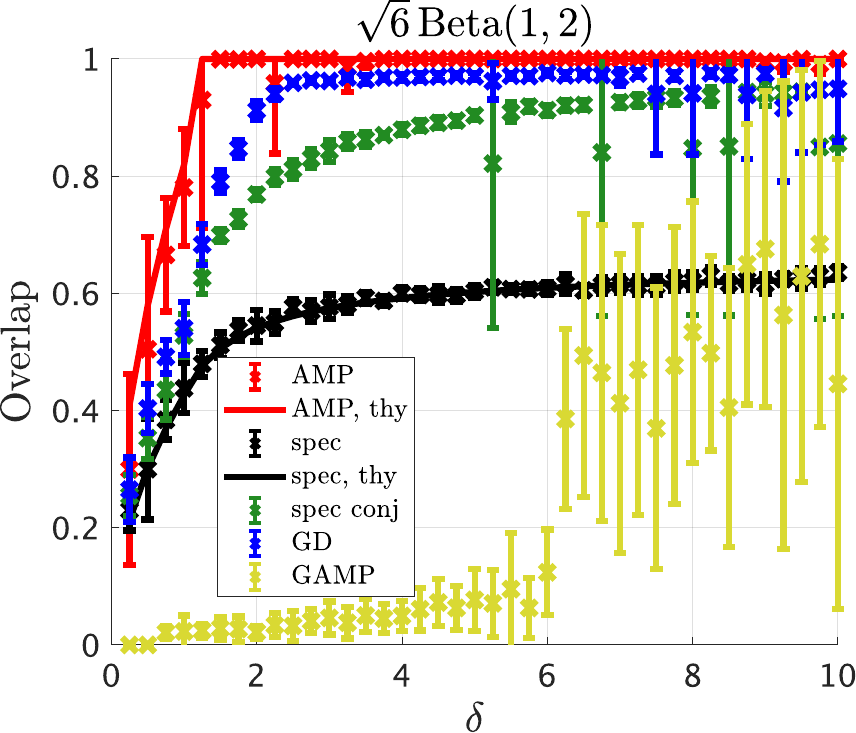}
    \end{subfigure}
    \begin{subfigure}{0.32\textwidth}
        \includegraphics[width=\linewidth]{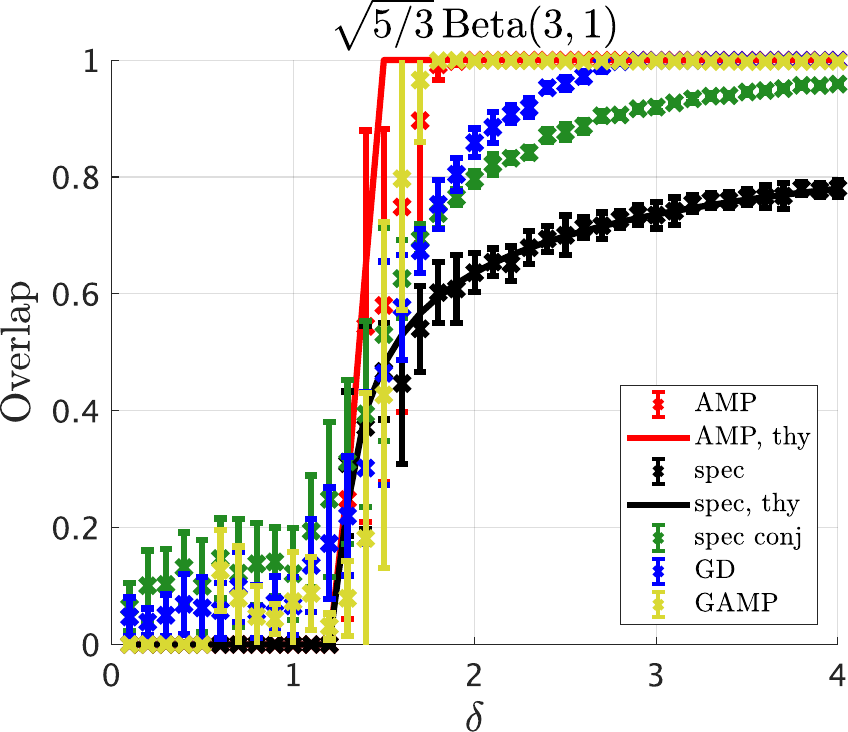}
    \end{subfigure}
    \caption{Overlaps achieved by $5$ algorithms on the phase retrieval model \Cref{eqn:phase_retrieval}. The same experiments are repeated for $3$ spectral distributions given in \Cref{eqn:eg_Lambda} of the design matrix $X$. }
    \label{fig:amp_pr}
\end{figure}

\begin{figure}[t]
    \centering
    \begin{subfigure}{0.32\textwidth}
        \includegraphics[width=\linewidth]{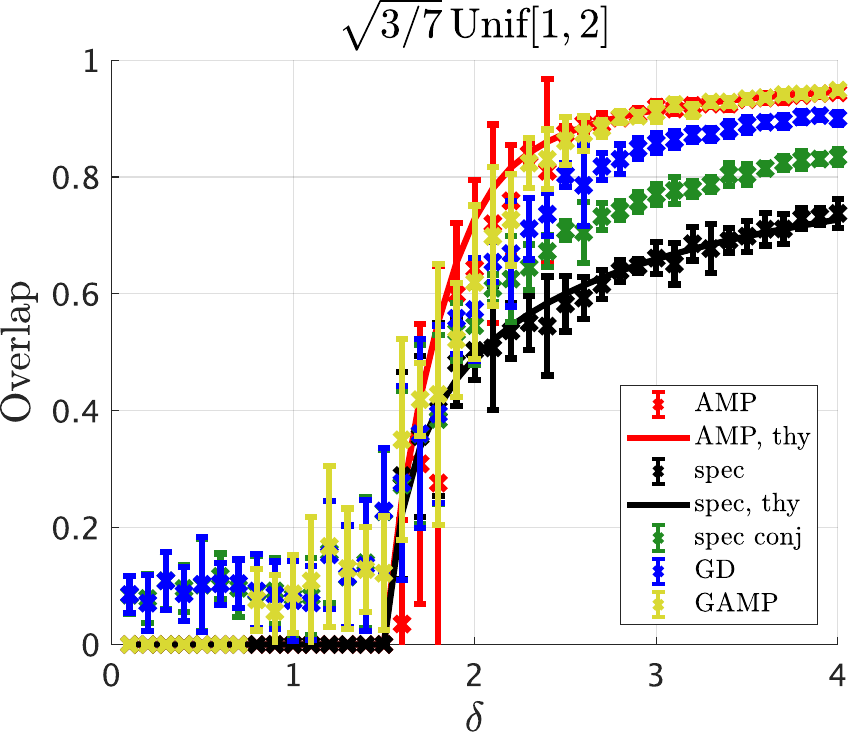}
    \end{subfigure}
    \begin{subfigure}{0.32\textwidth}
        \includegraphics[width=\linewidth]{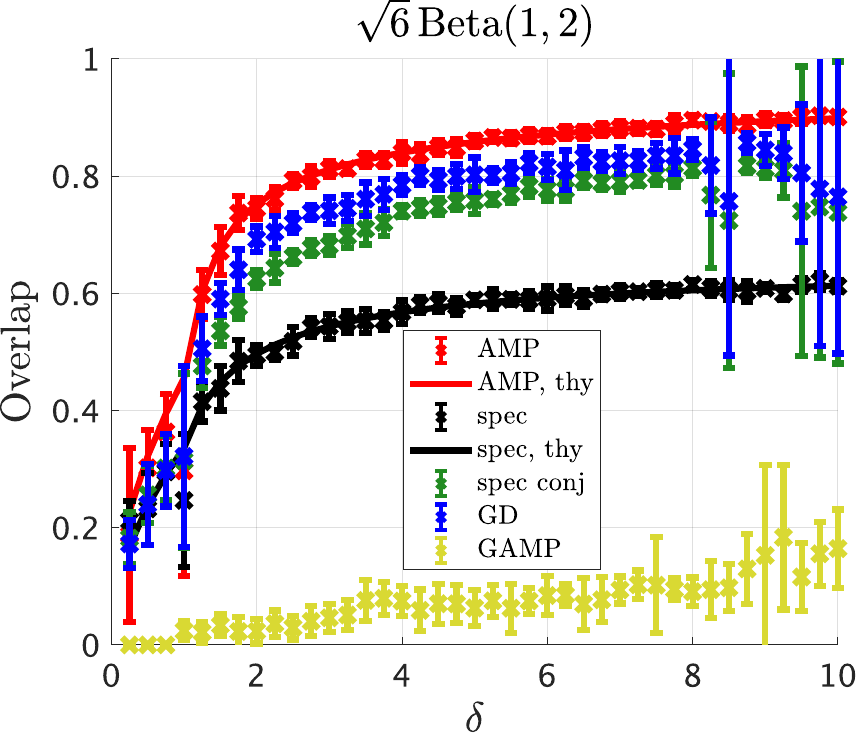}
    \end{subfigure}
    \begin{subfigure}{0.32\textwidth}
        \includegraphics[width=\linewidth]{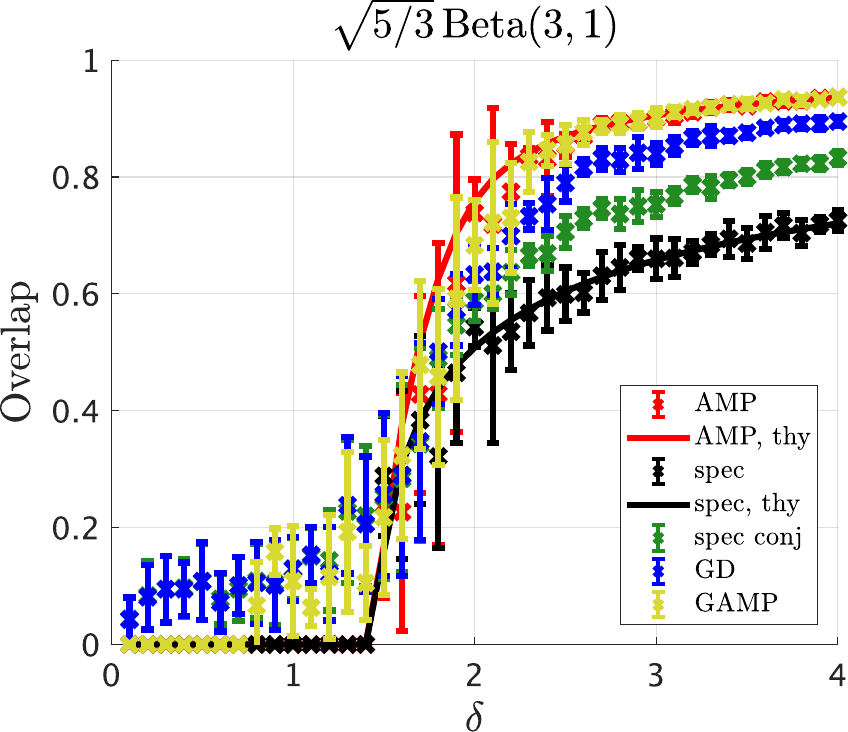}
    \end{subfigure}
    \caption{Overlaps achieved by $5$ algorithms on the Poisson regression model \Cref{eqn:pois_regr}. The same experiments are repeated for $3$ spectral distributions given in \Cref{eqn:eg_Lambda} of the design matrix $X$. }
    \label{fig:amp_pois}
\end{figure}

\Cref{fig:amp_pr,fig:amp_pois} compare the performance of five algorithms for phase retrieval and Poisson regression, respectively, for three choices of design spectral distribution. The algorithms are:
the spectral estimator in \Cref{thm:opt_thr} (`spec', shown in black in the plots), 
the spectral estimator in \cite[Equation (20)]{Luo_Alghamdi_Lu} (`spec conj' in green);
Bayes-GVAMP in \Cref{eqn:Bayes_GVAMP} (`AMP' in red), gradient descent (`GD' in blue), and the Generalized Approximate Message Passing algorithm developed for an i.i.d.\ Gaussian $X$ \cite[Section 3]{Mondelli_Venkataramanan} (`GAMP' in yellow). 

Bayes-GVAMP is initialized with `spec', whereas gradient descent and GVAMP are initialized with `spec conj'. 
The spectral estimator from \cite{Luo_Alghamdi_Lu} (`spec conj') takes the same form as \Cref{eqn:spec_est} but with the preprocessing function 
\begin{align}
    \cT^\circ(y) &= \frac{\ol{g}(y)}{\ol{g}(y) + 1} , \label{eqn:T_LAL} 
\end{align}
where $ \ol{g} $ is given by \Cref{eqn:olg_opt}. 
For i.i.d.\ Gaussian \cite{Luo_Alghamdi_Lu}, correlated Gaussian \cite{Zhang_Ji_Venkataramanan_Mondelli} and partial Haar designs \cite{Dudeja_Bakhshizadeh_Ma_Maleki}, the same $\cT^\circ$ in \Cref{eqn:T_LAL} is shown to be the optimal choice that simultaneously minimizes the weak recovery threshold and maximizes overlap within a large function class that surpasses the one defined by \Cref{asmp:preprocess}. 
Using non-rigorous methods from statistical physics, \cite[Conjecture 2]{Maillard_etal_construction} predicts that such optimality of $\cT^\circ$ holds for any orthogonally invariant design. 
The above $\cT^\circ$ does not satisfy \Cref{asmp:preprocess} and therefore our results in \Cref{thm:spec,thm:SE_Bayes_GVAMP} do not apply.
Here, we numerically test its performance and observe that the spectral estimator with respect to $\cT^\circ$ alone has inferior performance to Bayes-GVAMP initialized with our spectral estimator defined by \Cref{eqn:olg_opt} and \Cref{asmp:preprocess}. 

For gradient descent, we use the standard full-batch version
\begin{align}
    \theta^{t+1} &= \theta^t - \eta \nabla L(\theta^t) , \notag 
\end{align}
where the learning rate $\eta>0$ is fixed and the empirical loss 
\begin{align}
    L(\theta) &= \sum_{i = 1}^n \ell(\wh{q}(\inprod{x_i}{\theta}), y_i) \notag
\end{align}
is defined with respect to a postulated model $ \hat{q} \colon \bbR \to \bbR $ and a suitable loss function $\ell\colon\bbR^2 \to \bbR$ (both detailed below). 

We observe that GAMP \cite{Rangan,Mondelli_Venkataramanan}, which is believed to be computationally optimal for i.i.d.\ Gaussian designs \cite{Barbier_etal,Montanari_Wu}, is significantly suboptimal for the non-Gaussian designs in \Cref{fig:amp_pr,fig:amp_pois}. 
Indeed, in some cases, e.g.\ when $X$ has $ \Beta(1,2) $-distributed singular values (middle panels of \Cref{fig:amp_pr,fig:amp_pois}) or comes from GTEx datasets (\Cref{fig:GTEx}), GAMP does not even retain the performance of its spectral initializer \Cref{eqn:T_LAL}, and is also worse than general-purpose algorithms such as gradient descent. 

We observe that in all plots, the overlaps achieved by both the spectral estimator in \Cref{thm:opt_thr} and Bayes-GVAMP initialized with it are accurately predicted by our theories in \Cref{thm:SE_Bayes_GVAMP,thm:spec} (`spec, thy' and `AMP, thy' respectively plotted as solid curves with the corresponding color coding). 

\subsection{Synthetic data}
\label{sec:phase_synthetic}

In all numerical experiments with synthetic data, the design matrix is constructed according to \Cref{asmp:design} with $ \sfLambda $  taken to be 
\begin{align}
    \sfLambda &= \begin{cases}
        \wt{\sfLambda} , & \delta \le 1 \\
        \sqrt{\delta} \, \wt{\sfLambda} , & \delta > 1
    \end{cases} , \label{eqn:constr_Lambda} 
\end{align}
for some positive random variable $ \wt{\sfLambda} $ with compact support, second moment $1$ and positive variance. 
It is easy to verify that the above construction ensures $ \ol{\kappa}_2 = \ol{m}_2 = \expt{\sfLambda_n^2} = 1 $, a natural normalization consistent with the standard Gaussian design where $ X_{i,j} \iid \cN(0,1/d) $. 
All synthetic experiments are repeated for $3$ choices of $ \wt{\sfLambda} $: 
\begin{align}
&&
    &\sqrt{3/6} \unif([1,2]) , &
    &\sqrt{6} \Beta(1,2) , &
    &\sqrt{5/3} \Beta(3,1) , &
& \label{eqn:eg_Lambda}
\end{align}
whose densities are constant, linear and quadratic on their respective supports: 
\begin{align}
&&
    p_{\unif([1,2])}(x) &= \one_{[1,2]}(x) , & 
    p_{\Beta(1,2)}(x) &= 2(1-x) \one_{[0,1]}(x) , & 
    p_{\Beta(3,1)}(x) &= 3 x^2 \one_{[0,1]}(x) . & 
& \notag 
\end{align}
Throughout this section, we fix $ \rho = 1 $ and consider the Gaussian prior $ P_{\sfB_*} = \cN(0,1) $. 
This implies $ \sigma^2 = 1 $ by definition \Cref{eqn:sigma2}. 
We fix the number of regression coefficients to be $ d = 1000 $ and vary $\delta$ by varying $n$ accordingly. 
For each value of $\delta$, we report the overlap as per \Cref{eqn:overlap} of different algorithms averaged over $10$ i.i.d.\ trials with error bar at $1$ standard deviation. 

\paragraph{Phase retrieval.}
\label{sec:phase_retrieval}

Consider the (noiseless) phase retrieval model 
\begin{align}
    q(z, \eps) &= \abs{z} 
    \label{eqn:phase_retrieval}
\end{align}
for which $ Q\paren{y\mid z} = \delta_{\abs{z}}(y) $. 
For gradient descent, we take the postulated model $ \wh{q}(z) = \abs{z} $ to be the same as the true model $q$ in \Cref{eqn:phase_retrieval} and use the squared loss $ \ell(\wh{y},y) = \frac{1}{2} (\wh{y} - y)^2 $. 
The learning rate $\eta$ is taken to be $0.1,0.01,0.1$ respectively for the $3$ spectral distributions in \Cref{eqn:eg_Lambda}. 

\paragraph{Poisson regression.}
\label{sec:pois_regr}

Consider Poisson regression where 
\begin{align}
    y_i &\sim \pois(\inprod{x_i}{\beta_*}^2) . 
    \label{eqn:pois_regr}
\end{align}
That is, the conditional density $ Q $ is given by 
\begin{align}
&&
    Q\paren{y \mid z} &= \frac{z^{2y} e^{-z^2}}{y!} , & 
    \textnormal{ for } y\in\bbZ_{\ge0} . & 
& \notag 
\end{align}
For gradient descent, we take the postulated model $ \wh{q}(z) = z^2 $ to be the noiseless analogue of the true model \Cref{eqn:pois_regr} and use the absolute loss $ \ell(\wh{y},y) = \abs{\wh{y} - y} $. 
The learning rate is fixed to be $\eta = 0.001$ for all $3$ spectral distributions in \Cref{eqn:eg_Lambda}. 

\subsection{Real data}
\label{sec:phase_real}

For the phase retrieval model \Cref{eqn:phase_retrieval}, we also carry out numerical experiments in which the design matrices and/or regression coefficients arise from real data. 
In all such experiments, spectrally initialized Bayes-GVAMP in \Cref{eqn:Bayes_GVAMP} achieves the best performance against $4$ other competitors (\Cref{fig:CDP,fig:GTEx}). 
Unlike Bayes-GVAMP, methods such as GAMP \cite{Rangan,Mondelli_Venkataramanan} that are developed under i.i.d.\ assumptions do not adapt well to the data distribution. 
Moreover, the theoretical guarantees in \Cref{thm:spec,thm:SE_Bayes_GVAMP} respectively for the spectral estimator in \Cref{thm:opt_thr} and Bayes-GVAMP initialized with it remain accurate for real data that do not satisfy the distributional assumptions in \Cref{sec:model}. 
This suggests both the flexibility of the orthogonal invariance assumption and the universality of our theory with respect to the data distribution. 

\paragraph{Coded diffraction patterns.}
\label{sec:CDP}

Coded diffraction patterns (CDP) are a popular ensemble of designs that capture the physical measurements in phase retrieval. 
We consider the real analogue of binary and ternary complex CDP defined in \cite[Section 2.2]{candes_CDP}. 
More details are given in \Cref{sec:details_CDP}. 
The vector of regression coefficients is obtained by vectorizing, centering and normalizing the $3$ channels (R, G, B) of a $d=32\times32$ image of a truck from the CIFAR-10 dataset \cite{CIFAR10}. 
In both Bayes-GVAMP \Cref{eqn:Bayes_GVAMP} and GAMP \cite{Mondelli_Venkataramanan}, the vector denoisers are constructed assuming an i.i.d.\ Gaussian prior.  
The numerical results are shown in \Cref{fig:CDP}. 
Each data point in the plot is obtained by taking the average of $15$ trials with $5$ trials for each of the R, G, B channels. 

\paragraph{GTEx.}
\label{sec:GTEx}

We consider two GTEx datasets `Skin - Sun Exposed (Lower leg)' and `Muscle - Skeletal' (`skin' and `muscle' in respective plot titles) \cite{GTEx}. 
Upon standard preprocessing (see \Cref{sec:details_GTEx} for details), they become matrices of sizes $6461\times701$ and $7036\times803$ respectively. 
The design matrix $X$ is obtained by taking the first $n$ rows of these matrices. 
For Bayes-GVAMP, expectations involving $\sfLambda$ are replaced with averages with respect to the corresponding empirical spectral distribution of $X$. 
The regression coefficients are generated i.i.d.\ according to $ P_{\sfB_*} = \cN(0,1) $. 
The numerical results are shown in \Cref{fig:GTEx} where each data point is obtained by averaging over $100$ i.i.d.\ trials. 

\section{Discussion}
\label{sec:discuss}


\paragraph{Estimation of link function and prior distribution.}
Our general results for spectral estimators (\Cref{thm:spec}) and GVAMP (\Cref{thm:spec_GVAMP}) do not assume knowledge of the link function $q$ or prior distribution $P_{\sfB_*}$. 
However, the optimal choices of preprocessing functions for spectral estimators (\Cref{thm:opt_thr}) and denoisers for GVAMP (\Cref{thm:SE_Bayes_GVAMP}) do require such knowledge. 
When $ q $ and $ P_{\sfB_*} $ are parametrized by a fixed (relative to $n,d$) number of parameters, it is possible to consistently estimate them in the proportional regime (\Cref{asmp:proportional}) by e.g.\ moment methods. 
Then all our results continue to hold with these consistent estimates plugged in. 
In more general (semiparametric or nonparametric) scenarios, nuisance functionals such as $ q $ and $P_{\sfB_*}$ are no longer expected to be consistently estimable. 
Under both prior and link misspecification, the degradation of the posterior mean squared error is quantified in \cite{kabashima2008inference,takahashi2022macroscopic} using the non-rigorous replica method. 
Empirical Bayes methods \cite{Sawaya_Uematsu_Imaizumi,Mukherjee_Sen_Sen,Fan_Guan_Shen_Wu} have been proposed and analyzed for several Bayesian statistical models with i.i.d.\ design. 
We believe it will be fruitful to transfer similar ideas to our setting of orthogonally invariant GLM. 
A common drawback of these methods is that they require either sufficiently large sample size or sufficiently low signal-to-noise ratio. 
Developing data-driven estimation and debiasing procedures 
\cite{Zhong_Su_Fan,Li_Sur,Celentano_Montanari,Celentano_Wainwright} that retain the optimality of their oracle counterparts established in the present paper is a promising direction for future work.

\begin{figure}[t]
    \centering
    \begin{subfigure}{0.48\textwidth}
        \includegraphics[width=\linewidth]{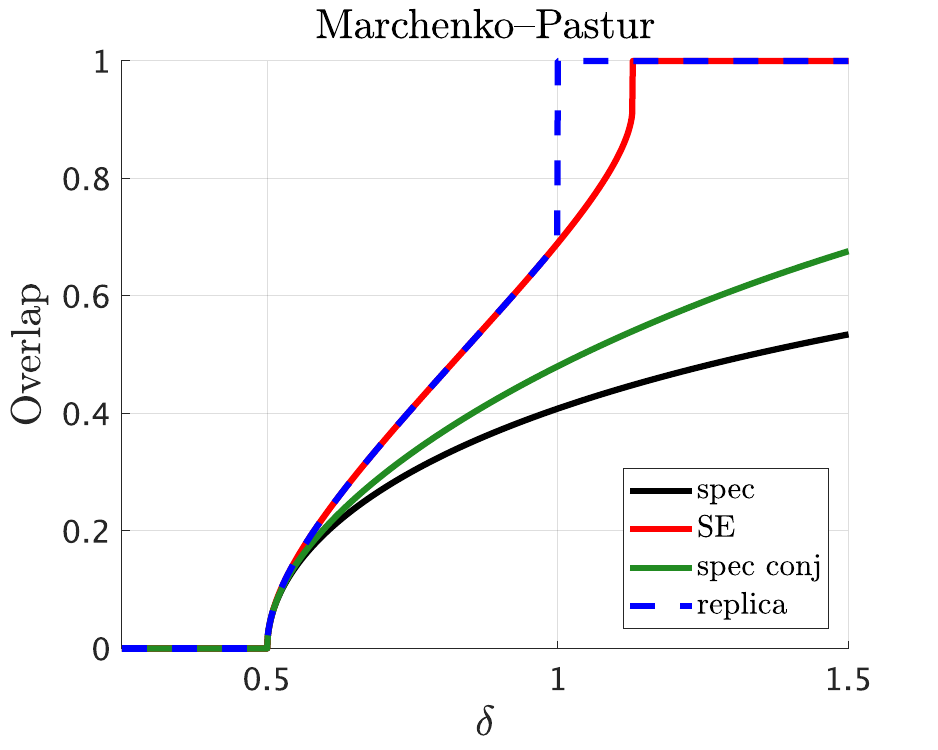}
    \end{subfigure}
    \hspace{1em}
    \begin{subfigure}{0.48\textwidth}
        \includegraphics[width=\linewidth]{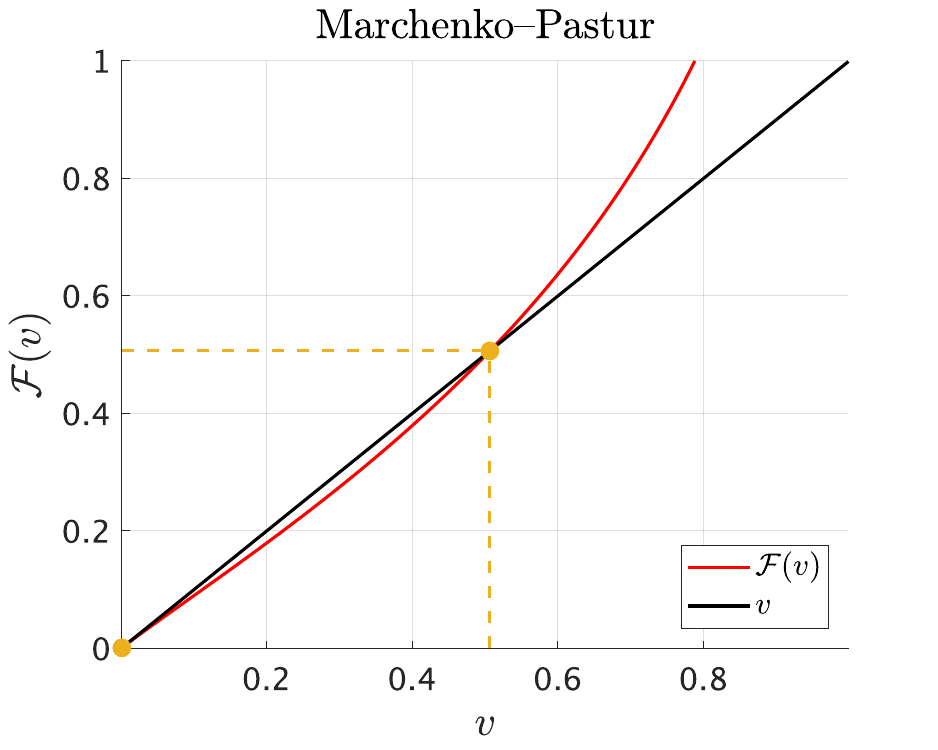}
    \end{subfigure}
    \caption{Statistical-computational tradeoffs in phase retrieval under the standard Gaussian design whose spectral distribution is given by the Marchenko--Pastur law. 
    The left panel shows $4$ curves obtained in the same way as in \Cref{fig:gap} with an additional one corresponding to the asymptotic overlap achieved by the spectral estimator given by \cite[Equation (20)]{Luo_Alghamdi_Lu}. 
    The right panel plots the function $ \cF $ at $\delta = 1.01$. }
    \label{fig:gauss}
\end{figure}

\paragraph{Statistical and computational optimality.}

For a generic orthogonally invariant GLM, our spectrally initialized Bayes-GVAMP stands as the state-of-the-art efficient estimator. 
Its performance matches the conjectured Bayes risk in some regimes while leaving gaps in other regimes, notably when the state evolution equations 
have multiple fixed points (see the discussion following \Cref{thm:replica}). 
It has been conjectured \cite{Maillard_etal_construction,Maillard_etal_threshold} that Bayes-GVAMP is optimal among all polynomial-time algorithms. 
However, we caution that the mere existence of multiple state evolution fixed points and the discrepancy between Bayes risk and the performance of spectrally initialized Bayes-GVAMP (see \Cref{fig:gauss}) do not necessarily entail a fundamental gap between statistical and computational limits. 
For instance, for noiseless phase retrieval \Cref{eqn:phase_retrieval} with an i.i.d.\ Gaussian design, there exist efficient algorithms \cite{Andoni_Hsu_Shi_Sun,Song_Zadik_Bruna} achieving overlap $1$ (i.e., perfect recovery) as long as $n\ge d+1$ (in particular, $ \delta \ge 1 $). 
In contrast, Bayes-GVAMP or GAMP requires $\delta \gtrsim 1.13$ to the same end. 
Understanding when and why Bayes-GVAMP is computationally optimal is an interesting avenue for future research.

\section*{Acknowledgement}
This work was partially done when Y.Z.\ was at the Institute of Science and Technology Austria, where he was funded by the European Union (ERC, INF$^2$, project number 101161364).
M.M.\ is funded by the European Union (ERC, INF$^2$, project number 101161364). Views and opinions expressed are however those of the author(s) only and do not necessarily reflect those of the European Union or the European Research Council Executive Agency. Neither the European Union nor the granting authority can be held responsible for them. Y.Z.\ thanks Rishabh Dudeja for an illuminating discussion during his visit to the Department of Statistics at UW-Madison in Feb 2024 which critically shifted the authors' attention from RI-GAMP \cite{VKM} to VAMP \cite{Dudeja_Lu_Sen,Dudeja_Sen_Lu}. 
Y.Z.\ also thanks Junjie Ma for a useful discussion of \cite{Dudeja_Liu_Ma} during the workshop ``Towards a theory of typical-case algorithmic hardness'' at the Les Houches Physics School in Feb 2025.

\renewcommand*{\bibfont}{\normalfont\small}
\printbibliography

\newpage 

\appendix

\section{Free cumulants and moments}
\label{sec:cumulant}

For a matrix $ X \in \bbR^{n\times d} $, denote by $ m_{2k} $ the $k$-th (where $ k \in \bbZ_{\ge0} $) moment of the empirical spectral distribution of $ XX^\top $. 
The $k$-th moment of the empirical spectral distribution of $ X^\top X $ is then given by 
\begin{align}
m'_{2k} &= \begin{cases}
1 , & k = 0 \\
\wh{\delta} m_{2k} , & k\ge1
\end{cases} , \notag 
\end{align}
where $ \wh{\delta} = n/d $. 
Denote by $ (\kappa_{2k})_{k\ge1} $ the rectangular free cumulants of $ XX^\top $; see, e.g., \cite[Section 2.4]{fan2020approximate} for the precise definition. 
The first four of them are explicitly given by: 
\begin{subequations}
\label{eqn:moment_cumulant}
\begin{align}
\kappa_2 &= m_2 = \frac{m_2'}{\wh{\delta}} , \\
\kappa_4 &= m_4 - (1 + \wh{\delta}) m_2^2
= \frac{m_4'}{\wh{\delta}} - \frac{1 + \wh{\delta}}{\wh{\delta}^2} m_2'^2 , \\
\kappa_6 &= m_6 - (3 + 3\wh{\delta}) m_4 m_2 + (2 + 3\wh{\delta} + 2\wh{\delta}^2) m_2^3
= \frac{m_6'}{\wh{\delta}} - \frac{3 + 3\wh{\delta}}{\wh{\delta}^2} m_4' m_2' + \frac{2 + 3\wh{\delta} + 2\wh{\delta}^2}{\wh{\delta}^3} m_2'^3 , \\
\kappa_8 &= m_8 - (4 + 4\wh{\delta}) m_6 m_2 - (2 + 2\wh{\delta}) m_4^2 + (10 + 16\wh{\delta} + 10\wh{\delta}^2) m_4 m_2^2 - (5 + 10\wh{\delta} + 10\wh{\delta}^2 + 5\wh{\delta}^3) m_2^4 \notag \\
&= \frac{m_8'}{\wh{\delta}} - \frac{4 + 4\wh{\delta}}{\wh{\delta}^2} m_6' m_2' - \frac{2 + 2\wh{\delta}}{\wh{\delta}^2} m_4'^2 + \frac{10 + 16\wh{\delta} + 10\wh{\delta}^2}{\wh{\delta}^3} m_4' m_2'^2 - \frac{5 + 10\wh{\delta} + 10\wh{\delta}^2 + 5\wh{\delta}^3}{\wh{\delta}^4} m_2'^4 . 
\end{align}
\end{subequations}
The rectangular free cumulants of $ X^\top X $ are denoted by $ (\kappa_{2k}')_{k\ge1} $ and satisfy $ \kappa_{2k}' = \wh{\delta} \kappa_{2k} $. 

If the empirical spectral distributions of $ XX^\top $ and $ X^\top X $ converge in the sense of \Cref{asmp:design}, then $ m_{2k} , m_{2k}' $ and $ \kappa_{2k}, \kappa_{2k}' $ also converge. 
We denote their limits by $ \ol{m}_{2k} , \ol{m}_{2k}' $ and $ \ol{\kappa}_{2k}, \ol{\kappa}_{2k}' $. 
Due to the one-to-one correspondence $ \ol{\kappa}_{2k} = \delta \ol{\kappa}_{2k}' $, we will always refer to $ \ol{\kappa}_{2k} $. 

\section{Proof of \Cref{thm:spec}}
\label{sec:lin_VAMP_SE}

As signposted in \Cref{sec:heuristics}, let us start by considering the following linearized GVAMP iteration: for $t\ge1$, 
\begin{subequations} 
\label{eqn:lin_GVAMP}
\begin{align}
&&
    v^{t} &= \ol{\kappa}_2^{-1} X^\top g^t , & 
    &&
& \\
&&
    u^{t} &= \paren{ \ol{\kappa}_2^{-1} X X^\top - I_n } g^t , &
    g^{t+1} &= \ol{G} u^t , &
& 
\end{align}
\end{subequations}
where $ \ol{G} = \diag(\ol{g}(y)) \in \bbR^{n\times n} $ for any function $ \ol{g} $ subject to \Cref{asmp:preprocess}. 
Suppose that the initializer $ g^1 $ is given by 
\begin{align}
    g^1 &= \ol{G} (\nu_0 z + \tau_0 n^0) , \label{eqn:lin_init} 
\end{align}
where $ n^0 \sim \cN(0_n, I_n) $ is independent of everything else and $ \nu_0, \tau_0 $ are defined as 
\begin{align}
&&
    \nu_0 &= \frac{1}{\sqrt{\rho + w_2}} , & 
    \tau_0 &= \sqrt{\frac{w_1}{\rho + w_2}} . & 
& \label{eqn:nu0_tau0}
\end{align}

\Cref{lem:SE_stay} below (proved in \Cref{sec:pf_spec_step1}) shows that under the initialization scheme \Cref{eqn:lin_init}, several state evolution parameters of \Cref{eqn:lin_GVAMP} stay at their fixed points across all iterates. 

\begin{lemma}[Stationary state evolution parameters]
\label{lem:SE_stay}
Consider the linearized GVAMP in \Cref{eqn:lin_GVAMP} initialized with \Cref{eqn:lin_init}. 
Let \Cref{eqn:thr} hold. 
Then for every $t\ge1$, there exist $ \chi_t,\nu_t \ne 0 $ and $ \omega_t,\tau_t>0 $ such that
\begin{align}
&&
    \matrix{v^t & \beta_*} &\xrightarrow{W_2} \matrix{\chi_t \sfB_* + \omega_t \sfM_t & \sfB_*} , & 
    \matrix{u^t & z & \eps} &\xrightarrow{W_2} \matrix{\nu_t \sfZ + \tau_t \sfN_t & \sfZ & \sfE} , & 
& \label{eqn:lin_conv} 
\end{align}
where $ (\sfB_*, \sfM_t) \sim P_{\sfB_*} \ot \cN(0,1) $ and $ (\sfZ, \sfM_t) \sim \cN(0,\sigma^2) \ot \cN(0,1) $. 
Moreover, 
\begin{align}
&&
    \frac{\tau_t^2}{\nu_t^2} &= w_1 , & 
    \frac{\omega_t^2}{\chi_t^2} &= w_2 , & 
    \frac{\nu_t}{\chi_t} &= w_3 , & 
    \frac{\chi_{t+1}}{\nu_t} &= w_4 , & 
    \frac{\nu_t}{\nu_{t-1}} &= \gamma , & 
& \label{eqn:SE_stay}
\end{align}
where the right-hand sides of all equalities above are positive. 
\end{lemma}

From \Cref{lem:SE_stay}, let us highlight a somewhat unusual feature of the linearized GVAMP algorithm in \Cref{eqn:lin_GVAMP}. 
Though the state evolution parameters are initialized with positive finite values in \Cref{eqn:nu0_tau0}, they may not have a positive finite limit as $t\to\infty$. 
According to \Cref{eqn:SE_stay}, the state evolution parameters form a geometric series indexed by $t$ of rate $\gamma$ and may converge to $0$ or $\infty$ as $t\to\infty$ unless $\gamma = 1$. 
This implies, for instance, that as $t\to\infty$, $ \lim_{n\to\infty} n^{-1} \normtwo{u^{t+1}}^2 $ may not converge to a positive finite value and $ \lim_{n\to\infty} n^{-1} \normtwo{u^{t+1} - u^t}^2 $ may not converge to $0$. 
However, \Cref{lem:SE_stay} shows that certain ratios of the state evolution parameters do have nontrivial limits. 
This can be used to show that as $t\to\infty$, $ \lim_{n\to\infty} \normtwo{u^{t+1}} / \normtwo{u^t} $ converges to $\gamma$ (see \Cref{eqn:gamma}) and $ \lim_{n\to\infty} \normtwo{\normtwo{u^{t+1}}^{-1} u^{t+1} - \normtwo{u^t}^{-1} u^t} $ converges to $0$ (see \Cref{eqn:whu_whu}). 
Identifying such ratios whose fixed points (in the large $t$ limit) exist is a critical step towards proving the alignment between $ v^{t+1} $ and $ v_1(D) $ in \Cref{lem:align}. 

\begin{lemma}[Right edge of the bulk]
\label{lem:edge}
Consider the spectral estimator defined through \Cref{eqn:spec_est} using a preprocessing function satisfying \Cref{asmp:preprocess}. 
Then we have 
\begin{align}
\lim_{d\to\infty} \lambda_2(D) &= \lambda^\circ , \notag 
\end{align}
where $ \lambda^\circ $ on the right-hand side is defined in \Cref{eqn:def_lambda2}. 
\end{lemma}

\Cref{lem:edge} is proved in \Cref{app:right_edge}.

The next lemma shows that $ v^t $ is aligned with the top eigenvector of $ D $ for sufficiently large $t$ and $d$. 
The proof is presented in \Cref{app:pf_lem:align}. 

\begin{lemma}
\label{lem:align}
Consider the linearized GVAMP in \Cref{eqn:lin_GVAMP} initialized with \Cref{eqn:lin_init}. 
Let \Cref{eqn:thr} hold. 
Then almost surely, 
\begin{align}
&&
    \lim_{t\to\infty} \lim_{d\to\infty} \frac{\abs{ \inprod{v^{t+1}}{v_1(D)} }}{\normtwo{v^{t+1}}} &= 1 , & 
    \lim_{d\to\infty} \lambda_1(D) &= \ol{\kappa}_2 , & 
& \notag 
\end{align}
where $D$ is defined in \Cref{eqn:D}. 
\end{lemma}

Now, \Cref{thm:spec} follows easily from the preceding $3$ lemmas. 

\begin{proof}[Proof of \Cref{thm:spec}]
The characterizations in \Cref{eqn:lambda2_D,eqn:lambda1_eta} of the first and second largest eigenvalues of $D$ have been obtained in \Cref{lem:edge,lem:align}. 
It remains to compute the performance of the spectral estimator. 
Due to \Cref{lem:align}, it suffices to compute the normalized squared overlap between $\beta_*$ and $v^t$. Using \Cref{lem:SE_stay}, we have
\begin{align}
\lim_{d\to\infty} \frac{\abs{ \inprod{\beta_*}{v_1(D)} }^2}{\normtwo{\beta_*}^2 }
&= \lim_{t\to\infty} \lim_{d\to\infty} \frac{\abs{\inprod{\beta_*}{v^t}}^2}{\normtwo{\beta_*}^2 \normtwo{v^t}^2}
= \lim_{t\to\infty} \frac{\chi_t^2 \rho}{ \chi_t^2 \rho + \omega_t^2 } 
= \lim_{t\to\infty} \frac{\rho}{ \rho + \frac{\omega_t^2}{\chi_t^2} } 
= \frac{\rho}{\rho + w_2}
= \eta , 
\notag 
\end{align}
where the last two equalities follow from the definitions of $ \eta, w_2 $ in \Cref{eqn:omega2/chi2,eqn:eta} and elementary calculations. 
Note that $ \eta \in (0,1) $ by $ w_2 \in (0,\infty) $ from \Cref{lem:SE_stay}. 
\end{proof}


\section{Proof of \Cref{lem:SE_stay}}
\label{sec:pf_spec_step1}

Consider the following iteration that is slightly more general than \Cref{eqn:lin_GVAMP}. 
For $t\ge1$, 
\begin{subequations} \label{eqn:GVAMP_lin} 
\begin{align}
&&
    v^{t} &= \ol{\kappa}_2^{-1} X^\top g^t , & 
    &&
& \\
&&
    u^{t} &= \paren{ \ol{\kappa}_2^{-1} X X^\top - I_n } g^t , &
    g^{t+1} &= g_{t+1}(u^t; z, \eps) , &
& 
\end{align}
\end{subequations}
initialized with $ g^1 = g_1(z, \eps) $.  
\Cref{eqn:GVAMP_lin} can be cast as an instance of GVAMP in \Cref{eqn:GVAMP} by setting
\begin{subequations}
\label{eqn:specialize_lin}
\begin{align}
&&
    \Phi_t(x) &= 0 , & 
    \wt{\Phi}_t(x) &= \ol{\kappa}_2^{-1} x , & 
    \Psi_t(x) &= \ol{\kappa}_2^{-1} x - 1 , & 
    \wt{\Psi}_t(x) &= 0 , & 
& 
\end{align}
and making the change of variables
\begin{align}
&&
    r^{t+1} &= v^t , & 
    \wt{r}^{t+1} &= 0_d , & 
    p^{t+1} &= u^t , & 
    \wt{p}^{t+1} &= g^{t+1} . & 
& 
\end{align}
\end{subequations}

Let us define the state evolution recursion for \Cref{eqn:GVAMP_lin}. 
\begin{subequations} \label{eqn:SE_lin}
Define random variables: for $ t\ge1 $, 
\begin{align}
\matrix{ \sfB_* & \sfW_1 & \cdots & \sfW_t }^\top &\sim P_{\sfB_*} \ot \cN(0_t, \Omega_t) , \\
\matrix{ \sfB_1 & \cdots & \sfB_t }^\top &= \matrix{\chi_1 & \cdots & \chi_t}^\top \sfB_* + \matrix{ \sfW_1 & \cdots & \sfW_t }^\top , \\
\matrix{\sfE & \sfZ & \sfZ_1 & \cdots & \sfZ_t }^\top &\sim P_{\sfE} \ot \cN(0_t, \Sigma_{t+1}) . 
\end{align}
The state evolution parameters $ \Omega_t \in \bbR^{t\times t} $, $ \chi_t \in \bbR $ and $ \Sigma_{t+1} \in \bbR^{(t+1)\times(t+1)} $ in the above display are defined as follows. 
We first recursively define semi-infinite arrays $ (\Omega_{i,j})_{(i,j)\in\bbZ_{\ge1}} $, $ (\Sigma_{i,j})_{(i,j)\in\bbZ_{\ge1}} $, and $ (\chi_i)_{i\in\bbZ_{\ge1}} $ initialized with
\begin{align}
&&
    \Sigma_{1,1} &= \sigma^2 , &
    \Omega_{1,1} &= \delta \frac{1}{\ol{\kappa}_2} \expt{g_1(\sfZ, \sfE)^2} + \delta \frac{\ol{\kappa}_4}{\ol{\kappa}_2^2} \rho \expt{\partial_z g_1(\sfZ, \sfE)}^2 , &
    \chi_1 &= \delta \expt{\partial_z g_1(\sfZ, \sfE)} . &
& \label{eqn:SE_lin_init}
\end{align}
Here and below, we use $ \partial_z g_t(z_{t-1}; z, e) $ to denote the partial derivative of $ g_t $ with respect to its second argument. 
For any $t\ge1$, $ 1\le i\le t$ and $0\le j\le t $, 
\begin{align}
\Sigma_{t+1,1} 
&= \ol{\kappa}_2 \rho \chi_t + \frac{\ol{\kappa}_4}{\ol{\kappa}_2} \rho \expt{\partial_z g_t(\sfZ_{t-1}; \sfZ, \sfE)} , 
\label{eqn:SE_lin1} \\
\Sigma_{t+1,i+1}
&= \ol{\kappa}_2 \paren{\chi_t \chi_{i} \rho + \Omega_{t,i}} + \frac{\ol{\kappa}_4}{\ol{\kappa}_2} \paren{\expt{\partial_z g_t(\sfZ_{t-1}; \sfZ, \sfE)} \rho \chi_{i} + \expt{\partial_z g_{i}(\sfZ_{i-1}; \sfZ, \sfE)} \rho \chi_t} \notag \\
&\quad + \frac{\ol{\kappa}_4}{\ol{\kappa}_2^2} \expt{g_t(\sfZ_{t-1}; \sfZ, \sfE) g_{i}(\sfZ_{i-1}; \sfZ, \sfE)} + \frac{\ol{\kappa}_6}{\ol{\kappa}_2^2} \rho \expt{\partial_z g_{i}(\sfZ_{i-1}; \sfZ, \sfE)} \expt{\partial_z g_t(\sfZ_{t-1}; \sfZ, \sfE)} , 
\label{eqn:SE_lin2} \\
\Omega_{t+1,j+1} 
&= \delta \paren{\frac{1}{\ol{\kappa}_2} \expt{g_{t+1}(\sfZ_t; \sfZ, \sfE) g_{j+1}(\sfZ_{j}; \sfZ, \sfE)} + \frac{\ol{\kappa}_4}{\ol{\kappa}_2^2} \rho \expt{\partial_z g_{t+1}(\sfZ_t; \sfZ, \sfE)} \expt{\partial_z g_{j+1}(\sfZ_{j}; \sfZ, \sfE)}} , 
\label{eqn:SE_lin3} \\
\chi_{t+1} 
&= \delta \expt{\partial_z g_{t+1}(\sfZ_t; \sfZ, \sfE)} , 
\label{eqn:SE_lin4}
\end{align}
and $ \Sigma_{1,t+1} = \Sigma_{t+1,1}, \Sigma_{i+1,t+1} = \Sigma_{t+1,i+1}, \Omega_{j+1,t+1} = \Omega_{t+1,j+1} $. 

Finally, $ \Omega_t = (\Omega_{i,j})_{(i,j)\in[t]^2} $ and $ \Sigma_{t+1} = (\Sigma_{i,j})_{(i,j)\in[t+1]^2} $ are obtained by taking the leading principal minors of the above semi-infinite arrays. 
\end{subequations}

The validity of the state evolution result, \Cref{prop:SE} below, crucially relies on the divergence-free property of the denoisers $ (g_{t+1})_{t\ge1} $. 

\begin{enumerate}[label=(A\arabic*)]
\setcounter{enumi}{\value{asmpctr}}

    \item\label[asmp]{asmp:div_free_lin} 
    For every $ t\ge0 $, $ g_t $ is Lipschitz in all of its arguments. 
    For every $ t \ge 1 $, 
    \begin{align}
    &&
        \expt{ g_{t+1}'(\sfZ_{t}; \sfZ, \sfE) } &= 0 , &
    & 
    \end{align}
    where the distribution of $ (\sfZ_t, \sfZ, \sfE) $ is specified in \Cref{eqn:SE_lin}. 
    Here and below, $g_{t+1}'(z_{t}; z, e)$ denotes the partial derivative of $g_{t+1}$ with respect to its first argument. 
    Moreover, $g_{t+1}'(z_t; z, \eps)$ and $ \partial_z g_{t+1}(z_t; z, \eps) $ are continuous with probability $1$ with respect to the law of $ (\sfZ_t, \sfZ, \sfE) $. 
    
\setcounter{asmpctr}{\value{enumi}}
\end{enumerate}

\begin{lemma}[State evolution]
\label{prop:SE}
Consider the linearized GVAMP algorithm in \Cref{eqn:GVAMP_lin} and let \Cref{asmp:div_free_lin} hold. 
For any $ t\ge1 $, 
\begin{align}
    \matrix{\beta_* & v^1 & \cdots & v^t} & \xrightarrow{W_2} \matrix{\sfB_* & \sfB_1 & \cdots & \sfB_t} , \notag \\
    \matrix{\eps & z & u^1 & \cdots & u^t & g^1 & g^2 & \cdots & g^{t+1}} & \xrightarrow{W_2} \matrix{\sfE & \sfZ & \sfZ_1 & \cdots & \sfZ_t & g_1(\sfZ, \sfE) & g_2(\sfZ_1; \sfZ, \sfE) & \cdots & g_{t+1}(\sfZ_t; \sfZ, \sfE)} , \notag 
\end{align}
where the random variables on the right-hand sides are defined in \Cref{eqn:SE_lin}. 
\end{lemma}

\Cref{prop:SE} is proved in \Cref{app:pf_prop:SE}. 

Next, we present an equivalent representation of the random variables $ \sfB_t, \sfZ_t $; see \Cref{app:pf_prop:BZ_equiv} for a proof. 
\begin{lemma}[Alternative state evolution]
\label{prop:BZ_equiv}
The random variables $ (\sfB_i)_{i\in[t]}, (\sfZ_i)_{i\in[t]} $ defined in \Cref{eqn:SE_lin} can be equivalently written as: 
\begin{align}
&&
    \sfB_i &= \chi_i \sfB_* + \omega_i \sfM_i , &
    \sfZ_i &= \nu_i \sfZ + \tau_i \sfN_i , &
    & \forall i\in[t] . & 
& \notag 
\end{align}
In the above display, $ (\omega_i \sfM_i)_{i\in[t]} $ is a centered Gaussian vector independent of $ \sfB_* \sim P_{\sfB_*} $ with $ \sfM_i \sim \cN(0,1), \forall i\in[t] $; 
$ (\tau_i \sfN_i)_{i\in[t]} $ is a centered Gaussian vector independent of $ \sfZ \sim \cN(0,\sigma^2) $ with $ \sfN_i \sim \cN(0,1), \forall i\in[t] $. 
The coefficients $ (\omega_i)_{i\in[t]}, (\tau_i)_{i\in[t]} \in \bbR_{>0}^t $ and the correlations $ (\expt{\omega_r \sfM_r \cdot \omega_s \sfM_s})_{(r,s)\in[t]^2} $, $ (\expt{\tau_r \sfN_r \cdot \nu_s \sfN_s})_{(r,s)\in[t]^2} $ follow the recursions below.  
For $t=1$, 
\begin{align}
&& 
    \chi_1 &= \delta \expt{\partial_z g_1(\sfZ, \sfE)} , & 
    \omega_1 &= \paren{ \frac{\delta}{\ol{\kappa}_2} \expt{g_1(\sfZ, \sfE)^2} + \delta \frac{\ol{\kappa}_4}{\ol{\kappa}_2^2} \rho \expt{\partial_z g_1(\sfZ, \sfE)}^2 }^{1/2} . & 
& \label{eqn:chi1_omega1} 
\end{align}
For $ t\ge1 $, 
\begin{align}
\nu_t &= \chi_t + \frac{\ol{\kappa}_4}{\ol{\kappa}_2^2} \expt{\partial_z g_t(\sfZ_{t-1}; \sfZ, \sfE)} , \notag \\
\tau_t^2 &= \ol{\kappa}_2 \omega_t^2 + \frac{\ol{\kappa}_4}{\ol{\kappa}_2^2} \expt{g_t(\sfZ_{t-1}; \sfZ, \sfE)^2} + \paren{\frac{\ol{\kappa}_6}{\kappa_2^2} - \frac{\ol{\kappa}_4^2}{\ol{\kappa}_2^3}} \rho \expt{\partial_z g_t(\sfZ_{t-1}; \sfZ, \sfE)}^2 , \notag \\
\chi_{t+1} &= \delta \expt{\partial_z g_{t+1}(\sfZ_t; \sfZ, \sfE)} , \notag \\
\omega_{t+1}^2 &= \frac{\delta}{\ol{\kappa}_2} \expt{g_{t+1}(\sfZ_{t}; \sfZ, \sfE)^2} + \delta \frac{\ol{\kappa}_4}{\ol{\kappa}_2^2} \rho \expt{\partial_z g_{t+1}(\sfZ_t; \sfZ, \sfE)}^2 . \notag
\end{align}
For $ r,s \ge 1 $, 
\begin{align}
    \expt{\tau_r \sfN_r \cdot \tau_s \sfN_s}
    &= \ol{\kappa}_2 \expt{\omega_r \sfM_r \cdot \omega_s \sfM_s} 
    + \frac{\ol{\kappa}_4}{\ol{\kappa}_2^2} \expt{g_r(\sfZ_{r-1}; \sfZ, \sfE) g_{s}(\sfZ_{s-1}; \sfZ, \sfE)} \notag \\
    &\quad + \paren{\frac{\ol{\kappa}_6}{\ol{\kappa}_2^2} - \frac{\ol{\kappa}_4^2}{\ol{\kappa}_2^3}} \rho \expt{\partial_z g_{r}(\sfZ_{r-1}; \sfZ, \sfE)} \expt{\partial_z g_s(\sfZ_{s-1}; \sfZ, \sfE)} , \notag \\
    \expt{\omega_r \sfM_r \cdot \omega_s \sfM_s}
    &= \delta \paren{\frac{1}{\ol{\kappa}_2} \expt{g_{r}(\sfZ_{r-1}; \sfZ, \sfE) g_{s}(\sfZ_{s-1}; \sfZ, \sfE)} + \frac{\ol{\kappa}_4}{\ol{\kappa}_2^2} \rho \expt{\partial_z g_{r}(\sfZ_{r-1}; \sfZ, \sfE)} \expt{\partial_z g_{s}(\sfZ_{s-1}; \sfZ, \sfE)}} . \notag 
\end{align}
\end{lemma}

As suggested by the heuristics in \Cref{sec:heuristics}, we will consider $ g_{t+1} $ of the following form which satisfies \Cref{asmp:div_free_lin}. 
\begin{enumerate}[label=(A\arabic*)]
\setcounter{enumi}{\value{asmpctr}}

    \item\label[asmp]{asmp:g(y)} \begin{subequations} \label{eqn:g_choice}
    The function $ g_{t+1} $ in \Cref{eqn:GVAMP_lin} takes the form
    \begin{align}
        g_{t+1}(u^t; z, \eps) &= \ol{g}_{t+1}(q(z, \eps)) u^t \label{eqn:g_sep}
    \end{align}
    for some function $ \ol{g}_{t+1} $ satisfying:
    \begin{align}
    &&
        \expt{ \ol{g}_{t+1}(\sfY) } &= 0 , & 
        \expt{\ol{g}_{t+1}(\sfY) \ol{\sfZ}^2} &> 0 . &
    & \label{eqn:mean0}
    \end{align}
    \end{subequations}
    
\setcounter{asmpctr}{\value{enumi}}
\end{enumerate}

Specializing \Cref{prop:SE} to this case, we have the following \Cref{prop:SE_lin} whose proof can be found in \Cref{app:pf_prop:SE_lin}. 

\begin{lemma}[State evolution, signal-to-noise ratios]
\label{prop:SE_lin}
Consider the linearized GVAMP in  \Cref{eqn:GVAMP_lin}. 
Let \Cref{asmp:g(y)} and \Cref{eqn:thr} hold. 
If $ \chi_1 \ne 0 $, then $ \nu_t, \chi_{t+1} \ne 0 $ and $ \tau_t, \omega_{t} > 0 $ for every $t\ge1$. 
Furthermore, the following recursions can be derived from the state evolution in \Cref{eqn:SE_lin}. 
For any $t\ge1$, $1\le i\le t$ and $0\le j\le t $, 
\begin{align}
\frac{\tau_t^2}{\nu_t^2} &= 
\frac{\paren{ \frac{\ol{\kappa}_4}{\ol{\kappa}_2^2} + \delta } \expt{\ol{g}_t(\sfY)^2}}{ \paren{ \frac{\ol{\kappa}_4}{\ol{\kappa}_2^2} + \delta }^2 \expt{\ol{g}_t(\sfY) \ol{\sfZ}^2}^2 } \frac{\tau_{t-1}^2}{\nu_{t-1}^2} 
+ \frac{\rho \brace{ \paren{ \frac{\ol{\kappa}_4}{\ol{\kappa}_2} + \delta \ol{\kappa}_2 } \expt{\ol{g}_t(\sfY)^2 \ol{\sfZ}^2} + \paren{\frac{\ol{\kappa}_6}{\ol{\kappa}_2^2} - \frac{\ol{\kappa}_4^2}{\ol{\kappa}_2^3} + \delta \frac{\ol{\kappa}_4}{\ol{\kappa}_2}} \expt{\ol{g}_t(\sfY) \ol{\sfZ}^2}^2 }}{ \paren{ \frac{\ol{\kappa}_4}{\ol{\kappa}_2^2} + \delta }^2 \expt{\ol{g}_t(\sfY) \ol{\sfZ}^2}^2 } .  
\label{eqn:tau2/nu2} \\
\frac{\omega_{t+1}^2}{\chi_{t+1}^2} &=
\frac{ \ol{\kappa}_2^{-1} \expt{\ol{g}_{t+1}(\sfY)^2} }{ \delta \expt{\ol{g}_{t+1}(\sfY) \ol{\sfZ}^2}^2 } \frac{\tau_t^2}{\nu_t^2}
+ \frac{ \rho \brace{ \expt{\ol{g}_{t+1}(\sfY)^2 \ol{\sfZ}^2} + \frac{\ol{\kappa}_4}{\ol{\kappa}_2^2} \expt{\ol{g}_{t+1}(\sfY) \ol{\sfZ}^2}^2 } }{ \delta \expt{\ol{g}_{t+1}(\sfY) \ol{\sfZ}^2}^2 } , 
\label{eqn:omega2/chi2} \\
\Sigma_{t+1,i+1} &= \ol{\kappa}_2 \paren{\chi_t \chi_{i} \rho + \Omega_{t,i}} + \frac{\ol{\kappa}_4}{\ol{\kappa}_2} \rho \paren{\expt{\ol{g}_t(\sfY) \ol{\sfZ}^2} \nu_{t-1} \chi_{i} + \expt{\ol{g}_i(\sfY) \ol{\sfZ}^2} \nu_{i-1} \chi_t} \notag \\
&\quad + \frac{\ol{\kappa}_4}{\ol{\kappa}_2^2} \expt{\ol{g}_t(\sfY) \ol{g}_i(\sfY)} \expt{\tau_{t-1} \sfN_{t-1} \cdot \tau_{i-1} \sfN_{i-1}} \notag \\
&\quad + \paren{
    \frac{\ol{\kappa}_4}{\ol{\kappa}_2} \expt{\ol{g}_t(\sfY) \ol{g}_i(\sfY) \ol{\sfZ}^2} 
    + \frac{\ol{\kappa}_6}{\ol{\kappa}_2^2} \expt{\ol{g}_t(\sfY) \ol{\sfZ}^2} \expt{\ol{g}_i(\sfY) \ol{\sfZ}^2}
} \rho \nu_{t-1} \nu_{i-1} , 
\label{eqn:Sigmati} \\
\Omega_{t+1,j+1} &= \delta \brace{
    \rho \paren{\expt{\ol{g}_{t+1}(\sfY) \ol{g}_{j+1}(\sfY) \ol{\sfZ}^2} + \frac{\ol{\kappa}_4}{\ol{\kappa}_2^2} \expt{\ol{g}_{t+1}(\sfY) \ol{\sfZ}^2} \expt{\ol{g}_{j+1}(\sfY) \ol{\sfZ}^2}} \nu_{t} \nu_{j} 
    + \frac{1}{\ol{\kappa}_2} \expt{\ol{g}_{t+1}(\sfY) \ol{g}_{j+1}(\sfY)} \expt{\tau_{t} \sfN_{t} \cdot \tau_{j} \sfN_{j}}
} . 
\label{eqn:Omegati}
\end{align}
\end{lemma}

If we further choose $ \ol{g}_{t+1} = \ol{g} $ to be time-invariant, then the state evolution has the following fixed points.  
See \Cref{app:pf_prop:SE_FP} for a proof. 

\begin{lemma}[State evolution fixed points]
\label{prop:SE_FP}
Take $g_{t+1}$ of the form \Cref{eqn:g_sep} with $ \ol{g}_{t+1} = \ol{g} $ satisfying \Cref{eqn:mean0}. 
Denote the fixed points of the state evolution parameters in \Cref{eqn:tau2/nu2,eqn:omega2/chi2} by the same notation with time indices removed. 
Assuming that \eqref{eqn:thr} holds, 
the fixed points satisfy: 
\begin{align}
&&
    \frac{\ol{\tau}^2}{\ol{\nu}^2} &= w_1 > 0 , &
    \frac{\ol{\omega}^2}{\ol{\chi}^2} &= w_2 > 0 . &
& \notag 
\end{align}
\end{lemma}

Finally, the above ingredients allow us to prove \Cref{lem:SE_stay} in \Cref{app:pf_lem:SE_stay}. 

\subsection{Proof of \Cref{prop:SE}}
\label{app:pf_prop:SE}

According to the specialization in \Cref{eqn:specialize_lin}, the linearized GVAMP in \Cref{eqn:GVAMP_lin} is an instance of the generic GVAMP in  \Cref{eqn:GVAMP}. 
Therefore its state evolution \Cref{eqn:SE_lin} can be derived by specializing \Cref{eqn:SE_GVAMP}. 
However, we will derive its state evolution from the following Rotationally Invariant Generalized Approximate Message Passing (RI-GAMP) algorithm studied in \cite[Equations (4) and (5)]{VKM}. 
This is for the technical convenience of making the appearance of free cumulants in \Cref{eqn:SE_lin} more transparent. 

The RI-GAMP iteration reads: for $t\ge1$, 
\begin{subequations}
\label{eqn:RIGAMP}
\begin{align}
&&
    b^t &= X^\top h^t - \sum_{i = 1}^{t-1} \beta_{t,i} f^i , &
    f^t &= f_t(b^1, \cdots, b^t), &
& \\
&&
    z^t &= X f^t - \sum_{i = 1}^t \alpha_{t,i} h^i , &
    h^{t+1} &= h_{t+1}(z^1, \cdots, z^{t}; z, \eps), &
&
\end{align}
\end{subequations}
initialized with $ h^1 = h_1(z, \eps) $. 

\begin{proposition}[{\cite[Theorem 1]{VKM}}]
\label{prop:SE_VKM}
For any fixed $t\ge1$, there exist random variables $ \sfZ, \sfE, (\sfB_i)_{i\in[t]}, (\sfZ_i)_{i\in[t]} $ (defined in \cite[Equations (11) to (22)]{VKM}) such that if the following conditions hold
\begin{enumerate}[label=(A\arabic*)]
\setcounter{enumi}{\value{asmpctr}}
    \item \label[asmp]{asmp:lin_denoiser} for every $1\le i\le t$, $f_i(b_1, \cdots, b_i)$ and $ h_i(z_1, \cdots, z_{i-1}; z, \eps) $ are Lipschitz in each of their arguments; 

    \item \label[asmp]{asmp:lin_deriv} for every $1\le i\le t$, $ 1\le j\le i $ and $ 1\le k\le i-1 $, the partial derivatives $ \partial_{x_j} f_i(x_1, \cdots, x_i) $, $ \partial_{z_k} h_i(z_1, \cdots, z_{i-1}, z, \eps) $ and $ \partial_z h_i(z_1, \cdots, z_{i-1}, z, \eps) $ are all continuous except on a measure $0$ set with respect to the joint law of $ \sfZ, \sfE, (\sfB_i)_{i\in[t]}, (\sfZ_i)_{i\in[t]} $, 
\setcounter{asmpctr}{\value{enumi}}
\end{enumerate}
then 
\begin{align}
    \matrix{\beta_* & v^1 & \cdots & v^t & f^1 & \cdots & f^t} & \xrightarrow{W_2} \matrix{\sfB_* & \sfB_1 & \cdots & \sfB_t & \sfF_1 & \cdots & \sfF_t} , \notag \\
    \matrix{\eps & z & u^1 & \cdots & u^t & h^1 & \cdots & h^t} & \xrightarrow{W_2} \matrix{\sfE & \sfZ & \sfZ_1 & \cdots & \sfZ_t & \sfH_1 & \cdots & \sfH_t} , \notag 
    \notag 
\end{align}
where $ \sfF_i \coloneqq f_i(\sfB_1, \cdots, \sfB_i) $ and $ \sfH_i \coloneqq h_{i}(\sfZ_1, \cdots, \sfZ_{i-1}; \sfZ, \sfE) $. 
\end{proposition}

Consider a special case of the above iteration where 
\begin{align}
&&
    f_t(b^t) &= b^t , &
    h_{t+1}(z^{t}; z, \eps) &= \ol{\kappa}_2^{-1} g_{t+1}(z^{t}; z, \eps) , &
    h_1(z, \eps) &= \ol{\kappa}_2^{-1} g_1(z, \eps) . &
& \label{eqn:fh} 
\end{align} 
\Cref{asmp:div_free_lin} guarantees that functions in \Cref{eqn:fh} satisfy \Cref{asmp:lin_denoiser,asmp:lin_deriv}. 
We can write \Cref{eqn:RIGAMP} as
\begin{subequations}
\label{eqn:RIGAMP_lin} 
\begin{align}
&&
    b^t &= X^\top h^t - \eps^t , &
    &&
& \label{eqn:bb} \\
&&
    z^t &= X b^t - \ol{\kappa}_2 h^t - \delta^t , &
    h^{t+1} &= h_{t+1}(z^{t}; z, \eps) , &
& \label{eqn:zz}
\end{align}
\end{subequations}
where 
\begin{align}
&&
    \eps^t &\coloneqq \sum_{i = 1}^{t-1} \beta_{t,i} b^i , &
    \delta^t &\coloneqq \sum_{i = 1}^{t-1} \alpha_{t,i} h^i + \paren{ \alpha_{t,t} - \ol{\kappa}_2 } h^t . &
& \notag 
\end{align}
Using the choice of $h_{t+1}$ in \Cref{eqn:fh} and combining \Cref{eqn:bb,eqn:zz}, we have
\begin{align}
&&
    b^t &= \ol{\kappa}_2^{-1} X^\top g^t - \eps^t , &
    &&
& \notag \\
&&
    z^t &= \paren{ X X^\top - \ol{\kappa}_2 I_n } h^t - X \eps^t - \delta^t &
    g^{t+1} &= g_{t+1}(z^t; z, \eps) . &
& \notag \\
&&
    &= \paren{ \ol{\kappa}_2^{-1} X X^\top - I_n  } g^t - X \eps^t - \delta^t , &
    &&
& \notag 
\end{align}
Comparing the above iteration with \Cref{eqn:GVAMP_lin}, we recognize that 
\begin{align}
&&
    v^t &= b^t + \eps^t , & 
    u^t &= z^t + X \eps^t + \delta^t . & 
& \label{eqn:eps_zeta}
\end{align}
Furthermore, evaluating the debiasing coefficients $ (\alpha_{t,i})_{i\in[t]} , (\beta_{t,i})_{i\in[t-1]} $ according to \cite[Equations (6) to (10)]{VKM}, we have, 
\begin{align}
&&
    \lim_{n\to\infty} \alpha_{t,t} &= \ol{\kappa}_2 , &
    \lim_{n\to\infty} \alpha_{t,i} &= \lim_{n\to\infty} \beta_{t,i} = 0 , &
    &\forall i\in[t-1] . &
& \notag 
\end{align}
This step crucially requires the divergence-free property of $ g_{t+1} $ in \Cref{asmp:div_free_lin} which causes the coefficients of all memory terms in \Cref{eqn:RIGAMP} to vanish and therefore makes the resulting state evolution recursion more tractable. 
Since $ \lim_{n\to\infty} \normtwo{X} < \infty $ by \Cref{asmp:design}, 
\begin{align}
\lim_{n\to\infty} n^{-1/2} \normtwo{X \eps^t} 
= \lim_{d\to\infty} d^{-1/2} \normtwo{\eps^t} 
&= \lim_{n\to\infty} n^{-1/2} \normtwo{\delta^t} 
= 0 , \label{eqn:eps_delta_0} 
\end{align}
hence it also holds that $ \lim_{n\to\infty} n^{-1/2} \normtwo{X \eps^t + \delta^t} = 0 $. 
We have shown that the iteration \Cref{eqn:GVAMP_lin} can be cast as an instance of \Cref{eqn:RIGAMP} up to vanishing perturbations in the iterates $ b^t, z^t $. 
We claim that the same state evolution result in \cite[Theorem 1]{VKM} applies to both iterations \Cref{eqn:GVAMP_lin,eqn:RIGAMP}. 
To see this, take any test function $ f \in \PL^2_{(t+1)\to1} $. 
Then there exists $C>0$ such that
\begin{align}
    & \abs{\frac{1}{d} \sum_{i = 1}^d f(\beta_{*,i}, v^1_i, \cdots, v^t_i) - \frac{1}{d} \sum_{i = 1}^d f(\beta_{*,i}, b^1_i, \cdots, b^t_i)} 
    \le \frac{1}{d} \sum_{i = 1}^d \abs{f(\beta_{*,i}, v^1_i, \cdots, v^t_i) - f(\beta_{*,i}, b^1_i, \cdots, b^t_i)} \label{eqn:CS1} \\
    &\le \frac{C}{d} \sum_{i = 1}^d \normtwo{(v^j_i)_{j\in[t]} - (b^j_i)_{j\in[t]}} \paren{1 + \normtwo{(v^j_i)_{j\in[t]}} + \normtwo{(b^j_i)_{j\in[t]}}} \label{eqn:use_PL} \\
    &= \frac{C}{d} \sum_{i = 1}^d \normtwo{(\eps^j_i)_{j\in[t]}} \paren{1 + \normtwo{(v^j_i)_{j\in[t]}} + \normtwo{(b^j_i)_{j\in[t]}}} \label{eqn:use_eps} \\
    &\le \frac{C}{d} \sqrt{\sum_{i = 1}^d \normtwo{(\eps_i^j)_{j\in[t]}}^2} \sqrt{\sum_{i = 1}^d \paren{1 + \normtwo{(v^j_i)_{j\in[t]}} + \normtwo{(b^j_i)_{j\in[t]}}}^2} \label{eqn:CS2} \\
    &\le \frac{\sqrt{3}C}{d} \sqrt{\sum_{i = 1}^d \normtwo{(\eps_i^j)_{j\in[t]}}^2} \sqrt{\sum_{i = 1}^d \paren{1 + \normtwo{(v^j_i)_{j\in[t]}}^2 + \normtwo{(b^j_i)_{j\in[t]}}^2}} \label{eqn:abc_ineq} \\
    &= \sqrt{3}C \sqrt{\sum_{j = 1}^t d^{-1} \normtwo{\eps^j}^2} \sqrt{1 + \sum_{j = 1}^t d^{-1} \normtwo{v^j}^2 + \sum_{j = 1}^t d^{-1} \normtwo{b^j}^2} \notag \\
    &\le \sqrt{3}C \sqrt{\sum_{j = 1}^t d^{-1} \normtwo{\eps^j}^2} \sqrt{1 + 2 \sum_{j = 1}^t d^{-1} \normtwo{\eps^j}^2 + 3 \sum_{j = 1}^t d^{-1} \normtwo{b^j}^2} . \label{eqn:ab_ineq}
\end{align}
In the above chain of (in)equalities, 
\Cref{eqn:use_PL} follows from the $ \PL^2 $ property of $f$ (see \Cref{def:pseudo_lip}); 
\Cref{eqn:use_eps} is by \Cref{eqn:eps_zeta}; 
\Cref{eqn:CS1,eqn:CS2} are both by Cauchy--Schwarz; 
\Cref{eqn:abc_ineq} uses the elementary inequality $ (a+b+c)^2 \le 3 (a^2 + b^2 + c^2) $ for any $ a,b,c \in \bbR $; 
\Cref{eqn:ab_ineq} follows since $ \normtwo{u+v}^2 \le 2(\normtwo{u}^2 + \normtwo{v}^2) $. 
The state evolution result \cite[Theorem 1]{VKM} applies to $ (b^j)_{j\in[t]} $ and implies in particular $ \lim_{d\to\infty} d^{-1/2} \normtwo{b^j} < \infty $ for every $j\in[t]$. 
Since $ t $ is fixed, this combined with \Cref{eqn:eps_delta_0} shows 
\begin{align}
    \lim_{d\to\infty} \abs{\frac{1}{d} \sum_{i = 1}^d f(\beta_{*,i}, v^1_i, \cdots, v^t_i) - \frac{1}{d} \sum_{i = 1}^d f(\beta_{*,i}, b^1_i, \cdots, b^t_i)}
    &= 0 . \notag 
\end{align}
A similar argument shows that
\begin{align}
    \lim_{n\to\infty} \abs{\frac{1}{n} \sum_{i = 1}^n g(\eps_i, z_i, u^1_i, \cdots, u^t_i, g^1_i, \cdots, g^t_i) - \frac{1}{n} \sum_{i = 1}^n g(\eps_i, z_i, z^1_i, \cdots, z^t_i, \ol{\kappa}_2 h^1_i, \cdots, \ol{\kappa}_2 h^t_i)}
    &= 0 , \notag 
\end{align}
for any $ g \in \PL^2_{2(t+1)\to1} $. 
Therefore, the Wasserstein-$2$ limits of $ \matrix{\beta_* & v^1 & \cdots & v^t} $ and $  \matrix{\eps & z & u^1 & \cdots & u^t & g^1 & \cdots g^t} $ are the same as those of $ \matrix{\beta_* & b^1 & \cdots & b^t} $ and $  \matrix{\eps & z & z^1 & \cdots & z^t & \ol{\kappa}_2 h^1 & \cdots & \ol{\kappa}_2 h^t} $, respectively, to the latter of which the state evolution result \cite[Theorem 1]{VKM} applies. 

Next, we specialize the state evolution recursion in \cite[Equations (11) to (22)]{VKM} to \Cref{eqn:RIGAMP_lin} which is also valid for \Cref{eqn:GVAMP_lin}. 
Define random variables: 
\begin{subequations}
\label{eqn:SE_lin_pf}
\begin{align}
&\matrix{\sfE & \sfZ & \sfZ_1 & \cdots & \sfZ_{t-1} }^\top  \sim P_{\sfE} \ot \cN(0_t, \Sigma_t) , \\
&\, \wt{\sfZ}_t = h_t(\sfZ_{t-1}; \sfZ, \sfE) , \\
&\matrix{ \sfB_* & \sfW_1 & \cdots & \sfW_t }^\top \sim P_{\sfB_*} \ot \cN(0_t, \Omega_t) , \\
&\matrix{ \sfB_1 & \cdots & \sfB_t }^\top = \matrix{\chi_1 & \cdots & \chi_t}^\top \sfB_* + \matrix{ \sfW_1 & \cdots & \sfW_t }^\top . 
\end{align}
The state evolution parameters $ \Sigma_t\in\bbR^{t\times t}, \Omega_t\in\bbR^{t\times t}, \chi_t \in \bbR $ are defined recursively as follows. 
For $ t=1 $, 
\begin{align}
&&
    \Sigma_1 &= \sigma^2 , &
    \Omega_1 &= \delta \ol{\kappa}_2 \expt{h_1(\sfZ, \sfE)^2} + \delta \ol{\kappa}_4 \rho \expt{\partial_z h_1(\sfZ, \sfE)}^2 , &
    \mu_1 &= \delta \ol{\kappa}_2 \expt{\partial_z h_1(\sfZ, \sfE)} , &
& 
\end{align}
where $ (\sfZ, \sfE) \sim \cN\paren{ 0, \sigma^2 } \ot P_{\sfE} $.

For $t\ge1$, we now describe how $ \Sigma_{t+1}, \Omega_{t+1}, \chi_{t+1} $ can be computed from $ \Sigma_t, \Omega_t, (\chi_i)_{i\in[t]} $. 
Define
\begin{align}
\phi_t &\coloneqq \matrix{ \expt{\partial_z h_1(\sfZ, \sfE)} & \expt{\partial_z h_2(\sfZ_1; \sfZ, \sfE)} & \cdots & \expt{\partial_z h_t(\sfZ_{t-1}; \sfZ, \sfE)} }^\top \in \bbR^t , \notag \\
\gamma_t &\coloneqq \matrix{ \expt{\sfB_* \sfB_1} & \cdots & \expt{\sfB_* \sfB_t} }^\top \in \bbR^t , \notag \\
\Gamma_t &\coloneqq \matrix{
    \expt{\sfB_i \sfB_j}
}_{(i,j)\in[t]^2} \in \bbR^{t\times t} , \notag \\
\Delta_t &\coloneqq \matrix{ \expt{\wt{\sfZ}_i \wt{\sfZ}_j} }_{(i,j)\in[t]^2} \in \bbR^{t\times t} . \notag 
\end{align}
Then $ \Sigma_{t+1} $ is updated as: 
\begin{align}
\Sigma_{t+1} &= \ol{\kappa}_2 \matrix{
    \rho & \gamma_t^\top \\
    \gamma_t & \Gamma_t
}
+ \ol{\kappa}_4 \paren{ \matrix{
    0 & \rho \phi_t^\top \\
    \rho \phi_t & \phi_t \gamma_t^\top + \gamma_t \phi_t^\top 
} + \matrix{
    0 & 0_t^\top \\
    0_t & \Delta_t
} }
+ \ol{\kappa}_6  \matrix{
    0 & 0_t^\top \\
    0_t & \rho \phi_t \phi_t^\top
} . \notag 
\end{align}
The matrix $ \Omega_{t+1} $ is updated as: 
\begin{align}
\Omega_{t+1} &\coloneqq \delta \paren{
    \ol{\kappa}_2 \Delta_{t+1}
    + \ol{\kappa}_4 \rho \phi_{t+1} \phi_{t+1}^\top
} . \notag 
\end{align}
The scalar $ \chi_{t+1} $ is updated as: 
\begin{align}
\chi_{t+1} &= \delta \ol{\kappa}_2 \phi_{t+1,t+1} , \notag 
\end{align}
where $ \phi_{t+1,t+1} $ denotes the last element of the vector $ \phi_{t+1} \in \bbR^{t+1} $. 

Note that $ \Sigma_t, \Omega_t $ are the leading principal minors (of order $t$) of $ \Sigma_{t+1}, \Omega_{t+1} $, respectively.
So we explicitly write the update rules for the last row and column of $ \Sigma_{t+1}, \Omega_{t+1} $ and $ \chi_{t+1} $. 
For $1\le i\le t$ and $0\le j\le t$, 
\begin{align}
\Sigma_{t+1,1} 
&= \ol{\kappa}_2 \expt{\sfB_* \sfB_t} + \ol{\kappa}_4 \rho \expt{\partial_z h_t(\sfZ_{t-1}; \sfZ, \sfE)} 
= \ol{\kappa}_2 \rho \chi_t + \frac{\ol{\kappa}_4}{\ol{\kappa}_2} \rho \expt{\partial_z g_t(\sfZ_{t-1}; \sfZ, \sfE)} , 
\\
\Sigma_{t+1,i+1}
&= \ol{\kappa}_2 \expt{\sfB_t \sfB_{i}} + \ol{\kappa}_4 \paren{\expt{\partial_z h_t(\sfZ_{t-1}; \sfZ, \sfE)} \expt{\sfB_* \sfB_{i}} + \expt{\partial_z h_{i}(\sfZ_{i-1}; \sfZ, \sfE)} \expt{\sfB_* \sfB_t}} \notag \\
&\quad + \ol{\kappa}_4 \expt{\wt{\sfZ}_t \wt{\sfZ}_{i}} + \ol{\kappa}_6 \rho \expt{\partial_z h_{i}(\sfZ_{i-1}; \sfZ, \sfE)} \expt{\partial_z h_t(\sfZ_{t-1}; \sfZ, \sfE)} \notag \\
&= \ol{\kappa}_2 \paren{\chi_t \chi_{i} \rho + \Omega_{t,i}} + \frac{\ol{\kappa}_4}{\ol{\kappa}_2} \paren{\expt{\partial_z g_t(\sfZ_{t-1}; \sfZ, \sfE)} \rho \chi_{i} + \expt{\partial_z g_{i}(\sfZ_{i-1}; \sfZ, \sfE)} \rho \chi_t} \notag \\
&\quad + \frac{\ol{\kappa}_4}{\ol{\kappa}_2^2} \expt{g_t(\sfZ_{t-1}; \sfZ, \sfE) g_{i}(\sfZ_{i-1}; \sfZ, \sfE)} + \frac{\ol{\kappa}_6}{\ol{\kappa}_2^2} \rho \expt{\partial_z g_{i}(\sfZ_{i-1}; \sfZ, \sfE)} \expt{\partial_z g_t(\sfZ_{t-1}; \sfZ, \sfE)} , 
\\
\Omega_{t+1,j+1} 
&= \delta \paren{\ol{\kappa}_2 \expt{\wt{\sfZ}_{t+1} \wt{\sfZ}_{j+1}} + \ol{\kappa}_4 \rho \expt{\partial_z h_{t+1}(\sfZ_t; \sfZ, \sfE)} \expt{\partial_z h_{j+1}(\sfZ_{j}; \sfZ, \sfE)}} \notag \\
&= \delta \paren{\frac{1}{\ol{\kappa}_2} \expt{g_{t+1}(\sfZ_t; \sfZ, \sfE) g_{j+1}(\sfZ_{j}; \sfZ, \sfE)} + \frac{\ol{\kappa}_4}{\ol{\kappa}_2^2} \rho \expt{\partial_z g_{t+1}(\sfZ_t; \sfZ, \sfE)} \expt{\partial_z g_{j+1}(\sfZ_{j}; \sfZ, \sfE)}} , 
\\
\chi_{t+1} &= \delta \ol{\kappa}_2 \expt{\partial_z h_{t+1}(\sfZ_t; \sfZ, \sfE)} 
= \delta \expt{\partial_z g_{t+1}(\sfZ_t; \sfZ, \sfE)} , 
\end{align}
\end{subequations}
where we recall the choice of $ h_t $ in \Cref{eqn:fh}. 

We conclude the proof by noting that the state evolution recursion \Cref{eqn:SE_lin_pf} is the same as \Cref{eqn:SE_lin}. 

\subsection{Proof of \Cref{prop:BZ_equiv}}
\label{app:pf_prop:BZ_equiv}

According to \Cref{eqn:SE_lin}, $ (\sfZ, \sfZ_1, \cdots ,\sfZ_t) $ is a centered Gaussian vector. 
So we can equivalently write 
\begin{align}
    \sfZ_i &= \nu_i \sfZ + \tau_i \sfN_i , \qquad i\in[t] \notag 
\end{align}
where $ (\sfN_i)_{i\in[t]} \sim \cN(0_t, I_t) $ is a standard Gaussian vector independent of $ \sfZ $. 
The coefficients $ (\nu_i)_{i\in[t]}, (\tau_i)_{i\in[t]} $ and the correlations $ (\expt{\tau_r \sfN_r \cdot \tau_s \sfN_s})_{(r,s)\in[t]^2} $ can be derived from the covariance matrix $ \Sigma_{t+1} $. 

To this end, we first note that the covariance structure of $ (\sfZ , \sfZ_t) $ is given by 
\begin{align}
    \expt{\sfZ^2} &= \sigma^2 , \notag \\
    \expt{\sfZ \sfZ_t} &= \ol{\kappa}_2 \rho \chi_t + \frac{\ol{\kappa}_4}{\ol{\kappa}_2} \rho \expt{\partial_z g_t(\sfZ_{t-1}; \sfZ, \sfE)} , \notag \\
    \expt{\sfZ_t^2} &= \ol{\kappa}_2 \paren{\chi_t^2 \rho + \Omega_{t,t}} + 2 \frac{\ol{\kappa}_4}{\ol{\kappa}_2} \expt{\partial_z g_t(\sfZ_{t-1}; \sfZ, \sfE)} \rho \chi_t + \frac{\ol{\kappa}_4}{\ol{\kappa}_2^2} \expt{g_t(\sfZ_{t-1}; \sfZ, \sfE)^2} + \frac{\ol{\kappa}_6}{\ol{\kappa}_2^2} \rho \expt{\partial_z g_t(\sfZ_{t-1}; \sfZ, \sfE)}^2 . \notag 
\end{align}
The coefficients $ \nu_t, \tau_t $ are therefore given by
\begin{align}
\nu_t &= \frac{\expt{\sfZ \sfZ_t}}{\sigma^2} 
= \chi_t + \frac{\ol{\kappa}_4}{\ol{\kappa}_2^2} \expt{\partial_z g_t(\sfZ_{t-1}; \sfZ, \sfE)} , \notag \\
\tau_t^2 &= \expt{\sfZ_t^2} - \nu_t^2 \sigma^2 
= \ol{\kappa}_2 \paren{\chi_t^2 \rho + \Omega_{t,t}} + 2 \frac{\ol{\kappa}_4}{\ol{\kappa}_2} \expt{\partial_z g_t(\sfZ_{t-1}; \sfZ, \sfE)} \rho \chi_t + \frac{\ol{\kappa}_4}{\ol{\kappa}_2^2} \expt{g_t(\sfZ_{t-1}; \sfZ, \sfE)^2} + \frac{\ol{\kappa}_6}{\ol{\kappa}_2^2} \rho \expt{\partial_z g_t(\sfZ_{t-1}; \sfZ, \sfE)}^2 \notag \\
&\quad - \paren{\chi_t + \frac{\ol{\kappa}_4}{\ol{\kappa}_2^2} \expt{\partial_z g_t(\sfZ_{t-1}; \sfZ, \sfE)}}^2 \sigma^2, \notag \\
&= \ol{\kappa}_2 \Omega_{t,t} + \frac{\ol{\kappa}_4}{\ol{\kappa}_2^2} \expt{g_t(\sfZ_{t-1}; \sfZ, \sfE)^2} + \paren{\frac{\ol{\kappa}_6}{\kappa_2^2} - \frac{\ol{\kappa}_4^2}{\ol{\kappa}_2^3}} \rho \expt{\partial_z g_t(\sfZ_{t-1}; \sfZ, \sfE)}^2 . \label{eqn:tau2} 
\end{align}

Similarly, the correlation $ \expt{\tau_r \sfN_r \cdot \tau_s \sfN_s} $ is given by 
\begin{align}
    & \expt{\tau_r \sfN_r \cdot \tau_s \sfN_s}
    = \Sigma_{r+1, s+1} - \nu_r \nu_s \sigma^2 \notag \\
    &= \ol{\kappa}_2 \paren{\chi_r \chi_s \rho + \Omega_{r,s}} + \frac{\ol{\kappa}_4}{\ol{\kappa}_2} \paren{\expt{\partial_z g_r(\sfZ_{r-1}; \sfZ, \sfE)} \rho \chi_{s} + \expt{\partial_z g_{s}(\sfZ_{s-1}; \sfZ, \sfE)} \rho \chi_r} \notag \\
    &\quad + \frac{\ol{\kappa}_4}{\ol{\kappa}_2^2} \expt{g_r(\sfZ_{r-1}; \sfZ, \sfE) g_{s}(\sfZ_{s-1}; \sfZ, \sfE)} + \frac{\ol{\kappa}_6}{\ol{\kappa}_2^2} \rho \expt{\partial_z g_{s}(\sfZ_{s-1}; \sfZ, \sfE)} \expt{\partial_z g_r(\sfZ_{r-1}; \sfZ, \sfE)} \notag \\
    &\quad - \paren{\chi_r + \frac{\ol{\kappa}_4}{\ol{\kappa}_2^2} \expt{\partial_z g_r(\sfZ_{r-1}; \sfZ, \sfE)}} \paren{\chi_s + \frac{\ol{\kappa}_4}{\ol{\kappa}_2^2} \expt{\partial_z g_s(\sfZ_{s-1}; \sfZ, \sfE)}} \sigma^2 \notag \\
    &= \ol{\kappa}_2 \Omega_{r,s} 
    + \frac{\ol{\kappa}_4}{\ol{\kappa}_2^2} \expt{g_r(\sfZ_{r-1}; \sfZ, \sfE) g_{s}(\sfZ_{s-1}; \sfZ, \sfE)} + \paren{\frac{\ol{\kappa}_6}{\ol{\kappa}_2^2} - \frac{\ol{\kappa}_4^2}{\ol{\kappa}_2^3}} \rho \expt{\partial_z g_{r}(\sfZ_{r-1}; \sfZ, \sfE)} \expt{\partial_z g_s(\sfZ_{s-1}; \sfZ, \sfE)} . \notag 
\end{align}
Note that for $ r=s=t $, the above expression coincides with $ \tau_t^2 $ in \Cref{eqn:tau2}. 

Turning to $ \sfB_{t+1} $, the representation $ \sfB_{t+1} = \chi_{t+1} \sfB_* + \omega_{t+1} \sfM_{t+1} $ (where $ (\sfB_*, \sfM_{t+1}) \sim P_{\sfB_*} \ot \cN(0,1) $) can be obtained by simply redefining $ \sfM_{t+1} = \sfW_{t+1} / \sqrt{\Omega_{t+1,t+1}} $ and $ \omega_{t+1} = \sqrt{\Omega_{t+1,t+1}} $ where the right-hand sides are given in \Cref{eqn:SE_lin}. 

\subsection{Proof of \Cref{prop:SE_lin}}
\label{app:pf_prop:SE_lin}

First off, let us use Stein's lemma and chain rule of derivatives to derive an alternative expression for $ \expt{\partial_z g_t(\sfZ_{t-1}; \sfZ, \sfE)} $ (which is used in \Cref{prop:BZ_equiv}): 
\begin{align}
    \expt{\partial_z g_t(\sfZ_{t-1}; \sfZ, \sfE)} &= \sigma^{-2} \expt{\sfZ g_t(\sfZ_{t-1}; \sfZ, \sfE)} - \nu_{t-1} \expt{g_t'(\sfZ_{t-1}; \sfZ, \sfE)} . \notag 
\end{align}
When $ g_t $ takes the form of \Cref{eqn:g_choice}, this identity becomes
\begin{align}
    \expt{\partial_z g_t(\sfZ_{t-1}; \sfZ, \sfE)} &= \sigma^{-2} \expt{\sfZ \ol{g}_t(\sfY) \sfZ_{t-1}} - \nu_{t-1} \expt{\ol{g}_t(\sfY)} 
    = \sigma^{-2} \nu_{t-1} \expt{\ol{g}_t(\sfY) \sfZ^2} 
    = \nu_{t-1} \expt{\ol{g}_t(\sfY) \ol{\sfZ}^2} , \label{eqn:g_stein} 
\end{align}
where we recall $ \ol{\sfZ} = \sigma^{-1} \sfZ $ defined in \Cref{sec:model}. 

Now, specializing \Cref{prop:BZ_equiv} to the case of \Cref{eqn:g_choice}, we have 
\begin{align}
\nu_t &= \paren{ 1 + \frac{\ol{\kappa}_4}{\delta \ol{\kappa}_2^2} } \chi_t
, \label{eqn:nu_chi} \\
\tau_t^2 &= \ol{\kappa}_2 \omega_t^2 + \frac{\ol{\kappa}_4}{\ol{\kappa}_2^2} \expt{\ol{g}_t(\sfY)^2 \sfZ_{t-1}^2} + \paren{\frac{\ol{\kappa}_6}{\kappa_2^2} - \frac{\ol{\kappa}_4^2}{\ol{\kappa}_2^3}} \rho \expt{\partial_z g_t(\sfZ_{t-1}; \sfZ, \sfE)}^2 \notag \\
&= \ol{\kappa}_2 \omega_t^2 + \frac{\ol{\kappa}_4}{\ol{\kappa}_2^2} \paren{ \nu_{t-1}^2 \expt{\ol{g}_t(\sfY)^2 \sfZ^2} + \tau_{t-1}^2 \expt{\ol{g}_t(\sfY)^2} } 
+ \paren{\frac{\ol{\kappa}_6}{\kappa_2^2} - \frac{\ol{\kappa}_4^2}{\ol{\kappa}_2^3}} \rho \expt{\ol{g}_t(\sfY) \ol{\sfZ}^2}^2 \nu_{t-1}^2 \label{eqn:stein1} \\
&= \ol{\kappa}_2 \omega_t^2 + \frac{\ol{\kappa}_4}{\ol{\kappa}_2^2} \expt{\ol{g}_t(\sfY)^2} \tau_{t-1}^2 
+ \rho \brace{ \frac{\ol{\kappa}_4}{\ol{\kappa}_2} \expt{\ol{g}_t(\sfY)^2 \ol{\sfZ}^2} + \paren{\frac{\ol{\kappa}_6}{\kappa_2^2} - \frac{\ol{\kappa}_4^2}{\ol{\kappa}_2^3}} \expt{\ol{g}_t(\sfY) \ol{\sfZ}^2}^2 } \nu_{t-1}^2 , \label{eqn:tau2_recur}\\
\chi_{t+1} &= \delta \expt{\ol{g}_{t+1}(\sfY) \ol{\sfZ}^2} \nu_t , \label{eqn:stein2}\\
\omega_{t+1}^2 &= \frac{\delta}{\ol{\kappa}_2} \expt{\ol{g}_{t+1}(\sfY)^2 \sfZ_t^2} + \delta \frac{\ol{\kappa}_4}{\ol{\kappa}_2^2} \rho \expt{\partial_z g_{t+1}(\sfZ_t; \sfZ, \sfE)}^2 \notag \\
&= \frac{\delta}{\ol{\kappa}_2} \paren{ \nu_t^2 \expt{\ol{g}_{t+1}(\sfY)^2 \sfZ^2} + \tau_t^2 \expt{\ol{g}_{t+1}(\sfY)^2} }
+ \delta \frac{\ol{\kappa}_4}{\ol{\kappa}_2^2} \rho \expt{\ol{g}_{t+1}(\sfY) \ol{\sfZ}^2}^2 \nu_t^2 \label{eqn:stein3} \\
&= \delta \ol{\kappa}_2^{-1} \expt{\ol{g}_{t+1}(\sfY)^2} \tau_t^2
+ \delta \rho \brace{ \expt{\ol{g}_{t+1}(\sfY)^2 \ol{\sfZ}^2} + \frac{\ol{\kappa}_4}{\ol{\kappa}_2^2} \expt{\ol{g}_{t+1}(\sfY) \ol{\sfZ}^2}^2 } \nu_t^2 . \label{eqn:omega_tau} 
\end{align}
The equalities \Cref{eqn:stein1,eqn:stein2,eqn:stein3} all follow from \Cref{eqn:g_stein}. 

By the assumption $ \chi_1 \ne 0 $, the recursion \Cref{eqn:nu_chi} and the observation \Cref{eqn:k4-by-k2+delta}, we have $ \nu_1 \ne 0 $. 
Eliminating $ \chi_t $ from \Cref{eqn:nu_chi,eqn:stein2}, we obtain the following recursion for $ \nu_t $: 
\begin{align}
    \nu_t &= \paren{ \frac{\ol{\kappa}_4}{\ol{\kappa}_2^2} + \delta } \expt{\ol{g}_t(\sfY) \ol{\sfZ}^2} \nu_{t-1} . \label{eqn:nu_nu} 
\end{align}
Using \Cref{eqn:k4-by-k2+delta} along with \Cref{asmp:g(y)} in \Cref{eqn:nu_nu} and inducting on $t$, we further obtain $ \nu_t \ne 0 $ for every $t\ge1$. 
Using \Cref{asmp:g(y)} again in \Cref{eqn:stein2}, this in turn implies $ \chi_{t+1} \ne 0 $ for every $t\ge1$. 

Next, let us recall the expression of $ \omega_1 $ from \Cref{eqn:chi1_omega1} and make two observations. 
Stein's lemma implies
\begin{align}
    \expt{\partial_z g_1(\sfZ, \sfE)} &= \sigma^{-2} \expt{\sfZ g_1(\sfZ, \sfE)}
    = \sigma^{-1} \expt{\ol{\sfZ} g_1(\sfZ, \sfE)} . \label{eqn:g1_stein}
\end{align}
Cauchy--Schwarz inequality implies
\begin{align}
    \expt{g_1(\sfZ, \sfE)^2} &\ge \expt{\ol{\sfZ} g_1(\sfZ, \sfE)}^2 . \label{eqn:g1_CS}
\end{align}
By the above two observations and recalling the definition \Cref{eqn:sigma2} of $ \sigma^2 $, 
\begin{align}
    \omega_1^2 &\ge \frac{\delta}{\ol{\kappa}_2} \sigma^2 \expt{\partial_z g_1(\sfZ, \sfE)}^2 + \delta \frac{\ol{\kappa}_4}{\ol{\kappa}_2^2} \rho \expt{\partial_z g_1(\sfZ, \sfE)}^2 
    = \delta \rho \expt{\partial_z g_1(\sfZ, \sfE)}^2 \paren{ 1 + \frac{\ol{\kappa}_4}{\ol{\kappa}_2^2} } . \label{eqn:omega1_factor}
\end{align}
By $ \chi_1 \ne 0 $ and \Cref{asmp:signal}, the first factor of \Cref{eqn:omega1_factor} is positive. 
We claim that the second factor surrounded by parentheses in \Cref{eqn:omega1_factor} is also positive. 
By the moment-cumulant relation (see \Cref{sec:cumulant}), 
\begin{align}
    1 + \frac{\ol{\kappa}_4}{\ol{\kappa}_2^2}
    &= \frac{\ol{m}_4}{\ol{m}_2^2} - 1 - \delta . \notag 
\end{align}
By definition, $ \ol{m}_{2k} = \expt{\sfLambda_n^{2k}} $ where $ \sfLambda_n $ is given in \Cref{eqn:sfLambda_nd}. 
Therefore 
\begin{align}
    1 + \frac{\ol{\kappa}_4}{\ol{\kappa}_2^2}
    &= \frac{\expt{\sfLambda_n^4}}{\expt{\sfLambda_n^2}^2} - \delta
    = \begin{cases}
        \frac{\expt{\sfLambda^4}}{\expt{\sfLambda^2}^2} - \delta , & \delta \le 1 \\ 
        \delta \frac{\expt{\sfLambda^4}}{\expt{\sfLambda^2}^2} - \delta , & \delta > 1
    \end{cases} . \label{eqn:1+k4-by-k2} 
\end{align}
In both cases, the right-hand side is positive, by Cauchy--Schwarz and \Cref{asmp:design}. 
This confirms the positivity of $ \omega_1 $. 

We then show $ \tau_1^2 > 0 $. 
By \Cref{prop:BZ_equiv}, 
\begin{align}
    \tau_1^2 &= \ol{\kappa}_2 \omega_1^2 + \frac{\ol{\kappa}_4}{\ol{\kappa}_2^2} \expt{g_1(\sfZ, \sfE)^2} + \paren{\frac{\ol{\kappa}_6}{\kappa_2^2} - \frac{\ol{\kappa}_4^2}{\ol{\kappa}_2^3}} \rho \expt{\partial_z g_1(\sfZ, \sfE)}^2 \notag \\
    &= \paren{\delta + \frac{\ol{\kappa}_4}{\ol{\kappa}_2^2}} \expt{g_1(\sfZ, \sfE)^2} + \paren{\delta \frac{\ol{\kappa}_4}{\ol{\kappa}_2} + \frac{\ol{\kappa}_6}{\kappa_2^2} - \frac{\ol{\kappa}_4^2}{\ol{\kappa}_2^3}} \rho \expt{\partial_z g_1(\sfZ, \sfE)}^2 , \label{eqn:tau2_TODO}
\end{align}
where the second step uses the expression of $ \omega_1 $. 
Using \Cref{eqn:g1_stein,eqn:g1_CS} in \Cref{eqn:tau2_TODO}, we obtain
\begin{align}
    \tau_1^2 &\ge \paren{\delta + \frac{\ol{\kappa}_4}{\ol{\kappa}_2^2} + \delta \frac{\ol{\kappa}_4}{\ol{\kappa}_2^2} + \frac{\ol{\kappa}_6}{\kappa_2^3} - \frac{\ol{\kappa}_4^2}{\ol{\kappa}_2^4}} \expt{\ol{\sfZ} g_1(\sfZ, \sfE)}^2 . \notag 
\end{align}
By the moment-cumulant relation (see \Cref{sec:cumulant}), 
\begin{align}
    \delta + \frac{\ol{\kappa}_4}{\ol{\kappa}_2^2} + \delta \frac{\ol{\kappa}_4}{\ol{\kappa}_2^2} + \frac{\ol{\kappa}_6}{\kappa_2^3} - \frac{\ol{\kappa}_4^2}{\ol{\kappa}_2^4}
    &= \frac{\ol{m}_2 \ol{m}_6 - \ol{m}_4^2}{\ol{m}_2^4} . \notag 
\end{align}
By the definition of $ \ol{m}_{2k} $ (see \Cref{sec:cumulant}), the numerator of the right-hand side above is
\begin{align}
    \ol{m}_2 \ol{m}_6 - \ol{m}_4^2
    &= \expt{\sfLambda_n^2} \expt{\sfLambda_n^6} - \expt{\sfLambda_n^4}^2 > 0 , \label{eqn:m_positive} 
\end{align}
by \holder's inequality applied to the non-negative random variable $\sfLambda_n \cdot \sfLambda_n^3$ (see \Cref{asmp:design}). 
Therefore $ \tau_1^2>0 $. 

Eliminating $ \omega_t $ from \Cref{eqn:tau2_recur,eqn:omega_tau}, we obtain the following recursion for $ \tau_t $: 
\begin{align}
\tau_t^2 &= 
\paren{ \frac{\ol{\kappa}_4}{\ol{\kappa}_2^2} + \delta } \expt{\ol{g}_t(\sfY)^2} \tau_{t-1}^2 
+ \rho \brace{ \paren{ \frac{\ol{\kappa}_4}{\ol{\kappa}_2} + \delta \ol{\kappa}_2 } \expt{\ol{g}_t(\sfY)^2 \ol{\sfZ}^2} + \paren{\frac{\ol{\kappa}_6}{\ol{\kappa}_2^2} - \frac{\ol{\kappa}_4^2}{\ol{\kappa}_2^3} + \delta \frac{\ol{\kappa}_4}{\ol{\kappa}_2}} \expt{\ol{g}_t(\sfY) \ol{\sfZ}^2}^2 } \nu_{t-1}^2 . \label{eqn:tau_tau}
\end{align}
Note that 
\begin{align}
    \expt{\ol{g}_t(\sfY)^2} &> 0 
    \label{eqn:g>0}
\end{align}
since otherwise by \Cref{asmp:g(y)}, Cauchy--Schwarz and $ \ol{\sfZ} \sim \cN(0,1) $, 
\begin{align}
    0 < \expt{\ol{g}_{t}(\sfY) \ol{\sfZ}^2} 
    &\le \sqrt{\expt{\ol{g}_{t}(\sfY)^2}} \cdot \sqrt{3} = 0 , \notag 
\end{align}
which is a contradiction. 
We claim that the multiplicative coefficient in front of $ \nu_{t-1}^2 $ on the right of \Cref{eqn:tau_tau} is also positive. 
Indeed, by Cauchy--Schwarz, 
\begin{align}
    \expt{\ol{g}_t(\sfY) \ol{\sfZ}^2}^2 &\le \expt{\ol{g}_t(\sfY)^2 \ol{\sfZ}^2} \expt{\ol{\sfZ}^2} = \expt{\ol{g}_t(\sfY)^2 \ol{\sfZ}^2} , \label{eqn:CS_gZ2} 
\end{align}
which can be used to lower bound the coefficient of interest as
\begin{align}
    \paren{ \frac{\ol{\kappa}_4}{\ol{\kappa}_2} + \delta \ol{\kappa}_2 } \expt{\ol{g}_t(\sfY)^2 \ol{\sfZ}^2} + \paren{\frac{\ol{\kappa}_6}{\ol{\kappa}_2^2} - \frac{\ol{\kappa}_4^2}{\ol{\kappa}_2^3} + \delta \frac{\ol{\kappa}_4}{\ol{\kappa}_2}} \expt{\ol{g}_t(\sfY) \ol{\sfZ}^2}^2 
    &\ge \paren{\frac{\ol{\kappa}_4}{\ol{\kappa}_2} + \delta \ol{\kappa}_2 + \frac{\ol{\kappa}_6}{\ol{\kappa}_2^2} - \frac{\ol{\kappa}_4^2}{\ol{\kappa}_2^3} + \delta \frac{\ol{\kappa}_4}{\ol{\kappa}_2}} \expt{\ol{g}_t(\sfY) \ol{\sfZ}^2}^2 . \label{eqn:coeff_pos} 
\end{align}
The second factor on the right is positive by \Cref{asmp:g(y)} and the first factor is also positive, as we have already seen in \Cref{eqn:m_positive}. 
Combining the above results with \Cref{eqn:k4-by-k2+delta}, recalling the positivity of $ (\nu_t^2)_{t\ge1} $ shown below \Cref{eqn:nu_nu} and inducting on $t$ in \Cref{eqn:tau_tau}, we conclude the positivity of $ \tau_t $ for every $t\ge1$. 

Moving to $\omega_t$, we consider the recursion \Cref{eqn:omega_tau}. 
The coefficient in front of $\nu_t^2$ on the right of \Cref{eqn:omega_tau} can be similarly seen positive: 
\begin{align}
    \expt{\ol{g}_{t+1}(\sfY)^2 \ol{\sfZ}^2} + \frac{\ol{\kappa}_4}{\ol{\kappa}_2^2} \expt{\ol{g}_{t+1}(\sfY) \ol{\sfZ}^2}^2
    &\ge \expt{\ol{g}_{t+1}(\sfY) \ol{\sfZ}^2}^2 \paren{1 + \frac{\ol{\kappa}_4}{\ol{\kappa}_2^2}} > 0 , \label{eqn:coeff_pos_omega} 
\end{align}
where we have used \Cref{eqn:CS_gZ2,eqn:1+k4-by-k2}. 
Combining this with \Cref{eqn:g>0} and the positivity of $ (\tau_t^2)_{t\ge1}, (\nu_t^2)_{t\ge1} $ just shown and inducting on $t$, we conclude the positivity of $ \omega_t $ for every $t\ge1$. 

Consider the recursions \Cref{eqn:nu_nu,eqn:tau_tau} for $ \nu_t, \tau_t $. 
Dividing both sides of \Cref{eqn:tau_tau} by $ \nu_t^2 $ and using \Cref{eqn:nu_nu} lead us to a recursion for the ratio $ \tau_t^2 / \nu_{t}^2 $: 
\begin{align}
    \frac{\tau_t^2}{\nu_t^2} &= \paren{ \frac{\ol{\kappa}_4}{\ol{\kappa}_2^2} + \delta } \expt{\ol{g}_t(\sfY)^2} \frac{\tau_{t-1}^2}{\nu_{t-1}^2} \frac{\nu_{t-1}^2}{\nu_t^2} 
    + \rho \brace{ \paren{ \frac{\ol{\kappa}_4}{\ol{\kappa}_2} + \delta \ol{\kappa}_2 } \expt{\ol{g}_t(\sfY)^2 \ol{\sfZ}^2} + \paren{\frac{\ol{\kappa}_6}{\ol{\kappa}_2^2} - \frac{\ol{\kappa}_4^2}{\ol{\kappa}_2^3} + \delta \frac{\ol{\kappa}_4}{\ol{\kappa}_2}} \expt{\ol{g}_t(\sfY) \ol{\sfZ}^2}^2 } \frac{\nu_{t-1}^2}{\nu_t^2} \notag \\
    &= \frac{\paren{ \frac{\ol{\kappa}_4}{\ol{\kappa}_2^2} + \delta } \expt{\ol{g}_t(\sfY)^2}}{\paren{ \frac{\ol{\kappa}_4}{\ol{\kappa}_2^2} + \delta }^2 \expt{\ol{g}_t(\sfY) \ol{\sfZ}^2}^2} \frac{\tau_{t-1}^2}{\nu_{t-1}^2}
    + \frac{\rho \brace{ \paren{ \frac{\ol{\kappa}_4}{\ol{\kappa}_2} + \delta \ol{\kappa}_2 } \expt{\ol{g}_t(\sfY)^2 \ol{\sfZ}^2} + \paren{\frac{\ol{\kappa}_6}{\ol{\kappa}_2^2} - \frac{\ol{\kappa}_4^2}{\ol{\kappa}_2^3} + \delta \frac{\ol{\kappa}_4}{\ol{\kappa}_2}} \expt{\ol{g}_t(\sfY) \ol{\sfZ}^2}^2 }}{\paren{ \frac{\ol{\kappa}_4}{\ol{\kappa}_2^2} + \delta }^2 \expt{\ol{g}_t(\sfY) \ol{\sfZ}^2}^2} , \notag 
\end{align}
matching \Cref{eqn:tau2/nu2}, as desired. 
Similarly, dividing both sides of \Cref{eqn:omega_tau} by $ \chi_{t+1}^2 $ and using \Cref{eqn:stein2}, we get \Cref{eqn:omega2/chi2}. 

Finally, we also specialize \Cref{eqn:SE_lin1,eqn:SE_lin2,eqn:SE_lin3} to the case subject to \Cref{asmp:g(y)}. 
Again using \Cref{eqn:g_stein}, we have: for any $ 1\le i\le t, 0\le j\le t $, 
\begin{align}
\Sigma_{t+1,i+1}
&= \ol{\kappa}_2 \paren{\chi_t \chi_{i} \rho + \Omega_{t,i}} + \frac{\ol{\kappa}_4}{\ol{\kappa}_2} \rho \paren{\expt{\ol{g}_t(\sfY) \ol{\sfZ}^2} \nu_{t-1} \chi_{i} + \expt{\ol{g}_i(\sfY) \ol{\sfZ}^2} \nu_{i-1} \chi_t} \notag \\
&\quad + \frac{\ol{\kappa}_4}{\ol{\kappa}_2^2} \expt{\ol{g}_t(\sfY) \ol{g}_i(\sfY) \sfZ_{t-1} \sfZ_{i-1}} + \frac{\ol{\kappa}_6}{\ol{\kappa}_2^2} \rho \expt{\ol{g}_t(\sfY) \ol{\sfZ}^2} \expt{\ol{g}_i(\sfY) \ol{\sfZ}^2} \nu_{t-1} \nu_{i-1} 
\notag \\
&= \ol{\kappa}_2 \paren{\chi_t \chi_{i} \rho + \Omega_{t,i}} + \frac{\ol{\kappa}_4}{\ol{\kappa}_2} \rho \paren{\expt{\ol{g}_t(\sfY) \ol{\sfZ}^2} \nu_{t-1} \chi_{i} + \expt{\ol{g}_i(\sfY) \ol{\sfZ}^2} \nu_{i-1} \chi_t} \notag \\
&\quad + \frac{\ol{\kappa}_4}{\ol{\kappa}_2^2} \paren{\expt{\ol{g}_t(\sfY) \ol{g}_i(\sfY) \ol{\sfZ}^2} \sigma^2 \nu_{t-1} \nu_{i-1} + \expt{\ol{g}_t(\sfY) \ol{g}_i(\sfY)} \expt{\tau_{t-1} \sfN_{t-1} \cdot \tau_{i-1} \sfN_{i-1}}} \notag \\
&\quad + \frac{\ol{\kappa}_6}{\ol{\kappa}_2^2} \rho \expt{\ol{g}_t(\sfY) \ol{\sfZ}^2} \expt{\ol{g}_i(\sfY) \ol{\sfZ}^2} \nu_{t-1} \nu_{i-1} \notag \\
&= \ol{\kappa}_2 \paren{\chi_t \chi_{i} \rho + \Omega_{t,i}} + \frac{\ol{\kappa}_4}{\ol{\kappa}_2} \rho \paren{\expt{\ol{g}_t(\sfY) \ol{\sfZ}^2} \nu_{t-1} \chi_{i} + \expt{\ol{g}_i(\sfY) \ol{\sfZ}^2} \nu_{i-1} \chi_t} \notag \\
&\quad + \frac{\ol{\kappa}_4}{\ol{\kappa}_2^2} \expt{\ol{g}_t(\sfY) \ol{g}_i(\sfY)} \expt{\tau_{t-1} \sfN_{t-1} \cdot \tau_{i-1} \sfN_{i-1}} \notag \\
&\quad + \paren{
    \frac{\ol{\kappa}_4}{\ol{\kappa}_2} \expt{\ol{g}_t(\sfY) \ol{g}_i(\sfY) \ol{\sfZ}^2} 
    + \frac{\ol{\kappa}_6}{\ol{\kappa}_2^2} \expt{\ol{g}_t(\sfY) \ol{\sfZ}^2} \expt{\ol{g}_i(\sfY) \ol{\sfZ}^2}
} \rho \nu_{t-1} \nu_{i-1} , 
\notag \\
\Omega_{t+1,j+1} 
&= \delta \paren{\frac{1}{\ol{\kappa}_2} \expt{\ol{g}_{t+1}(\sfY) \ol{g}_{j+1}(\sfY) \sfZ_t \sfZ_j} + \frac{\ol{\kappa}_4}{\ol{\kappa}_2^2} \rho \expt{\ol{g}_{t+1}(\sfY) \ol{\sfZ}^2} \expt{\ol{g}_{j+1}(\sfY) \ol{\sfZ}^2} \nu_{t} \nu_{j}} 
\notag \\
&= \delta \frac{1}{\ol{\kappa}_2} \paren{\expt{\ol{g}_{t+1}(\sfY) \ol{g}_{j+1}(\sfY) \ol{\sfZ}^2} \sigma^2 \nu_{t} \nu_{j} + \expt{\ol{g}_{t+1}(\sfY) \ol{g}_{j+1}(\sfY)} \expt{\tau_{t} \sfN_{t} \cdot \tau_{j} \sfN_{j}}} \notag \\
&\quad + \delta \frac{\ol{\kappa}_4}{\ol{\kappa}_2^2} \rho \expt{\ol{g}_{t+1}(\sfY) \ol{\sfZ}^2} \expt{\ol{g}_{j+1}(\sfY) \ol{\sfZ}^2} \nu_{t} \nu_{j} \notag \\
&= \delta \brace{
    \rho \paren{\expt{\ol{g}_{t+1}(\sfY) \ol{g}_{j+1}(\sfY) \ol{\sfZ}^2} + \frac{\ol{\kappa}_4}{\ol{\kappa}_2^2} \expt{\ol{g}_{t+1}(\sfY) \ol{\sfZ}^2} \expt{\ol{g}_{j+1}(\sfY) \ol{\sfZ}^2}} \nu_{t} \nu_{j}
    + \frac{1}{\ol{\kappa}_2} \expt{\ol{g}_{t+1}(\sfY) \ol{g}_{j+1}(\sfY)} \expt{\tau_{t} \sfN_{t} \cdot \tau_{j} \sfN_{j}}
} , 
\notag 
\end{align}
as claimed in \Cref{eqn:Sigmati,eqn:Omegati}. 

\subsection{Proof of \Cref{prop:SE_FP}}
\label{app:pf_prop:SE_FP}

The fixed point equation \Cref{eqn:tau2/nu2_FP} for $ \ol{\tau}^2/\ol{\nu}^2 $ follows by removing time indices from \Cref{eqn:tau2/nu2} and solving for $ \ol{\tau}^2/\ol{\nu}^2 $. 
The fixed point equation \Cref{eqn:omega2/chi2_FP} for $ \ol{\omega}^2/\ol{\chi}^2 $ follows by substituting  the expression in \Cref{eqn:tau2/nu2_FP} for $ \ol{\tau}^2/\ol{\nu}^2 $ into \Cref{eqn:omega2/chi2} (with time indices removed) and performing elementary algebraic manipulations.

Next, we show positivity of both $ \ol{\tau}^2/\ol{\nu}^2 $ and $ \ol{\omega}^2/\ol{\chi}^2 $. 
First consider $ \ol{\tau}^2/\ol{\nu}^2 $. 
The denominator of \Cref{eqn:tau2/nu2_FP} is positive by assumption \Cref{eqn:thr}. 
The positivity of the numerator has been shown in \Cref{eqn:coeff_pos}. 
Therefore, $ \ol{\tau}^2/\ol{\nu}^2 > 0 $. 
Using this together with \Cref{asmp:g(y)}, \Cref{eqn:g>0,eqn:coeff_pos_omega} in \Cref{eqn:omega2/chi2}, we obtain $ \ol{\omega}^2/\ol{\chi}^2 > 0 $. 

\subsection{Proof of \Cref{lem:SE_stay}}
\label{app:pf_lem:SE_stay}

The convergence result \Cref{eqn:lin_conv} has been guaranteed by \Cref{prop:SE}. 
It remains to verify the stationarity of the ratios in \Cref{eqn:SE_stay}. 
To this end, first note that by the choice of $ \nu_0, \tau_0 $ in \Cref{eqn:nu0_tau0}, 
\begin{align}
    \frac{\tau_0^2}{\nu_0^2} &= w_1 = \frac{\ol{\tau}^2}{\ol{\nu}^2} , \label{eqn:tau0/nu0}
\end{align}
where the last equality is by \Cref{prop:SE_FP}. 

By the choice \Cref{eqn:lin_init} of $ g_1 $ and Stein's lemma, 
\begin{align}
    \expt{\partial_z g_1(\sfZ, \sfE)} &= \sigma^{-2} \expt{\sfZ \cdot \ol{g}(\sfY) (\nu_0 \sfZ + \tau_0 \sfN_0)} = \nu_0 \expt{\ol{g}(\sfY) \ol{\sfZ}^2} , \notag \\
    \expt{g_1(\sfZ, \sfE)^2} &= \expt{\ol{g}(\sfY)^2 (\nu_0 \sfZ + \tau_0 \sfN_0)^2} = \sigma^2 \nu_0^2 \expt{\ol{g}(\sfY)^2 \ol{\sfZ}^2} + \tau_0^2 \expt{\ol{g}(\sfY)^2} , \notag 
\end{align}
where $ \sfN_0 \sim \cN(0,1) $ is independent of $ \sfZ, \sfE $. 
Using these formulas, we evaluate $ \chi_1, \omega_1 $ according to \Cref{eqn:SE_lin_init}: 
\begin{align}
    \chi_1 &= \delta \nu_0 \expt{\ol{g}(\sfY) \ol{\sfZ}^2} , \notag \\
    \omega_1^2 &= \delta \frac{1}{\ol{\kappa}_2} \paren{ \sigma^2 \nu_0^2 \expt{\ol{g}(\sfY)^2 \ol{\sfZ}^2} + \tau_0^2 \expt{\ol{g}(\sfY)^2} } + \delta \frac{\ol{\kappa}_4}{\ol{\kappa}_2^2} \rho \nu_0^2 \expt{\ol{g}(\sfY) \ol{\sfZ}^2}^2 \notag \\
    &= \delta \rho \nu_0^2 \expt{\ol{g}(\sfY)^2 \ol{\sfZ}^2} + \delta \frac{\tau_0^2}{\ol{\kappa}_2} \expt{\ol{g}(\sfY)^2} + \delta \rho \nu_0^2 \frac{\ol{\kappa}_4}{\ol{\kappa}_2^2} \expt{\ol{g}(\sfY) \ol{\sfZ}^2}^2 . \notag 
\end{align}
Since $ \nu_0 $ is positive, so is $ \chi_1 $. 
Then 
\begin{align}
    \frac{\omega_1^2}{\chi_1^2} &= \frac{\rho}{\delta} \frac{\expt{\ol{g}(\sfY)^2 \ol{\sfZ}^2} + \frac{\ol{\kappa}_4}{\ol{\kappa}_2^2} \expt{\ol{g}(\sfY) \ol{\sfZ}^2}^2}{\expt{\ol{g}(\sfY) \ol{\sfZ}^2}^2}
    + \frac{\expt{\ol{g}(\sfY)^2}}{\delta \ol{\kappa}_2 \expt{\ol{g}(\sfY) \ol{\sfZ}^2}^2} \frac{\tau_0^2}{\nu_0^2} \notag \\
    &\explain{\Cref{eqn:tau0/nu0}} \frac{\rho}{\delta} \frac{\expt{\ol{g}(\sfY)^2 \ol{\sfZ}^2} + \frac{\ol{\kappa}_4}{\ol{\kappa}_2^2} \expt{\ol{g}(\sfY) \ol{\sfZ}^2}^2}{\expt{\ol{g}(\sfY) \ol{\sfZ}^2}^2}
    + \frac{\expt{\ol{g}(\sfY)^2}}{\delta \ol{\kappa}_2 \expt{\ol{g}(\sfY) \ol{\sfZ}^2}^2} \frac{\ol{\tau}^2}{\ol{\nu}^2}
    \explain{\Cref{eqn:omega2/chi2}} \frac{\ol{\omega}^2}{\ol{\chi}^2}
    = w_2 , \notag 
\end{align}
where the last equality is by \Cref{prop:SE_FP}. 
This implies, by induction on $t$, the first two statements in \Cref{eqn:SE_stay}. 

The remaining statements in \Cref{eqn:SE_stay} follow immediately from \Cref{eqn:nu_chi,eqn:stein2,eqn:nu_nu}. 

\section{Proof of \Cref{lem:edge}}
\label{app:right_edge} 

Recall $ z = X \beta_* = O \Lambda Q^\top \beta_* $, where the second equality is by the SVD of $X$ (see \Cref{asmp:design}). 
Let $ w \coloneqq \Lambda Q^\top \beta_* $ such that $ z = O w $. 
Define 
\begin{align}
\wc{D} &\coloneqq Q \Lambda^\top \Pi_{w^\perp} \wt{O}^\top \Pi_{z^\perp}^\top T \Pi_{z^\perp} \wt{O} \Pi_{w^\perp}^\top \Lambda Q^\top , \label{eqn:whD}
\end{align}
where $ \wt{O} \sim \haar(\bbO(n-1)) $ is independent of everything else. 

\begin{lemma}
\label{lem:conditioning}
It holds almost surely that
\begin{align}
\lambda_3(\wc{D}) &\le \lambda_2(D) \le \lambda_1(\wc{D}) . \label{eqn:lambda321} 
\end{align}
\end{lemma}

\begin{proof}
We write the matrix $D$ as 
\begin{align}
D &= X^\top T X = Q \Lambda^\top O^\top T O \Lambda Q^\top . \notag 
\end{align}
It suffices to study the eigenvalues of $ \Lambda^\top O^\top T O \Lambda $ since conjugation by orthogonal matrices preserves the spectrum. 
By \Cref{prop:haar_cond}, we have
\begin{align}
O \mid \brace{ Ow = z }
&\eqqlaw \frac{z w^\top}{\normtwo{z}^2} + \Pi_{z^\perp} \wt{O} \Pi_{w^\perp}^\top , \notag 
\end{align}
where $ \wt{O} \sim \haar(\bbO(n - 1)) $ is independent of everything else. 
Then 
\begin{align}
\Lambda^\top O^\top T O \Lambda \mid \brace{ Ow = z } 
&\eqqlaw \Lambda^\top 
\paren{ \frac{w z^\top}{\normtwo{z}^2} + \Pi_{w^\perp} \wt{O}^\top \Pi_{z^\perp}^\top }
T
\paren{ \frac{z w^\top}{\normtwo{z}^2} + \Pi_{z^\perp} \wt{O} \Pi_{w^\perp}^\top }
\Lambda . \notag 
\end{align}
By the Courant--Fischer theorem, 
\begin{align}
\left. \lambda_2\paren{ \Lambda^\top O^\top T O \Lambda } \,|\, \brace{ Ow = z } \right.
&\eqqlaw \min_{\substack{\cV\subset\bbR^d \\ \dim(\cV) = d-1}} \max_{v\in\cV\cap\bbS^{d-1}} v^\top  
\paren{ \frac{\Lambda^\top w z^\top}{\normtwo{z}^2} + \Lambda^\top \Pi_{w^\perp} \wt{O}^\top \Pi_{z^\perp}^\top }
T
\paren{ \frac{z w^\top \Lambda}{\normtwo{z}^2} + \Pi_{z^\perp} \wt{O} \Pi_{w^\perp}^\top \Lambda }
 v \notag \\
&\le \max_{\substack{v\in\bbS^{d-1} \\ \inprod{v}{\frac{\Lambda^\top w}{\normtwo{z}}} = 0}} v^\top  
\paren{ \frac{\Lambda^\top w z^\top}{\normtwo{z}^2} + \Lambda^\top \Pi_{w^\perp} \wt{O}^\top \Pi_{z^\perp}^\top }
T
\paren{ \frac{z w^\top \Lambda}{\normtwo{z}^2} + \Pi_{z^\perp} \wt{O} \Pi_{w^\perp}^\top \Lambda }
 v \notag \\
&\le \max_{v\in\bbS^{d-1}} v^\top \Lambda^\top \Pi_{w^\perp} \wt{O}^\top \Pi_{z^\perp}^\top T \Pi_{z^\perp} \wt{O} \Pi_{w^\perp}^\top \Lambda v \notag \\
&= \lambda_1 \paren{ \Lambda^\top \Pi_{w^\perp} \wt{O}^\top \Pi_{z^\perp}^\top T \Pi_{z^\perp} \wt{O} \Pi_{w^\perp}^\top \Lambda } . \notag 
\end{align}
This proves the left-most inequality in \Cref{eqn:lambda321}. 

Again by the Courant--Fischer theorem, 
\begin{align}
\left. \lambda_2\paren{ \Lambda^\top O^\top T O \Lambda } \,|\, \brace{ Ow = z } \right.
&\eqqlaw \min_{\substack{\cV\subset\bbR^d \\ \dim(\cV) = d-1}} \max_{v\in\cV\cap\bbS^{d-1}} v^\top  
\paren{ \frac{\Lambda^\top w z^\top}{\normtwo{z}^2} + \Lambda^\top \Pi_{w^\perp} \wt{O}^\top \Pi_{z^\perp}^\top }
T
\paren{ \frac{z w^\top \Lambda}{\normtwo{z}^2} + \Pi_{z^\perp} \wt{O} \Pi_{w^\perp}^\top \Lambda }
 v \notag \\
&= \max_{\substack{v\in\bbS^{d-1} \\ \inprod{v}{v_*} = 0}} v^\top  
\paren{ \frac{\Lambda^\top w z^\top}{\normtwo{z}^2} + \Lambda^\top \Pi_{w^\perp} \wt{O}^\top \Pi_{z^\perp}^\top }
T
\paren{ \frac{z w^\top \Lambda}{\normtwo{z}^2} + \Pi_{z^\perp} \wt{O} \Pi_{w^\perp}^\top \Lambda }
 v \notag \\
&\ge \max_{\substack{v\in\bbS^{d-1} \\ \inprod{v}{v_*} = 0 \\ \inprod{v}{\frac{\Lambda^\top w}{\normtwo{z}}} = 0}} v^\top  
\paren{ \frac{\Lambda^\top w z^\top}{\normtwo{z}^2} + \Lambda^\top \Pi_{w^\perp} \wt{O}^\top \Pi_{z^\perp}^\top }
T
\paren{ \frac{z w^\top \Lambda}{\normtwo{z}^2} + \Pi_{z^\perp} \wt{O} \Pi_{w^\perp}^\top \Lambda }
 v \notag \\
&= \max_{\substack{v\in\bbS^{d-1} \\ \inprod{v}{v_*} = 0 \\ \inprod{v}{\frac{\Lambda^\top w}{\normtwo{z}}} = 0}} v^\top \Lambda^\top \Pi_{w^\perp} \wt{O}^\top \Pi_{z^\perp}^\top T \Pi_{z^\perp} \wt{O} \Pi_{w^\perp}^\top \Lambda v \notag \\
&\ge \min_{\substack{\cU\subset\bbR^d \\ \dim(\cU) = d-2}} \max_{v\in\cU\cap\bbS^{d-1}} v^\top \Lambda^\top \Pi_{w^\perp} \wt{O}^\top \Pi_{z^\perp}^\top T \Pi_{z^\perp} \wt{O} \Pi_{w^\perp}^\top \Lambda v \notag \\
&= \lambda_3\paren{ \Lambda^\top \Pi_{w^\perp} \wt{O}^\top \Pi_{z^\perp}^\top T \Pi_{z^\perp} \wt{O} \Pi_{w^\perp}^\top \Lambda } , \notag 
\end{align}
where we identify the minimizing subspace $ \cV_* $ of dimension $ d-1 $ with its normal vector $ v_* \in \bbS^{d-1} $, i.e., $ \cV_* = \brace{ v\in\bbR^d : \inprod{v}{v_*} = 0 } $. 
This proves the right-most inequality in \Cref{eqn:lambda321}. 
\end{proof}

It remains to study the extremal eigenvalues of the matrix $\wc{D}$ in \Cref{eqn:whD}. 
Due to potential interest from a  random matrix theory viewpoint, we present in \Cref{sec:edge1,sec:edge2,sec:edge3,sec:edge4} self-contained results for a slightly more general version of $ \wc{D} $ given in \Cref{eqn:def_Dhat}. 

\newtheorem{assump}[theorem]{Assumption}
\newcommand{\deq}{\coloneqq}
\newcommand{\e}[1]{\mathrm{e}^{#1}}
\newcommand{\R} {\mathbb{R}}
\newcommand{\C} {\mathbb{C}}
\newcommand{\N} {\mathbb{N}}
\newcommand{\E} {\mathbb{E}}
\newcommand{\adj}{^{*}} 
\newcommand{\tp}{^{\top}}
\newcommand{\caC}{{\mathcal C}}
\newcommand{\caH}{{\mathcal H}}
\newcommand{\caM}{{\mathcal M}}
\newcommand{\caN}{{\mathcal N}}
\newcommand{\caS}{{\mathcal S}}
\newcommand{\caT}{{\mathcal T}}

\newcommand{\bsa}{{\bm a}}
\newcommand{\bsp}{{\bm p}}
\newcommand{\bss}{{\bm s}}
\newcommand{\bst}{{\bm t}}
\newcommand{\bsu}{{\bm u}}
\newcommand{\bsL}{{\bm L}}

\newcommand{\lone}{\mathbbm{1}} 

\newcommand{\dd}{\mathrm{d}}
\newcommand{\ii}{\mathrm{i}}


\newcommand{\bssig}{\bm{\sigma}}
\newcommand\expct[1]{\mathbb{E}[#1]}
\newcommand\Absv[1]{\left\vert#1\right\vert}
\newcommand\absv[1]{\vert#1\vert}
\newcommand\brkt[1]{\langle#1\rangle}
\newcommand\Brkt[1]{\left\langle#1\right\rangle}


\newcommand{\Cr}{\color{red}}
\newcommand{\nc}{\normalcolor}



	
	\subsection{Setup} \label{sec:edge1}
	Let $\Sigma\in\mathbb{R}^{(d\wedge n)\times (d\wedge n)}$ and $T\in\mathbb{R}^{n\times n}$ be symmetric matrices such that $ \Sigma $ is strictly positive definite. Consider the symmetric random matrix
	\begin{align}
		\wt{D} \deq \Lambda\tp O\tp T O  \Lambda\in\mathbb{R}^{d\times d},
        \notag 
	\end{align}
	where $O\sim\haar(\mathbb{O}(n))$ and $\Lambda\in\R^{d\times n}$ is defined by
	\begin{align}
		\Lambda\deq\begin{cases}
			\begin{bmatrix}
				\Sigma^{1/2}	\\	0_{(n-d)\times d}
			\end{bmatrix}	&	\text{if }n\geq d,\\
		\begin{bmatrix}
			\Sigma^{1/2}	&	0_{n \times (d-n)}
		\end{bmatrix}	&	\text{if }n<d.
		\end{cases}
        \notag 
	\end{align}

	We assume the following on $\Sigma$ and $T$. 
    \begin{enumerate}[label=(A\arabic*)]
    \setcounter{enumi}{\value{asmpctr}}
		\item \label[asmp]{asmp:bulk_nd} $n/d\to\delta\in(0,\infty)$ as $n\to\infty$.
		\item \label[asmp]{asmp:bulk_norm} $\normtwo{\Sigma}$, $\normtwo{\Sigma^{-1}}$, and $\normtwo{T}$ are uniformly bounded over $n$.
		\item \label[asmp]{asmp:bulk_mu} As $n\to\infty$, the empirical spectral distributions $\mu_{T}$ and $\mu_{\Sigma}$ of $T$ and $\Sigma$ converge respectively to $\mu_{\sfT}$ and $\mu_{\sfSigma}$, with $\mu_{\sfT},\mu_{\sfSigma}\neq\delta_{0}$.
		\item \label[asmp]{asmp:strong} For each fixed $\epsilon>0$, we have
		\begin{align}
        &&
			\spec \Sigma&\subset \supp\mu_{\sfSigma}+(-\epsilon,\epsilon),&
			\spec T&\subset \supp\mu_{\sfT}+(-\epsilon,\epsilon) &
        & \notag 
		\end{align}
		for sufficiently large $n$, where $\spec(\cdot)$ denotes the set of eigenvalues.
		\item \label[asmp]{asmp:bulk_edge} The measures $\mu_{\sfSigma}$ and $\mu_{\sfT}$ do not decay too fast around the edges in the following sense:
		\begin{align}
        \begin{split}
				\lim_{z\nearrow\inf\supp\mu_{\sfSigma}}\int_{\R}\frac{x}{x-z}\dd \mu_{\sfSigma}(x)
				&=\infty 
				=\lim_{z\searrow\sup\supp\mu_{\sfSigma}}\int_{\R}\frac{x}{z-x}\dd \mu_{\sfSigma}(x),	\\
				\lim_{z\searrow\sup\supp\mu_{\sfT}}\int_{\R}\frac{x}{z-x}\mu_{\sfT}(x)&=\infty.
		\end{split}\notag 
        \end{align}
        \item \label[asmp]{asmp:nonneg} $\sup\supp\mu_{\sfT}>0$.
	\setcounter{asmpctr}{\value{enumi}}
    \end{enumerate}

    We comment on the above assumptions. 
    \Cref{asmp:bulk_nd} is a restatement of \Cref{asmp:proportional}. 
    Note that in the setting of \Cref{sec:model}, $ \Sigma^{1/2} $ corresponds to $ \diag(\lambda) $ where $ \lambda $ is defined in \Cref{asmp:design}. 
    Therefore, it can be checked that all other assumptions have been covered in \Cref{asmp:design,asmp:preprocess}. 
    In particular, the uniform boundedness of $\normtwo{T}$ in \Cref{asmp:bulk_norm} follows since $ \cT(y) = 1 - \frac{\gamma}{\ol{g}(y) + \gamma} $ is bounded from above due to positivity of both $\gamma$ and $ \ol{g}(y) + \gamma $ given by \Cref{asmp:preprocess}. 
    Also, $ \ol{g}(\sfY) $ is prevented from being constantly $0$ by \Cref{eqn:gbar}, implying $ \mu_{\sfT} \ne \delta_0 $ in \Cref{asmp:bulk_norm}. 
    Finally, since $ \ol{g}(\sfY) $ is larger than $-\gamma$, has mean $0$ and cannot be constantly $0$, the random variable $ \sfT = \frac{\ol{g}(\sfY)}{\ol{g}(\sfY) + \gamma} $ assigns positive mass to $ \bbR_{>0} $, therefore \Cref{asmp:nonneg} holds. 
	
	\subsection{Convergence of bulk edge} \label{sec:edge2}
	Consider the following variant of $\wt{D}$: 
	\begin{align}
		\wh{D}=\wt{\Lambda}\tp  \wt{O}\tp \wt{T} \wt{O} \wt{\Lambda}\in \R^{d\times d}, \label{eqn:def_Dhat} 
	\end{align}
	where now $\wt{O}$ is Haar distributed on $\bbO(n-1)$ independent of $\wt{T}$ and $\wt{\Lambda}$, $\wt{T}$ is a random minor of $T$, and $\wt{\Lambda}$ satisfies
	\begin{align}
		\wt{\Lambda}=\begin{cases}
			\begin{bmatrix}
				\wt{\Sigma}^{1/2}	&	0_{d-1}	\\
				0_{(n-d)\times (d-1)}	&	0
			\end{bmatrix}	&	\text{if }n\geq d, \\
			\begin{bmatrix}
				\wt{\Sigma}^{1/2}	&	0_{(n-1)\times (d-n+1)}	
			\end{bmatrix}	&	\text{if }n<d,
		\end{cases}
        \notag 
	\end{align}
	for a random minor $\wt{\Sigma}$ of $\Sigma$ with size $n\wedge d-1$. 

    Before proceeding, we verify that $\wc{D}$ in \Cref{eqn:whD} is a special case of $\wh{D}$ in \Cref{eqn:def_Dhat}. 
    Consequently, all results concerning $\wh{D}$ apply to $\wc{D}$. 
    Indeed, the spectrum of $\wc{D}$ is the same as that of
    \begin{align}
        & \Lambda^\top \Pi_{w^\perp} \wt{O}^\top \Pi_{z^\perp}^\top T \Pi_{z^\perp} \wt{O} \Pi_{w^\perp}^\top \Lambda . 
        \label{eqn:wcD}
    \end{align}
    Identifying 
    \begin{align}
    &&
        \wt{T} &= \Pi_{z^\perp}^\top T \Pi_{z^\perp} \in \bbR^{(n-1)\times(n-1)} , & 
        \wt{\Lambda} &= \Pi_{w^\perp}^\top \Lambda \in \bbR^{(n-1)\times d} , & 
    & \notag 
    \end{align}
    we see that the matrix \Cref{eqn:wcD} is a special case of $ \wh{D} $. 

    For the matrix $\wh{D}$, our goal is to show the following.
	\begin{lemma}\label{lem:edge_goal}
        It holds almost surely as $n\to\infty$ that
		\begin{align}
		\lambda_{1}(\wh{D}),\lambda_{3}(\wh{D})\to \sup\supp( \mu_{\sfSigma_n}\boxtimes\mu_{\sfT}),
        \notag 
		\end{align}
		where
		\begin{align}
			\mu_{\sfSigma_n}=\begin{cases}
				\frac{1}{\delta}\mu_{\sfSigma}+\frac{1-\delta}{\delta}\delta_{0}	&	\text{if }\delta>1,\\
				\mu_{\sfSigma}	&	\text{otherwise}.
			\end{cases}
            \label{eqn:mu_sfLambdan2}
		\end{align}
	\end{lemma}
    Note that $ \sfSigma_n $ in \Cref{eqn:mu_sfLambdan2} coincides with $ \sfLambda_n^2 $ in \Cref{eqn:sfLambda_nd}. 
	We prove \Cref{lem:edge_goal} at the end of this section. 
	\subsubsection{Strong asymptotic freeness of Haar orthogonal and deterministic matrices}
    \label{sec:strong_free}
	The first ingredient of the proof of \Cref{lem:edge_goal} is the strong asymptotic freeness of Haar orthogonal and deterministic matrices, that is of independent interest.
	\begin{proposition}\label{prop:str_asymp}
		Let $U_{1},\ldots,U_{p}\in\R^{N\times N}$ be a collection of independent random matrices Haar-distributed on $\mathbb{SO}(N)$ and $Y_{1},\ldots,Y_{q}\in\R^{N\times N}$ be deterministic matrices. Suppose that $\max_{j=1,\ldots,q}\normtwo{Y_{j}}\leq C$ for a constant $C>0$. Then, for any self-adjoint $*$-polynomial $Q$ in $(p+q)$ variables and any constant $\eps>0$, for sufficiently large $N$ we have
		\begin{align}
			\spec(Q(U_{1},\ldots,U_{p},Y_{1},\ldots,Y_{p}))\subset\spec(Q(u_{1},\ldots,u_{p},Y_{1},\ldots,Y_{p}))+(-\eps,\eps),
            \notag 
		\end{align}
		where the binary operation on the right-hand side is Minkowski sum and $(u_{1},\ldots,u_{p})$ is a collection of $*$-free Haar unitary elements\footnote{A Haar unitary element is a unitary operator $u$ in a tracial von Neumann algebra $(\caM,\brkt{\cdot})$ whose spectral distribution is uniform on $\bbS^{1}$.} that are $*$-free of $(Y_{1},\ldots,Y_{q})$.
	\end{proposition}
    Note that \Cref{prop:str_asymp} applies directly to our model $\wh{D}$, since we may take $\wt{O}$ to be Haar-distributed on $\mathbb{SO}(n-1)$, instead of $\mathbb{O}(n-1)$, without changing the law of $\wh{D}$. 
	
	On a high level, our proof follows \cite[Section 3]{Collins_Male}; we aim at expressing $U_{i}$'s as a functional calculus of some other matrix for which the strong asymptotic freeness is known. We start the proof with the strong asymptotic freeness of GOEs\footnote{The Gaussian orthogonal ensemble $ \GOE(d) $ is the probability measure on $ d\times d $ real symmetric matrices $ G $ such that $ G_{i,i} \iid \cN(0,2/d) $ and $ G_{i,j} \iid \cN(0,1/d) $ for all $ 1\le i< j\le d $. } and deterministic matrices.
	\begin{lemma}[{\cite[Theorem 4.2]{Fan_Sun_Wang}}]\label{lem:str_asymp}
		\Cref{prop:str_asymp} holds if we replace $(U_{1},\ldots,U_{p})$ with independent GOEs $(W_{1},\ldots,W_{p})$ and $(u_{1},\ldots,u_{p})$ with free semi-circular elements $(s_{1},\ldots,s_{p})$. 
	\end{lemma}
	In the next input, we represent a Haar orthogonal matrix in a canonical form $U\Theta U\tp$ with a Haar orthogonal $U$, so that one can express it as a functional calculus of an orthogonally invariant matrix. In what follows we write
	\begin{align}
		R(\theta)=\begin{bmatrix}
			\cos \theta	&	-\sin \theta \\
			\sin\theta&\cos\theta
		\end{bmatrix}\in\R^{2\times 2},\qquad \theta\in\R.
        \notag 
	\end{align}
	\begin{lemma}[{\cite[Theorem 3.5]{Meckes}}]\label{lem:Haar}
		Let $\Theta\in\R^{N\times N}$ be a random block diagonal matrix distributed as
		\begin{align}
			\Theta=\begin{cases}\diag(R(\theta_{1}),\ldots,R(\theta_{m}))	& N=2m,	\\
				\diag(R(\theta_{1}),\ldots,R(\theta_{m}),1)	&	N=2m+1,
			\end{cases}
            \notag 
		\end{align}
		where $0\leq \theta_{1}\leq \cdots\leq \theta_{m}\leq \pi$ is sampled from the joint density (with suitable choices of $Z_{N}$ and $Z_{N}'$)
		\begin{align}
			\frac{1}{Z_{N}}\prod_{1\leq j<k\leq m}	(\cos\theta_{j}-\cos \theta_{k})^{2},	\qquad & \text{if }N=2m, \notag \\
			\frac{1}{Z_{N}'}\prod_{1\leq j\leq m}\sin^{2}\left(\frac{\theta_{j}}{2}\right)\prod_{1\leq k<\ell\leq m}(\cos\theta_{j}-\cos \theta_{k})^{2},\qquad &	\text{if }N=2m+1.
            \notag 
		\end{align}
		Let $V$ be a Haar orthogonal matrix independent of $\Theta$. Then the matrix $V\Theta V\tp$ is Haar distributed on $\mathbb{SO}(N)$.
	\end{lemma}
	
	As one can see, unlike the usual diagonalization, we find that the matrix $\Theta$ is a block diagonal matrix consisting of $(2\times 2)$ blocks. Indeed, the true eigenvectors of a Haar orthogonal matrix $U$ should be genuinely complex with the complex-conjugate symmetry. This makes it difficult to express $U$ as a functional calculus of GOE, whose eigenvectors are real. 
	
	Instead, we shall express $U$ as a functional calculus of a \emph{skew}-symmetric matrix; skew-symmetric matrices already have the same canonical form consisting of $(2\times 2)$ blocks, the only difference being that the blocks are of the form $xJ$ for some $x\in\bbR$ instead of $R(\theta)$ where
	\begin{align}
		J=\begin{bmatrix}
			0	&	-1	\\
			1	&	0
		\end{bmatrix}.
        \notag 
	\end{align}
    Note that the spectrum of a skew-symmetric matrix is purely imaginary and the eigenvalues of $J$ are $-\ii,\ii$. 
	More fundamentally, it is natural to work with skew-symmetric matrices as the Lie algebra $\mathfrak{so}(2m)$ associated with the group $\mathbb{SO}(2m)$ is exactly the space of skew-symmetric matrices. In the complex case, the Lie algebra $\mathfrak{u}(n)$ associated to $\mathbb{U}(n)$ is the space of skew-Hermitian matrices, which is nothing but Hermitian matrices rotated by a factor of $\ii$.
	
	Thus we work with the simplest orthogonally invariant skew-symmetric matrix, which can be written as $X-X\tp$ where $X\in\R^{N\times N}$ has independent $\caN(0,1/N)$ entry, i.e.\ Ginibre ensemble. To this end, we next show that we can take $W_{i}$'s in \Cref{lem:str_asymp} to be real Ginibre ensembles:
	\begin{lemma}\label{lem:str_asymp_Gin}
		\Cref{prop:str_asymp} holds if we replace $(U_{1},\ldots,U_{p})$ with independent real Ginibre ensembles $(X_{1},\ldots,X_{p})\in\R^{N\times N}$, and $(u_{1},\ldots,u_{p})$ with $*$-free circular elements $(x_{1},\ldots,x_{p})$.
	\end{lemma}
	\begin{proof}
        In what follows, we use 
        \begin{align}
            A \ot B &= \matrix{
                B_{1,1} A & \cdots & B_{1,q} A \\
                \vdots & \ddots & \vdots \\
                B_{p,1} A & \cdots & B_{p,q} A
            } \notag 
        \end{align}
        to denote the Kronecker product of an $m\times n$ matrix $ A $ and a $p\times q$ matrix $ B $. 
		Let $W_{1},\ldots,W_{p}\in\R^{2N\times 2N}$ be independent GOEs, and take projection matrices 
		\begin{align}
			P_{1}=\begin{bmatrix}
				I_{N}	&	0_{N\times N}
			\end{bmatrix},\qquad P_{2}=\begin{bmatrix}
				0_{N\times N}	&	I_{N}
			\end{bmatrix},\qquad 
			\Pi_{k\ell}=P_{k}\tp P_{\ell}=I_{N}\otimes E_{k\ell}\in\R^{2N\times 2N},
            \notag 
		\end{align}
        where $ E_{k\ell} = e_k e_\ell^\top \in \bbR^{2\times 2} $ for $ k,\ell\in\{1,2\} $.
		Then we find that the collection
		\begin{align}
			X_{i}=\sqrt{2}P_{1}W_{i}P_{2}\tp \notag 
		\end{align}
		consists of independent Ginibre ensembles. Now we find that
		\begin{align}
				Q(X_{1},\ldots,X_{p},Y_{1},\ldots,Y_{q})\otimes I_{2}
				&=Q(X_{1}\otimes I_{2},\ldots, X_{p}\otimes I_{2},Y_{1}\otimes I_{2},\ldots,Y_{q}\otimes I_{2})	\notag \\
				&=\wt{Q}(W_{1},\ldots,W_{p},\Pi_{11},\Pi_{12},\Pi_{21},\Pi_{22},Y_{1}\otimes I_{2},\ldots,Y_{q}\otimes I_{2}) \notag 
		\end{align}
		for another $*$-polynomial $\wt{Q}$ that is obtained from $Q$ via
		\begin{align}
			X_{i}\otimes I_{2}=\sqrt{2}\sum_{k=1}^{2}P_{k}\tp P_{1}W_{i}P_{2}\tp P_{k}=\sqrt{2}\sum_{k=1}^{2}\Pi_{k1}W_{i}\Pi_{2k}. \notag 
		\end{align}
		Now an application of \Cref{lem:str_asymp} gives, in the limit $N\to\infty$, that
		\begin{align}
				\spec(Q(X_{1},\ldots,X_{p},Y_{1},\ldots,Y_{q}))
				&=\spec(Q(X_{1},\ldots,X_{p},Y_{1},\ldots,Y_{q})\otimes I_{2}) \label{eq:xlift} \\
				&\subset\spec(\wt{Q}(s_{1},\ldots,s_{p},\Pi_{11},\Pi_{12},\Pi_{21},\Pi_{22},Y_{1}\otimes I_{2},\ldots,Y_{q}\otimes I_{2}))+(-\eps,\eps) \notag \\
				&=\spec(Q(\wh{x}_{1},\ldots,\wh{x}_{p},Y_{1}\otimes I_{2},\ldots,Y_{q}\otimes I_{2})+(-\eps,\eps), \notag 
		\end{align}
		where $(s_{i})_{1\leq i\leq p}$ are free semi-circular elements that are $*$-free of $\{Y_{j}\otimes I_{2},\Pi_{k\ell}\}$, both realized in the free product von Neumann algebra $\caM*\C^{2N\times 2N}$, and in the last line we abbreviated
		\begin{align}
			\wh{x}_{i}\deq \sqrt{2}\sum_{k=1}^{2}\Pi_{k1}s_{i}\Pi_{2k}.
            \notag 
		\end{align}
		
		We next show the $*$-distributional identity
		\begin{align}\label{eq:hatx}
			\left(\wh{x}_{1},\ldots,\wh{x}_{p},Y_{1}\otimes I_{2},\ldots,Y_{q}\otimes I_{2}\right)\overset{d}{=}(x_{1},\ldots,x_{p},Y_{1},\ldots,Y_{q}),
		\end{align}
		where $x_{1},\ldots,x_{p}$ is a collection of circular elements that are $*$-free within themselves and from $Y_{1},\ldots,Y_{q}$. We achieve this via the so-called \emph{meta-argument}. By the usual asymptotic freeness of GOEs, for any $*$-polynomial $P$, we have 
		\begin{align}
			\brkt{P\left(\wh{x}_{1},\ldots,\wh{x}_{p},Y_{1}\otimes I_{2},\ldots,Y_{q}\otimes I_{2}\right)}
			=\lim_{M\to\infty}\Brkt{P\left(\wh{X}_{1},\ldots,\wh{X}_{p}, I_{M}\otimes Y_{1}\otimes I_{2},\ldots,I_{M}\otimes Y_{q}\otimes I_{2}\right)},
            \notag 
		\end{align}
		where $\brace{\wh{X}_{i}\in\R^{(2NM\times 2NM)} : 1\le i\le p}$ are given by
		\begin{align}
			\wh{X}_{i}=\sqrt{2}\sum_{k=1}^{2}(I_{M}\otimes \Pi_{k1})\wt{W}_{i}(I_{M}\otimes\Pi_{2k}) \notag
		\end{align}
		with a collection $(\wt{W}_{i})_{1\leq i\leq p}$ of GOEs of size $(2NM\times 2NM)$. Next notice that each $\wh{X}_{i}$ can be written as
		\begin{align}
			\wh{X}_{i}=\begin{bmatrix}
				\wt{X}_{i}	&	0	\\	0	& \wt{X}_{i}
			\end{bmatrix}=\wt{X}_{i}\otimes I_{2},
            \notag
		\end{align}
		where $(\wt{X}_{i})_{1\leq i\leq p}$ is a collection of independent real Ginibre ensemble of size $(NM\times NM)$. 
		Thus
		\begin{align}
			\Brkt{P\left(\wh{X}_{1},\ldots,\wh{X}_{p},I_{M}\otimes Y_{1}\otimes I_{2},\ldots,I_{M}\otimes Y_{q}\otimes I_{2}\right)}
			=\Brkt{P\left(\wt{X}_{1},\ldots,\wt{X}_{p},I_{M}\otimes Y_{1},\ldots,I_{M}\otimes Y_{q}\right)},
            \notag
		\end{align}
		and taking the limit as $M\to\infty$ gives
		\begin{align}
			\brkt{P\left(\wh{x}_{1},\ldots,\wh{x}_{p},Y_{1}\otimes I_{2},\ldots,Y_{q}\otimes I_{2}\right)}
			=\brkt{P\left(x_{1},\ldots,x_{p},Y_{1},\ldots,Y_{q}\right)}. \notag
		\end{align}
		Since $P$ was generic, this shows \Cref{eq:hatx} so that
		\begin{align}
			\spec(Q(\wh{x}_{1},\ldots,\wh{x}_{p},Y_{1}\otimes I_{2},\ldots,Y_{q}\otimes I_{2}))=\spec(Q(x_{1},\ldots,x_{p},Y_{1},\ldots,Y_{q})).
            \notag
		\end{align}
		Plugging this in \Cref{eq:xlift} proves the lemma.
	\end{proof}
	
	\begin{proof}[Proof of \Cref{prop:str_asymp}]	
		We aim at approximately coupling a Haar orthogonal matrix of order $N$ with 
		\begin{align}
			F((X-X\tp)/\sqrt{2}) \notag
		\end{align}
		for some function $F$ applied according to functional calculus, where $X\in\R^{N\times N}$ is a real Ginibre ensemble. Indeed, since $X-X\tp$ is invariant under orthogonal conjugation and skew-symmetric, we can write its canonical form as in \Cref{lem:Haar} as
		\begin{align}\label{eq:skew}
			\frac{X-X\tp}{\sqrt{2}}=V\wt{\Theta}V\tp,\qquad \wt{\Theta}=\begin{cases}
				\diag(\lambda_{1}J,\ldots,\lambda_{m}J),	&	N=2m,	\\
				\diag(\lambda_{1}J,\ldots,\lambda_{m}J,0),	&	N=2m+1,
			\end{cases}\qquad 
			J=\begin{bmatrix}
				0	&	-1	\\1	&	0
			\end{bmatrix},
		\end{align}
		where $V$ is a Haar orthogonal matrix independent of $\wt{\Theta}$ and 
        \begin{equation*}
            -\lambda_{m}\leq \cdots \leq -\lambda_{1}\leq \lambda_{1}\leq\cdots\leq \lambda_{m}
        \end{equation*} are the ordered eigenvalues of $(X-X\tp)/(\sqrt{2}\ii)$, with an extra nullity when $N$ is odd (recall that skew-symmetric matrices of odd order are singular). By the usual functional calculus we get
		\begin{align}
			F((X-X\tp)/\sqrt{2})=VF(\wt{\Theta})V\tp,
            \notag 
		\end{align}
		and taking $F(\cdot)=\exp(f(\cdot))$ for an odd function $f:\ii\R\to\ii\R$, we get
		\begin{align}
			F(\wt{\Theta})=\begin{cases}
				\diag(F(\lambda_{1}J),\ldots,F(\lambda_{m}J)),	&	N=2m,	\\
				\diag(F(\lambda_{1}J),\ldots,F(\lambda_{m}J),\exp(f(0))),	&	N=2m+1,
			\end{cases}
            \notag 
		\end{align}
		and by spectral theorem we can compute each block as
		\begin{align}
				F(\lambda_{i}J)&=F\left(\wt{U}\begin{bmatrix}\ii\lambda_{i}	&	0	\\	0	&	-\ii\lambda_{i}\end{bmatrix}\wt{U}\adj\right)
				=\exp\left(\wt{U}\begin{bmatrix}f(\ii\lambda_{i})	&	0	\\	0	&-f(\ii\lambda_{i})\end{bmatrix}\wt{U}\adj\right) \notag \\
				&=\exp(-\ii f(\ii\lambda_{i})J)	
				=\begin{bmatrix}
					\cos(\ii f(\ii\lambda_{i}))	&	\sin (\ii f(\ii\lambda_{i}))	\\	
					-\sin(\ii f(\ii\lambda_{i}))	&	\cos(\ii f(\ii\lambda_{i}))
				\end{bmatrix}
                =R(-\ii f(\ii\lambda_{i})), \notag 
		\end{align}
		with a suitable unitary matrix $\wt{U}\in\C^{2\times 2}$ diagonalizing $J$.
		
		Next, we choose a suitable odd function $f$ so that for independently sampled $\wt{\Theta}$ and $\Theta$ (recall that $(\lambda_{i})_{1\leq i\leq m}$ and $(\theta_{i})_{1\leq i\leq m}$ are in increasing order)
		\begin{align}
			\max_{1\leq i\leq m}\norm{}{F(\lambda_{i}J)-R(\theta_{i})}
			\leq\max_{1\leq i\leq m}\absv{-\ii f(\ii\lambda_{i})-\theta_{i}}=o(1),
            \notag 
		\end{align}
        which would then give a coupling between $F((X-X\tp)/\sqrt{2})$ and a $\mathbb{SO}(N)$-Haar matrix $U$ as
		\begin{align}
			\norm{}{F((X-X\tp)/\sqrt{2})-U}=\norm{}{VF(\wt{\Theta})V\tp -V\Theta V\tp}=\max_{1\leq i\leq m}\norm{}{F(\lambda_{i}J)-R(\theta_{i})}=o(1),
            \notag 
		\end{align}
		where $V$, $\wt{\Theta}$, and $\Theta$ are all independently sampled.
		
		In order to construct the appropriate function $f$, take $F_{\mathrm{qc}}:\R\to[0,1]$ to be the odd extension of the c.d.f.\ of the quarter-circle distribution, i.e., 
		\begin{align}
			F_{\mathrm{qc}}(x)\deq\int_{0}^{x}\lone_{[-2,2]}(t)\frac{\sqrt{4-t^{2}}}{\pi}\dd t.
            \notag 
		\end{align}
		Now since $\ii (X-X\tp)/\sqrt{2}$ is a complex Wigner matrix, the empirical distribution of $\lambda_{i}$'s converge to the quarter circle distribution on $[0,2]$, which gives the uniform convergence of c.d.f.s on $\R$:
		\begin{align}
			\sup_{x\in [0,2]}\Absv{F_{\mathrm{qc}}(x)-\frac{\absv{\{i\in[m]:\lambda_{i}\leq x\}}}{\lceil N/2\rceil}}\to 0,\qquad \text{almost surely}.
            \notag 
		\end{align}
		In particular, restricting to $x=\lambda_{i}$ yields
		\begin{align}
			\max_{1\leq i\leq N/2}\Absv{F_{\mathrm{qc}}(\lambda_{i})-\frac{i+\lone(N\text{ is odd})}{\lceil N/2\rceil}}\to 0,\qquad \text{almost surely}.
            \notag 
		\end{align}
		Likewise, taking $\mu_{\theta}$ to be the empirical c.d.f.\ of $\theta_{i}$'s (on $[0,\pi]$) from \Cref{lem:Haar}, we have that $\mu_{\theta}$ converges to the uniform distribution on $[0,\pi]$ (see e.g. \cite[Theorem 4.13]{Meckes}), so  
		\begin{align}
			\max_{1\leq i\leq N/2}\Absv{\mu_{\theta}(\theta_{i})-\frac{\theta_{i}}{\pi}}=\max_{1\leq i\leq N/2}\Absv{\frac{i+\lone(N\text{ is odd})}{\lceil N/2\rceil}-\frac{\theta_{i}}{\pi}}\to 0.
            \notag 
		\end{align}
		Therefore we find that
		\begin{align}
			\max_{1\leq i\leq N/2}\Absv{\pi F_{\mathrm{qc}}(\lambda_{i})-\theta_{i}}\to 0 \qquad \text{almost surely}.
            \notag 
		\end{align}
		Thus we choose the odd function $f$ as
		\begin{align}
			f(\ii x)= \ii \pi F_{\mathrm{qc}}(x),
            \notag 
		\end{align}
		which then guarantees, for the function $F=\e{f}$,
		\begin{align}
			\norm{}{VF(\wt{\Theta})V\tp -V\Theta V\tp}\to 0, \qquad \text{almost surely}.
            \notag 
		\end{align}
		
		Using the complex version of the Stone--Weierstrass theorem, we next take a sequence of polynomials $(P_{n})_{n\in\N}$ so that 
		\begin{align}
			\lim_{n\to\infty}\sup_{x\in[-4,4]}\absv{P_{n}(x)-\exp(\ii \pi F_{\mathrm{qc}}(x))}=0.
            \notag 
		\end{align}
		Then we have
		\begin{align}
                & \norm{}{Q(U_1, \ldots, U_p, Y_1, \ldots, Y_q) - Q(F((X_1-X_1^\top)/\sqrt{2}), \ldots, F((X_p-X_p^\top)/\sqrt{2}), Y_1, \ldots, Y_q)} \notag \\
				&= \norm{}{Q(V_1 \Theta_1 V_1^\top,\ldots,V_p \Theta_p V_p^\top,Y_{1},\ldots,Y_{q})
					-Q(V_{1}F(\wt{\Theta}_{1})V_{1}\tp, \ldots, V_{p}F(\wt{\Theta}_{p})V_{p}\tp,Y_{1},\ldots,Y_{q})} \notag \\
				\qquad&\lesssim \max_{1\leq j\leq p}\norm{}{\Theta_{i}-F(\wt{\Theta}_{i})}\to 0,\qquad \text{almost surely},
                \label{eq:str_1}
		\end{align}
		where the implicit constant depends only on $Q$ and the upper bound $C$ for $\norm{}{Y_{j}}$'s, $(V_{i},\Theta_{i},\wt{\Theta}_{i})_{1\leq i\leq p}$ are all independent, $V_{i}$'s are Haar distributed on $\mathbb{O}(N)$, $\Theta_{i}$'s are distributed as in \Cref{lem:Haar}, and $\wt{\Theta}_{i}$'s are distributed as \Cref{eq:skew}.
		Now we use the polynomial $P_{n}$ to get
		\begin{align}
				&\norm{}{Q(V_{1}F(\wt{\Theta}_{1})V_{1}\tp,\ldots, V_{p}F(\wt{\Theta}_{p})V_{p}\tp,Y_{1},\ldots,Y_{q})
					-Q(V_{1}P_{n}(-\ii\wt{\Theta}_{1})V_{1}\tp,\ldots, V_{p}P_{n}(-\ii\wt{\Theta}_{p})V_{p}\tp,Y_{1},\ldots,Y_{q})} \notag \\
				&\qquad\qquad	\lesssim \sup_{\substack{ x\in\R \\ \absv{x}\leq \max\limits_{j=1,\ldots,p}\norm{}{\wt{\Theta}_{j}}}}\absv{P_{n}(x)-F(\ii x)}.
                \label{eq:str_2}
		\end{align}
		The last quantity converges to $0$ in the limit $\lim_{n\to\infty}\lim_{N\to\infty}$ in this order, since $\norm{}{\wt{\Theta}_{j}}\leq 2\norm{}{X_{j}}$ is bounded by $4$ almost surely as $N\to\infty$ (since $X\tp X$ is a Wishart ensemble) and
		\begin{align}
			\sup_{\absv{x}\leq 4} \absv{P_{n}(x)-F(\ii x)}=\sup_{\absv{x}\leq 4}\absv{P_{n}(x)-\exp(\ii \pi F_{\mathrm{qc}}(x))}\to 0 \qquad\text{as }n\to\infty.
            \notag 
		\end{align}
		We thus conclude that, in the limit $\lim_{n\to\infty}\lim_{N\to\infty}$,
		\begin{align}
				&\spec(Q(U_{1},\ldots,U_{p},Y_{1},\ldots,Y_{q}))	\notag \\
				& \subset \spec(Q(V_{1}P_{n}(-\ii \wt{\Theta}_{1})V_{1}\tp,\ldots, V_{p}P_{n}(-\ii\wt{\Theta}_{p})V_{p}\tp,Y_{1},\ldots,Y_{q}))+(-\eps/3,\eps/3) \notag \\
				&=\spec(Q(P_{n}(-\ii (X_{1}-X_{1}\tp)/\sqrt{2}),\ldots, P_{n}(-\ii(X_{p}-X_{p}\tp)/\sqrt{2}),Y_{1},\ldots,Y_{q}))+(-\eps/3,\eps/3) \notag \\
				&\subset \spec(Q(P_{n}(-\ii (x_{1}-x_{1}\adj)/\sqrt{2}),\ldots, P_{n}(-\ii (x_{p}-x_{p}\adj)/\sqrt{2}),Y_{1},\ldots,Y_{q}))+(-2\eps/3,2\eps/3) \notag \\
				&\subset \spec(Q(F((x_{1}-x_{1}\adj)/\sqrt{2}),\ldots, F((x_{p}-x_{p}\adj)/\sqrt{2}),Y_{1},\ldots,Y_{q}))+(-\eps,\eps),
                \notag 
		\end{align}
		where the second line follows from \Cref{eq:str_1,eq:str_2}, the fourth line from \Cref{lem:str_asymp_Gin}, and the last line again from \Cref{eq:str_2}. 
        We conclude using that $(F((x_{i}-x_{i}\adj)/\sqrt{2}))_{1\leq i\leq p}$ are $*$-free Haar unitary elements: Their $*$-freeness follows from those of $x_{i}$'s, and each of them is Haar unitary since $\ii(x_{i}-x_{i}\adj)/\sqrt{2}=\sqrt{2}\re[\ii x_{i}]$ is a semi-circular element, whose push forward under $x\mapsto F(-\ii x)$ has the uniform spectral distribution on $\mathbb{S}^{1}$.
	\end{proof}
	
	\subsubsection{Insensitivity to $\delta$}
	The following result is our second ingredient, which we use to deal with the discrepancy between $\Lambda$ and $\wt{\Lambda}$ when $\delta=1$.
	\begin{lemma}\label{lem:minor_shrink}
		Let $\mu$ and $\nu$ be compactly supported Borel probability measures on $(0,\infty)$ and $\R$, respectively. Suppose that $\sup\supp \nu>0$. Then, for all $\epsilon\in(0,1]$, we have
		\begin{align}
			\sup\supp\left(((1-\epsilon)\mu+\epsilon\delta_{0})\boxtimes \nu\right)\leq \sup\supp\left(\mu\boxtimes \nu\right).
            \notag 
		\end{align}
	\end{lemma}
	\begin{proof}
		Let $\bss$ and $\bst$ be realizations of $\mu$ and $\nu$ in $W\adj$-probability spaces $\caM_{1}$ and $\caM_{2}$, respectively. Then we easily find that the measure
		$(1-\eps)\mu + \eps\delta_{0}$ is the spectral measure of $\bss\otimes (1-\bsp_{\epsilon})$ in the (von Neumann algebra) tensor product $\caM_{1}\otimes L^{\infty}([0,1])$, where $\bsp_{\epsilon}=\lone_{[0,\epsilon]}\in L^{\infty}([0,1])$ is the projection. Now we have that the free convolution
		\begin{align}
			((1-\epsilon)\mu+\epsilon\delta_{0})\boxtimes \nu
            \notag 
		\end{align}
		is the spectral measure of $(\sqrt{\bss}\otimes (1-\bsp_{\epsilon}))\bst(\sqrt{\bss}\otimes (1-\bsp_{\epsilon}))$ in the (von Neumann algebra) free product $(\caM_{1}\otimes L^{\infty}([0,1]))*\caM_{2}$. 
		
		Then we find
		\begin{align}
				\sup\supp\left(((1-\epsilon)\mu+\epsilon\delta_{0})\boxtimes \nu\right)
				&=\sup_{v\in \caH, \norm{}{v}= 1}\brkt{v,(\sqrt{\bss}\otimes (1-\bsp_{\epsilon}))\bst(\sqrt{\bss}\otimes (1-\bsp_{\epsilon}))v} \notag \\
				&=\sup_{v\in \caH, \norm{}{v}=1}\brkt{(1\otimes (1-\bsp_{\epsilon}))v,(\sqrt{\bss}\otimes1)\bst(\sqrt{\bss}\otimes  1)(1\otimes (1-\bsp_{\epsilon}))v} \notag \\
				&\leq \sup_{w\in\mathrm{Range}(1\otimes(1-\bsp_{\epsilon})), \norm{}{w}\leq 1}\brkt{w,(\sqrt{\bss}\otimes1)\bst(\sqrt{\bss}\otimes  1)w} \notag \\
				&=0\vee \sup_{w\in\mathrm{Range}(1\otimes(1-\bsp_{\epsilon})),\norm{}{w}=1}\brkt{w,(\sqrt{\bss}\otimes1)\bst(\sqrt{\bss}\otimes  1)w} \notag \\
				&\leq 0\vee \sup\supp (\mu\boxtimes\nu), \notag 
		\end{align}
        where $\caH$ is the Gelfand--Naimark--Segal (GNS) Hilbert space of $(\caM_{1}\otimes L^{\infty}([0,1]))*\caM_{2}$.
        To conclude, we only need to notice that $\sup\supp(\mu\boxtimes\nu)\geq 0$: If not, the realization $\sqrt{\bss}\bst\sqrt{\bss}$ would be negative semidefinite, which contradicts $\sup\supp\nu>0$ as $\bss$ is invertible since $\mu$ is supported on $(0,\infty)$.
	\end{proof}

\subsubsection{Proof of \Cref{lem:edge_goal}}
\label{sec:pf_lem:edge_goal}
Notice that, by Cauchy interlacing, the empirical eigenvalue distributions of $\wt{\Lambda}\wt{\Lambda}^{\tp}$ and $\wt{T}$ converge respectively to $\mu_{\sfSigma_n}$ and $\mu_{\sfT}$. Hence from the usual asymptotic freeness it immediately follows that the empirical spectral distribution of $\wh{D}$ converges to $\mu_{\sfSigma_n}\boxtimes\mu_{\sfT}$, which then implies
\begin{align}
	\liminf_{n\to\infty}\lambda_{3}(\wh{D})\geq\sup\supp(\mu_{\sfSigma_n}\boxtimes\mu_{\sfT}), \qquad\text{almost surely},
    \notag 
\end{align}
since the left-hand side is bounded from below by any $O(1)$-quantile of $\mu_{\sfSigma_n}\boxtimes\mu_{\sfT}$. 

For the upper bound, we first claim that $\lambda_{1}(\wh{D})=\lambda_{1}(\sqrt{\wt{\Lambda}\wt{\Lambda}\tp}\wt{O}\tp \wt{T}\wt{O}\sqrt{\wt{\Lambda}\wt{\Lambda}\tp})$. To see this, note that 
\begin{align}
	\lambda_{1}(\wh{D})=\begin{cases}
	    \lambda_{1}(\sqrt{\wt{\Lambda}\wt{\Lambda}\tp}\wt{O}\tp \wt{T}\wt{O}\sqrt{\wt{\Lambda}\wt{\Lambda}\tp})    &   \text{if } n\geq d,  \\
        0\vee \lambda_{1}(\sqrt{\wt{\Lambda}\wt{\Lambda}\tp}\wt{O}\tp \wt{T}\wt{O}\sqrt{\wt{\Lambda}\wt{\Lambda}\tp})=0\vee\lambda_{1}(\wt{\Sigma}^{1/2}\wt{O}\tp \wt{T}\wt{O}\wt{\Sigma}^{1/2})   &   \text{if }n<d,
	\end{cases}
    \notag 
\end{align}
and that for $n<d$ the matrix $\wt{\Sigma}^{1/2}\wt{O}\tp\wt{T}\wt{O}\tp\wt{\Sigma}$ cannot be negative semidefinite due to \Cref{asmp:nonneg}. Then, since we may replace $\wt{O}$ by a matrix Haar-distributed from $\mathbb{SO}(n-1)$, we may directly apply \Cref{prop:str_asymp} to obtain
\begin{align}
	\lim_{n\to\infty}\lambda_{1}(\wh{D})= \sup\supp (\mu_{\wt{\Lambda}\wt{\Lambda}\tp}\boxtimes \mu_{\wt{T}})
    \notag 
\end{align}
almost surely as $n\to\infty$, where we used
\begin{align}
	\mu_{\sqrt{\wt{\Lambda}\wt{\Lambda}\tp}\bsu\tp \wt{T}\bsu\sqrt{\wt{\Lambda}\wt{\Lambda}\tp}}=\mu_{\wt{\Lambda}\wt{\Lambda}\tp}\boxtimes \mu_{\wt{T}},
    \notag 
\end{align}
for a Haar unitary element $\bsu$. 

To finish the proof, one must take the limit $n\to\infty$ in the matrices $\wt{\Lambda}\wt{\Lambda}\tp$ and $\wt{T}$. Indeed, using \cite[Corollary 2.1 and Theorem A.1]{Male}, we have that
\begin{align}
	\sup\supp\mu_{\wt{\Lambda}\wt{\Lambda}\tp}\boxtimes \mu_{\wt{T}}
	\to\sup\supp\mu_{\sfSigma_n}\boxtimes\mu_{\sfT}
    \notag 
\end{align}
almost surely as $n\to\infty$, provided the matrices converges `strongly' i.e.\ for each fixed $\epsilon>0$,
\begin{align}
	\spec(\wt{T})\subset \supp\mu_{\sfT}+(-\epsilon,\epsilon),	\label{eq:str_wt_T}\\
	\spec(\wt{\Lambda}\wt{\Lambda}\tp)\subset \supp\mu_{\sfSigma_n}+(-\epsilon,\epsilon).\label{eq:str_wt_S}
\end{align}
In fact, \Cref{eq:str_wt_T} is always guaranteed by \Cref{asmp:strong}. On the other hand for \Cref{eq:str_wt_S}, notice that 
\begin{align}
	\mu_{\wt{\Lambda}\wt{\Lambda}\tp}
	=\begin{cases}
		\frac{d-1}{n-1}\mu_{\wt{\Sigma}}+\frac{n-d}{n-1}\delta_{0}	&	\text{if }n\geq d, \\
		\mu_{\wt{\Sigma}}	&	\text{if }n<d.
	\end{cases}
    \notag 
\end{align}
Hence, as long as $\delta\neq 1$, the measures $\mu_{\wt{\Lambda}\wt{\Lambda}\tp}$ and $\mu_{\sfSigma_{n}}$ have point masses at the origin of comparable sizes (or both have none). In this case, again \Cref{asmp:strong} guarantees \Cref{eq:str_wt_S}. Therefore for $\delta\neq 1$ we have proved
\begin{align}
	\limsup_{n\to\infty}\lambda_{1}(\wh{D})=\limsup_{n\to\infty}\sup\supp(\mu_{\wt{\Lambda}\wt{\Lambda}\tp}\boxtimes\mu_{\wt{T}})=\sup\supp(\mu_{\sfSigma_n}\boxtimes\mu_{\sfT}).
    \notag 
\end{align}

However \Cref{eq:str_wt_S} may fail when $\delta=1$, since $\mu_{\sfSigma_n}=\mu_{\sfSigma}$ does not have a point mass at zero whereas $\wt{\Lambda}\wt{\Lambda}\tp$ (as well as $\Lambda\Lambda\tp$) might have $o(1)$ nullity. In order to deal with this, we use \Cref{lem:minor_shrink} to get
\begin{align}
	\sup\supp(\mu_{\wt{\Lambda}\wt{\Lambda}\tp}\boxtimes\mu_{\wt{T}})\leq \sup\supp(\mu_{\wt{\Sigma}}\boxtimes\mu_{\wt{T}}).
    \notag 
\end{align}
Now $\mu_{\wt{\Sigma}}$ converges `strongly' to $\mu_{\sfSigma}=\mu_{\sfSigma_n}$ in the sense of \Cref{eq:str_wt_S} by Cauchy interlacing, so that 
\begin{align}
	\sup\supp(\mu_{\wt{\Sigma}}\boxtimes\mu_{\wt{T}})\to \sup\supp (\mu_{\sfSigma_n}\boxtimes\mu_{\sfT})
    \notag 
\end{align} 
almost surely as $n\to\infty$. This completes the proof of \Cref{lem:edge_goal}.

\subsection{Spectral edge of free multiplicative convolution} \label{sec:edge3}

In this section we provide a characterization of the free multiplicative convolution. Since it is of independent interest, we replace $\Sigma$ and $T$ with general operators $\bss$ and $\bst$ that are self-adjoint, free elements in a $W\adj$-probability space $(\caM,\brkt{\cdot})$ such that $\bss\geq0$. 
	Assume the following on the spectral measures $\mu_{\bss}$ and $\mu_{\bst}$. 
	\begin{enumerate}[label=(A\arabic*)]
    \setcounter{enumi}{\value{asmpctr}}
        \item \label[asmp]{asmp:free_nondeg} $\mu_{\bss}$ and $\mu_{\bst}$ are non-degenerate, i.e., not point masses.
		\item \label[asmp]{asmp:free_supp} The map
		\begin{align}
			\C\setminus[0,\infty)\ni z\mapsto \int_{0}^{\infty}\frac{s}{s-z}\dd\mu_{\bss}(s)=\Brkt{\frac{\bss}{\bss-z}}\in\C
            \notag 
		\end{align}
		extends analytically through the origin, i.e.\ the spectral measure $\mu_{\bss}$ is supported on a strictly positive interval modulo a possible point mass at the origin. We write $\wt{\mu}_{\bss}=\mu_{\bss}-\mu_{\bss}(\{0\})\delta_{0}$.
		\item \label[asmp]{asmp:free_nonneg} $\sup\supp\mu_{\bst}>0$.
		\item \label[asmp]{asmp:free_edge} We have
		\begin{align}
			\lim_{z\searrow \sup\supp\mu_{\bss}}\Brkt{\frac{\bss}{\bss-z}}&=-\infty=\lim_{z\nearrow \inf \supp \wt{\mu}_{\bss}}\Brkt{\frac{\bss}{z-\bss}},	\notag \\
			\lim_{z\searrow \sup\supp \mu_{\bst}}\Brkt{\frac{\bst}{\bst-z}}&=-\infty. \notag 
		\end{align}
		
	\setcounter{asmpctr}{\value{enumi}}
    \end{enumerate}
\begin{lemma}\label{lem:a}\
	\begin{enumerate}
		\item\label{itm:finite} Define $\caS_{\bst},\caS'_{\bst}\subset\R$ by
		\begin{align}
			\caS_{\bst}&\deq \left\{w \in\R\setminus\supp\mu_{\bst}:\Brkt{\frac{\bst}{\bst-w}}=0\right\},\notag \\
			\caS'_{\bst}&\deq \left\{w\in\R\setminus\supp\mu_{\bst}:\Brkt{\frac{\bst}{\bst-w}}=1-\mu_{\bss}(\{0\})\right\}.  \notag     
		\end{align}
		Then $\caS_{\bst}\cap(\sup\supp\mu_{\bst},\infty)$ and $\caS_{\bst}'\cap(\sup\supp\mu_{\bst},\infty)$ are finite sets. 
		
		\item\label{itm:omega1} For each $w\in (\sup\supp\mu_{\bst},\infty)\setminus\caS$, there exists a unique $\omega(w)\in\R\setminus(\inf\supp\wt{\mu}_{\bss},\sup\supp\wt{\mu}_{\bss})$ such that
		\begin{align}\label{eq:omega}
		\Brkt{\frac{\bst}{\bst-w}}=\Brkt{\frac{\bss}{\bss-\omega(w)}}.
		\end{align}
		
		\item \label{itm:omega} The map $\omega:(\sup\supp\mu_{\bst},\infty)\setminus\caS_{\bst}\to \R$ defined in \Cref{itm:omega1} extends meromorphically to an open set in $\C$ containing $(\sup\supp\mu_{\bst},\infty)$. The extension is analytic at each $w\in (\sup\supp\mu_{\bst},\infty)\setminus\caS_{\bst}$, has a pole at each $w\in\caS_{\bst}$, and has a zero at each $w\in\caS'_{\bst}$.
		
		\item\label{itm:psi} Define the function $\psi:(\sup\supp\mu_{\bst},\infty)\setminus\caS_{\bst}\to\R$ as
		\begin{align}\label{eq:psi}
		\psi(w)=w\omega(w)\left(1-\Brkt{\frac{\bss}{\bss-\omega(w)}}^{-1}\right)^{-1}\equiv w\omega(w)\left(1-\Brkt{\frac{\bst}{\bst-w}}^{-1}\right)^{-1}.
		\end{align}
		Then $\psi$ extends analytically to an open set in $\C$ containing $(\sup\supp\mu_{\bst},\infty)$, and $\lim_{w\to+\infty}\psi'(w)=\brkt{\bss}$.
	\end{enumerate}
\end{lemma}
\begin{proof}
	For proofs of \Cref{itm:finite,itm:omega,itm:omega1}, see \cite[Items 1 to 3 of Lemma E.2]{Zhang_Ji_Venkataramanan_Mondelli}. For \Cref{itm:psi}, note that we can write $\psi$ as
	\begin{align}
		\psi(w)=w\Brkt{\frac{\bss}{\bss-\omega(w)}}\Brkt{\frac{1}{\bss-\omega(w)}}^{-1}.
        \notag 
	\end{align}
	Thus $\psi$ is a meromorphic function on a neighborhood of $(\sup\supp\mu_{\bst},\infty)$ by \Cref{itm:omega1}, with possible poles at $\caS_{\bst}$. Hence we only need to check that each $w\in\caS_{\bst}$ is a removable singularity for $\psi$. Recall that $\absv{\omega(z)}\to\infty$ as $z\to w\in\caS_{\bst}$, so that by dominated convergence
	\begin{align}\label{eq:remov}
		\psi(z)&=z\Brkt{\frac{\bss\omega(z)}{\bss-\omega(z)}}\Brkt{\frac{\omega(z)}{\bss-\omega(z)}}^{-1}\to\brkt{\bss}w,\qquad \text{as }z\to w.
	\end{align}%
	This shows that $w$ is a removable singularity of $\psi$. 
    For the last part on $\psi'$, straightforward computation yields that whenever $w>\sup\supp\mu_{\bst}$, 
    \begin{align}
            \psi'(w)&= \frac{\dd}{\dd w} \paren{ \omega(w) \bracket{\frac{\bt}{\bt - w}} \bracket{\frac{1}{\bt - w}}^{-1} } \notag \\
            &= \paren{\omega'(w) \bracket{\frac{\bt}{\bt - w}} + \omega(w) \bracket{\frac{\bt}{(\bt - w)^2}}} \bracket{\frac{1}{\bt - w}}^{-1} - \omega(w) \bracket{\frac{\bt}{\bt - w}} \bracket{\frac{1}{(\bt - w)^2}} \bracket{\frac{1}{\bt - w}}^{-2} \notag \\
            &= \paren{ \bracket{\frac{\bt}{(\bt - w)^2}} \bracket{\frac{\bs}{(\bs - \omega(w))^2}}^{-1} \bracket{\frac{\bt}{\bt - w}} + \omega(w) \bracket{\frac{\bt}{(\bt - w)^2}} } \bracket{\frac{1}{\bt - w}}^{-1} \label{eqn:omega'} \\
            &\qquad - \omega(w) \bracket{\frac{\bt}{\bt - w}} \bracket{\frac{1}{(\bt - w)^2}} \bracket{\frac{1}{\bt - w}}^{-2} \notag \\
            &= \paren{ \bracket{\frac{\bs}{(\bs - \omega(w))^2}}^{-1} \bracket{\frac{\bs}{\bs - \omega(w)}} + \omega(w) } \bracket{\frac{\bt}{(\bt - w)^2}} \bracket{\frac{1}{\bt - w}}^{-1} \notag \\
            &\qquad - \omega(w) \bracket{\frac{\bt}{\bt - w}} \bracket{\frac{1}{(\bt - w)^2}} \bracket{\frac{1}{\bt - w}}^{-2} \notag \\
            &=\Brkt{\frac{\bst}{(\bst-w)^{2}}}\Brkt{\frac{\bss}{(\bss-\omega(w))^{2}}}^{-1}\Brkt{\frac{\bss^{2}}{(\bss-\omega(w))^{2}}}\Brkt{\frac{1}{\bst-w}}^{-1} \notag \\
			&\qquad -\Brkt{\frac{\omega(w)\bss}{\bss-\omega(w)}}\Brkt{\frac{1}{\bst-w}}^{-2}\Brkt{\frac{1}{(\bst-w)^{2}}},\label{eq:psi'}
    \end{align}
    where in \Cref{eqn:omega'} we use the following expression of $ \omega'(w) $ derived from the defining equation \Cref{eq:omega} and the implicit function theorem: 
    \begin{align}
        \omega'(w) &= \bracket{\frac{\bt}{(\bt - w)^2}} \bracket{\frac{\bs}{(\bs - \omega(w))^2}}^{-1} . \notag 
    \end{align}
    Using that $\absv{\omega(w)}\to\infty$ as $w\to+\infty$, we easily check that the first term of \Cref{eq:psi'} vanishes and the second term converges to $\brkt{\bss}$ in the limit.
\end{proof}

\begin{theorem}\label{thm:bulk}
	We have 
	\begin{align}
		\sup\supp \mu_{\sqrt{\bss}\bst\sqrt{\bss}}=\psi(a^{\circ}),
        \notag 
	\end{align}
	where $a^{\circ}>\sup\supp\mu_{\bst}$ is the largest critical point of the analytic function $\psi:(\sup\supp\mu_{\bst},\infty)\to\R$ defined above in \Cref{itm:psi} of \Cref{lem:a}. 
\end{theorem}
We prove \Cref{thm:bulk} at the end of this section after collecting the inputs.
\subsubsection{Setup}

Denote $ \C_{+}\deq \{z\in\C:\im z>0\} $ and $ \C_{-}\deq \{z\in\C:\im z<0\} $. 

	\begin{definition} \label{def:KL}
		Let $\bsa\in\caM$ be a non-zero self-adjoint element. Define two analytic maps on $\C_{+}$ as
		\begin{align}
			K_{\bsa}(z)&\deq 1-\Brkt{\frac{\bsa}{\bsa-z}}^{-1},	\notag \\
			L_{\bsa}(z)&\deq \frac{K_{\bsa}(z)}{z}\equiv \Brkt{\frac{1}{\bsa-z}}\Brkt{\frac{\bsa}{\bsa-z}}^{-1}. \notag 
		\end{align}
	\end{definition}
	\begin{remark}\label{rem:S_a}
	Notice that $K_{\bsa}$ is indeed well-defined, that is, the map $\brkt{\bsa/(\bsa-z)}$ has no zero in $\C_{+}$. To see this, suppose on the contrary that there is a zero $z\in\C_{+}$, so that
	\begin{align}
		\int_{\R}\frac{x}{x-z}\dd\mu_{\bsa}(x)=\Brkt{\frac{\bsa}{\bsa-z}}=0.
        \notag 
	\end{align} 
	Then we have 
	\begin{align}
		\im\int_{\R}\frac{x}{x-z}\dd\mu_{\bsa}(x)=\im z\int_{\R}\frac{x}{\absv{x-z}^{2}}\dd\mu_{\bsa}(x)=0
        \notag 
	\end{align}
	and
	\begin{align}
		\re\int_{\R}\frac{x}{x-z}\dd\mu_{\bsa}(x)=\int_{\R}\frac{x^{2}-x\re z}{\absv{x-z}^{2}}\dd\mu_{\bsa}(x)=0.
        \notag 
	\end{align}
	Combining the two equations gives
	\begin{align}
		\int_{\R}\frac{x^{2}}{\absv{x-z}^{2}}\dd\mu_{\bsa}(x)=0,
        \notag 
	\end{align}
	which obviously is a contradiction.
\end{remark}

%
%
	\begin{lemma}\label{lem:self}
		Let $\bsa\in\caM$ be self-adjoint with a non-degenerate spectral measure.
        Let $K_{\bsa}$ and $L_{\bsa}$ be given in \Cref{def:KL}. 
		\begin{enumerate}
			\item \label{itm:self1} $K_{\bsa}$ and $L_{\bsa}$ meromorphically extend to $\C\setminus\supp\mu_{\bsa}$ with $K_{\bsa}(\ol{z})=\ol{K_{\bsa}(z)}$ and $L_{\bsa}(\ol{z})=\ol{L_{\bsa}(z)}$. The two meromorphic extensions have the same set of poles in $\C\setminus\supp\mu_{\bsa}$, given by
			\begin{align}\label{eq:S_a}
				\caS_{\bsa}=\left\{z\in\R\setminus\supp\mu_{\bsa}:\Brkt{\frac{\bsa}{\bsa-z}}=0\right\}.
			\end{align}
			
			\item \label{itm:self2} $L_{\bsa}(\C_{+})\subset \C_{+}$. If, in addition, $\bsa\geq 0$, then $K_{\bsa}(\C_{+})\subset \C_{+}$. 
		\end{enumerate}
	\end{lemma}
	\begin{proof}
		\Cref{itm:self1} is an immediate consequence of the definitions of $K_{\bsa}$ and $L_{\bsa}$ and the fact that the map
		\begin{align}\label{eq:m}
			\C\setminus\supp\mu_{\bsa}\ni z\mapsto \Brkt{\frac{\bsa}{\bsa-z}}\ni\C
		\end{align}
		is analytic. The poles of $K_{\bsa}$ and $L_{\bsa}$ are exactly the zeros of this function; notice that the two maps
		\begin{align}
        &&
			& \Brkt{\frac{1}{\bsa-z}}, &
            & \Brkt{\frac{\bsa}{\bsa-z}} & 
        & \notag 
		\end{align}
		cannot have a common zero. Finally, the zeros of the function in \Cref{eq:m} are all real by \Cref{rem:S_a}, from which \Cref{eq:S_a} follows.
		
		For \Cref{itm:self2}, to check $L_{\bsa}(\C_{+})\subset\C_{+}$, we take $z\in\C_{+}$ and write
		\begin{align}
				\im L_{\bsa}(z)
                &= \im \bracket{\frac{1}{\ba-z}} \bracket{\frac{\ba}{\ba-z}}^{-1}
                = \abs{\bracket{\frac{\ba}{\ba-z}}}^{-2} \im\brack{\bracket{\frac{1}{\ba-z}} \paren{1 + \ol{z} \ol{\bracket{\frac{1}{\ba-z}}}}} \notag \\
				&=\im z\Absv{\int_{\R}\frac{x}{x-z}\dd\mu_{\bsa}(x)}^{-2}\left(\int_{\R}\frac{1}{\absv{x-z}^{2}}\dd\mu_{\bsa}(x)-\Absv{\int_{\R}\frac{1}{x-z}\dd\mu_{\bsa}(x)}^{2}\right)>0, \label{eq:imL} 
		\end{align}
		where we used Cauchy-Schwarz in the last inequality. For $\bsa\geq0$, the fact that $K_{\bsa}(\C_{+})\subset\C_{+}$ follows from 
        observing that
		\begin{align}
			\Brkt{\frac{\bsa}{\bsa-z}}=\int_{\R}\frac{1}{x-z}(x\dd\mu_{\bsa}(x))
            \notag 
		\end{align}
		is the Stieltjes transform of the positive measure $x\dd\mu_{\bsa}(x)$. 
	\end{proof}
\begin{lemma}\label{lem:K_inv}
	The map $K_{\bsa}$ is meromorphic around $\infty$ with a pole at $\infty$. The pole is simple if $\brkt{\bsa}\neq 0$ and double if $\brkt{\bsa}=0$. More precisely, there exist a constant $c>0$ and an analytic function $g$ on $\brace{z : \absv{z}<c}$ such that
	\begin{align}\label{eq:K_pole}
		\frac{1}{K_{\bsa}(1/z)}=\brkt{\bsa}z+(\brkt{\bsa^{2}}-\brkt{\bsa}^{2})z^{2}+z^{3}g(z),\qquad \absv{z}<c.
	\end{align}
	Consequently, there exists a constant $C>0$ such that following holds with the choices
	\begin{align}
		\caN\deq\{z:\absv{z}>C\},\qquad \caC\deq\{z:\absv{z}>C,\,\,\re z>C\absv{\im z}\}.
        \notag 
	\end{align}
	\begin{enumerate}
		\item \label{itm:K_inv1} If $\brkt{\bsa}>0$, there exists a unique analytic inverse $K_{\bsa}^{-1}$ of $K_{\bsa}$ on $\caN$. The map $K_{\bsa}^{-1}$ maps $\caN\cap\C_{+}$, $\caN\cap\C_{-}$, and $\caN\cap\R$ into $\C_{+},\C_{-}$, and $\R$, respectively.
		\item \label{itm:K_inv2} If $\brkt{\bsa}=0$, there exists a unique analytic map $h_{\bsa}$ on $\caN$ such that
		\begin{align}
				K_{\bsa}(h_{\bsa}(\sqrt{z}))=z=K_{\bsa}(h_{\bsa}(-\sqrt{z})),&\qquad z\in\caN,\label{eq:K_inv_h}\\
				\qquad h_{\bsa}(\sqrt{K_{\bsa}(w)})=w,&\qquad w\in\caC,\label{eq:K_inv_h2}
		\end{align}
		where we choose the branch of $\sqrt{z}$ to be $(-\infty,0)$, i.e.\ $\sqrt{r\e{\ii\theta}}=\sqrt{r}\e{\ii\theta/2}$ for $\theta\in(-\pi,\pi]$. The map $h_{\bsa}$ maps $\caC\cap\C_{+}$, $\caC\cap\C_{-}$, and $\caC\cap\R$ into $\C_{+},\C_{-},$ and $\R$, respectively. In this case we denote $K_{\bsa}^{-1}(z)\deq h_{\bsa}(\sqrt{z})$ on a suitable neighborhood of $\infty$. 
	\end{enumerate}
\end{lemma}
\begin{proof}
	Note first that \Cref{eq:K_pole} is a direct consequence of the definition:
	\begin{align}\label{tmp:1}\begin{split}
		1/K_{\bsa}(1/z)=z\Brkt{\frac{\bsa}{1-z\bsa}}\Brkt{\frac{1}{1-z\bsa}}^{-1}
		&=z\left(\brkt{\bsa}+z\brkt{\bsa^{2}}+O(\absv{z}^{2})\right)\left(1+\brkt{\bsa}z+O(\absv{z}^{2})\right)^{-1} \\
		&=\brkt{\bsa}z+(\brkt{\bsa^{2}}-\brkt{\bsa}^{2})z^{2}+O(\absv{z}^{3}).
	\end{split}\end{align}

	For \Cref{itm:K_inv1}, since the linear term in \Cref{eq:K_pole} is non-vanishing, we have a unique analytic inverse $f$ of $1/K_{\bsa}(1/z)$ around the origin. Then it is easy to check that $1/f(1/z)$ serves as the inverse function of $K_{\bsa}$. The last assertion follows from the fact that $\im K_{\bsa}(z)$ has the same sign as $\im z$ for large enough $\absv{z}$, which in turn follows from
	\begin{align}\begin{split}
		\im K_{\bsa}(z)&=\im z\Absv{\Brkt{\frac{\bsa}{\bsa-z}}}^{-2}\Brkt{\frac{\bsa}{\absv{\bsa-z}^{2}}}	
		=\im z\frac{\brkt{\bsa}\absv{z}^{-2}+O(\absv{z}^{-3})}{\absv{\brkt{\bsa}z^{-1}+O(\absv{z}^{-2})}^{2}}	
		=\im z \left(\frac{1}{\brkt{\bsa}}+O(\absv{z}^{-1})\right).
	\end{split}\notag\end{align}
	
	For \Cref{itm:K_inv2},  we can find a unique analytic function $p_{\bsa}(z)$ in a neighborhood of the origin so that
	\begin{align}
    &&
		\frac{1}{K_{\bsa}(1/z)}&=p_{\bsa}(z)^{2},&
        p_{\bsa}(0)&=0,&
        p_{\bsa}'(0)&=\sqrt{\brkt{\bsa^{2}}}>0.&
    & \label{eqn:p} 
	\end{align}
	By the same reasoning as above, the map $p_{\bsa}$ has an inverse $p_{\bsa}^{-1}$ around the origin. Then we choose $h_{\bsa}(z)=1/p_{\bsa}^{-1}(1/z)$ in a neighborhood of $\infty$. Now we have  the first equality of \Cref{eq:K_inv_h} as
	\begin{align}
		K_{\bsa}(h_{\bsa}(\sqrt{z}))=K_{\bsa}(1/p_{\bsa}^{-1}(1/\sqrt{z}))=\frac{1}{p_{\bsa}(p_{\bsa}^{-1}(1/\sqrt{z}))^{2}}=z.
        \notag 
	\end{align}
	Likewise we can prove the second equality of \Cref{eq:K_inv_h}.
	
	To prove \Cref{eq:K_inv_h2}, notice that $p_{\bsa}'(0)>0$ implies that there is a constant $c>0$ such that $z \mapsto p_{\bsa}(z)^{2}$ maps the cone $\caC'\deq \{\absv{w}<c,\absv{\im w}<c\re w\}$ into the slit-plane $\C\setminus(-\infty,0]$. On this cone we have $\sqrt{p_{\bsa}(w)^{2}}=p_{\bsa}(w)$, so that for $w^{-1}\in \caC'$ we have
	\begin{align}
		h_{\bsa}(\sqrt{K_{\bsa}(w)})=\frac{1}{p_{\bsa}^{-1}(1/\sqrt{K_{\bsa}(w)})}=\frac{1}{p_{\bsa}^{-1}(\sqrt{1/K_{\bsa}(w)})}=\frac{1}{p_{\bsa}^{-1}(\sqrt{p_{\bsa}(1/w)^{2}})}=w.
        \notag 
	\end{align}
	Noticing that the inverse-cone $\{w^{-1}:w\in\caC'\}$ contains another a set of the form $\{w:\absv{w}>C, \re w>C\absv{\im w}\}$, we have proved \Cref{eq:K_inv_h2}. The uniqueness follows from those of $p_{\bsa}$ and its inverse. Finally, for the last assertion on the parity of $h_{\bsa}(z)$, we use $\brkt{\bsa}=0$ to write
	\begin{align}\begin{split}\label{eq:K_inv_h3}
		\im K_{\bsa}(w)=\im\left[\Brkt{\frac{\bsa}{\bsa-w}+\frac{\bsa}{w}}^{-1}\right]
		=\im w\Absv{\Brkt{\frac{\bsa^{2}}{w(\bsa-w)}}}^{-2}\Brkt{\frac{\bsa^{2}(2\re w-\bsa)}{\absv{w(\bsa-w)}^{2}}}.
	\end{split}\end{align}
	Plugging in $w=K_{\bsa}^{-1}(z)$ and taking $z$ in a cone
	\begin{align}\label{eq:cone2}
		\{z:\absv{z}>C',\re z>C'\absv{\im z}\},
	\end{align}
	we have
	\begin{align}
		\re K_{\bsa}^{-1}(z)=\re h_{\bsa}(\sqrt{z})> \re[\sqrt{\brkt{\bsa}^{2}}\sqrt{z}(1+O(\absv{z}^{-1/2}))]>\norm{}{\bsa},
        \notag 
	\end{align}
	where we used that, by \Cref{eqn:p}, $h_{\bsa}$ admits the expansion
	 \begin{align}
	 	h_{\bsa}(z)=\frac{1}{p_{\bsa}^{-1}(1/z)}=\left(\frac{1}{p_{\bsa}'(0)z}+O(\absv{z}^{-2})\right)^{-1}=\sqrt{\brkt{\bsa^{2}}}z(1+O(\absv{z}^{-1})).
        \notag 
	 \end{align} 
 	This proves the second factor of \Cref{eq:K_inv_h3} is positive for $w=K_{\bsa}^{-1}(z)$, so that $\im K_{\bsa}^{-1}(z)=\im h_{\bsa}(\sqrt{z})$ has the same sign as $\im z$, hence as $\im \sqrt{z}$.
\end{proof}

\subsubsection{Analytic subordination}

\begin{proposition}[{\cite[Theorem 1.3]{Arizmendi_Hasebe_Kitagawa}}] \label{prop:subor}
	There is a unique pair of analytic functions $\Omega_{1},\Omega_{2}:\C_{+}\to\C$ that satisfies the following for each $z\in\C_{+}$:
		\begin{align}
			\Omega_{1}(z)/z\in\C_{+},\qquad\Omega_{2}(z)\in\C_{+}, \notag \\
			K_{\bss}(\Omega_{1}(z))=K_{\bst}(\Omega_{2}(z))=\frac{\Omega_{1}(z)\Omega_{2}(z)}{z},\label{eq:subor}	\\
			\lim_{r\to+\infty}\absv{\Omega_{1}(rz)}=\infty=\lim_{r\to+\infty}\absv{\Omega_{2}(rz)}. \notag 
		\end{align}
	Furthermore, the three quantities in \Cref{eq:subor} are all equal to
	\begin{align}
		K_{\bss^{1/2}\bst\bss^{1/2}}(z). 
        \notag 
	\end{align}
\end{proposition}

\subsubsection{Consequences of analytic subordination}
In what follows, we write
\begin{equation}
    \lambda^{\circ}\deq \sup\supp\mu_{\sqrt{\bss}\bst\sqrt{\bss}},
    \notag 
\end{equation}
and denote by $\ol{\R}$ the one-point compactification $\R\cup\{\infty\}\cong\bbS^{1}$ of $\R$.
\begin{corollary}\label{cor:subor_ext}
	 The maps $\Omega_{1}$ and (resp.\ $\Omega_{2}$) extends meromorphically (resp.\ analytically) through a neighborhood of the half-line $(\lambda^{\circ},\infty)$. Their images of $(\lambda^{\circ},\infty)$ satisfy $\Omega_{1}((\lambda^{\circ},\infty))\subset\ol{\R}\setminus\supp\mu_{\sfSigma}$, $\Omega_{2}((\lambda^{\circ},\infty))\subset (\sup\supp\mu_{\sfT},\infty)$. Finally, the map $\Omega_{2}$ is strictly increasing in $(\lambda^{\circ},\infty)$.
	 
\end{corollary}
\begin{proof}
	Applying \Cref{itm:omega} of \Cref{lem:a} with $\bst$ replaced by $\bss\bst$, we find that the map $K_{\bss}^{-1}\circ K_{\bss\bst}:(\lambda^{\circ},\infty)\to \ol{\R}\setminus\supp\mu_{\sfSigma}$ is well-defined and meromorphically extends to a neighborhood of $(\lambda^{\circ},\infty)$. On the other hand, \Cref{prop:subor}  shows that $K_{\bss}^{-1}\circ K_{\bss\bst}$ matches with $\Omega_{1}$ in a open set close to $\infty$. Thus $K_{\bss}^{-1}\circ K_{\bss\bst}$ provides a meromorphic extension of $\Omega_{1}$ to a neighborhood of $(\lambda^{\circ},\infty)$. Then we extend $\Omega_{2}$ to the same neighborhood via \Cref{eq:subor}:
	\begin{align}
		\Omega_{2}(z)=\frac{zK_{\bss\bst}(z)}{\Omega_{1}(z)}.
        \notag 
	\end{align}
	It also immediately follows that $\Omega_{1}((\lambda^{\circ},\infty))\subset\ol{\R}\setminus\supp\mu_{\sfSigma}$, and that $\Omega_{2}((\lambda^{\circ},\infty))\subset\ol{\R}$.
	
	We next show that $\Omega_{2}$ has no poles in $(\lambda^{\circ},\infty)$. Assume the contrary, so that there is a pole $x>\lambda^{\circ}$. Then we have
	\begin{align}
	\lim_{z\to x,z\in\C_{+}}\absv{L_{\bss}(\Omega_{1}(z))}=\lim_{z\to x,z\in\C_{+}}\frac{\absv{\Omega_{2}(z)}}{\absv{z}}=\infty,
    \notag 
	\end{align}
	which is a contradiction since $\Omega_{1}(x)\in\ol{\R}\setminus\supp\mu_{\sfSigma}$. Thus $\Omega_{2}$ has no poles in $(\lambda^{\circ},\infty)$ hence analytically continues through the line $(\lambda^{\circ},\infty)$.
	
	Next, we claim that there is a finite Borel measure $\wt{\mu}$ so that
	\begin{align}\label{eq:Pick}
		\Omega_{2}(z)=\frac{z}{\brkt{\bss}}-\frac{\brkt{\bss^{2}}-\brkt{\bss}^{2}}{\brkt{\bss}^{2}}\brkt{\bst}+\int_{\R}\frac{1}{x-z}\dd\wt{\mu}(x),\qquad z\in\C_{+},
	\end{align}
	which is in turn implied by, due to Nevalinna-Pick representation,
	\begin{align}\label{eq:Pick_pf}
		\limsup_{y\to\infty}y\Absv{\Omega_{2}(\ii y)-\frac{1}{\brkt{\bss}}\ii y+\frac{\brkt{\bss^{2}}-\brkt{\bss}^{2}}{\brkt{\bss}^{2}}\brkt{\bst}}<\infty.
	\end{align}
	To see \Cref{eq:Pick_pf}, for large enough $y$ so that $\Omega_{1}(\ii y)=K_{\bss}^{-1}(K_{\bss\bst}(\ii y))$ and $\Omega_{2}(\ii y)=K_{\bst}^{-1}(K_{\bss\bst}(\ii y))$, we notice that
	\begin{align}
		&\Omega_{2}(\ii y)-\frac{1}{\brkt{\bss}}\ii y
        =\ii y \left(L_{\bss}(\Omega_{1}(\ii y))-\frac{1}{\brkt{\bss}}\right) \label{eq:NP_1}.
        \end{align}
       We then use \Cref{eq:K_pole} to get
       \begin{equation*}
        L_{\bss}(\Omega_{1}(\ii y))-\frac{1}{\brkt{\bss}}=\frac{1}{\brkt{\bss}+(\brkt{\bss^{2}}-\brkt{\bss}^{2})\Omega_{1}(\ii y)^{-1}+O(\absv{\Omega_{1}(\ii y)}^{-2})}-\frac{1}{\brkt{\bss}}
        =-\frac{\brkt{\bss^{2}}-\brkt{\bss}^{2}}{\brkt{\bss}^{2}\Omega_{1}(\ii y)}+O(\absv{\Omega_{1}(\ii y)}^{-2}),
        \end{equation*}
        so that
        \begin{align}
        y\Absv{\Omega_{2}(\ii y)-\frac{1}{\brkt{\bss}}\ii y+\frac{\brkt{\bss^{2}}-\brkt{\bss}^{2}}{\brkt{\bss}^{2}}\brkt{\bst}}
        &\leq \frac{y(\brkt{\bss^{2}}-\brkt{\bss}^{2})}{\brkt{\bss}^{2}}\Absv{-\frac{\ii y}{\Omega_{1}(\ii y)}+\brkt{\bst}}+O\left(\frac{y^{2}}{\absv{\Omega_{1}(\ii y)}^{2}}\right) \notag \\
        &\lesssim y\Absv{-\frac{1}{L_{\bst}(\Omega_{2}(\ii y))}+\brkt{\bst}}+\absv{L_{\bst}(\Omega_{2}(\ii y))}^{-2} \notag \\
        &\lesssim \frac{y}{\absv{\Omega_{2}(\ii y)}}+\absv{L_{\bst}(\Omega_{2}(\ii y))}^{-2} \notag \\
        &=\absv{L_{\bss}(\Omega_{1}(\ii y)}^{-1}+\absv{L_{\bst}(\Omega_{2}(\ii y))}^{-2}, \notag 
        \end{align}
        where the third line follows from \Cref{eq:K_pole} applied to $\bst$ in place of $\bsa$. 
    This proves \Cref{eq:Pick_pf} as
	\begin{align}
		\lim_{y\to\infty}L_{\bss}(\Omega_{1}(\ii y))^{-1}=\lim_{\absv{\omega}\to\infty} L_{\bss}(\omega)^{-1}=\brkt{\bss}, \quad\text{and}\quad
        \lim_{y\to\infty}L_{\bst}(\Omega_{2}(\ii y))^{-1}=\lim_{\absv{\omega_{2}}\to\infty} L_{\bst}(\omega_{2})^{-1}=\brkt{\bst}.
        \notag 
	\end{align}
	
	Applying Stieltjes inversion to \Cref{eq:Pick}, since $\Omega_{2}$ extends analytically through $(\lambda^{\circ},\infty)$ and is real-valued therein, we find that $\supp\wt{\mu}\subset (-\infty,\lambda^{\circ}]$. Then it also follows that $\Omega_{2}$ is strictly increasing in $(\lambda^{\circ},\infty)$. If $\Omega_{2}(x)=\sup\supp\mu_{\sfT}$ for some $x>\lambda^{\circ}$, we would have
	\begin{align}
		-\infty=\lim_{w\searrow \sup\supp\mu_{\sfT}}\Brkt{\frac{\bst}{\bst-w}}=\lim_{z\searrow x}\Brkt{\frac{\bst}{\bst-\Omega_{2}(z)}}=\lim_{z\searrow x}\Brkt{\frac{\bss}{\bss-\Omega_{1}(z)}},
        \notag 
	\end{align}
	which in turn implies $\lim_{z\searrow x}\Omega_{1}(z)=\sup\supp\mu_{\sfSigma}$. But then we have a contradiction as
	\begin{align}
		0&=\lim_{z\searrow x}K_{\bss\bst}(z)=\lim_{z\searrow x}\frac{\Omega_{1}(z)\Omega_{2}(z)}{z}=\frac{(\sup\supp\mu_{\sfT})(\sup\supp\mu_{\sfSigma})}{x}.
        \notag 
        \qedhere
	\end{align}
\end{proof}

\subsubsection{Proof of \Cref{thm:bulk}} 
\Cref{thm:bulk} follows from the following lemma. 
\begin{lemma}\label{lem:bulk_main}
	Let $w>\sup\supp\mu_{\sfT}$. Recall the map $\psi$ defined in \Cref{itm:psi} of \Cref{lem:a}.
	\begin{enumerate}
		\item\label{itm:bulk_upper} If $\psi(\wt{w})>\lambda^{\circ}$ for all $\wt{w}\geq w$, then $\psi'(w)>0$.
		\item\label{itm:bulk_lower} If $\psi'(\wt{w})>0$ for all $\wt{w}\geq w$, then $\psi(\wt{w})\geq\lambda^{\circ}$.
	\end{enumerate}
\end{lemma}
We omit the proof of \Cref{thm:bulk} given \Cref{lem:bulk_main} since it is exactly the same as \cite[proof of Lemma 5.5 given Lemma E.1]{Zhang_Ji_Venkataramanan_Mondelli}.

\begin{proof}[Proof of \Cref{lem:bulk_main}]
	Recall from \Cref{lem:a} the triplet $(\psi(\wt{w}),\omega(\wt{w}),\wt{w})$ of meromorphic functions in $\wt{w}\in(\sup\supp\mu_{\sfT},\infty)$. It also follows from the definitions of $\psi$ and $\omega$ that the triplet satisfies the same subordination equation as in \Cref{eq:subor}, i.e.
	\begin{align}\label{eq:subor_1}
		K_{\bss}(\omega(\wt{w}))=K_{\bst}(\wt{w})=\frac{\wt{w}\omega(\wt{w})}{\psi(\wt{w})},\qquad \wt{w}>\sup\supp\mu_{\sfT}.
	\end{align}

	Recall from the proof of \Cref{lem:a} that $\psi'(\wt{w})$ is positive for large enough $\wt{w}>0$. We next claim that the same holds for $(\omega/\psi)'$. Note first that
	\begin{align}\begin{split}
		\frac{\omega(\wt{w})}{\psi(\wt{w})}=\frac{\omega(\wt{w})}{\wt{w}}L_{\bss}(\omega(\wt{w}))
		=\frac{1}{\wt{w}}K_{\bss}(\omega(\wt{w}))=\frac{1}{\wt{w}}K_{\bst}(\wt{w})=L_{\bst}(\wt{w}),
	\end{split}\notag\end{align}
	where we used the definition of $\omega$ in the last equality. Then we can compute the derivative easily using \Cref{eq:imL} via Nevanlinna derivative:
	\begin{align}\begin{split}
		\frac{\dd}{\dd \wt{w}}L_{\bst}(\wt{w})=\lim_{\eta\searrow 0}\frac{\im L_{\bst}(\wt{w}+\ii\eta)}{\eta}	
		=\Brkt{\frac{\bst}{\bst-\wt{w}}}^{-2}\left(\Brkt{\frac{1}{(\bst-\wt{w})^{2}}}-\Brkt{\frac{1}{\bst-\wt{w}}}^{2}\right)>0.
	\end{split}\notag\end{align}
	Now by taking large $\wt{w}>0$ at which derivatives of $\psi$ and $\omega/\psi$ are both positive, we find that there is a small, open sector in $\C_{+}$ based at $\wt{w}$ on which all of $(\psi(\wt{w}),\omega(\wt{w})/\psi(\wt{w}),\wt{w})$ are in $\C_{+}$. By analytic continuation the subordination equations \Cref{eq:subor_1} hold on the sector, so that by the uniqueness in \Cref{prop:subor} we have
	\begin{align}
		\Omega_{2}(\psi(\wt{w}))&=\wt{w},\label{eq:bulk_upper} \\ \Omega_{1}(\psi(\wt{w}))&=\omega(\wt{w}), \notag
	\end{align}
	on a neighborhood of $(\sup\supp\mu_{\sfT},\infty)$.
	
	We start with the proof of \Cref{itm:bulk_upper}. Taking the derivative of \Cref{eq:bulk_upper} and noticing from \Cref{cor:subor_ext} that $\Omega_{2}$ is strictly increasing in $(\lambda^{\circ},\infty)$ proves that $\psi'(\wt{w})>0$ for all $\wt{w}\geq w$.
	
	For \Cref{itm:bulk_lower}, we consider the analytic inverse function $\psi^{-1}:[\psi(w),\infty)\to (w,\infty)$. By \Cref{eq:bulk_upper} this inverse function must match with $\Omega_{2}$, so that $\Omega_{2}$ extends analytically to a neighborhood of $[\psi(w),\infty)$ and is real-valued on the interval. If we assume on the contrary that $\lambda^{\circ}>\psi(w)$, then for $x\in(\psi(w),\lambda^{\circ})\setminus \{0\}$ we would have
	\begin{align}
		\lim_{\eta\searrow 0}\im\Brkt{\frac{1}{\bss\bst-x-\ii\eta}}&=\frac{1}{x}\lim_{\eta\searrow0}\im\Brkt{\frac{\bss\bst}{\bss\bst-x-\ii\eta}} 
		=\frac{1}{x}\lim_{\eta\searrow 0}\im\Brkt{\frac{\bst}{\bst-\Omega_{2}(x+\ii\eta)}}=\frac{1}{x}\im\Brkt{\frac{\bst}{\bst-\Omega_{2}(x)}}=0, \notag
	\end{align}
	where the last two equalities follow from \Cref{cor:subor_ext}. By Stieltjes inversion this would imply that $(\psi(w),\lambda^{\circ})\setminus\{0\}$ does not intersect with $\supp\mu_{\bss\bst}$, contradicting the definition $\lambda^{\circ}=\sup\supp\mu_{\bss\bst}$.
\end{proof}

\subsection{Specialization to our setting} \label{sec:edge4}
Our choice of $T$ shall be given by \Cref{eqn:D} with the function $\cT$ therein satisfying \Cref{asmp:preprocess}. 
This implies $\sup\supp \mu_{\sfT}<1$, and also that
\begin{align}
	\expt{\frac{\sfT}{\sfT-1}}=\frac{1}{\gamma}\expt{\ol{g}(\sfY)}=0,
\end{align}
where $ \sfT = \cT(\sfY) $ (recall $\sfY$ from \Cref{eqn:BZEY}). 
    In what follows, in light of \Cref{lem:edge_goal}, we apply \Cref{thm:bulk} with the choices 
    \begin{align}
    &&
        \mu_{\bss}&=\mu_{\sfSigma_{n}},&
        \mu_{\bst}&=\mu_{\sfT}.&
    & \label{eq:s_t_choices}
    \end{align}
    Note that under these choices, all of \Cref{asmp:free_nondeg,asmp:free_supp,asmp:free_nonneg,asmp:free_edge} required in \Cref{sec:edge3} have been covered by \Cref{asmp:design,asmp:preprocess}.  
    In particular, by recalling $\sfSigma_n$ from \Cref{eqn:mu_sfLambdan2} and its relation to $ \sfLambda_n $ in \Cref{eqn:sfLambda_nd}, one may identify $ \wt{\mu}_{\bs} $ in \Cref{asmp:free_supp} as 
    \begin{align}
        \wt{\mu}_{\bs} &= \begin{cases}
            \law(\sfLambda^2) , & \delta\le1 \\
            \delta^{-1} \law(\sfLambda^2) , & \delta > 1
        \end{cases} . \notag 
    \end{align}

\begin{theorem}\label{thm:equiv}
    Let $ \ol{\psi}(a) = \psi\paren{\max\brace{a,a^\circ}} $. 
    Then the following three conditions are equivalent. 
    \begin{enumerate}
        \item \label{itm:cond1} $ a^\circ < 1 $; 
        \item \label{itm:cond2} $ \ol{\psi}(a^\circ) < \ol{\psi}(1) $, or more explicitly, $ \lambda^\circ < \ol{\kappa}_2 $; 
        \item \label{itm:cond3} $ \ol{\psi}'(1) > 0 $, or more explicitly, 
        \begin{align}
            \paren{\frac{\ol{\kappa}_4}{\ol{\kappa}_2^2} + \delta} \expt{\frac{\sfT^2}{(\sfT - 1)^2}} &< 1 . \label{eqn:cond3_equiv}
        \end{align}
    \end{enumerate}
\end{theorem}


\begin{proof}
	Notice that $\psi$ is eventually increasing by \Cref{itm:psi} of \Cref{lem:a}. 
	On the other hand since $a^{\circ}$ is the largest critical point, it follows that $\psi$ is increasing in $[a^{\circ},\infty)$. 
    Hence the equivalence of $ a^\circ < 1 $, $ \ol{\psi}(a^\circ) < \ol{\psi}(1) $ and $ \ol{\psi}'(1) > 0 $ holds. 

    For the explicit form of \Cref{itm:cond2}, we note that when $ a^\circ < 1 $, by the definition of $ \ol{\psi} $, it holds that $ \ol{\psi}(1) = \psi(1) $.
	We next compute the value of $\psi(1)$. Since 
	\begin{align}
		\Brkt{\frac{\bst}{\bst-1}}=(1\wedge\delta)\expt{\frac{\sfT}{\sfT-1}}=0,
        \notag 
	\end{align}
	we have $1\in\caS_{\bst}$ hence $\omega$ has a pole at $1$. Thus, by \Cref{eq:remov} we have $\psi(1)=\brkt{\bss} = \expt{\sfSigma_n} = \expt{\sfLambda_n^2} = \ol{m}_2 = \ol{\kappa}_2$.

    For the explicit form of \Cref{itm:cond3}, note that when $ a^\circ < 1 $, it holds that $ \ol{\psi}'(1) = \psi'(1) $. 
	We finally compute $\psi'(1)$ by taking the limit $a\to 1$ in \Cref{eq:psi'}. For the first term, we use $\absv{\omega(a)}\to\infty$ to get
	\begin{align}
			& \lim_{a\to 1}\Brkt{\frac{\bst}{(\bst-a)^{2}}}\Brkt{\frac{\bss}{(\bss-\omega(a))^{2}}}^{-1}\Brkt{\frac{\bss^{2}}{(\bss-\omega(a))^{2}}}\Brkt{\frac{1}{\bst-a}}^{-1} \notag \\
			&=\Brkt{\frac{\bst}{(\bst-1)^{2}}}\Brkt{\frac{1}{\bst-1}}^{-1}\lim_{a\to 1}\Brkt{\frac{\bss\omega(a)^{2}}{(\bss-\omega(a))^{2}}}^{-1}\Brkt{\frac{\bss^{2}\omega(a)^{2}}{(\bss-\omega(a))^{2}}}	
			=\Brkt{\frac{\bst}{(\bst-1)^{2}}}\Brkt{\frac{1}{\bst-1}}^{-1}\frac{\brkt{\bss^{2}}}{\brkt{\bss}}.
            \notag 
        \end{align}
	For the second term, we have
	\begin{align}
		\lim_{a\to 1}\Brkt{\frac{\omega(a)\bss}{\bss-\omega(a)}}\Brkt{\frac{1}{\bst-a}}^{-2}\Brkt{\frac{1}{(\bst-a)^{2}}}
		=-\brkt{\bss}\Brkt{\frac{1}{\bst-1}}^{-2}\Brkt{\frac{1}{(\bst-1)^{2}}}.
        \notag 
	\end{align} 
	Thus we conclude
	\begin{align}
		\psi'(1)=\brkt{\bss}\Brkt{\frac{1}{\bst-1}}^{-2}\left(\frac{\brkt{\bss^{2}}}{\brkt{\bss}^{2}}\Brkt{\frac{\bst}{(\bst-1)^{2}}}\Brkt{\frac{1}{\bst-1}}+\Brkt{\frac{1}{(\bst-1)^{2}}}\right) . 
        \notag 
	\end{align}
	Note that $\brkt{\bst/(\bst-1)}=0$ implies 
	\begin{align}
    &&
		\Brkt{\frac{1}{\bst-1}}&=-1,&
		\Brkt{\frac{1}{(\bst-1)^{2}}}&=1+\Brkt{\frac{\bst}{(\bst-1)^{2}}}=1+\Brkt{\frac{\bst^{2}}{(\bst-1)^{2}}},&
    & \notag 
	\end{align}
	from which we can further simplify the last expression as
	\begin{align}
			\psi'(1)=\brkt{\bss}\left(\left(-\frac{\brkt{\bss^{2}}}{\brkt{\bss}^{2}}+1\right)\Brkt{\frac{\bst^{2}}{(\bst-1)^{2}}}+1\right).
            \notag 
	\end{align}
	Recalling the definitions of $\bss$ and $\bst$ from \Cref{eq:s_t_choices}, we compute
	\begin{align}
		\psi'(1)
        &=\expt{\sfSigma_{n}}\left(\left(-\frac{\expt{\sfSigma_{n}^{2}}}{\expt{\sfSigma_{n}}^{2}}+1\right)\expt{\frac{\sfT^{2}}{(\sfT-1)^{2}}}+1\right) \notag \\
        &=\ol{m}_{2}\left(1-\left(\frac{\ol{m}_{4}-\ol{m}_{2}^{2}}{\ol{m}_{2}^{2}}\right)\expt{\frac{\sfT^{2}}{(\sfT-1)^{2}}}\right) \notag \\
        &=\ol{\kappa}_{2}\left(1-\left(\frac{\ol{\kappa}_{4}}{\ol{\kappa}_{2}^{2}}+\delta\right)\expt{\frac{\sfT^{2}}{(\sfT-1)^{2}}}\right).
        \notag 
 	\end{align}
    Therefore, $ \psi'(1) > 0 $ is equivalent to \Cref{eqn:cond3_equiv}. 
\end{proof}


\section{Proof of \Cref{lem:align}}
\label{app:pf_lem:align}

Recall from \Cref{eqn:GVAMP_lin} the linearized GVAMP iteration
\begin{align}
u^{t} &= \paren{ \ol{\kappa}_2^{-1} X X^\top - I_n } \paren{ \wt{G} - I_n } u^{t-1} , \label{eqn:recur1}
\end{align}
where 
\begin{align}
    \wt{G} \coloneqq \diag(\ol{g}(y)) + I_n . \label{eqn:wtG}
\end{align}
With the notation 
\begin{align}
&&
    \beta_t &\coloneqq \normtwo{u^t} , &
    \wh{u}^t &\coloneqq \frac{u^t}{\beta_t} , &
& \label{eqn:beta_whu} 
\end{align}
we can write \Cref{eqn:recur1} as 
\begin{align}
\frac{\beta_t}{\beta_{t-1}} \wh{u}^t &= 
\paren{ \ol{\kappa}_2^{-1} X X^\top - I_n } \paren{ \wt{G} - I_n } \wh{u}^{t-1} . \label{eqn:hut_hut1}
\end{align}

\begin{lemma}
Let $ \beta_t $ be defined in \Cref{eqn:beta_whu}. 
The following limits hold: 
\label{lem:technical}
\begin{align}
&&
    \lim_{t\to\infty} \lim_{n\to\infty} \frac{\inprod{u^t}{u^{t-1}}}{\normtwo{u^{t-1}}^2} 
    = \lim_{t\to\infty} \lim_{n\to\infty} \frac{\beta_{t}}{\beta_{t-1}}
    &= \gamma, &
& \label{eqn:gamma}
\end{align}
where $\gamma$ is defined in \Cref{eqn:T}. 
\end{lemma}
The proof of \Cref{lem:technical} is deferred to \Cref{app:pf_lem:technical}. 

Adding $ (\wt{G} - I_n) \wh{u}^t $ to  both sides of \Cref{eqn:hut_hut1}, we have
\begin{align}
\brack{ \paren{ \frac{\beta_t}{\beta_{t-1}} - 1 } I_n + \wt{G} } \wh{u}^t &= \ol{\kappa}_2^{-1} X X^\top \paren{ \wt{G} - I_n } \wh{u}^{t-1} + \paren{ \wt{G} - I_n } \paren{ \wh{u}^t - \wh{u}^{t-1} } . \notag 
\end{align}
Before proceeding, recall from the definitions \Cref{eqn:wtG} of $ \wt{G} $ that 
\begin{align}
    \paren{ \frac{\beta_t}{\beta_{t-1}} - 1 } I_n + \wt{G}
    &= \frac{\beta_t}{\beta_{t-1}} I_n + \diag(\ol{g}(y)) . \label{eqn:G_invertible} 
\end{align}
By the second identity of \Cref{eqn:gamma} and the first assumption on $\ol{g}$ in \Cref{eqn:gbar}, the above matrix is strictly positive definite with probability $1$ for all sufficiently large $n$. 
Now, let us multiply both sides by $ X^\top \paren{ \wt{G} - I_n } \brack{ \paren{ \frac{\beta_t}{\beta_{t-1}} - 1 } I_n + \wt{G} }^{-1} $: 
\begin{align}
X^\top \paren{ \wt{G} - I_n } \wh{u}^t &= \ol{\kappa}_2^{-1} X^\top \paren{ \wt{G} - I_n } \brack{ \paren{ \frac{\beta_t}{\beta_{t-1}} - 1 } I_n + \wt{G} }^{-1} X X^\top \paren{ \wt{G} - I_n } \wh{u}^{t-1} \notag \\
&\quad + X^\top \paren{ \wt{G} - I_n }^2 \brack{ \paren{ \frac{\beta_t}{\beta_{t-1}} - 1 } I_n + \wt{G} }^{-1} \paren{ \wh{u}^t - \wh{u}^{t-1} } . \notag 
\end{align}
Further define
\begin{align}
\wh{v}^{t+1} &\coloneqq X^\top \paren{ \wt{G} - I_n } \wh{u}^t , \label{eqn:vhat} 
\end{align}
and rewrite the above iteration as: 
\begin{align}
\wh{v}^{t+1} &= \ol{\kappa}_2^{-1} X^\top \paren{ \wt{G} - I_n } \brack{ \paren{ \gamma - 1 } I_n + \wt{G} }^{-1} X \wh{v}^t \notag \\
&\quad + \underbrace{ \ol{\kappa}_2^{-1} X^\top \paren{ \wt{G} - I_n } \brace{ \brack{ \paren{ \frac{\beta_t}{\beta_{t-1}} - 1 } I_n + \wt{G} }^{-1} - \brack{ \paren{ \gamma - 1 } I_n + \wt{G} }^{-1} } X \wh{v}^t }_{\wh{e}_2^t} \label{eqn:e2} \\
&\quad + \underbrace{ X^\top \paren{ \wt{G} - I_n }^2 \brack{ \paren{ \frac{\beta_t}{\beta_{t-1}} - 1 } I_n + \wt{G} }^{-1} \paren{ \wh{u}^t - \wh{u}^{t-1} } }_{\wh{e}_1^t} . \label{eqn:e1}
\end{align}
Recalling the definition of $D$ in \Cref{eqn:D}, the form of $\cT$ in \Cref{eqn:T} and the notation $\wt{G}$ in \Cref{eqn:wtG}, we recognize  
\begin{align}
\ol{\kappa}_2^{-1} X^\top \paren{ \wt{G} - I_n } \brack{ \paren{ \gamma - 1 } I_n + \wt{G} }^{-1} X
&= \ol{\kappa}_2^{-1} D . \notag 
\end{align}
So we arrive at 
\begin{align}
\wh{v}^{t+1} &= \ol{\kappa}_2^{-1} D \wh{v}^t + \wh{e}_2^t + \wh{e}_1^t . \notag 
\end{align}
We shift the spectrum of $\ol{\kappa}_2^{-1} D$ to the right so that the resulting matrix is positive definite (with probability $1$ for all sufficiently large $n$). 
Take $ \ell > 0 $ to be a sufficiently large constant such that $ \prob{\lambda_d(\ol{\kappa}_2^{-1} D + \ell I_d) > 0} = 1 $ for all sufficiently large $d$. 
Then, the iteration can be re-written as 
\begin{align}
\wh{v}^{t+1} &= \underbrace{\frac{\ol{\kappa}_2^{-1} D + \ell I_d}{1 + \ell}}_{\wt{D}} \wh{v}^t + \underbrace{\frac{\ell}{1 + \ell} \paren{ \wh{v}^{t+1} - \wh{v}^t } + \frac{1}{1 + \ell} \paren{ \wh{e}_2^t + \wh{e}_1^t } }_{\wh{e}_3^t} . \label{eqn:Dtilde_e3} 
\end{align}
We note that $D$ and $\wt{D}$ have the same eigenvectors.
Unrolling the above iteration for $ s $ steps, we get
\begin{align}
\wh{v}^{t+s} &= 
\wt{D}^s \wh{v}^t + \wh{e}^{t,s} , \label{eqn:norm_TODO} 
\end{align}
where 
\begin{align}
\wh{e}^{t,s} &\coloneqq \sum_{i = 0}^{s-1} 
\wt{D}^{s-1-i} \wh{e}_3^{t+i} . \label{eqn:ets} 
\end{align}

Let us examine the limiting squared norm of both sides of \Cref{eqn:norm_TODO}. 
For the left-hand side, we have  
\begin{align}
&\lim_{s\to\infty} \lim_{t\to\infty} \lim_{d\to\infty} \normtwo{\wh{v}^{t+s}}^2
=  \lim_{t\to\infty} \lim_{d\to\infty} \normtwo{\wh{v}^{t+1}}^2
=\lim_{t\to\infty} \lim_{d\to\infty} \frac{\normtwo{X^\top \paren{\wt{G} - I_n} u^t}^2}{\normtwo{u^t}^2}
\stackrel{(a)}{=} \lim_{t\to\infty} \lim_{d\to\infty} \frac{\ol{\kappa}_2^2 \normtwo{v^{t+1}}^2}{\normtwo{u^t}^2}\notag \\
& \stackrel{(b)}{=} \lim_{t\to\infty} \frac{\ol{\kappa}_2^2 \paren{ \chi_{t+1}^2 \rho + \omega_{t+1}^2 }}{\nu_t^2 \sigma^2 + \tau_t^2} 
= \lim_{t\to\infty} \frac{\chi_{t+1}^2}{\nu_t^2} \frac{\ol{\kappa}_2^2 \paren{ \rho + \frac{\omega_{t+1}^2}{\chi_{t+1}^2} }}{ \sigma^2 + \frac{\tau_t^2}{\nu_t^2}} 
\stackrel{(c)}{=} w_4^2 \frac{\ol{\kappa}_2^2 \paren{\rho + w_2}}{\sigma^2 + w_1} 
\in (0,\infty) , \label{eqn:l} 
\end{align}
where the equality $(a)$ is obtained by recalling the update rules for $ v^t $ in the linearized GVAMP dynamics \Cref{eqn:GVAMP_lin}, $(b)$ follows from \Cref{prop:SE,prop:BZ_equiv}, and $(c)$ from \Cref{eqn:SE_stay}. 

We then turn to the right-hand side of \Cref{eqn:norm_TODO}. 
The following lemma, whose proof is deferred to \Cref{app:error_estimates},  shows that the error term $ \wh{e}^{t,s} $ is asymptotically negligible. 
\begin{lemma}
\label{lem:whe_ts}
The error vector $ \wh{e}^{t,s} $ defined in \Cref{eqn:ets} satisfies
\begin{align}
\lim_{s\to\infty} \lim_{t\to\infty} \lim_{n\to\infty} \normtwo{ \wh{e}^{t,s} }^2
= 0 . \notag 
\end{align}
\end{lemma}

By \Cref{lem:whe_ts} and the spectral decomposition of $\wt{D}$, we have
\begin{align}
& \lim_{s\to\infty} \lim_{t\to\infty} \lim_{d\to\infty} \normtwo{ \wt{D}^s \wh{v}^t + \wh{e}^{t,s} }^2
= \lim_{s\to\infty} \lim_{t\to\infty} \lim_{d\to\infty} \normtwo{ \wt{D}^s \wh{v}^t }^2 \notag \\
&= \lim_{s\to\infty} \lim_{t\to\infty} \lim_{d\to\infty} \lambda_1(\wt{D})^{2s} \inprod{v_1(D)}{\wh{v}^t}^2 
+ \sum_{i = 2}^d \lambda_i(\wt{D})^{2s} \inprod{v_i(D)}{\wh{v}^t}^2 
\notag \\
&= \lim_{s\to\infty} \lim_{t\to\infty} \lim_{d\to\infty} \normtwo{ \wt{D}^s \Pi \wh{v}^t }^2 + \normtwo{ \wt{D}^s \Pi^\perp \wh{v}^t }^2 , \label{eqn:term12} 
\end{align}
where 
\begin{align}
&&
    \Pi &\coloneqq v_1(D) v_1(D)^\top , & 
    \Pi^\perp &\coloneqq I_d - \Pi . & 
& \notag 
\end{align}
The second term in \Cref{eqn:term12} can be upper bounded as 
\begin{align}
\normtwo{ \wt{D}^s \Pi^\perp \wh{v}^t }^2
&= \normtwo{\paren{\wt{D} \Pi^\perp}^s \wh{v}^t}^2
\le \normtwo{\wh{v}^t}^2 \sigma_1\paren{\paren{\wt{D} \Pi^\perp}^s}^{2}
= \normtwo{\wh{v}^t}^2 \lambda_1\paren{\wt{D} \Pi^\perp}^{2s}
= \normtwo{\wh{v}^t}^2 \lambda_2(\wt{D})^{2s} , \label{eqn:Dperp} 
\end{align}
where we have used the positive definiteness of $ \wt{D} \Pi^\perp $ in the penultimate equality. 

By \Cref{lem:edge},
\begin{align}
\lim_{d\to\infty} \lambda_2(D) &= \lambda^\circ < \ol{\kappa}_2 , \notag 
\end{align}
where the strict inequality holds under the assumption \Cref{eqn:thr}.
In view of the definition of $\wt{D}$ in \Cref{eqn:Dtilde_e3}, this implies
\begin{align}
\lim_{d\to\infty} \lambda_2(\wt{D}) &= \frac{\ol{\kappa}_2^{-1} \lambda^\circ + \ell}{1 + \ell} < 1 . \notag 
\end{align}
Getting back to \Cref{eqn:Dperp}, using the above bound and the boundedness of $ \lim_{t\to\infty} \lim_{d\to\infty} \normtwo{\wh{v}^t} $ shown in \Cref{eqn:l}, we obtain
\begin{align}
\lim_{s\to\infty} \lim_{t\to\infty} \lim_{d\to\infty} \normtwo{ \wt{D}^s \Pi^\perp \wh{v}^t }^2
\le \paren{ \lim_{t\to\infty} \lim_{d\to\infty} \normtwo{\wh{v}^t}^2 } \paren{\lim_{s\to\infty} \lim_{d\to\infty} \lambda_2(\wt{D})^{2s}} &= 0 . \label{eqn:term2} 
\end{align}
Using \Cref{eqn:term2} in \Cref{eqn:term12} yields
\begin{align}
\lim_{s\to\infty} \lim_{t\to\infty} \lim_{d\to\infty} \normtwo{\wt{D}^s \wh{v}^t + \wh{e}^{t,s}}^2
&= \lim_{s\to\infty} \lim_{t\to\infty} \lim_{d\to\infty} \lambda_1(\wt{D})^{2s} \inprod{v_1(D)}{\wh{v}^t}^2 . \label{eqn:r} 
\end{align}

Finally, combining \Cref{eqn:r,eqn:l} brings us to 
\begin{align}
1 &=  \frac{\displaystyle \lim_{s\to\infty} \lim_{t\to\infty} \lim_{d\to\infty} \lambda_1(\wt{D})^{2s}\inprod{v_1(D)}{\wh{v}^t}^2}{\displaystyle \lim_{t\to\infty} \lim_{d\to\infty} \normtwo{\wh{v}^t}^2} ,  \notag 
\end{align}
implying
\begin{align}
&&
    \lim_{d\to\infty} \lambda_1(\wt{D}) &= 1 , & 
    \lim_{t\to\infty} \lim_{d\to\infty} \frac{\inprod{v_1(D)}{\wh{v}^t}^2}{\normtwo{\wh{v}^t}^2} &= 1 . &
& \notag 
\end{align}
Recalling the definitions of $ \wt{D} $ in \Cref{eqn:Dtilde_e3} and noting from \Cref{eqn:vhat} that $ \wh{v}^t = \ol{\kappa}_2 v^t $, we conclude: 
\begin{align}
&&
    \lim_{d\to\infty} \lambda_1(D) &= \ol{\kappa}_2 , & 
    \lim_{t\to\infty} \lim_{d\to\infty} \frac{\inprod{v_1(D)}{v^t}^2}{\normtwo{v^t}^2} &= 1 , & 
& \notag
\end{align}
as claimed in \Cref{lem:align}. 

\subsection{Proof of \Cref{lem:technical}}
\label{app:pf_lem:technical}

Recalling the state evolution representation in \Cref{prop:BZ_equiv}, let
\begin{align}
n_t &\coloneqq \expt{ \sfN_t \sfN_{t-1} }
= \frac{1}{\tau_t \tau_{t-1}} \paren{ \Sigma_{t+1,t} - \nu_t \nu_{t-1} \sigma^2 } . \label{eqn:def_nt}  
\end{align}
We first show 
\begin{align}
    \lim_{t\to\infty}n_t &= 1 . \label{eqn:lim_nt}
\end{align}

To this end, we begin by deriving a recursion for $ n_t $. 
From the state evolution recursions \Cref{eqn:Sigmati,eqn:Omegati}, 
\begin{align}
\Sigma_{t+1,t} &= \ol{\kappa}_2 \paren{\chi_t \chi_{t-1} \rho + \Omega_{t,t-1}} + \frac{\ol{\kappa}_4}{\ol{\kappa}_2} \rho \expt{\ol{g}(\sfY) \ol{\sfZ}^2} \paren{\nu_{t-1} \chi_{t-1} + \nu_{t-2} \chi_t} \notag \\
&\quad + \frac{\ol{\kappa}_4}{\ol{\kappa}_2^2} \expt{\ol{g}(\sfY)^2} \expt{\tau_{t-1} \sfN_{t-1} \cdot \tau_{t-2} \sfN_{t-2}} 
+ \paren{
    \frac{\ol{\kappa}_4}{\ol{\kappa}_2} \expt{\ol{g}(\sfY)^2 \ol{\sfZ}^2} 
    + \frac{\ol{\kappa}_6}{\ol{\kappa}_2^2} \expt{\ol{g}(\sfY) \ol{\sfZ}^2}^2
} \rho \nu_{t-1} \nu_{t-2} , \notag \\
\Omega_{t+1,t} &= \delta \brace{
    \rho \paren{\expt{\ol{g}(\sfY)^2 \ol{\sfZ}^2} + \frac{\ol{\kappa}_4}{\ol{\kappa}_2^2} \expt{\ol{g}(\sfY) \ol{\sfZ}^2}^2} \nu_{t} \nu_{t-1} 
    + \frac{1}{\ol{\kappa}_2} \expt{\ol{g}(\sfY)^2} \expt{\tau_{t} \sfN_{t} \cdot \tau_{t-1} \sfN_{t-1}}
} . \notag 
\end{align}
Using the definition \Cref{eqn:def_nt} of $ n_t $, we obtain
\begin{align}
n_t 
&= \frac{1}{\tau_t \tau_{t-1}} \bigg\{ \ol{\kappa}_2 \paren{\chi_t \chi_{t-1} \rho + \Omega_{t,t-1}} + \frac{\ol{\kappa}_4}{\ol{\kappa}_2} \rho \expt{\ol{g}(\sfY) \ol{\sfZ}^2} \paren{\nu_{t-1} \chi_{t-1} + \nu_{t-2} \chi_t} \notag \\
&\quad + \frac{\ol{\kappa}_4}{\ol{\kappa}_2^2} \expt{\ol{g}(\sfY)^2} \tau_{t-1} \tau_{t-2} n_{t-1} 
+ \paren{
    \frac{\ol{\kappa}_4}{\ol{\kappa}_2} \expt{\ol{g}(\sfY)^2 \ol{\sfZ}^2} 
    + \frac{\ol{\kappa}_6}{\ol{\kappa}_2^2} \expt{\ol{g}(\sfY) \ol{\sfZ}^2}^2
} \rho \nu_{t-1} \nu_{t-2}
- \sigma^2 \nu_{t-1} \nu_t \bigg\} , \notag \\
\Omega_{t+1,t} &= \delta \brace{
    \rho \paren{\expt{\ol{g}(\sfY)^2 \ol{\sfZ}^2} + \frac{\ol{\kappa}_4}{\ol{\kappa}_2^2} \expt{\ol{g}(\sfY) \ol{\sfZ}^2}^2} \nu_{t} \nu_{t-1} 
    + \frac{1}{\ol{\kappa}_2} \expt{\ol{g}(\sfY)^2} \tau_{t} \tau_{t-1} n_t
} . \notag 
\end{align}
Eliminating $ \Omega_{t,t-1} $, we get 
\begin{align}
n_t &= \frac{1}{\tau_t \tau_{t-1}} \bigg\{ \ol{\kappa}_2 \chi_t \chi_{t-1} \rho
+ \delta \ol{\kappa}_2 \rho \paren{\expt{\ol{g}(\sfY)^2 \ol{\sfZ}^2} + \frac{\ol{\kappa}_4}{\ol{\kappa}_2^2} \expt{\ol{g}(\sfY) \ol{\sfZ}^2}^2} \nu_{t-1} \nu_{t-2} 
+ \delta \expt{\ol{g}(\sfY)^2} \tau_{t-1} \tau_{t-2} n_{t-1} \notag \\
&\quad + \frac{\ol{\kappa}_4}{\ol{\kappa}_2} \rho \expt{\ol{g}(\sfY) \ol{\sfZ}^2} \paren{\nu_{t-1} \chi_{t-1} + \nu_{t-2} \chi_t} 
+ \frac{\ol{\kappa}_4}{\ol{\kappa}_2^2} \expt{\ol{g}(\sfY)^2} \tau_{t-1} \tau_{t-2} n_{t-1} \notag \\
&\quad + \paren{
    \frac{\ol{\kappa}_4}{\ol{\kappa}_2} \expt{\ol{g}(\sfY)^2 \ol{\sfZ}^2} 
    + \frac{\ol{\kappa}_6}{\ol{\kappa}_2^2} \expt{\ol{g}(\sfY) \ol{\sfZ}^2}^2
} \rho \nu_{t-1} \nu_{t-2}
- \sigma^2 \nu_{t-1} \nu_t \bigg\}
\notag \\
&= \frac{1}{\tau_t \tau_{t-1}} \bigg\{ \ol{\kappa}_2 \chi_t \chi_{t-1} \rho 
+ \rho \paren{\paren{\delta \ol{\kappa}_2 + \frac{\ol{\kappa}_4}{\ol{\kappa}_2}} \expt{\ol{g}(\sfY)^2 \ol{\sfZ}^2} + \paren{\delta \frac{\ol{\kappa}_4}{\ol{\kappa}_2} + \frac{\ol{\kappa}_6}{\ol{\kappa}_2^2}} \expt{\ol{g}(\sfY) \ol{\sfZ}^2}^2} \nu_{t-1} \nu_{t-2} 
\notag \\
&\quad + \paren{\delta + \frac{\ol{\kappa}_4}{\ol{\kappa}_2^2}} \expt{\ol{g}(\sfY)^2} \tau_{t-1} \tau_{t-2} n_{t-1} 
+ \frac{\ol{\kappa}_4}{\ol{\kappa}_2} \rho \expt{\ol{g}(\sfY) \ol{\sfZ}^2} \paren{\nu_{t-1} \chi_{t-1} + \nu_{t-2} \chi_t} 
- \sigma^2 \nu_{t-1} \nu_t \bigg\} . 
\label{eqn:nt} 
\end{align}

Recalling $w_1, w_2, w_3, w_4$, and $\gamma$ from \Cref{lem:SE_stay} and noting that $ w_3 w_4 = \gamma $, we have the following identities for every sufficiently large $t$: 
\begin{align}
\frac{\chi_{t-1} \chi_t}{\tau_{t-1} \tau_t} &=  \frac{\chi_{t-1} \chi_t}{\nu_{t-1} \nu_t} \frac{\nu_{t-1} \nu_t}{\tau_{t-1} \tau_t} = w_3^{-2} w_1^{-1} , \notag \\
\frac{\nu_{t-2} \nu_{t-1}}{\tau_{t-1} \tau_t} &=  \frac{\nu_{t-2} \nu_{t-1}}{\nu_{t-1} \nu_t} \frac{\nu_{t-1} \nu_t}{\tau_{t-1} \tau_t} = \gamma^{-2} w_1^{-1} , \notag \\
\frac{\tau_{t-2} }{ \tau_t} &=  \frac{\tau_{t-2} }{\nu_{t-2} } \frac{\nu_{t-2}}{\nu_{t-1}} \frac{\nu_{t-1}}{\nu_t} \frac{\nu_t}{\tau_t} = \gamma^{-2} , \notag \\
\frac{\chi_{t-1} \nu_{t-1}}{\tau_{t-1} \tau_t} &=  \frac{\chi_{t-1}}{\nu_{t-1}} \frac{\nu_{t-1}}{\tau_{t-1}} \frac{\nu_{t-1}}{\nu_t} \frac{\nu_t}{\tau_t} = w_3^{-1} w_1^{-1} \gamma^{-1} , \notag \\
\frac{\chi_t \nu_{t-2}}{\tau_{t-1} \tau_t} &=  \frac{\chi_t}{\nu_{t-1}} \frac{\nu_{t-1}}{\tau_{t-1}} \frac{\nu_{t-2}}{\nu_{t-1}} \frac{\nu_{t-1}}{\nu_t} \frac{\nu_t}{\tau_t} = w_4 w_1^{-1} \gamma^{-2} , \notag \\
\frac{\nu_{t-1} \nu_t}{\tau_{t-1} \tau_t} &= w_1^{-1}. \notag 
\end{align}
Dividing both sides of \Cref{eqn:nt} by $ \tau_t \tau_{t-1} $ and using the above identities, we have 
\begin{align}
n_t &= \rho \ol{\kappa}_2 w_3^{-2} w_1^{-1} 
+ \rho \brace{ \paren{ \delta \ol{\kappa}_2 + \frac{\ol{\kappa}_4}{\ol{\kappa}_2} } \expt{\ol{g}(\sfY)^2 \ol{\sfZ}^2} + \paren{ \delta \frac{\ol{\kappa}_4}{\ol{\kappa}_2} + \frac{\ol{\kappa}_6}{\ol{\kappa}_2^2} } \expt{\ol{g}(\sfY) \ol{\sfZ}^2}^2 } \gamma^{-2} w_1^{-1} \notag \\
&\quad + \paren{ \delta + \frac{\ol{\kappa}_4}{\ol{\kappa}_2^2} } \expt{\ol{g}(\sfY)^2} \gamma^{-2} n_{t-1}
+ \frac{\ol{\kappa}_4}{\ol{\kappa}_2} \rho \expt{\ol{g}(\sfY) \ol{\sfZ}^2} \paren{ w_4 w_1^{-1} \gamma^{-2} + w_3^{-1} w_1^{-1} \gamma^{-1} } 
- \sigma^2 w_1^{-1} , \notag \\
&= A n_{t-1} + B , \label{eqn:nt_AB} 
\end{align}
where 
\begin{align}
    A &\coloneqq \paren{ \delta + \frac{\ol{\kappa}_4}{\ol{\kappa}_2^2} } \expt{\ol{g}(\sfY)^2} \gamma^{-2} = \frac{\expt{\ol{g}(\sfY)^2}}{\paren{\delta + \frac{\ol{\kappa}_4}{\ol{\kappa}_2^2}} \expt{\ol{g}(\sfY) \ol{\sfZ}^2}^2} , \notag \\
    B &\coloneqq \rho \ol{\kappa}_2 w_3^{-2} w_1^{-1} 
    + \rho \brace{ \paren{ \delta \ol{\kappa}_2 + \frac{\ol{\kappa}_4}{\ol{\kappa}_2} } \expt{\ol{g}(\sfY)^2 \ol{\sfZ}^2} + \paren{ \delta \frac{\ol{\kappa}_4}{\ol{\kappa}_2} + \frac{\ol{\kappa}_6}{\ol{\kappa}_2^2} } \expt{\ol{g}(\sfY) \ol{\sfZ}^2}^2 } \gamma^{-2} w_1^{-1} \notag \\
    &\quad + \frac{\ol{\kappa}_4}{\ol{\kappa}_2} \rho \expt{\ol{g}(\sfY) \ol{\sfZ}^2} \paren{ w_4 w_1^{-1} \gamma^{-2} + w_3^{-1} w_1^{-1} \gamma^{-1} } 
    - \sigma^2 w_1^{-1} . \notag 
\end{align}
Note that $ A \in (0,1) $ since the criticality condition \Cref{eqn:thr} holds. 

Recall that we would like to prove \Cref{eqn:lim_nt}. We will show that:
\begin{align}
    1 \le \ul{n} &\coloneqq \liminf_{t\to\infty} n_t \le \limsup_{t\to\infty} n_t \eqqcolon \ol{n} \le 1 . \notag 
\end{align}
To prove the right-most inequality, we take $\limsup$ on both sides of \Cref{eqn:nt_AB} and apply sub-additivity of $\limsup$: $ \ol{n} \le A \ol{n} + B $, implying $ \ol{n} \le B/(1-A) $. 
Using the definitions of $A,B$, it can be verified with a lengthy calculation that $ B/(1-A) = 1 $. 
Therefore $ \ol{n} \le 1 $. 
Similarly, super-additivity of $ \liminf $ allows us to show the left-most inequality $ \ul{n} \ge 1 $. 
Therefore \Cref{eqn:lim_nt} holds.

Now  consider the first limit in \Cref{eqn:gamma}. 
Using the definition \Cref{eqn:def_nt} of $n_t$, its limit in \Cref{eqn:lim_nt} and \Cref{lem:SE_stay}, we have
\begin{align}
    \lim_{t\to\infty} \lim_{n\to\infty} \frac{\inprod{u^t}{u^{t-1}}}{\normtwo{u^{t-1}}^2} 
    &= \lim_{t\to\infty} \frac{\Sigma_{t+1,t}}{\Sigma_{t,t}} = \lim_{t\to\infty} \frac{\tau_t \tau_{t-1} n_t + \nu_t \nu_{t-1} \sigma^2}{\tau_{t-1}^2 + \nu_{t-1}^2 \sigma^2} 
    = \lim_{t\to\infty} \frac{\frac{\tau_t \tau_{t-1}}{\nu_t \nu_{t-1}} n_t + \sigma^2}{\frac{\nu_{t-1}}{\nu_t} \frac{\tau_{t-1}^2}{\nu_{t-1}^2} + \frac{\nu_{t-1}}{\nu_t} \sigma^2}
    = \frac{w_1 + \sigma^2}{\gamma^{-1} w_1 + \gamma^{-1} \sigma^2} 
    = \gamma , \notag 
\end{align}
as claimed. 

Turning to the second limit in \Cref{eqn:gamma}, again using \Cref{lem:SE_stay}, we have
\begin{align}
\lim_{t\to\infty} \lim_{n\to\infty} \frac{\normtwo{u^t}^2}{\normtwo{u^{t-1}}^2}
&= \lim_{t\to\infty} \frac{\nu_t^2 \sigma^2 + \tau_t^2}{\nu_{t-1}^2 \sigma^2 + \tau_{t-1}^2}
= \lim_{t\to\infty} \frac{\nu_t^2}{\nu_{t-1}^2}\frac{\sigma^2 + \frac{\tau_t^2}{\nu_t^2}}{\sigma^2 + \frac{\tau_{t-1}^2}{\nu_{t-1}^2}} 
= \gamma^2 \frac{\sigma^2 + w_1}{\sigma^2 + w_1}
= \gamma^2 , 
\notag
\end{align}
implying, by the positivity of $\gamma$ (see \Cref{asmp:preprocess}) that
\begin{align}
    \lim_{t\to\infty} \lim_{n\to\infty} \frac{\normtwo{u^t}}{\normtwo{u^{t-1}}} &= \gamma , \notag 
\end{align}
as claimed in the second equality of \Cref{eqn:gamma}. 

\subsection{Proof of \Cref{lem:whe_ts}}
\label{app:error_estimates}

First, by the definition \Cref{eqn:e2}, we can write $ \wh{e}_2^t $ as
\begin{align}
\wh{e}_2^t &= \ol{\kappa}_2^{-1} \paren{ \gamma - \frac{\beta_t}{\beta_{t-1}} } X^\top \paren{ \wt{G} - I_n } \brack{ \paren{ \frac{\beta_t}{\beta_{t-1}} - 1 } I_n + \wt{G} }^{-1} \brack{ \paren{ \gamma - 1 } I_n + \wt{G} }^{-1} X \wh{v}^t . \notag 
\end{align}
Then we immediately have
\begin{align}
\lim_{t\to\infty} \lim_{d\to\infty} \normtwo{\wh{e}_2^t} &= 0 ,\label{eqn:e20} 
\end{align}
since by \Cref{lem:technical}, 
\begin{align}
\lim_{t\to\infty} \lim_{d\to\infty} \frac{\beta_t}{\beta_{t-1}} &= \gamma , \notag 
\end{align}
and by \Cref{eqn:l} and the invertibility of \Cref{eqn:G_invertible}, 
\begin{align}
\lim_{t\to\infty} \lim_{d\to\infty} \normtwo{ \ol{\kappa}_2^{-1} X^\top \paren{ \wt{G} - I_n } \brace{ \brack{ \paren{ \frac{\beta_t}{\beta_{t-1}} - 1 } I_n + \wt{G} }^{-1} \brack{ \paren{ \gamma - 1 } I_n + \wt{G} }^{-1} } X \wh{v}^t } &< \infty . \notag 
\end{align}

Next, we claim that 
\begin{align}
\lim_{t\to\infty} \lim_{d\to\infty} \normtwo{ \wh{u}^t - \wh{u}^{t-1} }^2 &= 0 . \label{eqn:whu_whu}
\end{align}
Expanding the square on the left-hand side, we have
\begin{align}
\lim_{t\to\infty} \lim_{d\to\infty} \normtwo{ \wh{u}^t - \wh{u}^{t-1} }^2
&= 2 - 2 \lim_{t\to\infty} \lim_{d\to\infty} \frac{\inprod{u^t}{u^{t-1}}}{\normtwo{u^t} \normtwo{u^{t-1}}} , \label{eqn:uhat_err} 
\end{align}
since $ \wh{u}^t \in \bbS^{n-1} $ by the definition \Cref{eqn:beta_whu}. 
Invoking \Cref{lem:technical}, we have
\begin{align}
    \lim_{t\to\infty} \lim_{d\to\infty} \frac{\inprod{u^t}{u^{t-1}}}{\normtwo{u^t} \normtwo{u^{t-1}}}
    &= \lim_{t\to\infty} \lim_{d\to\infty} \frac{\beta_{t-1}}{\beta_t} \frac{\inprod{u^t}{u^{t-1}}}{\normtwo{u^{t-1}}^2}
    = 1 , \notag 
\end{align}
which implies \Cref{eqn:whu_whu}. 

In view of the definition of $ \wh{e}_1^t $ in \Cref{eqn:e1} and $ \wh{v}^t $ in \Cref{eqn:vhat}, we have
\begin{align}
&&
    \lim_{t\to\infty} \lim_{d\to\infty} \normtwo{\wh{e}_1^t} &= 0 , & 
    \lim_{t\to\infty} \lim_{d\to\infty} \normtwo{\wh{v}^{t+1} - \wh{v}^t} &= 0 . &
& \label{eqn:e1_vhat_0}
\end{align}
Recalling the definition of $ \wh{e}_3^t $ in \Cref{eqn:Dtilde_e3}, by \Cref{eqn:e1_vhat_0,eqn:e20}, we then have
\begin{align}
\lim_{t\to\infty} \lim_{d\to\infty} \normtwo{\wh{e}_3^t} &= 0 . \label{eqn:e30}
\end{align}
Finally, consider $ \wh{e}^{t,s} $ defined in \Cref{eqn:ets}. 
Using \Cref{eqn:e30} and the triangle inequality, 
\begin{align}
\lim_{s\to\infty} \lim_{t\to\infty} \lim_{d\to\infty} \normtwo{\wh{e}^{t,s}} &= 0 . \notag 
\end{align}

\section{Proof of \Cref{thm:opt_thr}}
\label{app:pf_thm:opt_thr}

\paragraph{Proof of \Cref{itm:conv}.}

Consider any spectral estimator identified with a function $ \ol{g} $ satisfying \Cref{eqn:gbar} such that the criticality condition \Cref{eqn:thr} holds. 
By \Cref{eqn:gbar} and the Cauchy--Schwarz inequality, we have 
\begin{align}
\expt{\ol{g}(\sfY) \ol{\sfZ}^2} &= \expt{\ol{g}(\sfY) \paren{\ol{\sfZ}^2 - 1}} \notag \\
&= \iint p_\sfZ(z) Q(y\,|\,z) \ol{g}(y) \paren{\sigma^{-2} z^2 - 1} \diff z\diff y \notag \\
&= \int \frac{\expt{Q(y\,|\,\sfZ) \paren{\ol{\sfZ}^2 - 1}}}{\sqrt{\expt{Q(y\,|\,\sfZ)}}} \ol{g}(y) \sqrt{\expt{Q(y\,|\,\sfZ)}} \diff y \notag \\
&\le \paren{ \int \frac{\expt{Q(y\,|\,\sfZ) \paren{\ol{\sfZ}^2 - 1}}^2}{\expt{Q(y\,|\,\sfZ)}} \diff y }^{1/2} \paren{\int \ol{g}(y)^2 \expt{Q(y\,|\,\sfZ)} \diff y}^{1/2} \notag \\
&= \paren{ \int \frac{\expt{Q(y\,|\,\sfZ) \paren{\ol{\sfZ}^2 - 1}}^2}{\expt{Q(y\,|\,\sfZ)}} \diff y }^{1/2} \expt{\ol{g}(\sfY)^2}^{1/2} . \notag 
\end{align}
Using this in \Cref{eqn:thr}, we deduce the following inequality:
\begin{align}
1 &< \paren{ \frac{\ol{\kappa}_4}{\ol{\kappa}_2^2} + \delta } \frac{\expt{\ol{g}(\sfY) \ol{\sfZ}^2}^2}{\expt{\ol{g}(\sfY)^2}} 
\le \paren{ \frac{\ol{\kappa}_4}{\ol{\kappa}_2^2} + \delta } \int \frac{\expt{Q(y\,|\,\sfZ) \paren{\ol{\sfZ}^2 - 1}}^2}{\expt{Q(y\,|\,\sfZ)}} \diff y , \notag 
\end{align}
which is equivalent to \Cref{eqn:opt_thr}, as desired. 

\paragraph{Proof of \Cref{itm:ach}.}

Assume that \Cref{eqn:opt_thr} holds and consider $ \ol{g} $ given by \Cref{eqn:olg_opt}. 
We can compute:
\begin{align}
\expt{\ol{g}(\sfY) \ol{\sfZ}^2}
&= \iint p_\sfZ(z) Q(y \,|\, z) \frac{\int p_\sfZ(z') \paren{\frac{z'^2}{\sigma^2} - 1} Q(y \,|\, z') \diff z'}{\int p_\sfZ(z'') Q(y \,|\, z'') \diff z''} \frac{z^2}{\sigma^2} \diff y \diff z \notag \\
&= \int \frac{\int p_\sfZ(z') \paren{\frac{z'^2}{\sigma^2} - 1} Q(y \,|\, z') \diff z'}{\int p_\sfZ(z'') Q(y \,|\, z'') \diff z''} \brack{ \int p_\sfZ(z) Q(y \,|\, z) \frac{z^2}{\sigma^2} \diff z } \diff y , \notag \\
&= \int \frac{\int p_\sfZ(z') \paren{\frac{z'^2}{\sigma^2} - 1} Q(y \,|\, z') \diff z'}{\int p_\sfZ(z'') Q(y \,|\, z'') \diff z''} \brack{ \int p_\sfZ(z) Q(y \,|\, z) \paren{ \frac{z^2}{\sigma^2} - 1 } \diff z + \int p_\sfZ(z) Q(y \,|\, z) \diff z } \diff y , \notag \\
&= \int \frac{\brack{ \int p_\sfZ(z') \paren{\frac{z'^2}{\sigma^2} - 1} Q(y \,|\, z') \diff z' }^2}{ \int p_\sfZ(z'') Q(y \,|\, z'') \diff z'' } \diff y
+ \iint p_\sfZ(z') \paren{\frac{z'^2}{\sigma^2} - 1} Q(y \,|\, z') \diff z' \diff y \notag \\
&= \int \frac{\expt{\paren{\ol{\sfZ}^2 - 1} Q(y \,|\, \sfZ)}^2}{ \expt{Q(y \,|\, \sfZ)} } \diff y , \label{eqn:use_Z}
\end{align}
where \Cref{eqn:use_Z} follows since $ \expt{\sfZ^2} = \sigma^2 $, therefore the second term in the previous step vanishes.  We also have:
\begin{align}
\expt{\ol{g}(\sfY)^2} &= \int p_\sfY(y) \frac{\brack{ \int p_\sfZ(z') \paren{\frac{z'^2}{\sigma^2} - 1} Q(y \,|\, z') \diff z' }^2}{\brack{ \int p_\sfZ(z'') Q(y \,|\, z'') \diff z'' }^2} \diff y \notag \\
&= \int \frac{\brack{ \int p_\sfZ(z') \paren{\frac{z'^2}{\sigma^2} - 1} Q(y \,|\, z') \diff z' }^2}{ \int p_\sfZ(z'') Q(y \,|\, z'') \diff z'' } \diff y 
= \int \frac{\expt{\paren{\ol{\sfZ}^2 - 1} Q(y \,|\, \sfZ)}^2}{ \expt{Q(y \,|\, \sfZ)} } \diff y . \label{eqn:Eg2}
\end{align}
By the above expressions, the condition \Cref{eqn:opt_thr} is equivalent to \Cref{eqn:thr}. 

Next, we verify the validity of \Cref{eqn:gbar,eqn:T_edge} for the choice of $ \ol{g} $ in \Cref{eqn:olg_opt}. 
The middle identity of \Cref{eqn:gbar} holds by the tower property of conditional expectation: 
\begin{align}
    \expt{\ol{g}(\sfY)} &= \expt{\expt{\ol{\sfZ}^2 - 1 \mid \sfY}} = \expt{\ol{\sfZ}^2} - 1 = 0 . \notag 
\end{align}

For the third inequality of \Cref{eqn:gbar}, recall from \Cref{eqn:use_Z} that 
\begin{align}
    \expt{\ol{g}(\sfY) \ol{\sfZ}^2} &= \expt{\expt{\ol{\sfZ}^2 - 1 \mid \sfY}^2} , \notag 
\end{align}
which is positive by the assumption \Cref{eqn:nontrivial}. 

For the first inequality in \Cref{eqn:gbar}, observe that $ \ol{g}(y) = \expt{\ol{\sfZ}^2 \mid \sfY = y} - 1 \ge -1 $ for any $y$. 
Since $ \expt{\ol{g}(\sfY) \ol{\sfZ}^2} = \expt{\ol{g}(\sfY)^2} $ by \Cref{eqn:use_Z,eqn:Eg2}, we have 
\begin{align}
    \ol{g}(y) &\ge - \frac{\expt{\ol{g}(\sfY)^2}}{\expt{\ol{g}(\sfY) \ol{\sfZ}^2}}
    > - \paren{\frac{\ol{\kappa}_4}{\ol{\kappa}_2^2} + \delta} \expt{\ol{g}(\sfY) \ol{\sfZ}^2}
    = - \gamma , \notag 
\end{align}
where the strict inequality is by \Cref{eqn:thr}. 

Finally, we move to \Cref{eqn:T_edge}. 
We claim that $ \supp\mu_{\sfT} $ is finite. 
Indeed, for any $y\in\supp\sfY$, we have $ \cT(y) = 1 - \frac{\gamma}{\ol{g}(y) + \gamma} \le 1 $ since the second term is positive. 
If $ \sfT $ has a point mass at $ \supp\mu_{\sfT} $, then \Cref{eqn:T_edge} holds. 
Otherwise, we can perturb $\sfT$ by adding a small point mass at $ \supp\mu_{\sfT} $ while preserving other properties required by \Cref{asmp:preprocess}. 
As long as the perturbation is sufficiently small, the perturbed preprocessing function still satisfies \Cref{eqn:thr} and therefore achieves positive limiting overlap. 
This type of argument has been employed in multiple prior works (see e.g.\ \cite[Proof of Theorem 2]{Mondelli_Montanari} or \cite[Section C]{Zhang_Ji_Venkataramanan_Mondelli}) and we omit the details. 

\section{Proof of \Cref{thm:spec_GVAMP}}
\label{app:pf_SE_GVAMP}

The proof of the state evolution result for spectrally initialized GVAMP (\Cref{thm:spec_GVAMP}) consists of $3$ major steps. 
First, \Cref{thm:SE_GVAMP} below provides a state evolution result for GVAMP with an `oracle' initialization that is independent of the data $(X,y)$. 
Some technical assumptions required by this result are then relaxed in \Cref{cor:degenerate}. 
Note that both theorems above exclude spectral initialization which correlates with $(X,y)$. 
This issue is resolved in \Cref{app:pf_thm:spec_GVAMP} by simulating spectrally initialized GVAMP with a two-phase GVAMP algorithm. 

During the finalization of this work, we became aware of the recent independent work \cite{Chen_Liu_Ma} which contains in Appendix B a state evolution result for a general Orthogonal Approximate Message Passing (OAMP) algorithm. 
A general OAMP takes a similar shape to \Cref{eqn:GVAMP}, with a distinction that it does not assume any underlying statistical model such as \Cref{eqn:model}. 
On the other hand, a general OAMP in \cite{Chen_Liu_Ma} allows the vector denoisers $f_{t}, g_{t}$ to depend on all iterates up to time $t$ (subject to divergence-free conditions with respect to all arguments). 
Due to this flexibility, we acknowledge the possibility of proving \Cref{thm:SE_GVAMP} by reducing GVAMP to general OAMP. 
However, treating spectral initialization promised by \Cref{thm:spec_GVAMP} still requires techniques beyond those in \cite{Chen_Liu_Ma} for spiked matrix models. 
We therefore opt to proceed with our proof from first principles. 
Furthermore, it is straightforward to extend our results to allow memory in the vector denoisers as in \cite{Chen_Liu_Ma}. 
Since such extensions are not necessary for \Cref{sec:BVAMP}, we do not pursue them so as to keep the notation more succinct. 

Let us start with \Cref{thm:SE_GVAMP} which requires additional assumptions listed below. 

\begin{enumerate}[label=(A\arabic*)]
\setcounter{enumi}{\value{asmpctr}}

    \item \label[asmp]{asmp:nondegenerate} For every $t\ge0$, both $ \max\brace{\expt{\Phi_t(\sfLambda_d^2)^2}, \var{\wt{\Phi}_t(\sfLambda_d)}} > 0 $ and $ \max\brace{\expt{\Psi_t(\sfLambda_n^2)^2}, \var{\wt{\Psi}_t(\sfLambda_n)}} > 0 $ hold. 

    \item \label[asmp]{asmp:nondegenerate_fg} 
    $ \wt{\sfR}_0 $ is not equal to $ \pm\sfB_* $ almost surely and $ \wt{\sfP}_0 $ is not equal to $ \pm\sfZ $ almost surely. 
    For every $t\ge0$, there do not exist constants $ \beta_0, \cdots, \beta_t, \alpha_0, \cdots, \alpha_t \in \bbR $ such that 
    \begin{align}
    &&
        &\textnormal{either}&
        \wt{\sfR}_{t+1} - \wt{b}_{t+1} \sfB_* &= \sum_{s = 0}^t \beta_s \paren{\wt{\sfR}_{s} - \wt{b}_{s} \sfB_*} & 
        &\textnormal{or}&
        \wt{\sfP}_{t+1} - \wt{a}_{t+1} \sfZ &= \sum_{s = 0}^t \alpha_s \paren{\wt{\sfP}_{s} - \wt{a}_{s} \sfZ} & 
    & \notag 
    \end{align}
    holds almost surely. 
    
    \item \label[asmp]{asmp:regularity} For every $ t\ge1 $, the function $ f_t \in \PG_{(2+d')\to1} $ is weakly differentiable in the first argument and continuous in other arguments. 
    The weak derivative of $ r \mapsto f_t(r; b, \theta) $, denoted by $ f_t'(r;b, \theta) $, is also in $ \PG_{(2+d')\to1} $, and $ f_t'(\sfR_t; \sfB_*, \sfTheta) $ is continuous with probability $1$. 

    For every $ t\ge1 $, the function $ g_t \in \PG_{(3+n')\to1} $ is weakly differentiable in the first argument and continuous in other arguments. 
    The weak derivative of $ p \mapsto g_t(p; z, e, \xi) $, denoted by $ g_t'(p;z,e, \xi) $, is also in $ \PG_{(3+n')\to1} $, and $ g_t'(\sfP_t; \sfZ, \sfE, \sfXi) $ is continuous with probability $1$. 

\setcounter{asmpctr}{\value{enumi}}
\end{enumerate}

The initializers and side information are required to satisfy the following assumption.

\begin{enumerate}[label=(A\arabic*)]
\setcounter{enumi}{\value{asmpctr}}

    \item\label[asmp]{asmp:init} $ \Theta, \Xi $ are independent of each other and of $ O,Q $ from \Cref{asmp:design}. 
    The initializers $ \wt{r}^0, \wt{p}^0 $ take the form
    \begin{align}
    &&
        \wt{r}^0 &= f_0(\beta_*, \Theta) , &
        \wt{p}^0 &= g_0(z, \eps, \Xi) & 
    & \notag 
    \end{align}
    for some continuous functions $ f_0 \in \PG_{(1+d')\to1} $ and $ g_0 \in \PG_{(2+n')\to1} $, and satisfy
    \begin{align}
    && 
        \matrix{\wt{r}^0 & \beta_* & \Theta} &\wto \matrix{\wt{\sfR}_0 & \sfB_* & \sfTheta^\top} , & 
        \matrix{\wt{p}^0 & z & \eps & \Xi} &\wto \matrix{\wt{\sfP}_0 & \sfZ & \sfE & \sfXi^\top} , &
    & \label{eqn:GVAMP_init} 
    \end{align}
    where the random variables on the right sides are independent of $\sfLambda_n, \sfLambda_d$. 

\setcounter{asmpctr}{\value{enumi}}
\end{enumerate}

\begin{theorem}
\label{thm:SE_GVAMP}
Consider the GVAMP iteration in \Cref{eqn:GVAMP} with initialization \Cref{eqn:GVAMP_init}. 
Let \Cref{asmp:init,asmp:tr_free,asmp:div_free,asmp:nondegenerate,asmp:nondegenerate_fg,asmp:regularity} hold. 
Then for any fixed $t\ge0$, 
\begin{subequations}
\label{eqn:conv_rp}
\begin{align}
\matrix{
    r^1 & \cdots & r^{t+1} & 
    \wt{r}^0 & \cdots & \wt{r}^{t+1} & 
    \beta_* & \Theta 
} 
&\wto \matrix{
    \sfR_1 & \cdots & \sfR_{t+1} & 
    \wt{\sfR}_0 & \cdots & \wt{\sfR}_{t+1} & 
    \sfB_* & \sfTheta^\top
} , \label{eqn:SE_G_r} \\
\matrix{
    p^1 & \cdots & p^{t+1} & 
    \wt{p}^0 & \cdots & \wt{p}^{t+1} & 
    z & \eps & \Xi
} 
&\wto \matrix{
    \sfP_1 & \cdots & \sfP_{t+1} & 
    \wt{\sfP}_0 & \cdots & \wt{\sfP}_{t+1} & 
    \sfZ & \sfE & \sfXi^\top
} , \label{eqn:SE_G_p}
\end{align}
\end{subequations}
where the random variables on the right-hand sides are defined in \Cref{eqn:SE_GVAMP}. 
Furthermore, the covariance matrices of $ (\sfJ_{s+1})_{0\le s\le t}, (\sfK_{s+1})_{0\le s\le t} $ are invertible for any $t\ge0$. 
\end{theorem}

The full proof of \Cref{thm:SE_GVAMP} is given in \Cref{sec:pf_thm:SE_GVAMP}. 

The next result relaxes \Cref{asmp:nondegenerate,asmp:regularity} required in \Cref{thm:SE_GVAMP} to \Cref{asmp:remove}, at the cost of concluding Wasserstein convergence of a lower order. 

\begin{theorem}
\label{cor:degenerate}
Consider the GVAMP iteration in \Cref{eqn:GVAMP} with initialization \Cref{eqn:GVAMP_init}. 
Let \Cref{asmp:init,asmp:tr_free,asmp:div_free,asmp:remove,asmp:remove_PhiPsi} hold. 
Then for any fixed $t\ge0$, 
\begin{subequations}
\label{eqn:conv_rp_remove}
\begin{align}
\matrix{
    r^1 & \cdots & r^{t+1} & 
    \wt{r}^0 & \cdots & \wt{r}^{t+1} & 
    \beta_* & \Theta 
} 
&\xrightarrow{W_2} \matrix{
    \sfR_1 & \cdots & \sfR_{t+1} & 
    \wt{\sfR}_0 & \cdots & \wt{\sfR}_{t+1} & 
    \sfB_* & \sfTheta^\top
} , \label{eqn:SE_G_r_remove} \\
\matrix{
    p^1 & \cdots & p^{t+1} & 
    \wt{p}^0 & \cdots & \wt{p}^{t+1} & 
    z & \eps & \Xi
} 
&\xrightarrow{W_2} \matrix{
    \sfP_1 & \cdots & \sfP_{t+1} & 
    \wt{\sfP}_0 & \cdots & \wt{\sfP}_{t+1} & 
    \sfZ & \sfE & \sfXi^\top
} , \label{eqn:SE_G_p_remove}
\end{align}
\end{subequations}
where the random variables on the right-hand sides are defined in \Cref{eqn:SE_GVAMP}. 
\end{theorem}

See \Cref{app:pf_cor:degenerate} for a proof. 

To conclude \Cref{thm:spec_GVAMP}, it remains to treat spectral initialization which is detailed in \Cref{app:pf_thm:spec_GVAMP}. 

\subsection{Proof of \Cref{thm:SE_GVAMP}}
\label{sec:pf_thm:SE_GVAMP}

To derive the state evolution result in \Cref{thm:SE_GVAMP}, we reparametrize the GVAMP iteration in \Cref{eqn:GVAMP} into a form of an abstract VAMP algorithm presented in \Cref{eqn:abs_VAMP} below. 
We then prove a state evolution result for this instance of the abstract VAMP, \Cref{thm:SE_abs_VAMP},  which implies \Cref{thm:SE_GVAMP}.

Let us introduce an abstract iteration driven by $ (O, Q) \sim \haar(\bbO(n)) \ot \haar(\bbO(d)) $. 
We will soon configure this iteration  to implement the original GVAMP algorithm in \Cref{eqn:GVAMP}. 
Let $ r_\inn, r_\out, s_\inn, s_\out, k_\inn, k_\out, \ell_\inn, \ell_\out $ be non-negative integers that are fixed relative to $n,d$. 
Let $ w_{p,\inn}\in\bbR^{d\times r_\inn},w_{p,\out}\in\bbR^{n\times r_\out},w_{q,\inn}\in\bbR^{d\times s_\inn},w_{q,\out}\in\bbR^{n\times s_\out} $ be side information that is independent of $O,Q$. 
Suppose that we are given the following functions for each $ t\ge0 $: 
\begin{subequations}
\label{eqn:fpq}
\begin{align}
F_{q,\inn}^t &\colon \bbR^{\ell_\inn} \times \bbR^{\ell_\out} \times \bbR^{s_\inn} \to \bbR^{k_\inn} , & 
F_{q,\out}^t &\colon \bbR^{\ell_\inn} \times \bbR^{\ell_\out} \times \bbR^{s_\out} \to \bbR^{k_\out} , & \\
F_{p,\inn}^{t+1} &\colon \bbR^{k_\inn} \times \bbR^{k_\out} \times \bbR^{r_\inn} \to \bbR^{\ell_\inn} , & 
F_{p,\out}^{t+1} &\colon \bbR^{k_\inn} \times \bbR^{k_\out} \times \bbR^{r_\out} \to \bbR^{\ell_\out} . &
\end{align}
\end{subequations}
Consider the iteration: for $ t\ge0 $, 
\begin{subequations}
\label{eqn:abs_VAMP}    
\begin{align}
q_\inn^t &= Q^\top f_{p,\inn}^t \in \bbR^{d\times\ell_\inn} , & 
q_\out^t &= O^\top f_{p,\out}^t \in \bbR^{n\times\ell_\out} , & \label{eqn:OQ1} \\
f_{q,\inn}^t &= F_{q,\inn}^t(q_\inn^t, q_\out^t, w_{q,\inn}) \in \bbR^{d\times k_\inn} , & 
f_{q,\out}^t &= F_{q,\out}^t(q_\inn^t, q_\out^t, w_{q,\out}) \in \bbR^{n\times k_\out} , & \label{eqn:OQ2} \\
p_\inn^{t+1} &= Q f_{q,\inn}^t \in \bbR^{d\times k_\inn} , & 
p_\out^{t+1} &= O f_{q,\out}^t \in \bbR^{n\times k_\out} , & \label{eqn:OQ3} \\
f_{p,\inn}^{t+1} &= F_{p,\inn}^{t+1}(p_\inn^{t+1}, p_\out^{t+1}, w_{p,\inn}) \in \bbR^{d\times\ell_\inn} , & 
f_{p,\out}^{t+1} &= F_{p,\out}^{t+1}(p_\inn^{t+1}, p_\out^{t+1}, w_{p,\out}) \in \bbR^{n\times\ell_\out} , & \label{eqn:OQ4} 
\end{align}
\end{subequations}
initialized with $ f_{p,\inn}^0 \in \bbR^{d\times\ell_\inn}, f_{p,\out}^0 \in \bbR^{n\times\ell_\out} $ independent of $ O,Q $. 

The functions in \Cref{eqn:fpq} are understood to apply to their matrix inputs row-wise in the updates \Cref{eqn:OQ2,eqn:OQ4}. 
Note that these functions may take input matrices with a varying number of rows. 
If the number of rows of the output matrix is larger than that of some of the input matrices,  the function is understood to apply to the shorter input matrices appended with additional zero rows. 
Otherwise, the taller input matrices are truncated before applying the function row-wise. 
These operations ensure that in both cases the column dimension of all input matrices equals that of the output matrix. 
For instance, consider the first equation in \Cref{eqn:OQ2}. 
For any $ i\in[d] $, according to the above convention, we have that if $ d > n $, 
\begin{align}
    f_{q,\inn}^t &= \matrix{
        F_{q,\inn}^t( (q_{\inn}^t)_{1,:}, (q_\out^t)_{1,:}, (w_{q,\inn})_{1,:} )^\top \\
        \vdots \\
        F_{q,\inn}^t( (q_{\inn}^t)_{n,:}, (q_\out^t)_{n,:}, (w_{q,\inn})_{n,:} )^\top \\
        F_{q,\inn}^t( (q_{\inn}^t)_{n+1,:}, 0_{\ell_\out\times1}, (w_{q,\inn})_{n+1,:} )^\top \\
        \vdots \\
        F_{q,\inn}^t( (q_{\inn}^t)_{d,:}, 0_{\ell_\out\times1}, (w_{q,\inn})_{d,:} )^\top 
    } \in \bbR^{d\times k_\inn} , \notag 
\end{align}
and otherwise if $ d \le n $, 
\begin{align}
    f_{q,\inn}^t &= \matrix{
        F_{q,\inn}^t( (q_{\inn}^t)_{1,:}, (q_\out^t)_{1,:}, (w_{q,\inn})_{1,:} )^\top \\
        \vdots \\
        F_{q,\inn}^t( (q_{\inn}^t)_{d,:}, (q_\out^t)_{d,:}, (w_{q,\inn})_{d,:} )^\top 
    } \in \bbR^{d\times k_\inn} . \notag 
\end{align}

We now design the functions $ F_{q,\inn}^t, F_{q,\out}^t, F_{p,\inn}^{t+1}, F_{p,\out}^{t+1} $ suitably so that the iteration \Cref{eqn:abs_VAMP} implements (a slight generalization of) \Cref{eqn:GVAMP}. 
To this end, let us first fix any $ \Gamma \ge 0 $ as an upper bound on the time index. 
Recall from \Cref{asmp:design} that $ O\in\bbR^{n\times n},Q\in\bbR^{d\times d} $ are independent Haar matrices collecting the left, right singular vectors of $X$ and $ \lambda\in\bbR^{n\wedge d} $ is an entry-wise non-negative vector collecting the singular values of $X$.  
We define the vectors $\lambda_n \in \bbR^n_{\ge0}$ and 
$ \lambda_d \in \bbR^d_{\ge0}$ as:
\begin{align}
&&
    \lambda_n &= \begin{cases}
        \lambda , & n\le d \\
        \matrix{\lambda^\top & 0_{n-d}^\top}^\top , & n > d
    \end{cases} , & 
    \lambda_d &= \begin{cases}
        \matrix{\lambda^\top & 0_{d-n}^\top}^\top , & n < d \\
        \lambda , & n \ge d
    \end{cases} . & 
& \label{eqn:lambda_nd} 
\end{align}
These vectors have Wasserstein limits $ \lambda_n \wto \sfLambda_n, \lambda_d \wto \sfLambda_d $ where the limiting random variables are defined in \Cref{eqn:sfLambda_nd}. Consider $4$ sequences (indexed by $ (n,d) $) of vectors $ (\phi_t)_{0\le t\le\Gamma}\subset\bbR^{d} , (\wt{\phi}_t)_{0\le t\le\Gamma}\subset\bbR_{\ge0}^{n\wedge d} , (\psi_t)_{0\le t\le\Gamma}\subset\bbR^{n} , (\wt{\psi}_t)_{0\le t\le\Gamma}\subset\bbR_{\ge0}^{n\wedge d} $, with $ \wt{\phi}_{t,d}, \wt{\psi}_{t,n} $ denoting vectors of length $d,n$, respectively, obtained by truncating or expanding $\wt{\phi}_{t}$ and $\wt{\psi}_{t}$.
These sequences of vectors are subject to the following assumption.
\begin{enumerate}[label=(A\arabic*)]
\setcounter{enumi}{\value{asmpctr}}
    \item\label[asmp]{asmp:4matrices} $ (\phi_t)_{t\ge0} , (\wt{\phi}_t)_{t\ge0} , (\psi_t)_{t\ge0} , (\wt{\psi}_t)_{t\ge0} $ are independent of $ O,Q,\eps,\beta_*,\Theta,\Xi $. 
    For any fixed $ \Gamma \ge 0 $, 
    \begin{align}
    \begin{split}
        \matrix{\lambda_d & \phi_0 & \cdots & \phi_\Gamma & \wt{\phi}_{0,d} & \cdots & \wt{\phi}_{\Gamma,d}} &\wto \matrix{\sfLambda_d & \sfPhi_0 & \cdots & \sfPhi_\Gamma & \wt{\sfPhi}_{0} & \cdots & \wt{\sfPhi}_{\Gamma}} , \\
        \matrix{\lambda_n & \psi_0 & \cdots & \psi_\Gamma & \wt{\psi}_{0,n} & \cdots & \wt{\psi}_{\Gamma,n}} &\wto \matrix{\sfLambda_n & \sfPsi_0 & \cdots & \sfPsi_\Gamma & \wt{\sfPsi}_{0} & \cdots & \wt{\sfPsi}_{\Gamma}} . 
    \end{split}
    \label{eqn:conv_phi_psi}
    \end{align}
    Furthermore, for each $ t\ge0 $, the following conditions hold: $ \expt{\sfPhi_t} = \expt{\sfPsi_t} = 0 $, $ \max\brace{\expt{\sfPhi_t^2}, \var{\wt{\sfPhi}_t}} > 0 $, $ \max\brace{\expt{\sfPsi_t^2}, \var{\wt{\sfPsi}_t}} > 0 $. 
\setcounter{asmpctr}{\value{enumi}}
\end{enumerate}

The above sequences of vectors along with the pair of independent Haar matrices $(O,Q)$ can be used to define an iteration that slightly generalizes the GVAMP iteration in \Cref{eqn:GVAMP}: 
\begin{subequations}
\label{eqn:GVAMP_gen}
\begin{align}
&&
    r^{t+1} &= Q \diag(\phi_t) Q^\top \wt{r}^t + Q \diag_{d\times n}(\wt{\phi}_t) O^\top \wt{p}^t , & 
    \wt{r}^{t+1} &= f_{t+1}(r^{t+1}; \beta_*, \Theta) , & 
& \label{eqn:GVAMP_gen1} \\
&&
    p^{t+1} &= O \diag(\psi_t) O^\top \wt{p}^t + O \diag_{n\times d}(\wt{\psi}_t) Q^\top \wt{r}^t , & 
    \wt{p}^{t+1} &= g_{t+1}(p^{t+1}; z, \eps, \Xi) , & 
& \label{eqn:GVAMP_gen2}
\end{align}
\end{subequations}
initialized with $ \wt{r}^0 , \wt{p}^0 $ from \Cref{asmp:init}.  
The above iteration is driven by $4$ sequences (indexed by $t\ge0$) of random matrices $ Q \diag(\phi_t) Q^\top$, $O \diag_{n\times d}(\wt{\phi}_t) Q^\top$, $O \diag(\psi_t) O^\top$, $O \diag_{n\times d}(\wt{\psi}_t) Q^\top $. 
These matrices share common eigen-/singular spaces but \Cref{asmp:4matrices} allows them to have arbitrarily correlated spectra across time. 
Clearly, the matrices $ \Phi_t(X^\top X), \wt{\Phi}_t(X), \Psi_t(XX^\top), \wt{\Psi}_t(X) $ in \Cref{eqn:GVAMP} are a special case, hence the GVAMP iteration in \Cref{eqn:GVAMP} is a special case of \Cref{eqn:GVAMP_gen}. Indeed,  if we set
\begin{align}
&&
    \phi_t &= \Phi_t(\lambda_d \circ \lambda_d) , & 
    \wt{\phi}_t &= \wt{\Phi}_t(\lambda) , & 
    \psi_t &= \Psi_t(\lambda_n \circ \lambda_n) , & 
    \wt{\psi}_t &= \wt{\Psi}_t(\lambda) , & 
& \label{eqn:4functions} 
\end{align}
where all functions are applied component-wise, \Cref{eqn:GVAMP_gen} reduces to \Cref{eqn:GVAMP}. 
The additional flexibility offered by \Cref{asmp:4matrices} will be useful for technical purposes in \Cref{cor:degenerate}. 

Before diving into further specifications of \Cref{eqn:abs_VAMP}, we condition on a sequence (indexed by $ (n,d) $) of $ \lambda$, $(\phi_t)_{0\le t\le \Gamma} $, $(\wt{\phi}_t)_{0\le t\le \Gamma} $, $(\psi_t)_{0\le t\le \Gamma} $, $(\wt{\psi}_t)_{0\le t\le \Gamma} $ satisfying \Cref{asmp:design,asmp:4matrices} throughout the rest of the this section. 
\Cref{thm:SE_abs_VAMP} below is a state evolution result for $\Gamma$ steps of \Cref{eqn:abs_VAMP} conditioned on the above vectors.  Choosing these vectors according to \Cref{eqn:4functions}  gives  \Cref{thm:SE_GVAMP}. 

Back to the abstract VAMP iteration \Cref{eqn:abs_VAMP}, we fix the dimensions to be: 
\begin{align}
&&
    k_\inn &= 1, &
    k_\out &= 2, &
    \ell_\inn &= 2, &
    \ell_\out &= 1, &
    r_\inn &= 1 + d', &
    r_\out &= 1 + n', &
    s_\inn &= 1 + 2(\Gamma+1), &
    s_\out &= 1 + 2(\Gamma+1). &
& \notag 
\end{align}
Let the side information be given by 
\begin{align}
\begin{split}
    w_{q,\inn} &= \matrix{\lambda_d & \phi_0 & \cdots & \phi_\Gamma & \wt{\phi}_{0,d} & \cdots & \wt{\phi}_{\Gamma,d}} \in \bbR^{d\times(2(\Gamma+1)+1)} , \\
    w_{q,\out} &= \matrix{\lambda_n & \psi_0 & \cdots & \psi_\Gamma & \wt{\psi}_{0,n} & \cdots & \wt{\psi}_{\Gamma,n}} \in \bbR^{n\times(2(\Gamma+1)+1)} , \\
    w_{p,\inn} &= \matrix{\beta_* & \Theta} \in \bbR^{d\times(1 + d')} , \qquad  
    w_{p,\out} = \matrix{\eps & \Xi} \in \bbR^{n\times(1 + n')} . 
\end{split}
\label{eqn:side_info} 
\end{align}
We will use $ w_{q,\inn,\lambda}, w_{q,\inn,t}, w_{q,\inn,\wt{t}}, w_{q,\out,\lambda}, w_{q,\out,t}, w_{q,\out,\wt{t}} $ to respectively refer to $ \lambda_d, \phi_t, \wt{\phi}_{t,d}, \lambda_n, \wt{\psi}_t, \wt{\psi}_{t,n} $, for any $ 0\le t\le \Gamma $. 
We will also use the notation $ w_{p,\inn,1}, w_{p,\inn,-} $ to respectively refer to $ \beta_*, \Theta $, and $ w_{p,\out,1}, w_{p,\out,-} $ to respectively refer to $ \eps, \Xi $. 
The iterates are specialized as follows. 
For $ t = 0 $, let 
\begin{align}
&& 
    f_{p,\inn}^0 &= \matrix{ \wt{r}^0 - \wt{b}_0 \beta_* & \beta_* } \in \bbR^{d\times2} , & 
    f_{p,\out}^0 &= \wt{p}^0 - \wt{a}_0 z \in \bbR^n , &
& \label{eqn:abs_VAMP_init} 
\end{align}
where $ \wt{r}^0, \wt{p}^0 $ are the initializers of \Cref{eqn:GVAMP} given in \Cref{eqn:GVAMP_init} and $ \wt{b}_0, \wt{a}_0 $ are defined in \Cref{eqn:ab_tilde}. 
For $ t\ge0 $, we define the functions in \Cref{eqn:fpq}  as follows. For inputs
$ q_\inn = \matrix{ q_{\inn,1} & q_{\inn,2} }^\top \in \bbR^2 $, 
$ q_\out \in \bbR $, 
$ p_\inn \in \bbR $, 
$ p_\out = \matrix{ p_{\out,1} & p_{\out,2} }^\top \in \bbR^2 $, 
$ w_{q,\inn} = \matrix{w_{q,\inn,\lambda} & w_{q,\inn,0} & \cdots & w_{q,\inn,\Gamma} & \cdots & w_{q,\inn,\wt{0}} & \cdots & w_{q,\inn,\wt{\Gamma}} }^\top \in \bbR^{2(\Gamma+1)+1} $, 
$ w_{q,\out} = \matrix{w_{q,\out,\lambda} & w_{q,\out,0} & \cdots & w_{q,\out,\Gamma} & \cdots & w_{q,\out,\wt{0}} & \cdots & w_{q,\out,\wt{\Gamma}} }^\top \in \bbR^{2(\Gamma+1)+1} $, 
the outputs of the functions are given by
\begin{subequations}
\label{eqn:F_config}    
\begin{align}
    F_{q,\inn}^t(q_\inn, q_\out, w_{q,\inn}) &= w_{q,\inn,t} (q_{\inn,1} + \wt{\frb}_t q_{\inn,2}) + w_{q,\inn,\wt{t}} (q_\out + \wt{\fra}_t w_{q,\inn,\lambda} q_{\inn,2}) - \frb_{t+1} q_{\inn,2} \in \bbR , \\
    F_{q,\out}^t(q_\inn, q_\out, w_{q,\out}) &= \matrix{ 
        w_{q,\out,t} (q_\out + \wt{\fra}_t w_{q,\out,\lambda} q_{\inn,2}) + w_{q,\out,\wt{t}} (q_{\inn,1} + \wt{\frb}_t q_{\inn,2}) - \fra_{t+1} w_{q,\out,\lambda} q_{\inn,2} & 
        w_{q,\out,\lambda} q_{\inn,2}
    }^\top \in \bbR^2 , \\
    F_{p,\inn}^{t+1}(p_\inn, p_\out, w_{p,\inn}) &= \matrix{ f_{t+1}(p_\inn + \frb_{t+1} w_{p,\inn,1}; w_{p,\inn}) - \wt{\frb}_{t+1} w_{p,\inn,1} & w_{p,\inn,1} }^\top \in \bbR^2 , \\
    F_{p,\out}^{t+1}(p_\inn, p_\out, w_{p,\out}) &= g_{t+1}( p_{\out,1} + \fra_{t+1} p_{\out,2} ; p_{\out,2} , w_{p,\out} ) - \wt{\fra}_{t+1} p_{\out,2} \in \bbR , 
\end{align}
\end{subequations}
where the coefficients $ \wt{\fra}_t, \wt{\frb}_t, \fra_{t+1}, \frb_{t+1} $ are defined in \Cref{eqn:frab_tilde,eqn:frab} below. 
These coefficients are designed in such a way that for every $t\ge0$, 
\begin{align}
&&
    \frac{1}{n} \inprod{f_{p,\out}^{t}}{z} &\to 0 , & 
    \frac{1}{d} \inprod{f_{p,\inn,1}^{t}}{\beta_*} &\to 0 , 
& \notag \\
&&
    \frac{1}{n} \inprod{p_\out^{t+1}}{z} &\to 0 , & 
    \frac{1}{d} \inprod{p_\inn^{t+1}}{\beta_*} &\to 0 , & 
& \notag 
\end{align}
which will be formally justified in the proof of \Cref{thm:SE_abs_VAMP}. 

We now clarify the connection between the GVAMP iteration in \Cref{eqn:GVAMP} and the abstract VAMP iteration in \Cref{eqn:abs_VAMP} with initialization \Cref{eqn:abs_VAMP_init}, side information \Cref{eqn:side_info} and denoising functions in \Cref{eqn:fpq} specialized to \Cref{eqn:F_config}. 
As shall be seen from the proof of \Cref{thm:SE_GVAMP} below (in particular \Cref{itm:ind1} therein), 
under the specification in \Cref{eqn:4functions}, 
it holds that $ \wt{\fra}_t = \wt{a}_t, \wt{\frb}_t = \wt{b}_t, \fra_{t+1} = a_{t+1}, \frb_{t+1} = b_{t+1} $ for all $t\ge0$, where the quantities on the right are defined in \Cref{eqn:ab,eqn:ab_tilde}. 
Given this, by induction on $t$, one can easily verify that \Cref{eqn:abs_VAMP} reduces to the following iteration: 
\begin{align}
\begin{split}
q_\inn^t 
    &= \matrix{ Q^\top (\wt{r}^t - \wt{b}_t \beta_*) & Q^\top \beta_* } \in \bbR^{d\times2} ,  \\
q_\out^t 
    &= O^\top (\wt{p}^t - \wt{a}_t z) \in \bbR^{n} ,  \\
f_{q,\inn}^t &= \Phi_t(\Lambda^\top \Lambda) Q^\top \wt{r}^t + \wt{\Phi}_t(\Lambda)^\top O^\top \wt{p}^t - b_{t+1} Q^\top \beta_* \in \bbR^{d} ,  \\
f_{q,\out}^t &= \matrix{ \Psi_t(\Lambda \Lambda^\top) O^\top \wt{p}^t + \wt{\Psi}_t(\Lambda) Q^\top \wt{r}^t - a_{t+1} \Lambda Q^\top \beta_* & \Lambda Q^\top \beta_* } \in \bbR^{n\times2} ,  \\
p_\inn^{t+1} &= r^{t+1} - b_{t+1} \beta_* \in \bbR^{d} ,  \\
p_\out^{t+1} &= \matrix{ p^{t+1} - a_{t+1} z & z } \in \bbR^{n\times2} ,  \\
f_{p,\inn}^{t+1} &= \matrix{ \wt{r}^{t+1} - \wt{b}_{t+1} \beta_* & \beta_* } \in \bbR^{d\times2} ,  \\
f_{p,\out}^{t+1} &= \wt{p}^{t+1} - \wt{a}_{t+1} z \in \bbR^{n} ,  
\end{split}
\label{eqn:config}
\end{align}
from which the equivalence between \Cref{eqn:abs_VAMP,eqn:GVAMP} immediately follows. 

Let us denote by $ q_{\inn,1}^t, f_{q,\out,1}^t, p_{\out,1}^{t+1}, f_{p,\inn,1}^{t+1} $ the first columns of $ q_{\inn}^t, f_{q,\out}^t, p_{\out}^{t+1}, f_{p,\inn}^{t+1} $, and by $ q_{\inn,2}, f_{q,\out,2}, p_{\out,2}, f_{p,\inn,2} $ their second columns (since the latter terms are time invariant). 

At a high level, the state evolution result in \Cref{thm:SE_abs_VAMP} shows that the joint empirical distribution of the iterates converges to a multivariate Gaussian whose covariance matrix can be tracked recursively over time. 
For instance, the empirical distribution of the rows of $ \matrix{q_{\inn,1}^0 & \cdots & q_{\inn,1}^t} $ converges in Wasserstein to $ \cN(0_{t+1}, \Sigma_{q,\inn}^t) $. 
Similar results hold for $ (q_{\out}^s)_{0\le s\le t}, (p_{\inn}^s)_{1\le s\le t+1}, (p_{\out,1}^s)_{1\le s\le t+1} $ with corresponding covariance matrices $ \Sigma_{q,\out}^t, \Sigma_{p,\inn}^{t+1}, \Sigma_{p,\out}^{t+1} $. 
We now define these covariance matrices along with other state evolution parameters associated with the abstract VAMP iteration \Cref{eqn:abs_VAMP}. 
For $t=0$, they are identified with the initial conditions: 
\begin{subequations}
\label{eqn:SE_OQ}
\begin{align}
& \matrix{ f_{p,\inn,1}^0 & w_{q,\inn,\lambda} & w_{q,\inn,0} & \cdots & w_{q,\inn,\Gamma} & w_{q,\inn,\wt{0}} & \cdots & w_{q,\inn,\wt{\Gamma}} & w_{p,\inn,1} & w_{p,\inn,-} } \notag \\
&\qquad = \matrix{ \wt{r}^0 - \wt{b}_0 \beta_* & \lambda_d & \phi_0 & \cdots & \phi_\Gamma & \wt{\phi}_{0,d} & \cdots & \wt{\phi}_{\Gamma,d} & \beta_* & \Theta } \notag \\
&\qquad \wto \matrix{ \wt{\sfR}_0 - \wt{b}_0 \sfB_* & \sfLambda_d & \sfPhi_0 & \cdots & \sfPhi_\Gamma & \wt{\sfPhi}_{0} & \cdots & \wt{\sfPhi}_{\Gamma} & \sfB_* & \sfTheta^\top } \notag \\
&\qquad \eqqcolon \matrix{ \sfF_{p,\inn,1}^0 & \sfW_{q,\inn,\lambda} & \sfW_{q,\inn,0} & \cdots & \sfW_{q,\inn,\Gamma} & \sfW_{q,\inn,\wt{0}} & \cdots & \sfW_{q,\inn,\wt{\Gamma}} & \sfW_{p,\inn,1} & \sfW_{p,\inn,-}^\top } , \label{eqn:SE_OQ_init_in} \\
& \matrix{ f_{p,\out}^0 & w_{q,\out,\lambda} & w_{q,\out,0} & \cdots & w_{q,\out,\Gamma} & w_{q,\out,\wt{0}}, & \cdots & w_{q,\out,\wt{\Gamma}} & w_{p,\out,1} & w_{p,\out,-} } \notag \\
&\qquad = \matrix{ \wt{p}^0 - \wt{a}_0 z & \lambda_n & \psi_0 & \cdots & \psi_\Gamma & \wt{\psi}_{0,n} & \cdots & \wt{\psi}_{\Gamma,n} & \eps & \Xi } \notag \\
&\qquad \wto \matrix{ \wt{\sfP}_0 - \wt{a}_0 \sfZ & \sfLambda_n & \sfPsi_0 & \cdots & \sfPsi_\Gamma & \wt{\sfPsi}_{0} & \cdots & \wt{\sfPsi}_{\Gamma} & \sfE & \sfXi^\top } \notag \\
&\qquad \eqqcolon \matrix{ \sfF_{p,\out}^0 & \sfW_{q,\out,\lambda} & \sfW_{q,\out,0} & \cdots & \sfW_{q,\out,\Gamma} & \sfW_{q,\out,\wt{0}}, & \cdots & \sfW_{q,\out,\wt{\Gamma}} & \sfW_{p,\out,1} & \sfW_{p,\out,-}^\top } , \label{eqn:SE_OQ_init_out}
\end{align}
where the existence of the limiting random variables on the right is guaranteed by \Cref{asmp:init,asmp:4matrices}, and the constants $ \wt{a}_0, \wt{b}_0 $ are defined as per \Cref{eqn:ab_tilde}. 
We will also use the notation 
$ \sfW_{p,\inn} = \matrix{\sfW_{p,\inn,1} & \sfW_{p,\inn,-}^\top}^\top $, 
$ \sfW_{p,\out} = \matrix{\sfW_{p,\out,1} & \sfW_{p,\out,-}^\top}^\top $,
$ \sfW_{q,\inn} = \matrix{\sfW_{q,\inn,\lambda} & \sfW_{q,\inn,0} & \cdots & \sfW_{q,\inn,\Gamma} & \sfW_{q,\inn,\wt{0}} & \cdots & \sfW_{q,\inn,\wt{\Gamma}}}^\top $, 
and $ \sfW_{q,\out} = \matrix{\sfW_{q,\out,\lambda} & \sfW_{q,\out,0} & \cdots & \sfW_{q,\out,\Gamma} & \sfW_{q,\out,\wt{0}}, & \cdots & \sfW_{q,\out,\wt{\Gamma}}}^\top $. 
For notational convenience, let the random vector $ \sfA \in \bbR^{4(\Gamma+1)+3+d'+n'} $ collect: 
\begin{align}
    \sfA &= \matrix{\sfPhi_0 & \cdots & \sfPhi_\Gamma & \wt{\sfPhi}_0 & \cdots & \wt{\sfPhi}_\Gamma & \sfPsi_0 & \cdots & \sfPsi_\Gamma & \wt{\sfPsi}_0 & \cdots & \wt{\sfPsi}_\Gamma & \sfLambda & \sfB_* & \sfE & \sfTheta^\top & \sfXi^\top}^\top . \notag 
\end{align}
Using these random variables, we define  
\begin{align}
&& 
    \Sigma_{q,\inn}^0 &= \expt{(\wt{\sfR}_0 - \wt{b}_0 \sfB_*)^2} = \expt{(\sfF_{p,\inn,1}^0)^2} \in \bbR , & 
    \Sigma_{q,\out}^0 &= \expt{ (\wt{\sfP}_0 - \wt{a}_0 \sfZ)^2 } = \expt{(\sfF_{p,\out}^0)^2} \in \bbR . &
& \label{eqn:Sigmaq_inout0} 
\end{align}
Note that 
\begin{align}
&&
    \Sigma_{q,\inn}^0 &> 0 , &
    \Sigma_{q,\out}^0 &> 0 , &
& \label{eqn:Sigma0_pos}
\end{align}
since by the definitions of $ \wt{b}_0 $ in \Cref{eqn:ab_tilde},  
\begin{align}
    \expt{(\wt{\sfR}_0 - \wt{b}_0 \sfB_*)^2} &= \expt{\wt{\sfR}_0^2} - 2 \wt{b}_0 \expt{\wt{\sfR}_0 \sfB_*} + \wt{b}_0^2 \expt{\sfB_*^2}
    = \expt{\wt{\sfR}_0^2} - \wt{b}_0^2 \rho , \notag
\end{align}
and similarly, $\expt{ (\wt{\sfP}_0 - \wt{a}_0 \sfZ)^2 } = \expt{\wt{\sfP}_0^2} - \wt{a}_0 \sigma^2 $. 
By \Cref{asmp:init} and the equality condition of Cauchy--Schwarz, the random variables $ \wt{\sfR}_0, \wt{\sfP}_0 $ satisfy the following inequalities strictly: 
\begin{align}
&&
    \expt{ \wt{\sfR}_0 \sfB_* }^2 &< \expt{\wt{\sfR}_0^2} \expt{\sfB_*^2} = \rho \expt{\wt{\sfR}_0^2} , & 
    \expt{ \wt{\sfP}_0 \sfZ }^2 &< \expt{\wt{\sfP}_0^2} \expt{\sfZ^2} = \sigma^2 \expt{\wt{\sfP}_0^2} , &
& \notag 
\end{align}
implying the claimed positivity. 

Now, we recursively define the state evolution parameters $ \fra_{t+1}, \frb_{t+1} \in \bbR $, $ \Sigma_{p,\inn}^{t+1},\Sigma_{p,\out}^{t+1} \in \bbR^{(t+1) \times (t+1)} $, $\wt{\fra}_{t+1}, \wt{\frb}_{t+1} \in \bbR$ and $ \Sigma_{q,\inn}^{t+1}, \Sigma_{q,\out}^{t+1} \in \bbR^{(t+2) \times (t+2)} $. 
Let 
\begin{align}
&&
    \wt{\fra}_0 &= \wt{a}_0, &
    \wt{\frb}_0 &= \wt{b}_0. & 
& \label{eqn:frab_tilde0}
\end{align}
For $ t\ge 0 $, given $\Sigma_{q,\inn}^t, \Sigma_{q,\out}^t$, let 
\begin{align}
    \matrix{ \sfQ_{\inn,2} & \sfQ_{\inn,1}^0 & \cdots & \sfQ_{\inn,1}^t }^\top &\sim \cN\paren{0, \rho} \ot \cN(0_{t+1}, \Sigma_{q,\inn}^t) , \label{eqn:Qin}  \\
    \matrix{ \sfQ_\out^0 & \cdots & \sfQ_\out^t }^\top &\sim \cN(0_{t+1}, \Sigma_{q,\out}^t) \label{eqn:Qout} 
\end{align}
be two mutually independent Gaussian vectors that are also independent of $ \sfA $, $ (\sfP_\inn^{s+1})_{0\le s\le t-1} $, $ (\sfP_\out^{s+1})_{0\le s\le t-1} $. 
For every $ 0\le s\le t $, denote $ \sfQ_\inn^s = \matrix{ \sfQ_{\inn,1}^s & \sfQ_{\inn,2} }^\top $. 
Let 
\begin{align}
    \sfF_{q,\inn}^t &= F_{q,\inn}^t(\sfQ_\inn^t, \sfQ_\out^t, \sfLambda_d) \notag \\
    &= \sfPhi_t (\sfQ_{\inn,1}^t + \wt{\frb}_t \sfQ_{\inn,2}) + \wt{\sfPhi}_t (\sfQ_\out^t + \wt{\fra}_t \sfLambda_d \sfQ_{\inn,2}) - \frb_{t+1} \sfQ_{\inn,2} , \label{eqn:Fqin} \\
    \sfF_{q,\out,1}^t &= F_{q,\out,1}^t(\sfQ_\inn^t, \sfQ_\out^t, \sfLambda_n) \notag \\
    &= \sfPsi_t (\sfQ_\out^t + \wt{\fra}_t \sfLambda_n \sfQ_{\inn,2}) + \wt{\sfPsi}_t (\sfQ_{\inn,1}^t + \wt{\frb}_t \sfQ_{\inn,2}) - \fra_{t+1} \sfLambda_n \sfQ_{\inn,2} , \label{eqn:Fqout1} \\
    \sfF_{q,\out,2} &= F_{q,\out,2}^t(\sfQ_\inn^t, \sfQ_\out^t, \sfLambda_n)
    = \sfLambda_n \sfQ_{\inn,2} , \label{eqn:Fqout2} 
\end{align}
where the constants $ (\fra_{s+1})_{0\le s\le t}, (\frb_{s+1})_{0\le s\le t} $ are defined as
\begin{align}
&&
    \fra_{r+1} &= \frac{1}{\ol{\kappa}_2} \paren{ \expt{ \sfPsi_r \sfLambda_n^2 } \wt{\fra}_r + \expt{ \wt{\sfPsi}_r \sfLambda_n } \wt{\frb}_r } , & 
    \frb_{r+1} &= \expt{ \wt{\sfPhi}_r \sfLambda_d } \wt{\fra}_r . & 
& \label{eqn:frab}
\end{align}
Let $ \sfF_{q,\out}^t = \matrix{ \sfF_{q,\out,1}^t & \sfF_{q,\out,2} }^\top $. 
Using these random variables, define covariance matrices
\begin{align}
&&
    \Sigma_{p,\inn}^{t+1} &= \expt{ \sfF_{q,\inn}^{0:t} (\sfF_{q,\inn}^{0:t})^\top } \in \bbR^{(t+1)\times(t+1)} , &
    \Sigma_{p,\out}^{t+1} &= \expt{ \sfF_{q,\out,1}^{0:t} (\sfF_{q,\out,1}^{0:t})^\top } \in \bbR^{(t+1) \times (t+1)} . & 
& \label{eqn:cov_p}
\end{align}

Next, define two mutually independent Gaussian vectors with the above covariances: 
\begin{align}
\matrix{ \sfP_\inn^1 & \cdots, \sfP_\inn^{t+1} }^\top &\sim \cN\paren{ 0_{t+1}, \Sigma_{p,\inn}^{t+1} } , \label{eqn:Pin} \\
\matrix{ \sfZ & \sfP_{\out,1}^1 & \cdots & \sfP_{\out,1}^{t+1} }^\top &\sim \cN\paren{0, \sigma^2} \ot \cN\paren{ 0_{t+1}, \Sigma_{p,\out}^{t+1} } , \label{eqn:Pout}
\end{align}
that are independent of $ (\sfQ_\inn^s)_{0\le s\le t}, (\sfQ_\out^s)_{0\le s\le t}, \sfA $. 
We will also use the notation $ \sfP_{\out,2} = \sfZ $ and $ \sfP_\out^{s+1} = \matrix{ \sfP_{\out,1}^{s+1} & \sfP_{\out,2} }^\top $, for every $ 0\le s\le t $. 
Let 
\begin{align}
&&
    \sfF_{p,\inn,1}^{t+1} &= F_{p,\inn,1}^{t+1}(\sfP_\inn^{t+1}, \sfP_\out^{t+1}, \sfB_*, \sfTheta) = f_{t+1}(\sfP_\inn^{t+1} + \frb_{t+1} \sfB_*; \sfB_*, \sfTheta) - \wt{\frb}_{t+1} \sfB_* , \label{eqn:Fpin1} & 
&  \\
&& 
    \sfF_{p,\out}^{t+1} &= F_{p,\out}^{t+1}(\sfP_\inn^{t+1}, \sfP_\out^{t+1}, \sfE, \sfXi) = g_{t+1}(\sfP_{\out,1,}^{t+1} + \fra_{t+1} \sfP_{\out,2}; \sfP_{\out,2}, \sfE, \sfXi) - \wt{\fra}_{t+1} \sfP_{\out,2} , \label{eqn:Fpout} & 
&  
\end{align}
where the constants $ \wt{\fra}_{t+1}, \wt{\frb}_{t+1} $ are defined as 
\begin{align}
\begin{split}
    \wt{\fra}_{t+1} &= \frac{1}{\sigma^2} \expt{ g_{t+1}(\sfP_{\out,1}^{t+1} + \fra_{t+1} \sfP_{\out,2}; \sfP_{\out,2}, \sfE, \sfXi) \sfZ } , \\
    \wt{\frb}_{t+1} &= \frac{1}{\rho} \expt{ f_{t+1}(\sfP_\inn^{t+1} + \frb_{t+1} \sfB_*; \sfB_*, \sfTheta) \sfB_* } . 
\end{split}
\label{eqn:frab_tilde}
\end{align}
Let also $ \sfF_{p,\inn,2} = \sfB_* $ and $ \sfF_{p,\inn}^{t+1} = \matrix{ \sfF_{p,\inn,1}^{t+1} & \sfF_{p,\inn,2} }^\top $. 
Use these random variables to define covariance matrices: 
\begin{align}
&& 
    \Sigma_{q,\inn}^{t+1} &= \expt{ \sfF_{p,\inn,1}^{0:t+1} (\sfF_{p,\inn,1}^{0:t+1})^\top } \in \bbR^{(t+2) \times (t+2)} , &
    \Sigma_{q,\out}^{t+1} &= \expt{ \sfF_{p,\out}^{0:t+1} (\sfF_{p,\out}^{0:t+1})^\top } \in \bbR^{(t+2) \times (t+2)} . &
& \label{eqn:cov_q}
\end{align}
\end{subequations}

Note that by the definitions in \Cref{eqn:cov_p,eqn:cov_q}, for every $ t\ge0 $, $ \Sigma_{q,\inn}^t, \Sigma_{q,\out}^t, \Sigma_{p,\inn}^{t+1}, \Sigma_{p,\out}^{t+1} $ are the leading principal minors (of order $ t+1 $) of $ \Sigma_{q,\inn}^{t+1}, \Sigma_{q,\out}^{t+1}, \Sigma_{p,\inn}^{t+2}, \Sigma_{p,\out}^{t+2} $, respectively. 
As a notational convention, we index elements in $ \Sigma_{q,\inn}^t, \Sigma_{q,\out}^t $ with $ (i,j) \in \brace{0,1,\cdots,t}^2 $ and index elements in $ \Sigma_{p,\inn}^{t+1}, \Sigma_{p,\out}^{t+1} $ with $ (i,j) \in \brace{1,2,\cdots,t+1}^2 $, so that we conveniently have, e.g., $ (\Sigma_{q,\inn}^t)_{r,s} = \expt{\sfQ_\inn^r \sfQ_\inn^s} $ and $ (\Sigma_{p,\inn}^{t+1})_{r+1,s+1} = \expt{\sfP_\inn^{r+1} \sfP_\inn^{s+1}} $ for any $ 0\le r,s\le t $. 

To summarize, for any $ t\ge0 $, \Cref{eqn:SE_OQ} defines random variables $ \sfA $, $ (\sfQ_\inn^s)_{0\le s\le t}$, $ (\sfQ_\out^s)_{0\le s\le t}$, $ (\sfF_{q,\inn}^s)_{0\le s\le t}$, $ (\sfF_{q,\out}^s)_{0\le s\le t}$, $ (\sfP_\inn^{s+1})_{0\le s\le t}$, $ (\sfP_\out^{s+1})_{0\le s\le t}$, $ (\sfF_{p,\inn}^s)_{0\le s\le t+1}$, $ (\sfF_{p,\out}^s)_{0\le s\le t+1} $. 
The main technical result of this section, \Cref{thm:SE_abs_VAMP} below, asserts that these random variables identify the Wasserstein limits of the iterates in \Cref{eqn:abs_VAMP}. 

\begin{enumerate}[label=(A\arabic*)]
\setcounter{enumi}{\value{asmpctr}}
    \item \label[asmp]{asmp:nondegenerate_abs} 
    For every $t\ge0$, there do not exist constants $ \beta_0, \cdots, \beta_t, \alpha_0, \cdots, \alpha_t \in \bbR $ such that 
    \begin{align}
    &&
        &\textnormal{either}&
        \sfF_{p,\inn,1}^{t+1} &= \sum_{s = 0}^t \beta_s \sfF_{p,\inn,1}^s & 
        &\textnormal{or}&
        \sfF_{p,\out}^{t+1} &= \sum_{s = 0}^t \alpha_s \sfF_{p,\out}^s & 
    & \notag 
    \end{align}
    holds almost surely. 
\setcounter{asmpctr}{\value{enumi}}
\end{enumerate}

For any $ t\ge0 $, let $ \bq_{\inn,1}^t \coloneqq \matrix{ q_{\inn,1}^0 & \cdots & q_{\inn,1}^t } \in \bbR^{d\times(t+1)} $, $ \bm{\sfQ}_{\inn,1}^{t} \coloneqq \matrix{ \sfQ_{\inn,1}^0 & \cdots & \sfQ_{\inn,1}^t } \in \bbR^{1\times(t+1)} $ and other notation including $ \bf_{q,\inn}^t$, $\bp_\inn^{t+1}$, $\bf_{p,\inn,1}^{t+1}$, $\bq_\out^t$, $\bf_{q,\out,1}^t$, $\bp_{\out,1}^{t+1}$, $\bf_{p,\out}^{t+1} $, $ \bm{\sfF}_{q,\inn}^t$, $\bm{\sfP}_{\inn}^{t+1}$, $\bm{\sfF}_{p,\inn,1}^{t+1}$, $\bm{\sfQ}_{\out}^t$, $\bm{\sfF}_{q,\out,1}^t$, $\bm{\sfP}_{\out,1}^{t+1}$, $\bm{\sfF}_{p,\out}^{t+1} $ is similarly defined. 

\begin{theorem}[State evolution of abstract VAMP]
\label{thm:SE_abs_VAMP}
Consider the abstract VAMP iteration \Cref{eqn:abs_VAMP} with initialization \Cref{eqn:abs_VAMP_init}, side information \Cref{eqn:side_info} and denoising functions in \Cref{eqn:fpq} specialized to \Cref{eqn:F_config}. 
Let \Cref{asmp:init,asmp:4matrices,asmp:regularity,asmp:div_free,asmp:nondegenerate_abs} hold. 
Then for any fixed $ t\ge0 $, 
\begin{align}
& \matrix{
    q_{\inn,2} &
    \bq_{\inn,1}^t &
    \bf_{q,\inn}^t & 
    \bp_\inn^{t+1} & 
    f_{p,\inn,2} & \bf_{p,\inn,1}^{t+1} & 
    w_{p,\inn} & w_{q,\inn}
} \notag \\
&\qquad\wto \matrix{
    \sfQ_{\inn,2} &
    \bm{\sfQ}_{\inn,1}^{t} &
    \bm{\sfF}_{q,\inn}^{t} &
    \bm{\sfP}_\inn^{t+1} &
    \sfF_{p,\inn,2} & \bm{\sfF}_{p,\inn,1}^{t+1} &
    \sfW_{p,\inn}^\top & \sfW_{q,\inn}^\top
} , \label{eqn:concl1} \\
& \matrix{
    \bq_\out^t &
    f_{q,\out,2} & 
    \bf_{q,\out,1}^t & 
    p_{\out,2} &
    \bp_{\out,1}^{t+1} & 
    \bf_{p,\out}^{t+1} & 
    w_{p,\out} & w_{q,\out}
} \notag \\
&\qquad\wto \matrix{
    \bm\sfQ_\out^{t} &
    \sfF_{q,\out,2} &
    \bm\sfF_{q,\out,1}^{t} &
    \sfP_{\out,2} &
    \bm\sfP_{\out,1}^{t+1} &
    \bm\sfF_{p,\out}^{t+1} &
    \sfW_{p,\out}^\top & \sfW_{q,\out}^\top
} . \label{eqn:concl2} 
\end{align}
\end{theorem}


Since the abstract VAMP in \Cref{eqn:abs_VAMP} is configured to simulate the GVAMP dynamics in \Cref{eqn:GVAMP} (see \Cref{eqn:config}), the state evolution result in \Cref{thm:SE_abs_VAMP} for the former iteration can be translated to that for the latter, as stated in \Cref{thm:SE_GVAMP}. 

\begin{proof}[Proof of \Cref{thm:SE_GVAMP}]
We specialize the vectors $ \phi_t, \wt{\phi}_t, \psi_t, \wt{\psi}_t $ (for every $ t\ge0 $) according to \Cref{eqn:4functions}. 
Due to the regularity of the functions $ \Phi_t, \wt{\Phi}_t, \Psi_t, \wt{\Psi}_t $ guaranteed by \Cref{asmp:nondegenerate}, \Cref{prop:wto2} implies the Wasserstein convergence required in \Cref{asmp:4matrices}. 
In particular, we can identify for every $t\ge0$, 
\begin{align}
&&
    \sfPhi_t &= \Phi_t(\sfLambda_d^2) , & 
    \wt{\sfPhi}_t &= \wt{\Phi}_t(\sfLambda_d) , & 
    \sfPsi_t &= \Psi_t(\sfLambda_n^2) , & 
    \wt{\sfPsi}_t &= \wt{\Psi}_t(\sfLambda_n) . & 
& \notag 
\end{align}
The trace-free property of $ \Phi_t, \Psi_t $ in \Cref{asmp:tr_free} implies the zero mean condition of $ \sfPhi_t, \sfPsi_t $ required in \Cref{asmp:4matrices}. 
Therefore, \Cref{thm:SE_abs_VAMP} is applicable to the present setting. 

By induction on $t\ge0$, we will prove the following statements which allow us to translate the state evolution result in \Cref{thm:SE_abs_VAMP} to \Cref{thm:SE_GVAMP}. 
\begin{enumerate}
    \item \label{itm:ind1} $ \wt{\fra}_t = \wt{a}_t, \wt{\frb}_t = \wt{b}_t, \fra_{t+1} = a_{t+1}, \frb_{t+1} = b_{t+1} $. 

    \item \label{itm:ind2} The joint law of 
    \begin{align}
        \matrix{
            \sfR_1 & \cdots & \sfR_{t+1} & 
            \wt{\sfR}_0 & \cdots & \wt{\sfR}_{t+1} & 
            \sfB_* & \sfTheta^\top
        } \label{eqn:law_r} 
    \end{align}
    is equal to that of 
    \begin{align}
        \matrix{
            \frb_1 \sfF_{p,\inn,2} + \sfP_\inn^1 & \cdots & \frb_{t+1} \sfF_{p,\inn,2} + \sfP_\inn^{t+1} & 
            \sfF_{p,\inn,1}^0 + \wt{\frb}_0 \sfF_{p,\inn,2} & \cdots & \sfF_{p,\inn,1}^{t+1} + \wt{\frb}_{t+1} \sfF_{p,\inn,2} & 
            \sfW_{p,\inn,1} & \sfW_{p,\inn,-}^\top 
        } . \label{eqn:law_in}
    \end{align}
    The joint law of 
    \begin{align}
        \matrix{
            \sfP_1 & \cdots & \sfP_{t+1} & 
            \wt{\sfP}_0 & \cdots & \wt{\sfP}_{t+1} & 
            \sfZ & \sfE & \sfXi^\top
        } \label{eqn:law_p}
    \end{align}
    is equal to that of
    \begin{align}
        \matrix{
            \fra_1 \sfP_{\out,2} + \sfP_{\out,1}^1 & \cdots & \fra_{t+1} \sfP_{\out,2} + \sfP_{\out,1}^{t+1} & 
            \sfF_{p,\out}^0 + \wt{\fra}_0 \sfP_{\out,2} & \cdots & \sfF_{p,\out}^{t+1} + \wt{\fra}_{t+1} \sfP_{\out,2} & 
            \sfP_{\out,2} & \sfW_{p,\out,1} & \sfW_{p,\out,-}^\top
        } . \label{eqn:law_out}
    \end{align}
\end{enumerate}

Recall $ \sfP_{\out,2} = \sfZ , \sfF_{p,\inn,2} = \sfW_{p,\inn,1} = \sfB_* , \sfW_{p,\out,1} = \sfE , \sfW_{p,\inn,-} = \sfTheta , \sfW_{p,\out,-} = \sfXi $. 
We identify the random variables $ (\sfK_{s+1})_{0\le s\le t}, (\sfJ_{s+1})_{0\le s\le t} $ with the following change of variables: 
\begin{align}
&& 
    \sfK_{s+1} &= \sfP_\inn^{s+1} , & 
    \sfJ_{s+1} &= \sfP_{\out,1}^{s+1} . & 
& \label{eqn:KJ_to_PP} 
\end{align}
Using this and comparing the definitions \Cref{eqn:wtP_wtR} with \Cref{eqn:Fpin1,eqn:Fpout}, we further identify $ (\wt{\sfR}_s)_{0\le s\le t+1}, (\wt{\sfP}_s)_{0\le s\le t+1} $ as
\begin{align}
&&
    \wt{\sfR}_{s} &= \sfF_{p,\inn,1}^{s} + \wt{\frb}_{s} \sfB_* , & 
    \wt{\sfP}_{s} &= \sfF_{p,\out}^{s} + \wt{\fra}_{s} \sfZ . & 
& \label{eqn:relation} 
\end{align}
To conclude the convergence result in \Cref{eqn:conv_rp}, it remains to verify that the joint distribution of the random variables on the left-hand sides above is consistent with the descriptions in \Cref{eqn:SE_GVAMP}. 

\paragraph{\textsc{Base case.}}
For $t=0$, we have $ \wt{\fra}_0 = \wt{a}_0, \wt{\frb}_0 = \wt{b}_0 $ by definition \Cref{eqn:frab_tilde0}. 
Comparing \Cref{eqn:ab,eqn:frab} then gives $ \fra_1 = a_1, \frb_1 = \frb_1 $. 
Therefore \Cref{itm:ind1} holds at $t=0$. 
Given \Cref{itm:ind1} together with \Cref{eqn:SE_OQ_init_in,eqn:SE_OQ_init_out}, it remains to verify $ \sfR_1 = b_1 \sfB_* + \sfP_\inn^1, \sfP_1 = a_1 \sfZ + \sfP_{\out,1}^1 $ which, upon recalling \Cref{eqn:Fpin1,eqn:Fpout,eqn:relation} and comparing them with \Cref{eqn:wtP_wtR}, then implies \Cref{itm:ind2} for $t=0$. 
To this end, we note from \Cref{eqn:Pin,eqn:Pout} that $ \sfP_\inn^1, \sfP_{\out,1}^1 $ are centered Gaussians independent of $ \sfB_*, \sfZ $. 
So we only need to check that their covariances $ \Sigma_{p,\inn}^1, \Sigma_{p,\out}^1 $ equal $ \expt{\sfK_1^2}, \expt{\sfJ_1^2} $ given in \Cref{eqn:K_cov,eqn:J_cov}, respectively.
We evaluate the former quantities according the their definitions in \Cref{eqn:cov_p}. 
For $ \Sigma_{p,\inn}^1 $, we have 
\begin{align}
   &  \Sigma_{p,\inn}^1 = \expt{(\sfF_{q,\inn}^0)^2}
    = \expt{\brace{ \Phi_0(\sfLambda_d^2) (\sfQ_{\inn,1}^0 + \wt{\frb}_0 \sfQ_{\inn,2}) + \wt{\Phi}_0(\sfLambda_d) (\sfQ_\out^0 + \wt{\fra}_0 \sfLambda_d \sfQ_{\inn,2}) - \frb_1 \sfQ_{\inn,2} }^2} \notag \\
    &\explain{} \expt{\Phi_0(\sfLambda_d^2)^2} \Sigma_{q,\inn}^0 + \expt{\paren{ \wt{\frb}_0 \Phi_0(\sfLambda_d^2) + \wt{\fra}_0 \wt{\Phi}_0(\sfLambda_d) \sfLambda_d - \frb_1}^2} \rho + \expt{\wt{\Phi}_0(\sfLambda_d)^2} \Sigma_{q,\out}^0 \notag \\
    &= \expt{\Phi_0(\sfLambda_d^2)^2} \expt{\wt{\sfR}_0^2} - b_1^2 \rho + \expt{\wt{\Phi}_0(\sfLambda_d)^2} \paren{\expt{\wt{\sfP}_0^2} - \wt{a}_0^2 \sigma^2} + \expt{\wt{\Phi}_0(\sfLambda_d)^2 \sfLambda_d^2} \wt{a}_0^2 \rho + 2 \expt{\Phi_0(\sfLambda_d^2) \wt{\Phi}_0(\sfLambda_d) \sfLambda_d} \wt{a}_0 \wt{b}_0 \rho . \notag 
\end{align}
The penultimate line is taken from \Cref{eqn:Sigmapin1_calc} in the proof of \Cref{thm:SE_abs_VAMP}.
The last line is obtained by using \Cref{itm:ind1} for $t=0$ just shown, plugging in the definitions $ \Sigma_{q,\inn}^0 = \expt{\wt{\sfR}_0^2} - \wt{a}_0^2 \rho $ and $ \Sigma_{q,\out}^0 = \expt{\wt{\sfP}_0^2} - \wt{b}_0 \sigma^2 $ from \Cref{eqn:Sigmaq_inout0}, and simplifying. 
The result coincides with $ \expt{\sfK_1^2} $ in \Cref{eqn:K_cov}, as desired. 
For $ \Sigma_{p,\out}^1 $, similar calculations starting from \Cref{eqn:Sigmapout1_calc} lead to the claimed result $ \Sigma_{p,\out}^1 = \expt{\sfJ_1^2} $. 
This proves \Cref{itm:ind2} for $t=0$. 

\paragraph{\textsc{Induction step.}}

Assuming the validity of \Cref{itm:ind1,itm:ind2} for time up to $ t-1 $, we verify their validity for time $t$. 
Using the induction hypothesis in \Cref{itm:ind2} and comparing the definitions \Cref{eqn:ab_tilde,eqn:frab_tilde}, we immediately get $ \wt{\fra}_t = \wt{a}_t, \wt{\frb}_t = \wt{b}_t $. 
In view of \Cref{eqn:ab,eqn:frab}, this in turn implies $ \fra_{t+1} = a_{t+1}, \frb_{t+1} = b_{b+1} $, thereby justifying \Cref{itm:ind1} for time $t$. 
In the sequel, we will therefore replace all occurrences of $ \wt{\fra}_s, \wt{\frb}_s, \fra_{s+1}, \frb_{s+1} $ with $ \wt{a}_s, \wt{b}_s, a_{s+1}, b_{s+1} $, for any $ 0\le s \le t $. 

According to the identification \Cref{eqn:KJ_to_PP} and the definition \Cref{eqn:Pin,eqn:Pout} of the right-hand side thereof, we have that $ (\sfK_s)_{1\le s\le t+1} $ and $ (\sfJ_s)_{1\le s\le t+1} $ are two independent centered Gaussian processes that are also independent of $ \sfLambda_n, \sfLambda_d, \sfB_*, \sfZ, \sfE, \sfTheta, \sfXi $. 
So comparing the description \Cref{eqn:SE_GVAMP} of the laws of \Cref{eqn:law_r,eqn:law_p} with the description \Cref{eqn:SE_OQ} of the laws of \Cref{eqn:law_in,eqn:law_out}, we see that to prove \Cref{itm:ind2} for time $t$, we only need to examine the covariance structures of the Gaussian processes $ (\sfK_s)_{1\le s\le t+1} $ and $ (\sfJ_s)_{1\le s\le t+1} $. 
The covariance of $(\sfJ_s)_{1\le s\le t+1}$ is given as follows. 
For any $ 0\le r,s\le t $, 
\begin{align}
    & \expt{ \sfJ_{r+1} \sfJ_{s+1} } = \expt{ \sfP_{\out,1}^{r+1} \sfP_{\out,1}^{s+1} }
    = (\Sigma_{p,\out}^{t+1})_{r+1,s+1}
    = \expt{ \sfF_{q,\out,1}^r \sfF_{q,\out,1}^s } \notag \\
    &= \bbE\Big[
        \Big( \Psi_r(\sfLambda_n^2) (\sfQ_\out^r + \wt{a}_r \sfLambda_n \sfQ_{\inn,2}) + \wt{\Psi}_r(\sfLambda_n) (\sfQ_{\inn,1}^r + \wt{b}_r \sfQ_{\inn,2}) - a_{r+1} \sfLambda_n \sfQ_{\inn,2} \Big) \notag \\
        &\qquad \cdot \Big( \Psi_s(\sfLambda_n^2) (\sfQ_\out^s + \wt{a}_s \sfLambda_n \sfQ_{\inn,2}) + \wt{\Psi}_s(\sfLambda_n) (\sfQ_{\inn,1}^s + \wt{b}_s \sfQ_{\inn,2}) - a_{s+1} \sfLambda_n \sfQ_{\inn,2} \Big) 
    \Big] \notag \\
    &= \expt{ \Psi_r(\sfLambda_n^2) \Psi_s(\sfLambda_n^2) (\sfQ_\out^r + \wt{a}_r \sfLambda_n \sfQ_{\inn,2}) (\sfQ_\out^s + \wt{a}_s \sfLambda_n \sfQ_{\inn,2}) } \label{eqn:line1} \\
    &\quad + \expt{ \wt{\Psi}_r(\sfLambda_n) \wt{\Psi}_s(\sfLambda_n) } \expt{ (\sfQ_{\inn,1}^r + \wt{b}_r \sfQ_{\inn,2}) (\sfQ_{\inn,1}^s + \wt{b}_s \sfQ_{\inn,2}) }
    + a_{r+1} a_{s+1} \expt{ \sfLambda_n^2 } \expt{ \sfQ_{\inn,2}^2 } \label{eqn:line2} \\
    &\quad + \expt{ \Psi_r(\sfLambda_n^2) \wt{\Psi}_s(\sfLambda_n) (\sfQ_\out^r + \wt{a}_r \sfLambda_n \sfQ_{\inn,2}) (\sfQ_{\inn,1}^s + \wt{b}_s \sfQ_{\inn,2}) } \label{eqn:line3} \\
    &\quad + \expt{ \Psi_s(\sfLambda_n^2) \wt{\Psi}_r(\sfLambda_n) (\sfQ_\out^s + \wt{a}_s \sfLambda_n \sfQ_{\inn,2}) (\sfQ_{\inn,1}^r + \wt{b}_r \sfQ_{\inn,2}) } \label{eqn:line4} \\
    &\quad - a_{r+1} \paren{ \expt{\Psi_s(\sfLambda_n^2) \sfLambda_n (\sfQ_\out^s + \wt{a}_s \sfLambda_n \sfQ_{\inn,2}) \sfQ_{\inn,2}} + \expt{\wt{\Psi}_s(\sfLambda_n) \sfLambda_n} \expt{(\sfQ_{\inn,1}^s + \wt{b}_s \sfQ_{\inn,2}) \sfQ_{\inn,2}} } \label{eqn:line5} \\
    &\quad - a_{s+1} \paren{ \expt{\Psi_r(\sfLambda_n^2) \sfLambda_n (\sfQ_\out^r + \wt{a}_r \sfLambda_n \sfQ_{\inn,2}) \sfQ_{\inn,2}} + \expt{\wt{\Psi}_r(\sfLambda_n) \sfLambda_n} \expt{(\sfQ_{\inn,1}^r + \wt{b}_r \sfQ_{\inn,2}) \sfQ_{\inn,2}} } , \label{eqn:line6}
\end{align}
where the second line in the above display follows from the definition \Cref{eqn:Fqout1} of $ \sfF_{q,\out,1}^t $. 
We will separately examine the terms in \Cref{eqn:line1,eqn:line2,eqn:line3,eqn:line4,eqn:line5,eqn:line6}. 

We first note several facts that will be used in the  calculations: 
\begin{align}
    \expt{\sfQ_\out^r \sfQ_\out^s} 
    &= (\Sigma_{q,\out}^t)_{r,s} 
    = \expt{\sfF_{p,\out}^r \sfF_{p,\out}^s} 
    = \expt{\paren{\wt{\sfP}_r - \wt{a}_r \sfZ} \paren{\wt{\sfP}_s - \wt{a}_s \sfZ}} 
    = \expt{\wt{\sfP}_r \wt{\sfP}_s} - \wt{a}_r \wt{a}_s \sigma^2 , \label{eqn:obs1} \\
    \expt{\sfQ_{\inn,1}^r \sfQ_{\inn,1}^s} 
    &= (\Sigma_{q,\inn}^t)_{r,s} 
    = \expt{\sfF_{p,\inn,1}^r \sfF_{p,\inn,1}^s} 
    = \expt{\paren{\wt{\sfR}_r - \wt{b}_r \sfB_*} \paren{\wt{\sfR}_s - \wt{b}_s \sfB_*}} 
    = \expt{\wt{\sfR}_r \wt{\sfR}_s} - \wt{b}_r \wt{b}_s \rho , \label{eqn:obs2} \\
    \expt{\sfLambda_n^2} \expt{\sfQ_{\inn,2}^2}
    &= \ol{m}_2 \expt{\sfB_*^2} = \ol{\kappa}_2 \rho = \sigma^2 . \label{eqn:obs3}
\end{align}
In the above, we have used the
the relation \Cref{eqn:relation} and
the definitions of $ \wt{a}_t, \wt{b}_t $ in \Cref{eqn:ab_tilde}.

Now we simplify \Cref{eqn:line1} to: 
\begin{align}
    & \expt{ \Psi_r(\sfLambda_n^2) \Psi_s(\sfLambda_n^2) (\sfQ_\out^r + \wt{a}_r \sfLambda_n \sfQ_{\inn,2}) (\sfQ_\out^s + \wt{a}_s \sfLambda_n \sfQ_{\inn,2}) } \notag \\
    &= \expt{ \Psi_r(\sfLambda_n^2) \Psi_s(\sfLambda_n^2) } \expt{ \sfQ_\out^r \sfQ_\out^s }
    + \wt{a}_r \wt{a}_s \expt{ \Psi_r(\sfLambda_n^2) \Psi_s(\sfLambda_n^2) \sfLambda_n^2 } \expt{ \sfQ_{\inn,2}^2 } \notag \\
    &= \expt{ \Psi_r(\sfLambda_n^2) \Psi_s(\sfLambda_n^2) } \paren{ \expt{\wt{\sfP}_r \wt{\sfP}_s} - \wt{a}_r \wt{a}_s \sigma^2 }
    + \wt{a}_r \wt{a}_s \expt{ \Psi_r(\sfLambda_n^2) \Psi_s(\sfLambda_n^2) \sfLambda_n^2 } \rho , \label{eqn:E2} 
\end{align}
where we have used the independence between $ (\sfQ_\inn^s)_{0\le s\le t} $ and $ (\sfQ_\out^s)_{0\le s\le t} $, and the fact \Cref{eqn:obs1}.  

The second expectation in \Cref{eqn:line2} can be simplified as follows: 
\begin{align}
    & \expt{ (\sfQ_{\inn,1}^r + \wt{b}_r \sfQ_{\inn,2}) (\sfQ_{\inn,1}^s + \wt{b}_s \sfQ_{\inn,2}) }
    = \expt{\sfQ_{\inn,1}^r \sfQ_{\inn,1}^s} + \wt{b}_r \wt{b}_s \expt{\sfQ_{\inn,2}^2}
    = \expt{\wt{\sfR}_r \wt{\sfR}_s} , \label{eqn:E1} 
\end{align}
where we have used the independence between $ \sfQ_{\inn,2} $ and $ (\sfQ_{\inn,1}^s)_{0\le s\le t} $ (see \Cref{eqn:Qin}), and the fact \Cref{eqn:obs2}.

We then move to \Cref{eqn:line3}:  
\begin{align}
    & \expt{ \Psi_r(\sfLambda_n^2) \wt{\Psi}_s(\sfLambda_n) (\sfQ_\out^r + \wt{a}_r \sfLambda_n \sfQ_{\inn,2}) (\sfQ_{\inn,1}^s + \wt{b}_s \sfQ_{\inn,2}) } \notag \\
    &= \wt{a}_r \expt{ \Psi_r(\sfLambda_n^2) \wt{\Psi}_s(\sfLambda_n) \sfLambda_n } \expt{ \sfQ_{\inn,2} (\sfQ_{\inn,1}^s + \wt{b}_s \sfQ_{\inn,2}) }
    = \wt{a}_r \wt{b}_s \expt{ \Psi_r(\sfLambda_n^2) \wt{\Psi}_s(\sfLambda_n) \sfLambda_n } \rho . \label{eqn:E3} 
\end{align}
Similarly, \Cref{eqn:line4} evaluates to
\begin{align}
    \expt{ \Psi_s(\sfLambda_n^2) \wt{\Psi}_r(\sfLambda_n) (\sfQ_\out^s + \wt{a}_s \sfLambda_n \sfQ_{\inn,2}) (\sfQ_{\inn,1}^r + \wt{b}_r \sfQ_{\inn,2}) } 
    &= \wt{a}_s \wt{b}_r \expt{ \Psi_s(\sfLambda_n^2) \wt{\Psi}_r(\sfLambda_n) \sfLambda_n } \rho . \label{eqn:E4}
\end{align}

Finally, consider the terms in the parentheses in \Cref{eqn:line6}: 
\begin{align}
    & \expt{\Psi_r(\sfLambda_n^2) \sfLambda_n (\sfQ_\out^r + \wt{a}_r \sfLambda_n \sfQ_{\inn,2}) \sfQ_{\inn,2}} + \expt{\wt{\Psi}_r(\sfLambda_n) \sfLambda_n} \expt{(\sfQ_{\inn,1}^r + \wt{b}_r \sfQ_{\inn,2}) \sfQ_{\inn,2}} \notag \\
    &= \wt{a}_r \expt{\Psi_r(\sfLambda_n^2) \sfLambda_n^2} \expt{\sfQ_{\inn,2}^2} + \wt{b}_r \expt{\wt{\Psi}_r(\sfLambda_n) \sfLambda_n} \expt{\sfQ_{\inn,2}^2} \notag \\
    &= \frac{\rho}{\delta} \paren{ \wt{a}_r \expt{\Psi_r(\sfLambda_d^2) \sfLambda_d^2} + \wt{b}_r \expt{\wt{\Psi}_r(\sfLambda_d) \sfLambda_d} }
    = \ol{\kappa}_2 \rho a_{r+1}
    = \sigma^2 a_{r+1} , \notag 
\end{align}
where the second line follows since both $\sfQ_\out^r$ and $\sfQ_{\inn,1}^r$ are independent of $\sfQ_{\inn,2}$. 
Similarly, the terms in the parentheses in \Cref{eqn:line5} can be simplified to $\sigma^2 a_{s+1}$. 

Putting all the terms together, we conclude that $\expt{\sfJ_{r+1} \sfJ_{s+1}}$ evaluates to \Cref{eqn:J_cov}.  

We use similar steps to check the covariance structure of $(\sfK_{r+1})_{0\le r\le t}$. 
For any $ 0\le r,s\le t $, 
\begin{align}
    \expt{ \sfK_{r+1} \sfK_{s+1} }
    &= \expt{ \sfP_\inn^{r+1} \sfP_\inn^{s+1} }
    = (\Sigma_{p,\inn}^{t+1})_{r+1,s+1}
    = \expt{ \sfF_{q,\inn}^r \sfF_{q,\inn}^s } \notag \\
    &= \bbE\Big[ 
        \paren{ \Phi_r(\sfLambda_d^2) (\sfQ_{\inn,1}^r + \wt{b}_r \sfQ_{\inn,2}) + \wt{\Phi}_r(\sfLambda_d) (\sfQ_\out^r + \wt{a}_r \sfLambda_d \sfQ_{\inn,2}) - b_{r+1} \sfQ_{\inn,2} } \notag \\
    &\qquad \cdot \paren{ \Phi_s(\sfLambda_d^2) (\sfQ_{\inn,1}^s + \wt{b}_s \sfQ_{\inn,2}) + \wt{\Phi}_s(\sfLambda_d) (\sfQ_\out^s + \wt{a}_s \sfLambda_d \sfQ_{\inn,2}) - b_{s+1} \sfQ_{\inn,2} }
    \Big] \notag \\
    &= \expt{ \Phi_r(\sfLambda_d^2) \Phi_s(\sfLambda_d^2) } \expt{ (\sfQ_{\inn,1}^r + \wt{b}_r \sfQ_{\inn,2}) (\sfQ_{\inn,1}^s + \wt{b}_s \sfQ_{\inn,2}) } \label{eqn:ln1} \\
    &\quad + \expt{ \wt{\Phi}_r(\sfLambda_d) \wt{\Phi}_s(\sfLambda_d) (\sfQ_\out^r + \wt{a}_r \sfLambda_d \sfQ_{\inn,2}) (\sfQ_\out^s + \wt{a}_s \sfLambda_d \sfQ_{\inn,2}) }
    + b_{r+1} b_{s+1} \expt{ \sfQ_{\inn,2}^2 } \label{eqn:ln2} \\
    &\quad + \expt{ \Phi_r(\sfLambda_d^2) \wt{\Phi}_s(\sfLambda_d) (\sfQ_{\inn,1}^r + \wt{b}_r \sfQ_{\inn,2}) (\sfQ_\out^s + \wt{a}_s \sfLambda_d \sfQ_{\inn,2}) } \label{eqn:ln3} \\
    &\quad + \expt{ \Phi_s(\sfLambda_d^2) \wt{\Phi}_r(\sfLambda_d) (\sfQ_{\inn,1}^s + \wt{b}_s \sfQ_{\inn,2}) (\sfQ_\out^r + \wt{a}_r \sfLambda_d \sfQ_{\inn,2}) } \label{eqn:ln4} \\
    &\quad - b_{r+1} \expt{ \wt{\Phi}_s(\sfLambda_d) (\sfQ_\out^s + \wt{a}_s \sfLambda_d \sfQ_{\inn,2}) \sfQ_{\inn,2} } 
    - b_{s+1} \expt{ \wt{\Phi}_r(\sfLambda_d) (\sfQ_\out^r + \wt{a}_r \sfLambda_d \sfQ_{\inn,2}) \sfQ_{\inn,2} } , \label{eqn:ln5} 
\end{align}
where we have used the trace-free property of $ \Phi_t $ (see \Cref{asmp:tr_free}) in \Cref{eqn:ln5}. 

We examine each term in \Cref{eqn:ln1,eqn:ln2,eqn:ln3,eqn:ln4,eqn:ln5} separately. 
The second expectation in \Cref{eqn:ln1} has been evaluated in \Cref{eqn:E1}. 
The first term in \Cref{eqn:ln2} can be computed similarly to \Cref{eqn:E2}: 
\begin{align}
    & \expt{ \wt{\Phi}_r(\sfLambda_d) \wt{\Phi}_s(\sfLambda_d) (\sfQ_\out^r + \wt{a}_r \sfLambda_d \sfQ_{\inn,2}) (\sfQ_\out^s + \wt{a}_s \sfLambda_d \sfQ_{\inn,2}) } \notag \\
    &= \expt{ \wt{\Phi}_r(\sfLambda_d) \wt{\Phi}_s(\sfLambda_d) } \paren{ \expt{\wt{\sfP}_r \wt{\sfP}_s} - \wt{a}_r \wt{a}_s \sigma^2 }
    + \wt{a}_r \wt{a}_s \expt{ \wt{\Phi}_r(\sfLambda_d) \wt{\Phi}_s(\sfLambda_d) \sfLambda_d^2 } \rho . \notag 
\end{align}
The terms in \Cref{eqn:ln3,eqn:ln4} can be computed similarly to \Cref{eqn:E3,eqn:E4}: 
\begin{align}
    \expt{ \Phi_r(\sfLambda_d^2) \wt{\Phi}_s(\sfLambda_d) (\sfQ_{\inn,1}^r + \wt{b}_r \sfQ_{\inn,2}) (\sfQ_\out^s + \wt{a}_s \sfLambda_d \sfQ_{\inn,2}) } 
    &= \wt{a}_s \wt{b}_r \expt{ \Phi_r(\sfLambda_d^2) \wt{\Phi}_s(\sfLambda_d) \sfLambda_d } \rho , \notag \\
    \expt{ \Phi_s(\sfLambda_d^2) \wt{\Phi}_r(\sfLambda_d) (\sfQ_{\inn,1}^s + \wt{b}_s \sfQ_{\inn,2}) (\sfQ_\out^r + \wt{a}_r \sfLambda_d \sfQ_{\inn,2}) } 
    &= \wt{a}_r \wt{b}_s \expt{ \Phi_s(\sfLambda_d^2) \wt{\Phi}_r(\sfLambda_d) \sfLambda_d } \rho . \notag 
\end{align}
The first expectation in \Cref{eqn:ln5} evaluates to: 
\begin{align}
    \expt{ \wt{\Phi}_s(\sfLambda_d) (\sfQ_\out^s + \wt{a}_s \sfLambda_d \sfQ_{\inn,2}) \sfQ_{\inn,2} }
    &= \wt{a}_s \expt{ \wt{\Phi}_s(\sfLambda_d) \sfLambda_d } \rho 
    = b_{s+1} \rho , \notag 
\end{align}
where the last step is by the definition of $ b_{s+1} $ in \Cref{eqn:ab}. 
Similarly the second expectation in \Cref{eqn:ln5} is: 
\begin{align}
    \expt{ \wt{\Phi}_r(\sfLambda_d) (\sfQ_\out^r + \wt{a}_r \sfLambda_d \sfQ_{\inn,2}) \sfQ_{\inn,2} }
    &= b_{r+1} \rho . \notag 
\end{align}
Putting all the terms back, we obtain \Cref{eqn:K_cov} as the desired expression for $ \expt{\sfK_{r+1} \sfK_{s+1}} $.
\end{proof}


\begin{proof}[Proof of \Cref{thm:SE_abs_VAMP}]
The proof is by induction on $t\ge0$. 
We will inductively show the following statements which imply the conclusion of \Cref{thm:SE_abs_VAMP}: 
\begin{enumerate}
    \item\label{itm:hyp1} For every $ t\ge0 $, $ \expt{ \sfF_{p,\inn,2} \sfF_{p,\inn,1}^t } = \expt{ \sfF_{q,\out,2} \sfF_{q,\out,1}^t } = 0 $. 

    \item\label{itm:hyp3} For every $t\ge0$, $ \Sigma_{q,\inn}^t, \Sigma_{q,\out}^t, \Sigma_{p,\inn}^{t+1}, \Sigma_{p,\out}^{t+1} $ are non-singular. 

    \item\label{itm:hyp4} For every $ t\ge0 $, there exist Gaussian random variables $ \sfU_{q,\inn}^t, \sfU_{q,\out}^t $ with zero mean and strictly positive variances, independent of each other and of everything else such that 
    \begin{align}
    \begin{split}
        \sfQ_{\inn,1}^t &= \sum_{s = 0}^{t-1} \sfQ_{\inn,1}^s \paren{(\Sigma_{q,\inn}^{t-1})^{-1} (\Sigma_{q,\inn}^t)_{0:t-1,t}}_s + \sfU_{q,\inn}^t , \\
        \sfQ_{\out}^t &= \sum_{s = 0}^{t-1} \sfQ_{\out}^s \paren{(\Sigma_{q,\out}^{t-1})^{-1} (\Sigma_{q,\out}^t)_{0:t-1,t}}_s + \sfU_{q,\out}^t . 
    \end{split}
    \label{eqn:hyp4_in} 
    \end{align}
    There exist Gaussian random variables $ \sfU_{p,\inn}^{t+1}, \sfU_{p,\out}^{t+1} $ with zero mean and strictly positive variances, independent of each other and of everything else such that 
    \begin{align}
    \begin{split}
        \sfP_\inn^{t+1} &= \sum_{s = 0}^{t-1} \sfP_\inn^{s+1} \paren{(\Sigma_{p,\inn}^t)^{-1} (\Sigma_{p,\inn}^{t+1})_{1:t,t+1}}_{s+1} + \sfU_{p,\inn}^{t+1} , \\
        \sfP_{\out,1}^{t+1} &= \sum_{s = 0}^{t-1} \sfP_{\out,1}^{s+1} \paren{(\Sigma_{p,\out}^t)^{-1} (\Sigma_{p,\out}^{t+1})_{1:t,t+1}}_{s+1} + \sfU_{p,\out}^{t+1} . 
    \end{split}
    \label{eqn:hyp4_out}
    \end{align}

    \item\label{itm:hyp2} \Cref{eqn:concl1,eqn:concl2} hold. 
\end{enumerate}

\paragraph{\textsc{Base case.}}
We take $ t = 0 $ as the base case and verify the validity of \Cref{itm:hyp1,itm:hyp2,itm:hyp3,itm:hyp4} for $t=0$. 

By the choice of $ \wt{\frb}_0 $ in \Cref{eqn:frab_tilde0}, we have $ \expt{\sfF_{p,\inn,2} \sfF_{p,\inn,1}^0} = \expt{(\wt{\sfR}_0 - \wt{\frb}_0 \sfB_*) \sfB_*} = 0 $ which proves the first statement in \Cref{itm:hyp1}. 

The non-singularity of $ \Sigma_{q,\inn}^0, \Sigma_{q,\out}^0 \in \bbR $ in \Cref{itm:hyp3} reduces to positivity at $t=0$, which follows from \Cref{asmp:init}, as noted in \Cref{eqn:Sigma0_pos}. 

We now verify \Cref{eqn:concl1,eqn:concl2} in \Cref{itm:hyp2} for $t=0$. 
Recall that $ f_{p,\inn,2} = w_{p,\inn,1} = \beta_* $.
The joint convergence of 
\begin{align}
&&
    & \matrix{f_{p,\inn,2} & f_{p,\inn,1}^0 & w_{p,\inn} & w_{q,\inn}} & 
    & \textnormal{and} &
    & \matrix{f_{p,\out}^0 & w_{p,\out} & w_{q,\out}} & 
& \label{eqn:notime}
\end{align}
holds due to the initial condition \Cref{eqn:SE_OQ}. 
We then show joint convergence of the iterates involved in the left of \Cref{eqn:concl1,eqn:concl2} (for $t=0$) by examining the updates \Cref{eqn:OQ1,eqn:OQ2,eqn:OQ3,eqn:OQ4} sequentially.
The first update \Cref{eqn:OQ1} generates 
\begin{align}
&&
    q_\inn^0 &= \matrix{Q^\top (\wt{r}^0 - \wt{\frb}_0 \beta_*) & Q^\top \beta_*} &
    & \textnormal{and} &
    q_\out^0 &= O^\top (\wt{p}^0 - \wt{\fra}_0 z) = O^\top \wt{p}^0 - \wt{\fra}_0 \Lambda Q^\top \beta_* &
& \notag 
\end{align}
whose joint convergence with vectors in \Cref{eqn:notime} is guaranteed by \Cref{prop:wto4}: 
\begin{align}
\begin{split}
    \matrix{f_{p,\inn,2} & f_{p,\inn,1}^0 & w_{p,\inn} & w_{q,\inn} & q_{\inn,2} & q_{\inn,1}^0}
    &\wto \matrix{\sfF_{p,\inn,2} & \sfF_{p,\inn,1}^0 & \sfW_{p,\inn}^\top & \sfW_{q,\inn}^\top & \sfQ_{\inn,2} & \sfQ_{\inn,1}^0} , \\
    \matrix{f_{p,\out}^0 & w_{p,\out} & w_{q,\out} & q_\out^0}
    &\wto \matrix{\sfF_{p,\out}^0 & \sfW_{p,\out}^\top & \sfW_{q,\out}^\top & \sfQ_\out^0} , 
\end{split}
\label{eqn:t0_q} 
\end{align}
where 
\begin{align}
    \matrix{ \sfQ_{\inn,2} \\ \sfQ_{\inn,1}^0 } &\sim \cN\paren{\matrix{0 \\ 0}, \matrix{\rho & 0 \\ 0 & \Sigma_{q,\inn}^0}} \label{eqn:Qin0_cov} 
\end{align}
is independent of $ \sfF_{p,\inn,2}, \sfF_{p,\inn,1}^0, \sfW_{p,\inn}, \sfW_{q,\inn} $, and $ \sfQ_\out^0 \sim \cN(0, \Sigma_{q,\out}^0) $ is independent of $ \sfF_{p,\out}^0, \sfW_{p,\out}, \sfW_{q,\out} $. 
By positivity of $ \Sigma_{q,\inn}^0 $ and $ \rho $ (see \Cref{eqn:Sigma0_pos} and \Cref{asmp:signal}), the covariance matrix of $\matrix{ \sfQ_{\inn,2} & \sfQ_{\inn,1}^0 }^\top$ is strictly positive definite. 
Note that $ \matrix{\sfQ_{\inn,2} & \sfQ_{\inn,1}^0}^\top $ and $ \sfQ_\out^0 $ are also independent of each other. 
This proves \Cref{eqn:hyp4_in} in \Cref{itm:hyp4} for $t=0$. 

The second update \Cref{eqn:OQ2} generates 
\begin{align}
    f_{q,\inn}^0 &= \diag(\phi_0) Q^\top \wt{r}^0 + \diag_{n\times d}(\wt{\phi}_0)^\top O^\top \wt{p}^0 - \frb_1 Q^\top \beta_* , \notag \\
    f_{q,\out}^0 &= \matrix{\diag(\psi_0) O^\top \wt{p}^0 + \diag_{n\times d}(\wt{\psi}_0) Q^\top \wt{r}^0 - \fra_1 \Lambda Q^\top \beta_* & \Lambda Q^\top \beta_*}, \notag 
\end{align}
whose joint convergence with vectors in \Cref{eqn:t0_q} is guaranteed by \Cref{prop:wto4,prop:wto2}, 
\begin{subequations}
\label{eqn:t0_fq}     
\begin{align}
    & \matrix{f_{p,\inn,2} & f_{p,\inn,1}^0 & w_{p,\inn} & w_{q,\inn} & q_{\inn,2} & q_{\inn,1}^0 & f_{q,\inn}^0} \notag \\
    &\qquad \wto \matrix{\sfF_{p,\inn,2} & \sfF_{p,\inn,1}^0 & \sfW_{p,\inn}^\top & \sfW_{q,\inn}^\top & \sfQ_{\inn,2} & \sfQ_{\inn,1}^0 & \sfF_{q,\inn}^0} , \\
    & \matrix{f_{p,\out}^0 & w_{p,\out} & w_{q,\out} & q_\out^0 & f_{q,\out,2} & f_{q,\out,1}^0} \notag \\
    &\qquad \wto \matrix{\sfF_{p,\out}^0 & \sfW_{p,\out}^\top & \sfW_{q,\out}^\top & \sfQ_\out^0 & \sfF_{q,\out,2} & \sfF_{q,\out,1}^0} , 
\end{align}
\end{subequations}
where $ \sfF_{q,\inn}^0 , \sfF_{q,\out,1}^0, \sfF_{q,\out,2} $ are defined in \Cref{eqn:Fqin,eqn:Fqout1,eqn:Fqout2}, respectively. 
Now $ \Sigma_{p,\inn}^1, \Sigma_{p,\out}^1 \in \bbR $ can be computed as per \Cref{eqn:cov_p}. 
We argue that they must be strictly positive, thereby justifying \Cref{itm:hyp3} for $t=0$. 
For $ \Sigma_{p,\inn}^1 $, we write it as 
\begin{align}
    \Sigma_{p,\inn}^1 &= \expt{(\sfF_{q,\inn}^0)^2}
    = \expt{\brace{ \sfPhi_0 (\sfQ_{\inn,1}^0 + \wt{\frb}_0 \sfQ_{\inn,2}) + \wt{\sfPhi}_0 (\sfQ_\out^0 + \wt{\fra}_0 \sfLambda_d \sfQ_{\inn,2}) - \frb_1 \sfQ_{\inn,2} }^2} \notag \\
    &= \expt{\brace{ \sfPhi_0 \sfQ_{\inn,1}^0 + \paren{\wt{\frb}_0 \sfPhi_0 + \wt{\fra}_0 \wt{\sfPhi}_0 \sfLambda_d - \frb_1} \sfQ_{\inn,2} + \wt{\sfPhi}_0 \sfQ_\out^0 }^2} \notag \\
    &= \expt{\sfPhi_0^2} \Sigma_{q,\inn}^0 + \expt{\paren{ \wt{\frb}_0 \sfPhi_0 + \wt{\fra}_0 \wt{\sfPhi}_0 \sfLambda_d - \wt{\fra}_0 \expt{\wt{\sfPhi}_0 \sfLambda_d}}^2} \rho + \expt{\wt{\sfPhi}_0^2} \Sigma_{q,\out}^0 , \label{eqn:Sigmapin1_calc}
\end{align}
which can be seen positive since at least one of the first and last terms above is positive by \Cref{eqn:Sigma0_pos} and \Cref{asmp:4matrices}. 
Similarly, the positivity of 
\begin{align}
    \Sigma_{p,\out}^1
    &= \expt{\sfPsi_0^2} \Sigma_{q,\out}^0 + \expt{\paren{\wt{\fra}_0 \sfPsi_0 \sfLambda_n + \wt{\frb}_0 \wt{\sfPsi}_0 - \fra_1 \sfLambda_n}^2} \rho + \expt{\wt{\sfPsi}_0^2} \Sigma_{q,\inn}^0 \label{eqn:Sigmapout1_calc}
\end{align}
follows from the positivity of $ \Sigma_{q,\out}^0, \Sigma_{q,\inn}^0 $ and at least one of $ \expt{\sfPsi_0^2}, \expt{\wt{\sfPsi}_0^2} $. 

Moving on, the third update \Cref{eqn:OQ3} generates 
\begin{align}
    p_\inn^1 &= Q \diag(\phi_0) Q^\top \wt{r}^0 + Q \diag_{n\times d}(\wt{\phi}_0)^\top O^\top \wt{p}^0 - \frb_1 \beta_* = r^1 - \frb_1 \beta_* , \notag \\
    p_\out^1 &= \matrix{O \diag(\psi_0) O^\top \wt{p}^0 + O \diag_{n\times d}(\wt{\psi}_0) Q^\top \wt{r}^0 - \fra_1 O \Lambda Q^\top \beta_* & O \Lambda Q^\top \beta_*} = \matrix{p^1 - \fra_1 z & z} . \notag
\end{align}
To understand their joint convergence, we condition on all vectors on the left of \Cref{eqn:t0_fq}. 
The laws of $ Q, O $ are then conditioned on the events $ Q q_\inn^0 = f_{p,\inn}^0, O q_\out^0 = f_{p,\out}^0 $, respectively. 
The Wasserstein convergence of $ q_\inn^0 $ and $ q_\out^0 $ in \Cref{eqn:t0_q} implies $ d^{-1} (q_\inn^0)^\top q_\inn^0 \to \matrix{\Sigma_{q,\inn}^0 & 0 \\ 0 & \rho} \succ 0_{2\times2} $ (see \Cref{eqn:Qin0_cov}) and $ n^{-1} (q_\out^0)^\top q_\out^0 \to \Sigma_{q,\out}^0 > 0 $ (see \Cref{eqn:Sigma0_pos}). 
In particular, the left-hand sides are invertible for all sufficiently large $n,d$. 
Using \Cref{prop:haar_cond}, for all sufficiently large $n,d$, we can write the conditional law of $ p_\inn^1, p_\out^1 $ as
\begin{align}
\begin{split}
    \lr{p_\inn^1 \mid \brace{Q q_\inn^0 = f_{p,\inn}^0}} &= \lr{Q f_{q,\inn}^0 \mid \brace{Q q_\inn^0 = f_{p,\inn}^0}} \eqqlaw f_{p,\inn}^0 ((q_\inn^0)^\top q_\inn^0)^{-1} (q_\inn^0)^\top f_{q,\inn}^0 + \Pi_{(f_{p,\inn}^0)^\perp} \wt{Q} \Pi_{(q_\inn^0)^\perp}^\top f_{q,\inn}^0 , \\
    \lr{p_\out^1 \mid \brace{O q_\out^0 = f_{p,\out}^0}} &= \lr{O f_{q,\out}^0 \mid \brace{O q_\out^0 = f_{p,\out}^0}} \eqqlaw f_{p,\out}^0 ((q_\out^0)^\top q_\out^0)^{-1} (q_\out^0)^\top f_{q,\out}^0 + \Pi_{(f_{p,\out}^0)^\perp} \wt{O} \Pi_{(q_\out^0)^\perp}^\top f_{q,\out}^0 , 
\end{split}
\label{eqn:p1}
\end{align}
where $ (\wt{Q}, \wt{O}) \sim \haar(\bbO(d-2)) \ot \haar(\bbO(n-1)) $ are independent of everything else. 
We claim that 
\begin{align}
&&
    \frac{1}{d} (q_\inn^0)^\top f_{q,\inn}^0 &\to \matrix{0 \\ 0} , & 
    \frac{1}{n} (q_\out^0)^\top f_{q,\out}^0 &\to \matrix{0 & 0} , & 
& \label{eqn:q_fq_0}
\end{align}
entry-wise. 
To show the first claim in \Cref{eqn:q_fq_0}, using the joint convergence result \Cref{eqn:t0_fq}, zero mean of $ \sfPhi_0 $ and the definition of $ \frb_1$ in \Cref{eqn:frab}, we have:
\begin{align}
    \frac{1}{d} q_{\inn,2}^\top f_{q,\inn}^0 &\to \expt{\sfQ_{\inn,2} \brace{\sfPhi_0 (\sfQ_{\inn,1}^0 + \wt{\frb}_0 \sfQ_{\inn,2}) + \wt{\sfPhi}_0 (\sfQ_\out^0 + \wt{\fra}_0 \sfLambda_d \sfQ_{\inn,2}) - \frb_1 \sfQ_{\inn,2}}} \notag \\
    &= \wt{\frb}_0 \expt{\sfPhi_0} \expt{\sfQ_{\inn,2}^2} + \wt{\fra}_0 \expt{\wt{\sfPhi}_0 \sfLambda_d} \expt{\sfQ_{\inn,2}^2} - \frb_1 \expt{\sfQ_{\inn,2}^2} \notag \\
    &= \rho \brace{\wt{\fra}_0 \expt{\wt{\sfPhi}_0 \sfLambda_d} - \frb_1} = 0 , \label{eqn:qin2_fq} \\
    \frac{1}{d} (q_{\inn,1}^0)^\top f_{q,\inn}^0 &\to \expt{\sfQ_{\inn,1}^0 \brace{\sfPhi_0 (\sfQ_{\inn,1}^0 + \wt{\frb}_0 \sfQ_{\inn,2}) + \wt{\sfPhi}_0 (\sfQ_\out^0 + \wt{\fra}_0 \sfLambda_d \sfQ_{\inn,2}) - \frb_1 \sfQ_{\inn,2}}} \notag \\
    &= \expt{\sfPhi_0} \expt{(\sfQ_{\inn,1}^0)^2} = 0 . \label{eqn:qin1_fq} 
\end{align}
 The second claim in \Cref{eqn:q_fq_0} can be similarly proved using the joint convergence result \Cref{eqn:t0_fq} and zero mean of $ \sfPsi_0 $. 
Therefore, by \Cref{prop:wto3,prop:wto4} we have that $ \lr{p_\inn^1 \mid \brace{Q q_\inn^0 = f_{p,\inn}^0}} $ and $ \lr{p_\out^1 \mid \brace{O q_\out^0 = f_{p,\out}^0}} $ have Gaussian limits due to the second terms on the right of \Cref{eqn:p1}, denoted by $ \sfU_{p,\inn}^1 \in \bbR $ and $ (\sfU_{p,\out}^1)^\top = \matrix{\sfU_{p,\out,1}^1 & \sfU_{p,\out,2}} \in \bbR^{2} $, respectively, independent of everything else. 
In formulas, 
\begin{align}
\begin{split}
    & \matrix{f_{p,\inn,2} & f_{p,\inn,1}^0 & w_{p,\inn} & w_{q,\inn} & q_{\inn,2} & q_{\inn,1}^0 & f_{q,\inn}^0 & p_\inn^1} \\
    &\qquad \wto \matrix{\sfF_{p,\inn,2} & \sfF_{p,\inn,1}^0 & \sfW_{p,\inn}^\top & \sfW_{q,\inn}^\top & \sfQ_{\inn,2} & \sfQ_{\inn,1}^0 & \sfF_{q,\inn}^0 & \sfU_{p,\inn}^1} , \\
    & \matrix{f_{p,\out}^0 & w_{p,\out} & w_{q,\out} & q_\out^0 & f_{q,\out,2} & f_{q,\out,1}^0 & p_{\out,1}^1 & p_{\out,2}} \\
    &\qquad \wto \matrix{\sfF_{p,\out}^0 & \sfW_{p,\out}^\top & \sfW_{q,\out}^\top & \sfQ_\out^0 & \sfF_{q,\out,2} & \sfF_{q,\out,1}^0 & \sfU_{p,\out,1}^1 & \sfU_{p,\out,2}} . 
\end{split}
\label{eqn:t0_p} 
\end{align}
We then argue that $ \sfU_{p,\inn}^1, \sfU_{p,\out,1}^1 , \sfU_{p,\out,2} $ are precisely $ \sfP_\inn^1, \sfP_{\out,1}^1, \sfP_{\out,2} $ by verifying that the joint distributions between the former random variables with others on the right of \Cref{eqn:t0_p} are consistent with the definition \Cref{eqn:Pin,eqn:Pout}. 
By the first claim in \Cref{eqn:q_fq_0}, $ \sfU_{p,\inn}^1 $ is Gaussian with variance: 
\begin{align}
    \lim_{d\to\infty} \frac{1}{d} \normtwo{p_\inn^1}^2 &= \lim_{d\to\infty} \frac{1}{d} \normtwo{f_{q,\inn}^0}^2 = \Sigma_{p,\inn}^1 , \notag 
\end{align}
independent of all other random variables on the right of \Cref{eqn:t0_p}. 
Similarly, by the second claim in \Cref{eqn:q_fq_0} and the convergence result \Cref{eqn:t0_fq}, $ (\sfU_{p,\out,1}^1, \sfU_{p,\out,2}) $ are jointly Gaussian with covariance: 
\begin{align}
    \lim_{n\to\infty} \frac{1}{n} (p_\out^1)^\top p_\out^1 &= \lim_{n\to\infty} \frac{1}{n} (f_{q,\out}^0)^\top f_{q,\out}^0 = \matrix{\Sigma_{p,\out}^1 & 0 \\ 0 & \sigma^2} , \notag 
\end{align}
independent of all other random variables on the right of \Cref{eqn:t0_p}. 
The bottom-right entry of the covariance matrix above follows since $ \expt{(\sfF_{q,\out,2})^2} = \expt{\sfLambda_n^2} \expt{\sfQ_{\inn,2}^2} = \ol{\kappa}_2 \rho = \sigma^2 $. 
The off-diagonal entries follow from 
\begin{align}
    \expt{\sfF_{q,\out,2} \sfF_{q,\out,1}^0}
    &= \expt{\sfPsi_0 \sfLambda_n \sfQ_{\inn,2} (\sfQ_\out^0 + \wt{\fra}_0 \sfLambda_n \sfQ_{\inn,2})} 
    + \expt{\wt{\sfPsi}_0 \sfLambda_n \sfQ_{\inn,2} (\sfQ_{\inn,1}^0 + \wt{\frb}_0 \sfQ_{\inn,2})}
    - \fra_1 \expt{\sfLambda_n^2 \sfQ_{\inn,2}^2}\notag \\
    &= \wt{\fra}_0 \expt{\sfPsi_0 \sfLambda_n^2} \rho
    + \wt{\frb}_0 \expt{\wt{\sfPsi}_0 \sfLambda_n} \rho
    - \fra_1 \ol{\kappa}_2 \rho
    = 0 , \label{eqn:Fqout2_Fqout1_0} 
\end{align}
and the choices of $ \wt{\fra}_0, \wt{\frb}_0 $ in \Cref{eqn:frab_tilde0} and of $ \fra_1 $ in \Cref{eqn:frab}. 
This proves \Cref{eqn:hyp4_out} in \Cref{itm:hyp4} for $t=0$. 
Also, \Cref{eqn:Fqout2_Fqout1_0} above proves the second statement in \Cref{itm:hyp1} for $t=0$. 

Finally, the fourth update \Cref{eqn:OQ4} generates 
\begin{align}
&&
    f_{p,\inn}^1 &= \matrix{f_1(r^1; \beta_*, \Theta) - \wt{\frb}_1 \beta_* & \beta_*} &
    & \textnormal{and} &
    f_{p,\out}^1 &= g_1(p^1; z, \eps, \Xi) - \wt{\fra}_1 z, & 
& \notag 
\end{align}
whose joint convergence with all preceding iterates is guaranteed by \Cref{prop:wto2,prop:wto4}. 
In formulas, we have the desired result: 
\begin{align}
\begin{split}
    & \matrix{q_{\inn,2} & q_{\inn,1}^0 & f_{q,\inn}^0 & p_\inn^1 & f_{p,\inn,2} & f_{p,\inn,1}^0 & f_{p,\inn,1}^1 & w_{p,\inn} & w_{q,\inn}} \\
    &\qquad \wto \matrix{\sfQ_{\inn,2} & \sfQ_{\inn,1}^0 & \sfF_{q,\inn}^0 & \sfP_\inn^1 & \sfF_{p,\inn,2} & \sfF_{p,\inn,1}^0 & \sfF_{p,\inn,1}^1 & \sfW_{p,\inn}^\top & \sfW_{q,\inn}^\top} \\
    & \matrix{q_\out^0 & f_{q,\out,2} & f_{q,\out,1}^0 & p_{\out,2} & p_{\out,1}^1 & f_{p,\out}^0 & f_{p,\out}^1 & w_{p,\out} & w_{q,\out}} \\
    &\qquad \wto \matrix{\sfQ_\out^0 & \sfF_{q,\out,2} & \sfF_{q,\out,1}^0 & \sfP_{\out,2} & \sfP_{\out,1}^1 & \sfF_{p,\out}^0 & \sfF_{p,\out}^1 & \sfW_{p,\out}^\top & \sfW_{q,\out}^\top} . 
\end{split}
\label{eqn:goal0}
\end{align}

\paragraph{\textsc{Induction step.}}

Now assuming the validity of \Cref{itm:hyp1,itm:hyp2,itm:hyp3,itm:hyp4} for time up to $t\ge0$, we verify their validity for time $t+1$. 

The statement $ \expt{ \sfF_{p,\inn,1}^{t+1} \sfF_{p,\inn,2} } = 0 $ in \Cref{itm:hyp1} (for time $t+1$) 
holds since 
\begin{align}
    \expt{ \sfF_{p,\inn,1}^{t+1} \sfF_{p,\inn,2} } &= \expt{ f_{t+1}(\sfP_\inn^{t+1} + \frb_{t+1} \sfB_*; \sfB_*, \sfTheta) \sfB_* } - \wt{\frb}_{t+1} \rho = 0, \label{eqn:Fpin1_Fpin2} 
\end{align}
by the definition of $\sfF_{p,\inn,1}^{t+1}$ in \Cref{eqn:Fpin1} and the definition of $ \wt{\frb}_{t+1} $ in \Cref{eqn:frab_tilde}. 

After the updates for time $t$, the covariance matrices $ \Sigma_{p,\inn}^{t+1}, \Sigma_{p,\out}^{t+1}, \Sigma_{q,\inn}^{t+1}, \Sigma_{q,\out}^{t+1} $ have been computed. 
We write the first two in block form: 
\begin{align}
&& 
    \Sigma_{q,\inn}^{t+1} &= \matrix{
        \Sigma_{q,\inn}^{t} & \sigma_{q,\inn}^{t} \\
        (\sigma_{q,\inn}^{t})^\top & s_{q,\inn}^{t+1}
    } \in \bbR^{(t+2)\times(t+2)} , & 
    \Sigma_{q,\out}^{t+1} &= \matrix{
        \Sigma_{q,\out}^{t} & \sigma_{q,\out}^{t} \\
        (\sigma_{q,\out}^{t})^\top & s_{q,\out}^{t+1}
    } \in \bbR^{(t+2)\times(t+2)} , & 
& \label{eqn:Sigma_q>0} 
\end{align}
where we have introduced the notation $ \sigma_{q,\inn}^t, \sigma_{q,\out}^t \in\bbR^{t+1}$ and $s_{q,\inn}^{t+1}, s_{q,\out}^{t+1}\in\bbR $. 
Notationally, we index the elements of $ \sigma_{q,\inn}^t, \sigma_{q,\out}^t $ using the set $ \brace{0,1,\cdots,t} $, which is consistent with the indexing of $ \Sigma_{q,\inn}^{t+1}, \Sigma_{q,\out}^{t+1} $. 
We verify that $ \Sigma_{q,\inn}^{t+1} $ is strictly positive definite. 
Recalling the definition \Cref{eqn:cov_q}, given $ \Sigma_{q,\inn}^t $ being non-singular due to the induction hypothesis in \Cref{itm:hyp3}, $ \Sigma_{q,\inn}^{t+1} $ is singular if and only if there exist coefficients $ \alpha_0, \cdots, \alpha_t \in \bbR $, such that almost surely, 
\begin{align}
    \sfF_{p,\inn,1}^{t+1} &= \sum_{s = 0}^t \alpha_s \sfF_{p,\inn,1}^s . \notag
\end{align}
However, this is excluded by \Cref{asmp:nondegenerate_abs}. 
So $ \Sigma_{q,\inn}^{t+1} $ is invertible. 
A similar argument applies to $ \Sigma_{q,\out}^{t+1} $. 
This proves part of \Cref{itm:hyp3} for time $t+1$. 

\paragraph{Induction step for $ q_\inn^{t+1} $.}
We collect some past iterates into matrices $U$ and $V$: 
\begin{align}
    U &= \matrix{ q_{\inn,2} & q_{\inn,1}^0 & \cdots & q_{\inn,1}^t & f_{q,\inn}^0 & \cdots & f_{q,\inn}^t } \in \bbR^{d \times (1+2(t+1))} , \label{eqn:U1} \\
    V &= \matrix{ f_{p,\inn,2} & f_{p,\inn,1}^0 & \cdots f_{p,\inn,1}^t & p_\inn^1 & \cdots & p_\inn^{t+1} } \in \bbR^{d \times (1+2(t+1))} . \label{eqn:V1}
\end{align}
According to the update rules of the abstract VAMP iteration \Cref{eqn:abs_VAMP}, the distribution of $ q_\inn^{t+1} $ given $ w_{p,\inn}, w_{p,\out}, w_{q,\inn}, w_{q,\out} $, $ (q_\inn^s)_{0\le s\le t} $, $ (q_\out^s)_{0\le s\le t} $, $ (f_{q,\inn}^s)_{0\le s\le t} $, $ (f_{q,\out}^s)_{0\le s\le t} $, $ (p_\inn^{s+1})_{0\le s\le t} $, $ (p_\out^{s+1})_{0\le s\le t} $, $ (f_{p,\inn}^s)_{0\le s\le t+1} $, $ (f_{p,\out}^s)_{0\le s\le t+1} $ is the same as the law of $ q_\inn^{t+1} $ conditioned on the event $ U = Q^\top V $. 
Here we have used the orthogonality of $Q$ which implies that the condition $ p_\inn^{s+1} = Q f_{q,\inn}^s $ given by the first update in \Cref{eqn:OQ3} is identical to $ Q^\top p_\inn^{s+1} = f_{q,\inn}^s $ for all $ 0\le s \le t $. 
For all sufficiently large $d$, the latter conditional distribution equals
\begin{align}
    \left. q_{\inn}^{t+1} \mid \brace{ U = Q^\top V } \right. 
    &= \left. Q^\top f_{p,\inn}^{t+1} \mid \brace{ U = Q^\top V } \right. 
    \eqqlaw U (U^\top U)^{-1} V^\top f_{p,\inn}^{t+1} + \Pi_{U^\perp} \wt{Q}^\top \Pi_{V^\perp}^\top f_{p,\inn}^{t+1} \in \bbR^{d\times2} , \label{eqn:qin_TODO} 
\end{align}
where the last step follows from \Cref{prop:haar_cond} (which contains the definitions of the projectors $ \Pi_{U^\perp}, \Pi_{V^\perp} $) and $ \wt{Q} \sim \haar(\bbO(d - (1 + 2(t+1)))) $ is independent of everything else. 
We claim that $ U^\top U $ in the above display is invertible for all sufficiently large $d$ and will justify this claim in \Cref{eqn:UU_result}. 
(We remark that across different passages of the induction proof, the notation $ U,V,\wt{O},\wt{Q} $ will be overloaded.)
To analyze the first term on the right in \Cref{eqn:qin_TODO}, we  compute the limits of $ d^{-1} U^\top U $ and $ d^{-1} V^\top f_{p,\inn,1}^{t+1} $. 

Let us start with the computation of the limit of $ d^{-1} U^\top U $. 
Note that by the induction hypothesis \Cref{eqn:concl1}, for any $ 0\le r,s\le t $, 
\begin{subequations}
\label{eqn:diag}     
\begin{align}
&& 
    \frac{1}{d} \inprod{q_{\inn,2}}{q_{\inn,2}} &= \frac{1}{d} \inprod{\beta_*}{\beta_*} \to \expt{\sfB_*^2} = \rho , & 
    \frac{1}{d} \inprod{q_{\inn,2}}{q_{\inn,1}^s} &\to \expt{ \sfQ_{\inn,2} \sfQ_{\inn,1}^s } = 0 , & 
& \label{eqn:diag1} \\
&&
    \frac{1}{d} \inprod{q_{\inn,1}^r}{q_{\inn,1}^s} &\to \expt{ \sfQ_{\inn,1}^r \sfQ_{\inn,1}^s } = (\Sigma_{q,\inn}^t)_{r,s} , & 
    \frac{1}{d} \inprod{f_{q,\inn}^{r}}{f_{q,\inn}^{s}} &\to \expt{ \sfF_{q,\inn}^r \sfF_{q,\inn}^s } = (\Sigma_{p,\inn}^{t+1})_{r,s} . & 
& 
\end{align}
\end{subequations}
We will then show
\begin{align}
&& 
    \frac{1}{d} \inprod{q_{\inn,1}^r}{f_{q,\inn}^s} &\to 0 , &
    \frac{1}{d} \inprod{q_{\inn,2}}{f_{q,\inn}^s} &\to 0 . &
& \label{eqn:entry12} 
\end{align}
The first term is: 
\begin{align}
    \frac{1}{d} \inprod{q_{\inn,1}^r}{f_{q,\inn}^s}
    &= \frac{1}{d} \inprod{q_{\inn,1}^r}{\diag(\phi_s) (q_{\inn,1}^s + \wt{\frb}_s q_{\inn,2})} \label{eqn:rhs1} \\
    &\quad + \frac{1}{d} \inprod{q_{\inn,1}^r}{\diag_{n\times d}(\wt{\phi}_s)^\top (q_\out^s + \wt{\fra}_s \Lambda q_{\inn,2})} \label{eqn:rhs2} \\
    &\quad - \frac{\frb_{s+1}}{d} \inprod{q_{\inn,1}^r}{q_{\inn,2}} , \label{eqn:rhs3}
\end{align}
where, by the induction hypotheses \Cref{eqn:concl1,eqn:concl2} in \Cref{itm:hyp2}, the right-hand sides have the following limits: 
\begin{align}
    \Cref{eqn:rhs1} &\to \expt{ \sfQ_{\inn,1}^r (\sfQ_{\inn,1}^s + \wt{\frb}_s \sfQ_{\inn,2}) } \expt{ \sfPhi_s } = 0 , \label{eqn:lim1} \\
    \Cref{eqn:rhs2} 
    &= \frac{1}{d} \inprod{q_{\inn,1}^r}{\diag_{n\times d}(\wt{\phi}_s)^\top q_\out^s}
    + \frac{\wt{\fra}_s}{d} \inprod{q_{\inn,1}^r}{\diag_{n\times d}(\wt{\phi}_s)^\top \Lambda q_{\inn,2}} \notag \\
    &\to 0 + \wt{\fra}_s \expt{ \sfQ_{\inn,1}^r \sfQ_{\inn,2} } \expt{ \wt{\sfPhi}_s \sfLambda_d } = 0 , \label{eqn:lim2} \\
    \Cref{eqn:rhs3} &= - \frac{\frb_{s+1}}{d} \inprod{Q^\top f_{p,\inn}^1}{Q^\top \beta_*}
    = - \frac{\frb_{s+1}}{d} \inprod{f_r(p_\inn^r + \frb_r \beta_*; \beta_*, \Theta) - \wt{\frb}_r \beta_*}{\beta_*} 
    \to 0 . \label{eqn:use_btilde}
\end{align}
The equality in \Cref{eqn:lim1} follows from independence between $ (\sfQ_\inn^s)_{0\le s\le t}, \sfPhi_s $ and the  zero mean of $ \sfPhi_s $ in \Cref{asmp:4matrices}. 
The equality in \Cref{eqn:lim2} follows from independence between $ (\sfQ_{\inn,1}^s)_{0\le s\le t} $ and $ \sfQ_{\inn,2} , (\sfQ_{\out,1}^s)_{0\le s\le t} $. 
The equality in \Cref{eqn:use_btilde} follows from the definition of $ \wt{\frb}_r $ in \Cref{eqn:frab_tilde}. 
Therefore, the first statement in \Cref{eqn:entry12} holds. 

Turning to the second statement in \Cref{eqn:entry12}, we have 
\begin{align}
    \frac{1}{d} \inprod{q_{\inn,2}}{f_{q,\inn}^s}
    &= \frac{1}{d} \inprod{q_{\inn,2}}{\diag(\phi_s) (q_{\inn,1}^s + \wt{\frb} q_{\inn,2})}
    + \frac{1}{d} \inprod{q_{\inn,2}}{\diag_{n\times d}(\wt{\phi}_s)^\top (q_\out^s + \wt{\fra}_s \Lambda q_{\inn,2}) - \frb_{s+1} q_{\inn,2}} . \notag 
\end{align}
The first term has limit $0$ due to similar reasoning leading to \Cref{eqn:lim1}. 
For the second term, by manipulations similar to those leading to \Cref{eqn:lim2}, we have
\begin{align}
    & \frac{1}{d} \inprod{q_{\inn,2}}{\diag_{n\times d}(\wt{\phi}_s)^\top (q_\out^s + \wt{\fra}_s \Lambda q_{\inn,2})}
    - \frac{\frb_{s+1}}{d} \inprod{q_{\inn,2}}{q_{\inn,2}} \notag \\
    &= \frac{1}{d} \inprod{q_{\inn,2}}{\diag_{n\times d}(\wt{\phi}_s)^\top q_\out^s}
    + \frac{\wt{\fra}_s}{d} \inprod{q_{\inn,2}}{\diag_{n\times d}(\wt{\phi}_s)^\top \Lambda q_{\inn,2}}
    - \frac{\frb_{s+1}}{d} \inprod{q_{\inn,2}}{q_{\inn,2}} \notag \\
    &\to 0 + \wt{\fra}_s \expt{ \wt{\sfPhi}_s \sfLambda_d } \expt{ \sfQ_{\inn,2}^2 }
    - \frb_{s+1} \expt{ \sfQ_{\inn,2}^2 } \notag \\
    &= \wt{\fra}_s \expt{ \wt{\sfPhi}_s \sfLambda_d } \rho
    - \frb_{s+1} \rho 
    = 0 , \notag 
\end{align}
where the last step follows from the definition of $ \frb_{s+1} $ in \Cref{eqn:frab}. 
Therefore the second equality in \Cref{eqn:entry12} also holds. 
Taking \Cref{eqn:entry12,eqn:diag} collectively gives 
\begin{align}
    \frac{1}{d} U^\top U &\to \matrix{
        \rho & 0_{1\times(t+1)} & 0_{1\times(t+1)} \\
        0_{(t+1)\times1} & \Sigma_{q,\inn}^t & 0_{(t+1) \times (t+1)} \\
        0_{(t+1)\times1} & 0_{(t+1) \times (t+1)} & \Sigma_{p,\inn}^{t+1}
    } \in \bbR^{(1+2(t+1)) \times (1+2(t+1))} . \label{eqn:UU_result} 
\end{align}
Note that the matrix on the right is strictly positive definite by $ \rho > 0 $ in \Cref{asmp:signal} and the induction hypothesis in \Cref{itm:hyp3}. 

Next, consider $ d^{-1} V^\top f_{p,\inn}^{t+1} $.
By the induction hypothesis \Cref{eqn:concl1} in \Cref{itm:hyp2}, for any $ 0\le s \le t $, 
\begin{subequations}
\label{eqn:V_fpin_1} 
\begin{align}
    \frac{1}{d} \inprod{f_{p,\inn,2}}{f_{p,\inn,1}^{t+1}} &\to \expt{ \sfF_{p,\inn,2} \sfF_{p,\inn,1}^{t+1} } = 0 , \label{eqn:V_fpin_11} \\
    \frac{1}{d} \inprod{f_{p,\inn,1}^s}{f_{p,\inn,1}^{t+1}} &\to \expt{ \sfF_{p,\inn,1}^s \sfF_{p,\inn,1}^{t+1} } = (\Sigma_{q,\inn}^{t+1})_{s,t+1} = (\sigma_{q,\inn}^t)_s ,  \\
    \frac{1}{d} \inprod{f_{p,\inn,2}}{f_{p,\inn,2}} &\to \expt{ \sfF_{p,\inn,2}^2 } = \expt{\sfB_*^2} = \rho , \\
    \frac{1}{d} \inprod{f_{p,\inn,1}^s}{f_{p,\inn,2}} &\to \expt{ \sfF_{p,\inn,1}^s \sfF_{p,\inn,2} } = 0 . \label{eqn:V_fpin_14}
\end{align}
\end{subequations}
\Cref{eqn:V_fpin_11} follows from \Cref{eqn:Fpin1_Fpin2}, and 
\Cref{eqn:V_fpin_14} follows from the induction hypothesis in \Cref{itm:hyp1}. 
Furthermore, we study the limit of 
\begin{align}
    \frac{1}{d} (p_\inn^{s+1})^\top f_{p,\inn}^{t+1}
    &= \matrix{
        \frac{1}{d} \inprod{p_\inn^{s+1}}{f_{p,\inn,1}^{t+1}} 
        & \frac{1}{d} \inprod{p_\inn^{s+1}}{f_{p,\inn,2}}
    } . \label{eqn:entry12_fpin} 
\end{align}
The second entry of \Cref{eqn:entry12_fpin} can be easily seen vanishing: 
\begin{align}
    \frac{1}{d} \inprod{p_\inn^{s+1}}{\beta_*}
    &\to \expt{ \sfP_\inn^{s+1} } \expt{ \sfB_* } = 0 , \label{eqn:entry2_1} 
\end{align}
since $ \sfP_\inn^{s+1} $ and $ \sfB_* $ are independent. 
As for the first entry of \Cref{eqn:entry12_fpin}, we have
\begin{align}
    \frac{1}{d} \inprod{p_\inn^{s+1}}{f_{p,\inn,1}^{t+1}} 
    &= \frac{1}{d} \inprod{p_\inn^{s+1}}{f_{t+1}(p_\inn^{t+1} + \frb_{t+1} \beta_*; \beta_*, \Theta)} 
    - \frac{\wt{\frb}_{t+1}}{d} \inprod{p_\inn^{s+1}}{\beta_*} . \notag 
\end{align}
We have shown in \Cref{eqn:entry2_1} that the second term on the right vanishes. 
By \Cref{asmp:regularity} and \Cref{prop:wto2}, the first term converges to 
\begin{align}
    \frac{1}{d} \inprod{p_\inn^{s+1}}{f_{t+1}(p_\inn^{t+1} + \frb_{t+1} \beta_*; \beta_*, \Theta)}
    &\to \expt{ \sfP_\inn^{s+1} f_{t+1}(\sfP_\inn^{t+1} + \frb_{t+1} \sfB_*; \sfB_*, \sfTheta) } . \label{eqn:entry2_2} 
\end{align}
Since $ (\sfP_\inn^{s+1}, \sfP_\inn^{t+1}) $ are jointly Gaussian and independent of $ \sfB_* $, the above expression is zero by Stein's lemma \cite[Theorem 2.1]{shrinkage_estimation_book} and divergence-freeness of $ f_{t+1} $ (see \Cref{asmp:div_free}). 
Putting \Cref{eqn:entry2_1,eqn:entry2_2} together, we see that 
\begin{align}
    \frac{1}{d} (p_\inn^{s+1})^\top f_{p,\inn}^{t+1}
    &\to \matrix{0 & 0} . \label{eqn:V_fpin_2} 
\end{align}
Combining \Cref{eqn:V_fpin_1,eqn:V_fpin_2}, we arrive at
\begin{align}
    \frac{1}{d} V^\top f_{p,\inn}^{t+1} &\to \matrix{
        0 & \rho \\
        \sigma_{q,\inn}^t & 0_{t+1} \\
        0_{t+1} & 0_{t+1}
    } \in \bbR^{(1+2(t+1))\times2} . \label{eqn:V_fpin_result}
\end{align}

We are now in a position to compute the right-hand side of \Cref{eqn:qin_TODO}. 
Using \Cref{eqn:UU_result,eqn:V_fpin_result}, we have
\begin{align}
    (U^\top U)^{-1} V^\top f_{p,\inn}^{t+1}
    &= (d^{-1} U^\top U)^{-1} \cdot d^{-1} V^\top f_{p,\inn}^{t+1}
    \to \matrix{
        \rho^{-1} & 0_{1\times(t+1)} & 0_{1\times(t+1)} \\
        0_{(t+1)\times1} & (\Sigma_{q,\inn}^t)^{-1} & 0_{(t+1) \times (t+1)} \\
        0_{(t+1)\times1} & 0_{(t+1) \times (t+1)} & (\Sigma_{p,\inn}^{t+1})^{-1}
    } \matrix{
        0 & \rho \\
        \sigma_{q,\inn}^t & 0_{t+1} \\
        0_{t+1} & 0_{t+1}
    } \notag \\
    &= \matrix{
        0 & 1 \\
        (\Sigma_{q,\inn}^t)^{-1} \sigma_{q,\inn}^t & 0_{t+1} \\
        0_{t+1} & 0_{t+1}
    } \in \bbR^{(1+2(t+1))\times2} . \notag 
\end{align}
By \Cref{prop:wto4}, the second term in \Cref{eqn:qin_TODO} has a Gaussian limit in Wasserstein distance, denoted by $ \matrix{ \sfU_{q,\inn}^{t+1} & 0 } \in \bbR^{1\times2} $, independent of $ \sfA $, $ (\sfQ_\inn^s)_{0\le s\le t} $, $ (\sfQ_\out^s)_{0\le s\le t} $, $ (\sfF_{q,\inn}^s)_{0\le s\le t} $, $ (\sfF_{q,\out}^s)_{0\le s\le t} $, $ (\sfP_\inn^{s+1})_{0\le s\le t} $, $ (\sfP_{\out}^{s+1})_{0\le s\le t} $, $ (\sfF_{p,\inn}^s)_{0\le s\le t+1} $, $ (\sfF_{p,\out}^s)_{0\le s\le t+1} $.  
Note that the second entry is zero.
This follows since every column of $ \Pi_{V^\perp} $ is orthogonal to the first column of $V$ which equals $ f_{p,\inn,2} $, implying that $ \Pi_{V^\perp}^\top f_{p,\inn,2} $ vanishes. 
Using this together with the above display, recalling that  
\begin{align}
    U \wto \matrix{ \sfQ_{\inn,2} & \sfQ_{\inn,1}^0 & \cdots & \sfQ_{\inn,1}^t & \sfF_{q,\inn}^0 & \cdots & \sfF_{q,\inn}^t } \notag 
\end{align}
by the induction hypothesis \Cref{eqn:concl1} in \Cref{itm:hyp2}, and applying \Cref{prop:wto3}, we conclude
\begin{subequations}
\label{eqn:Qint+1}
\begin{align}
     \left. q_{\inn,1}^{t+1} \mid \brace{U = Q^\top V} \right. &\wto \matrix{ \sfQ_{\inn,1}^0 & \cdots & \sfQ_{\inn,1}^t } (\Sigma_{q,\inn}^t)^{-1} \sigma_{q,\inn}^t 
     + \sfU_{q,\inn}^{t+1} 
     = \sum_{s = 0}^t \sfQ_{\inn,1}^s \paren{ (\Sigma_{q,\inn}^t)^{-1} \sigma_{q,\inn}^t }_s + \sfU_{q,\inn}^{t+1} , \label{eqn:Qin1^t+1} \\
     \left. q_{\inn,2} \mid \brace{U = Q^\top V} \right. &\wto \sfQ_{\inn,2} . 
\end{align} 
\end{subequations}

We claim that  the right side of the two equations above is precisely $ \sfQ_\inn^{t+1} = \matrix{ \sfQ_{\inn,1}^{t+1} & \sfQ_{\inn,2} }^\top $. 
To justify this, note that $ (\sfQ_\inn^s)_{0\le s\le t+1} $ is jointly Gaussian and is independent of $ \sfA, (\sfP_\inn^{s+1})_{0\le s\le t}, (\sfP_\out^{s+1})_{0\le s\le t} $. 
So by the definition \Cref{eqn:Qin}, it remains to verify that the right-hand side of \Cref{eqn:Qin1^t+1}, denoted by $ \wt{\sfQ}_\inn^{t+1} $, has the desired covariance structure with $ (\sfQ_\inn^s)_{0\le s\le t} $. 
For any $ 0\le r\le t $, we have 
\begin{subequations}
\label{eqn:check_cov}
\begin{align}
    \expt{ \wt{\sfQ}_{\inn,1}^{t+1} \sfQ_{\inn,1}^r } &= \sum_{s = 0}^t \expt{ \sfQ_{\inn,1}^s \sfQ_{\inn,1}^r } \paren{ (\Sigma_{q,\inn}^t)^{-1} \sigma_{q,\inn}^t }_s + \expt{ \sfU_{q,\inn}^{t+1} \sfQ_{\inn,1}^r } \notag \\
    &= \sum_{s = 0}^t (\Sigma_{q,\inn}^t)_{r,s} \paren{ (\Sigma_{q,\inn}^t)^{-1} \sigma_{q,\inn}^t }_s + 0 
    = \paren{ \Sigma_{q,\inn}^t (\Sigma_{q,\inn}^t)^{-1} \sigma_{q,\inn}^t }_r 
    = (\sigma_{q,\inn}^t)_r . 
\end{align}
Note also that 
\begin{align}
    \expt{ (\wt{\sfQ}_{\inn,1}^{t+1})^2 } 
    &= \lim_{d\to\infty} \frac{1}{d} \normtwo{ q_{\inn,1}^{t+1} }^2 
    = \lim_{d\to\infty} \frac{1}{d} \normtwo{ f_{p,\inn,1}^{t+1} }^2
    = \expt{ (\sfF_{p,\inn,1}^{t+1})^2 }
    = s_{q,\inn}^{t+1} . \label{eqn:sqint+1}
\end{align}
Moreover, $ \wt{\sfQ}_{\inn,1}^{t+1} $ is independent of $ \sfQ_{\inn,2} $ since $ (\sfQ_{\inn,1}^s)_{0\le s\le t} $ and $ \sfU_{q,\inn}^{t+1} $ are all independent of $ \sfQ_{\inn,2} $. 
\end{subequations}
From this and \Cref{eqn:check_cov}, we see that $ \matrix{\wt{\sfQ}_{\inn,1}^{t+1} & \sfQ_{\inn,2}}^\top $ has the desired covariance structure with $ (\sfQ_{\inn}^s)_{0\le s\le t} $ which we identify as $\sfQ_\inn^{t+1}$, and therefore
\begin{align}
     & \matrix{
        q_{\inn,2} & \bq_{\inn,1}^{t+1} & 
        \bf_{q,\inn}^{t} & 
        \bp_\inn^{t+1} & 
        f_{p,\inn,2} & \bf_{p,\inn,1}^{t+1} & 
        w_{p,\inn} & w_{q,\inn}
    } \notag \\
    &\qquad\wto \matrix{
        \sfQ_{\inn,2} & \bm\sfQ_{\inn,1}^{t+1} &
        \bm\sfF_{q,\inn}^{t} &
        \bm\sfP_\inn^{t+1} &
        \sfF_{p,\inn,2} & \bm\sfF_{p,\inn,1}^{t+1} &
        \sfW_{p,\inn}^\top & \sfW_{q,\inn}^\top
    } . \label{eqn:qin_concl}
\end{align} 

We check that the variance of $ \sfU_{q,\inn}^{t+1} $ is strictly positive. 
Since $ (\sfQ_{\inn,1}^s)_{0\le s\le t} \sim \cN(0_{t+1}, \Sigma_{q,\inn}^t) $, we have 
\begin{align}
    \sum_{s = 0}^t \sfQ_{\inn,1}^s \paren{ (\Sigma_{q,\inn}^t)^{-1} \sigma_{q,\inn}^t }_s &\sim \cN\paren{ 0, \paren{(\Sigma_{q,\inn}^t)^{-1} \sigma_{q,\inn}^t}^\top \cdot \Sigma_{q,\inn}^t \cdot \paren{ (\Sigma_{q,\inn}^t)^{-1} \sigma_{q,\inn}^t } } \notag \\
    &= \cN\paren{0, (\sigma_{q,\inn}^t)^\top (\Sigma_{q,\inn}^t)^{-1} \sigma_{q,\inn}^t} . \label{eqn:var} 
\end{align}
According to \Cref{eqn:Qint+1}, the variance of $ \sfU_{q,\inn}^{t+1} $ is given by 
\begin{align}
    \expt{ (\sfU_{q,\inn}^{t+1})^2 } &= s_{q,\inn}^{t+1} - (\sigma_{q,\inn}^t)^\top (\Sigma_{q,\inn}^t)^{-1} \sigma_{q,\inn}^t , \notag 
\end{align}
where the two terms on the right come from \Cref{eqn:sqint+1,eqn:var}. 
By the Schur complement formula, 
\begin{align}
    \expt{ (\sfU_{q,\inn}^{t+1})^2 } &= \frac{1}{\brack{(\Sigma_{q,\inn}^{t+1})^{-1}}_{t+1,t+1}} = \paren{ e_{t+1}^\top (\Sigma_{q,\inn}^{t+1})^{-1} e_{t+1} }^{-1} > 0 , \label{eqn:Uqin^t+1>0} 
\end{align}
where the strict inequality follows from strict positive definiteness of $ \Sigma_{q,\inn}^{t+1} $ shown around \Cref{eqn:Sigma_q>0}.
This proves the first statement of \Cref{eqn:hyp4_in} in \Cref{itm:hyp4}. 

\paragraph{Induction step for $ q_\out^{t+1} $.}
Redefine the matrices $ U,V $ as 
\begin{align}
    U &= \matrix{ q_{\out}^0 & \cdots & q_{\out}^t & f_{q,\out,2} & f_{q,\out,1}^0 & \cdots & f_{q,\out,1}^t } \in \bbR^{n\times(1+2(t+1))} , \label{eqn:U2} \\
    V &= \matrix{ f_{p,\out}^0 & \cdots & f_{p,\out}^t & p_{\out,2} & p_{\out,1}^1 & \cdots & p_{\out,1}^{t+1} } \in \bbR^{n\times(1+2(t+1))} . \label{eqn:V2}
\end{align}
The law of $ q_\out^{t+1} $ given all previous iterates (including the initializers and side information) is the same as the law of $ q_\inn^{t+1} $ given the condition $ U = O^\top V $. 
By \Cref{prop:haar_cond}, the latter conditional law is 
\begin{align}
    \left. q_\out^{t+1} \mid \brace{ U = O^\top V } \right. 
    = \left. O^\top f_\out^{t+1} \mid \brace{ U = O^\top V } \right. 
    &\eqqlaw U (U^\top U)^{-1} V^\top f_{p,\out}^{t+1} + \Pi_{U^\top} \wt{O}^\top \Pi_{V^\perp}^\top f_{p,\out}^{t+1} , \label{eqn:qoutt+1_cond} 
\end{align}
where $ \wt{O} \sim \haar(\bbO(n - (1+2(t+1)))) $ is independent of everything else. 
To analyze the first term on the right, we compute the limits of $ n^{-1} U^\top U $ and $ n^{-1} V^\top f_{p,\out}^{t+1} $. 
In particular, \Cref{eqn:UU_qoutt+1} below shows that $ U^\top U $ on the right of \Cref{eqn:qoutt+1_cond} is invertible for every sufficiently large $n$. 

Starting with $ n^{-1} U^\top U $, for any $ 0\le r,s\le t $, by induction hypotheses \Cref{itm:hyp1,itm:hyp2}, we have 
\begin{subequations}
\label{eqn:qoutt+1_cov}
\begin{align}
&& 
    \frac{1}{n} \inprod{q_{\out}^r}{q_{\out}^s} &\to (\Sigma_{q,\out}^t)_{r,s} , & 
    \frac{1}{n} \inprod{f_{q,\out,2}}{f_{q,\out,2}} &= \frac{1}{d} \inprod{z}{z} \to \expt{\sfZ^2} = \sigma^2 , & 
& \label{eqn:qoutt+1_cov1} \\
&& 
    \frac{1}{n} \inprod{f_{q,\out,2}}{f_{q,\out,1}^s} &\to 0 , & 
    \frac{1}{n} \inprod{f_{q,\out,1}^r}{f_{q,\out,1}^s} &\to (\Sigma_{p,\out}^{t+1})_{r+1,s+1} . & 
&  
\end{align}
\end{subequations}
Note that the second equation in \Cref{eqn:qoutt+1_cov1} implies $ \expt{ \sfF_{q,\out,2}^2 } = \sigma^2 $. 
We then examine 
\begin{align}
&& 
    & \frac{1}{n} \inprod{q_{\out}^r}{f_{q,\out,1}^s} , & 
    & \frac{1}{n} \inprod{q_{\out}^r}{f_{q,\out,2}} . & 
& \label{eqn:qoutt+1_TODO} 
\end{align}
The second term is easily seen to be vanishing:
\begin{align}
    \frac{1}{n} \inprod{q_{\out}^r}{f_{q,\out,2}}
    &= \frac{1}{n} \inprod{O^\top f_{p,\out}^r}{\Lambda q_{\inn,2}}
    = \frac{1}{n} \inprod{g_r(p_{\out,1}^r + \fra_r z; z, \eps, \Xi) - \wt{\fra}_r z}{z}
    \to 0 , \label{eqn:use_wta} 
\end{align}
by the definition of $ \wt{\fra}_r $ in \Cref{eqn:frab_tilde}. 
Turning to the first term in \Cref{eqn:qoutt+1_TODO}, 
\begin{align}
    & \frac{1}{n} \inprod{q_{\out}^r}{f_{q,\out,1}^s}
    = \frac{1}{n} \inprod{q_{\out}^r}{\diag(\psi_s) (q_\out^s + \wt{\fra}_s \Lambda q_{\inn,2}) + \diag_{n\times d}(\wt{\psi}_s) (q_{\inn,1}^s + \wt{\frb}_s q_{\inn,2}) - \fra_{s+1} \Lambda q_{\inn,2}} \notag \\
    &= \frac{1}{n} \inprod{q_\out^r}{\diag(\psi_s) q_\out^s} 
    + \frac{\wt{\fra}_s}{n} \inprod{q_\out^r}{\diag(\psi_s) \Lambda q_{\inn,2}} 
    + \frac{1}{n} \inprod{q_{\out}^r}{\diag_{n\times d}(\wt{\psi}_s) (q_{\inn,1}^s + \wt{\frb}_s q_{\inn,2})} 
    - \frac{\fra_{s+1}}{n} \inprod{O^\top f_{p,\out}^r}{\Lambda q_{\inn,2}} . \notag 
\end{align}
The limit of the first term is zero since $ (\sfQ_\out^r, \sfQ_\out^s) $ is independent of $ \sfPsi_s $ which has zero mean (see \Cref{asmp:4matrices}). 
The limit of the second term is also zero since the Wasserstein limits $ \sfQ_\out^r, \sfQ_{\inn,2} $ of $ q_\out^r, q_{\inn,2} $ are independent of one another. 
The third term again converges to zero due to the independence between $ \sfQ_\out^r $ and $ \sfQ_\inn^s $. 
We have just shown in \Cref{eqn:use_wta} that the fourth term converges to zero. 
Therefore, $ \frac{1}{n} \inprod{q_{\out}^r}{f_{q,\out,1}^s} \to 0 $. 

We have shown that both terms in \Cref{eqn:qoutt+1_TODO} have limit zero. 
By this and \Cref{eqn:qoutt+1_cov}, 
\begin{align}
    \frac{1}{n} U^\top U &\to \matrix{
        \Sigma_{q,\out}^t & 0_{(t+1)\times1} & 0_{(t+1)\times(t+1)} \\
        0_{1\times(t+1)} & \sigma^2 & 0_{1\times(t+1)} \\
        0_{(t+1)\times(t+1)} & 0_{(t+1)\times1} & \Sigma_{p,\out}^{t+1}
    } \in \bbR^{(1+2(t+1))\times(1+2(t+1))} . \label{eqn:UU_qoutt+1} 
\end{align}
The matrix on the right is invertible by the induction hypothesis in \Cref{itm:hyp3} and the positivity of $ \sigma^2 = \ol{\kappa}_2 \rho $ implied by \Cref{asmp:design,asmp:signal}. 

Next, we compute the limit of $ n^{-1} V^\top f_{p,\out}^{t+1} $. 
By the induction hypothesis \Cref{eqn:concl2} in \Cref{itm:hyp2}, for any $ 0\le s\le t $, we have 
\begin{align}
    \frac{1}{n} \inprod{f_{p,\out}^s}{f_{p,\out}^{t+1}} &\to (\Sigma_{q,\out}^{t+1})_{s,t+1} = (\sigma_{q,\out}^t)_s . \label{eqn:Vfpout} 
\end{align}
We then examine 
\begin{align}
&& 
    & \frac{1}{n} \inprod{p_{\out,2}}{f_{p,\out}^{t+1}} , & 
    & \frac{1}{n} \inprod{p_{\out,1}^{s+1}}{f_{p,\out}^{t+1}} . &
& \label{eqn:Vfpout_12}
\end{align}
The first term equals $ n^{-1} \inprod{z}{g_{t+1}(p_{\out,1}^{t+1} + \fra_{t+1} z; z, \eps) - \wt{\fra}_{t+1} z} $ which converges to zero, by the definition of $ \wt{\fra}_{t+1} $ in \Cref{eqn:frab_tilde}. 
As for the second term in \Cref{eqn:Vfpout_12}, 
\begin{align}
    \frac{1}{n} \inprod{p_{\out,1}^{s+1}}{f_{p,\out}^{t+1}}
    &= \frac{1}{n} \inprod{p_{\out,1}^{s+1}}{g_{t+1}(p_{\out,1}^{t+1} + \fra_{t+1} p_{\out,2} ; p_{\out,2}, w_{p,\out})}
    - \frac{\wt{\fra}_{t+1}}{n} \inprod{p_{\out,1}^{s+1}}{p_{\out,2}} . \label{eqn:second12}
\end{align}
Applying \Cref{asmp:regularity} and \Cref{prop:wto2}, noting that $ (\sfP_{\out,1}^{s+1}, \sfP_{\out,1}^{t+1}, \sfP_{\out,2}) $ is jointly Gaussian independent of $ \sfW_{p,\out} $, using Stein's lemma \cite[Theorem 2.1]{shrinkage_estimation_book} and the divergence-freeness of $ g_{t+1} $ (see \Cref{asmp:div_free}), the first term converges to zero. 
We claim that the second term in \Cref{eqn:second12} also has limit zero. 
To show this, note that 
\begin{align}
    \frac{1}{n} \inprod{p_{\out,1}^{s+1}}{p_{\out,2}}
    &= \frac{1}{n} \inprod{f_{q,\out,1}^s}{f_{q,\out,2}}
    \to \expt{ \sfF_{q,\out,1}^s \sfF_{q,\out,2} }
    = 0 , \notag 
\end{align}
by the induction hypothesis in \Cref{itm:hyp1}, as claimed. 
Thus via \Cref{eqn:second12}, we see that the second term in \Cref{eqn:Vfpout_12} has limit zero. 

We have shown that both terms in \Cref{eqn:Vfpout_12} converge to zero. 
Using this and \Cref{eqn:Vfpout}, we have 
\begin{align}
    \frac{1}{n} V^\top f_{p,\out}^{t+1} &\to \matrix{
        \sigma_{q,\out}^t \\ 0 \\ 0_{t+1}
    } \in \bbR^{1+2(t+1)} . \notag 
\end{align}
Further combining the above display with \Cref{eqn:UU_qoutt+1}, and recalling that 
\begin{align}
    U &\wto \matrix{\sfQ_\out^0 & \cdots & \sfQ_\out^t & \sfF_{q,\out,2} & \sfF_{q,\out,1}^0 & \cdots & \sfF_{q,\out.1}^1} \notag 
\end{align}
by the induction hypothesis \Cref{eqn:concl2} in \Cref{itm:hyp2}, we obtain the Wasserstein limit of the first term on the right-hand side of \Cref{eqn:qoutt+1_cond}: 
\begin{align}
    & U \cdot (n^{-1} U^\top U)^{-1} \cdot n^{-1} V^\top f_{p,\out}^{t+1} \notag \\
    &\wto \matrix{ \sfQ_\out^0 & \cdots & \sfQ_\out^t & \sfF_{q,\out,2} & \sfF_{q,\out,1}^0 & \cdots & \sfF_{q,\out,1}^t }
    \matrix{
        \Sigma_{q,\out}^t & 0_{(t+1)\times1} & 0_{(t+1)\times(t+1)} \\
        0_{1\times(t+1)} & \sigma^2 & 0_{1\times(t+1)} \\
        0_{(t+1)\times(t+1)} & 0_{(t+1)\times1} & \Sigma_{p,\out}^{t+1}
    }^{-1} \matrix{
        \sigma_{q,\out}^t \\ 0 \\ 0_{t+1}
    } \notag \\
    &= \sum_{s = 0}^t \sfQ_\out^s \paren{ (\Sigma_{q,\out}^t)^{-1} \sigma_{q,\out}^t }_s . \notag 
\end{align}
Moreover, by \Cref{prop:wto4} the second term on the right side of \Cref{eqn:qoutt+1_cond} admits a Gaussian distributional limit, denoted by $ \sfU_{q,\out}^{t+1} \in \bbR $, independent of everything else. 
So overall, the conditional empirical distribution of $ q_\out^{t+1} $ converges to that of the random variable 
\begin{align}
    \wt{\sfQ}_\out^{t+1} = \sum_{s = 0}^t \sfQ_\out^s \paren{ (\Sigma_{q,\out}^t)^{-1} \sigma_{q,\out}^t }_s + \sfU_{q,\out}^{t+1} \label{eqn:Qout^t+1}
\end{align}
which is independent of $ \sfA $ and $ (\sfP_\inn^{s+1})_{0\le s\le t}, (\sfP_\out^{s+1})_{0\le s\le t} $. 
We will next verify that this random variable $ \wt{\sfQ}_\out^{t+1} $ is precisely $ \sfQ_\out^{t+1} $ by noting that $ \sfQ_{\out}^0, \cdots, \sfQ_{\out}^t, \wt{\sfQ}_{\out}^{t+1} $ are jointly Gaussian and checking their covariance structure. 

For any $0\le s\le t$, it is easy to compute: 
\begin{subequations}
\label{eqn:qout_cov}
\begin{align}
    \expt{ \wt{\sfQ}_\out^{t+1} \sfQ_\out^r } &= \sum_{s = 0}^t (\Sigma_{q,\out}^t)_{r,s} \paren{ (\Sigma_{q,\out}^t)^{-1} \sigma_{q,\out}^t }_s
    = (\sigma_{q,\out}^t)_r , \\
    \expt{ (\wt{\sfQ}_\out^{t+1})^2 } &= \lim_{n\to\infty} \frac{1}{n} \normtwo{q_\out^{t+1}}^2 
    = \lim_{n\to\infty} \frac{1}{n} \normtwo{f_{p,\out}^{t+1}}^2 = \expt{ (\sfF_{p,\out}^{t+1})^2 } 
    = s_{q,\out}^{t+1} , 
\end{align}
\end{subequations}
thereby verifying the desired covariance structure of $ \sfQ_{\out}^0, \cdots, \sfQ_{\out}^t, \wt{\sfQ}_{\out}^{t+1} $. 
This allows us to identify $ \wt{\sfQ}_\out^{t+1} $ with $ \sfQ_\out^{t+1} $ and conclude
\begin{align}
    & \matrix{
        \bq_\out^{t+1} & 
        f_{q,\out,2} & \bf_{q,\out,1}^{t} & 
        p_{\out,2} & \bp_{\out,1}^{t+1} & 
        \bf_{p,\out}^{t+1} & 
        w_{p,\out} & w_{q,\out}
    } \notag \\
    &\qquad\wto \matrix{
        \bm\sfQ_\out^{t+1} & 
        \sfF_{q,\out,2} & \bm\sfF_{q,\out,1}^{t} &
        \sfP_{\out,2} & \bm\sfP_{\out,1}^{t+1} &
        \bm\sfF_{p,\out}^{t+1} & 
        \sfW_{p,\out}^\top & \sfW_{q,\out}^\top
    } . \label{eqn:qout_concl}
\end{align}

We also check that the variance of $ \sfU_{q,\out}^{t+1} $ is strictly positive: 
\begin{align}
    \expt{ (\sfU_{q,\out}^{t+1})^2 } &= \expt{ (\sfQ_{\out}^{t+1})^2 } - \expt{ \paren{ \sum_{s = 0}^t \sfQ_\out^s \paren{ (\Sigma_{q,\out}^t)^{-1} \sigma_{q,\out}^t } }^2 } \notag \\
    &= s_{q,\out}^{t+1} - \paren{ (\Sigma_{q,\out}^t)^{-1} \sigma_{q,\out}^t }^\top \cdot \Sigma_{q,\out}^t \cdot \paren{ (\Sigma_{q,\out}^t)^{-1} \sigma_{q,\out}^t } \notag \\
    &= s_{q,\out}^{t+1} - (\sigma_{q,\out}^t)^\top (\Sigma_{q,\out}^t)^{-1} \sigma_{q,\out}^t 
    = \frac{1}{e_{t+1}^\top (\Sigma_{q,\out}^{t+1})^{-1} e_{t+1}} > 0 , \label{eqn:Uqout^t+1>0} 
\end{align}
by the Schur complement formula and strict positive definiteness of $ \Sigma_{q,\out}^{t+1} $ shown around \Cref{eqn:Sigma_q>0}. 
This proves the second statement of \Cref{eqn:hyp4_in} in \Cref{itm:hyp4}. 

\paragraph{Induction step for $ f_{q,\inn}^{t+1}, f_{q,\out}^{t+1} $.}
As a direct consequence of \Cref{eqn:qin_concl,eqn:qout_concl} shown in the previous two steps, the update rules in \Cref{eqn:OQ2} of the abstract VAMP algorithm, and \Cref{prop:wto2}, we have 
\begin{align}
    & \matrix{
        q_{\inn,2} & \bq_{\inn,1}^{t+1} & 
        \bf_{q,\inn}^{t+1} & 
        \bp_\inn^{t+1} & 
        f_{p,\inn,2} & \bf_{p,\inn,1}^{t+1} & 
        w_{p,\inn} & w_{q,\inn}
    } \notag \\
    &\qquad\wto \matrix{
        \sfQ_{\inn,2} & \bm\sfQ_{\inn,1}^{t+1} &
        \bm\sfF_{q,\inn}^{t+1} &
        \bm\sfP_\inn^{t+1} &
        \sfF_{p,\inn,2} & \bm\sfF_{p,\inn,1}^{t+1} &
        \sfW_{p,\inn}^\top & \sfW_{q,\inn}^\top
    } , \label{eqn:fqin_conv} \\
    & \matrix{
        \bq_\out^{t+1} & 
        f_{q,\out,2} & \bf_{q,\out,1}^{t+1} & 
        p_{\out,2} & \bp_{\out,1}^{t+1} & 
        \bf_{p,\out}^{t+1} & 
        w_{p,\out} & w_{q,\out}
    } \notag \\
    &\qquad\wto \matrix{
        \bm\sfQ_\out^{t+1} &
        \sfF_{q,\out,2} & \bm\sfF_{q,\out,1}^{t+1} &
        \sfP_{\out,2} & \bm\sfP_{\out,1}^{t+1} &
        \bm\sfF_{p,\out}^{t+1} &
        \sfW_{p,\out}^\top & \sfW_{q,\out}^\top
    } . \label{eqn:fqout_conv}
\end{align}

We verify that the induction hypothesis in \Cref{itm:hyp1} continues to hold for $ \sfF_{q,\out,1}^{t+1} $. 
By the convergence result \Cref{eqn:fqout_conv} above, 
\begin{align}
    & \expt{ \sfF_{q,\out,1}^{t+1} \sfF_{q,\out,2} }
    = \lim_{n\to\infty} \frac{1}{n} \inprod{f_{q,\out,1}^{t+1}}{f_{q,\out,2}} \notag \\
    &= \lim_{n\to\infty} \frac{1}{n} \inprod{\diag(\psi_{t+1}) (q_\out^{t+1} + \wt{\fra}_{t+1} \Lambda q_{\inn,2}) + \diag_{n\times d}(\wt{\psi}_{t+1}) (q_{\inn,1}^{t+1} + \wt{\frb}_{t+1} q_{\inn,2}) - \fra_{t+2} \Lambda q_{\inn,2}}{\Lambda q_{\inn,2}} \notag \\
    &= \lim_{n\to\infty} \frac{1}{n} \inprod{\diag(\psi_{t+1}) (q_\out^{t+1} + \wt{\fra}_{t+1} \Lambda q_{\inn,2})}{\Lambda q_{\inn,2}} 
    + \lim_{n\to\infty} \frac{1}{n} \inprod{\diag_{n\times d}(\wt{\psi}_{t+1}) (q_{\inn,1}^{t+1} + \wt{\frb}_{t+1} q_{\inn,2})}{\Lambda q_{\inn,2}} 
    - \fra_{t+2} \lim_{n\to\infty} \frac{1}{n} \inprod{z}{z}.
    \label{eqn:second123} 
\end{align}
The last term above equals $ - \fra_{t+1} \sigma^2 $. 
The first term can be computed as follows: 
\begin{subequations}
\label{eqn:Fqout_term12}    
\begin{align}
    & \frac{1}{n} \inprod{\diag(\psi_{t+1}) (q_\out^{t+1} + \wt{\fra}_{t+1} \Lambda q_{\inn,2})}{\Lambda q_{\inn,2}} \notag \\
    &= \frac{1}{n} \inprod{\Lambda^\top \diag(\psi_{t+1}) q_\out^{t+1}}{q_{\inn,2}}
    + \frac{\wt{\fra}_{t+1}}{n} \inprod{\Lambda^\top \diag(\psi_{t+1}) \Lambda q_{\inn,2}}{q_{\inn,2}} \notag \\
    &\to 0 + \frac{\wt{\fra}_{t+1}}{\delta} \expt{ \sfPsi_{t+1} \sfLambda_d^2 } \expt{ \sfQ_{\inn,2}^2 } 
    = \frac{\rho}{\delta} \expt{ \sfPsi_t \sfLambda_d^2 } \wt{\fra}_{t+1} ,  
\end{align}
where the last line follows from independence of $ \sfQ_\out^{t+1}, \sfQ_{\inn,2} $. 
Finally, we turn to the second term in \Cref{eqn:second123}: 
\begin{align}
    & \frac{1}{n} \inprod{\diag_{n\times d}(\wt{\psi}_{t+1}) (q_{\inn,1}^{t+1} + \wt{\frb}_{t+1} q_{\inn,2})}{\Lambda q_{\inn,2}} \notag \\
    &= \frac{1}{n} \inprod{\Lambda^\top \diag_{n\times d}(\wt{\psi}_{t+1}) q_{\inn,1}^{t+1}}{q_{\inn,2}}
    + \frac{\wt{\frb}_{t+1}}{n} \inprod{\Lambda^\top \diag_{n\times d}(\wt{\psi}_{t+1}) q_{\inn,2}}{q_{\inn,2}} \notag \\
    &\to 0 + \frac{\wt{\frb}_{t+1}}{\delta} \expt{ \sfLambda_d \wt{\sfPsi}_{t+1} } \expt{ \sfQ_{\inn,2}^2 }
    = \frac{\rho}{\delta} \expt{ \wt{\sfPsi}_{t+1} \sfLambda_d } \wt{\frb}_{t+1} , 
\end{align}
\end{subequations}
since $ \sfQ_{\inn,2} $ and $ \sfQ_{\inn,1}^{t+1} $ are independent. 
From \Cref{eqn:Fqout_term12}, we conclude that
\begin{align}
    \frac{1}{n} \inprod{f_{q,\out,1}^{t+1}}{f_{q,\out,2}}
    &\to \frac{\rho}{\delta} \paren{
        \expt{ \sfPsi_{t+1} \sfLambda_d^2 } \wt{\fra}_{t+1} 
        + \expt{ \wt{\sfPsi}_{t+1} \sfLambda_d } \wt{\frb}_{t+1}
    }
    - \fra_{t+2} \sigma^2
    = 0 , \label{eqn:pout12} 
\end{align}
by the definition of $ \fra_{t+2} $ in \Cref{eqn:frab}. 
Here it helps to recall $ \sfLambda_n, \sfLambda_d $ from \Cref{eqn:sfLambda_nd}. 
This establishes $ \expt{\sfF_{q,\out,1}^{t+1} \sfF_{q,\out,2}} = 0 $, which, together with \Cref{eqn:Fpin1_Fpin2}, confirms the validity of \Cref{itm:hyp1} for time $t+1$. 

Given the new random variables $ \sfF_{q,\inn}^{t+1}, \sfF_{q,\out}^{t+1} $, the covariance matrices $ \Sigma_{p,\inn}^{t+2}, \Sigma_{p,\out}^{t+2} $ can be computed according to \Cref{eqn:cov_p}. 
We write these matrices in block form: 
\begin{align}
&& 
    \Sigma_{p,\inn}^{t+2} &= \matrix{
        \Sigma_{p,\inn}^{t+1} & \sigma_{p,\inn}^{t+1} \\
        (\sigma_{p,\inn}^{t+1})^\top & s_{p,\inn}^{t+2}
    } \in \bbR^{(t+2)\times(t+2)} , & 
    \Sigma_{p,\out}^{t+2} &= \matrix{
        \Sigma_{p,\out}^{t+1} & \sigma_{p,\out}^{t+1} \\
        (\sigma_{p,\out}^{t+1})^\top & s_{p,\out}^{t+2}
    } \in \bbR^{(t+2)\times(t+2)} , & 
& \label{eqn:Sigmap>0} 
\end{align}
where we have introduced the notation $ \sigma_{p,\inn}^{t+1}, \sigma_{p,\out}^{t+1} \in \bbR^{t+1}$, $ s_{p,\inn}^{t+2}, s_{p,\out}^{t+2} \in \bbR $. 
Notationally, we index the elements of $ \sigma_{p,\inn}^{t+1}, \sigma_{p,\out}^{t+1} $ using the set $ \brace{1,2,\cdots,t+1} $, which is consistent with the indexing of $ \Sigma_{p,\inn}^{t+2}, \Sigma_{p,\out}^{t+2} $. 
We now verify the invertibility of $ \Sigma_{p,\inn}^{t+2} $. 
Recalling the definition \Cref{eqn:cov_p}, given the invertibility of $ \Sigma_{p,\inn}^{t+1} $ guaranteed by the induction hypothesis in \Cref{itm:hyp3}, $ \Sigma_{p,\inn}^{t+2} $ is singular if and only if there exist $ \alpha_0, \cdots, \alpha_t \in \bbR $, such that almost surely, 
\begin{align}
    \sfF_{q,\inn}^{t+1} &= \sum_{s = 0}^t \alpha_s \sfF_{q,\inn}^s . \notag
\end{align}
By the definition \Cref{eqn:Fqin} of $ \sfF_{q,\inn}^{t+1} $, the above condition is equivalent to 
\begin{align}
    \sfPhi_{t+1} \sfQ_{\inn,1}^{t+1} + \wt{\sfPhi}_{t+1} \sfQ_\out^{t+1} 
    &= \sum_{s = 0}^t \alpha_s \sfF_{q,\inn}^s - \wt{\frb}_{t+1} \sfPhi_{t+1} \sfQ_{\inn,2} -  \wt{\fra}_{t+1} \wt{\sfPhi}_{t+1} \sfLambda_d \sfQ_{\inn,2} + \frb_{t+2} \sfQ_{\inn,2} . \label{eqn:inv_cond}
\end{align}
\Cref{asmp:4matrices} ensures that at least one of the random variables $ \sfPhi_{t+1}, \wt{\sfPhi}_{t+1} $ is not almost surely a constant. 
Suppose that this is the case for $ \sfPhi_{t+1} $. 
Using the representations of $ \sfQ_{\inn,1}^{t+1} $ identified in \Cref{eqn:Qin1^t+1}, we can further write \Cref{eqn:inv_cond} as 
\begin{align}
\begin{split}
    \sfPhi_{t+1} \sfU_{q,\inn}^{t+1} &= \sum_{s = 0}^t \alpha_s \sfF_{q,\inn}^s - \wt{\frb}_{t+1} \sfPhi_{t+1} \sfQ_{\inn,2} -  \wt{\fra}_{t+1} \wt{\sfPhi}_{t+1} \sfLambda_d \sfQ_{\inn,2} + \frb_{t+2} \sfQ_{\inn,2} \\
    &\quad - \wt{\sfPhi}_{t+1} \sfQ_\out^{t+1} - \sfPhi_{t+1} \sum_{s = 0}^t \sfQ_{\inn,1}^s \paren{(\Sigma_{q,\inn}^t)^{-1} \sigma_{q,\inn}^t}_s . 
\end{split}
\label{eqn:RV}
\end{align}
Then we condition on a realization $ (\phi, f) $ of $ \sfPhi_{t+1} $ and the right-hand side of the above equation such that $ \phi \ne 0 $. 
Such $\phi$ exists since we assumed that $ \sfPhi_{t+1} $ is non-constant. 
The conditioning does not change the distribution of $ \sfU_{\inn,1}^{t+1} $ since it is independent of other random variables in \Cref{eqn:RV}. 
Since $ \sfU_{q,\inn}^{t+1} $ is Gaussian with strictly positive variance as shown in \Cref{eqn:Uqin^t+1>0}, we have $ \prob{\sfU_{q,\inn}^{t+1} \ne f/\phi} > 0 $. 
This implies that \Cref{eqn:RV} cannot hold almost surely and therefore $ \Sigma_{p,\inn}^{t+2} $ must be non-singular. 
On the other hand, if $ \sfPhi_{t+1} $ is constant, \Cref{asmp:4matrices} ensures that $ \wt{\sfPhi}_{t+1} $ is non-constant. 
We can run the same argument but now use the representation of $ \sfQ_\out^{t+1} $ identified in \Cref{eqn:Qout^t+1} and leave $ \wt{\sfPhi}_{t+1} \sfQ_\out^{t+1} $ on the left-hand side of \Cref{eqn:RV}. 
Positivity of the variance of $ \sfU_{q,\out}^{t+1} $ shown in \Cref{eqn:Uqout^t+1>0} allows us to draw the same conclusion of invertibility of $ \Sigma_{p,\inn}^{t+2} $. 
An analogous proof applies to $ \Sigma_{p,\out}^{t+2} $ whose details we omit. 

\paragraph{Induction step for $ p_\inn^{t+2} $.}
We redefine $ U,V $ to be 
\begin{align}
    U &= \matrix{
        q_{\inn,2} & q_{\inn,1}^0 & \cdots & q_{\inn,1}^{t+1} & f_{q,\inn}^0 & \cdots & f_{q,\inn}^t
    } \in \bbR^{d\times2(t+2)} , \label{eqn:U3} \\
    V &= \matrix{
        f_{p,\inn,2} & f_{p,\inn,1}^0 & \cdots & f_{p,\inn,1}^{t+1} & p_\inn^1 & \cdots & p_\inn^{t+1}
    } \in \bbR^{d\times2(t+2)} . \label{eqn:V3}
\end{align}
The conditional law of $ p_\inn^{t+2} $ given all past iterates and side information is the same as that of $ p_\inn^{t+2} $ given the condition $ V = Q U $. 
By \Cref{prop:haar_cond}, for every sufficiently large $d$, the latter conditional law can be written as
\begin{align}
    \left. p_\inn^{t+2} \mid \brace{ V = Q U } \right. 
    = \left. Q f_{q,\inn}^{t+1} \mid \brace{ V = Q U } \right. 
    &\eqqlaw V (U^\top U)^{-1} U^\top f_{q,\inn}^{t+1} + \Pi_{V^\perp} \wt{Q} \Pi_{U^\perp}^\top f_{q,\inn}^{t+1} , \label{eqn:step3_TODO} 
\end{align}
where $ \wt{Q} \sim \haar(\bbO(d - 2(t+2))) $ is independent of everything else. 
To proceed, we compute the limits of $ d^{-1} U^\top U $ and $ d^{-1} U^\top f_{q,\inn}^{t+1} $. 
In particular, \Cref{eqn:UU3} below shows that $ U^\top U $ on the right of \Cref{eqn:step3_TODO} is invertible for every sufficiently large $d$. 

The steps for computing the limit of $ d^{-1} U^\top U $ are similar to those for \Cref{eqn:UU_result}. 
Compared to \Cref{eqn:U1}, $U$ defined in \Cref{eqn:U3} above has one additional column $ q_{\inn,1}^{t+1} $. 
So we only need to check the inner product between $ q_{\inn,1}^{t+1} $ and other columns of $U$. 
By the convergence result \Cref{eqn:qin_concl} shown above, for $ 0\le s\le t $, we have 
\begin{align}
&&
    \frac{1}{d} \inprod{q_{\inn,1}^{t+1}}{q_{\inn,2}} &\to 0 , &
    \frac{1}{d} \inprod{q_{\inn,1}^{t+1}}{q_{\inn,1}^s} &\to (\sigma_{q,\inn}^t)_s , &
    \frac{1}{d} \inprod{q_{\inn,1}^{t+1}}{q_{\inn,1}^{t+1}} &\to s_{q,\inn}^{t+1} . &
& \notag 
\end{align}
Furthermore, by an analysis similar to that for the first term in \Cref{eqn:entry12}, we have 
\begin{align}
    \frac{1}{d} \inprod{q_{\inn,1}^{t+1}}{f_{q,\inn}^s} &\to 0 . \notag 
\end{align}
Combining these facts with \Cref{eqn:UU_result}, we obtain
\begin{align}
    \frac{1}{d} U^\top U &\to \matrix{
        \rho & 0_{1\times(t+2)} & 0_{1\times(t+1)} \\
        0_{(t+2)\times1} & \Sigma_{q,\inn}^{t+1} & 0_{(t+2) \times (t+1)} \\
        0_{(t+1)\times1} & 0_{(t+1) \times (t+2)} & \Sigma_{p,\inn}^{t+1}
    } \in \bbR^{2(t+2) \times 2(t+2)} . \label{eqn:UU3}  
\end{align}
The matrix on the right-hand side is strictly positive definite by the positivity of $\rho$ guaranteed by \Cref{asmp:signal}, the positive definiteness of $ \Sigma_{q,\inn}^{t+1} $ shown around \Cref{eqn:Sigma_q>0}, and the positive definiteness of $ \Sigma_{p,\inn}^{t+1} $ due to the induction hypothesis in \Cref{itm:hyp3}. 

Moving to $ d^{-1} U^\top f_{q,\inn}^{t+1} $, 
by the convergence result \Cref{eqn:fqin_conv} shown above and  analysis similar to that for the two terms in \Cref{eqn:entry12}, for every $ 0\le s \le t+1 $ and $ 0\le r\le t $, 
\begin{align}
&& 
    \frac{1}{d} \inprod{q_{\inn,2}}{f_{q,\inn}^{t+1}} &\to 0 , & 
    \frac{1}{d} \inprod{q_{\inn,1}^s}{f_{q,\inn}^{t+1}} &\to 0 , & 
    \frac{1}{d} \inprod{f_{q,\inn}^r}{f_{q,\inn}^{t+1}} &\to (\sigma_{p,\inn}^{t+1})_r . & 
& \notag 
\end{align}
Therefore,
\begin{align}
    \frac{1}{d} U^\top f_{q,\inn}^{t+1} &\to \matrix{
        0 \\ 0_{t+2} \\ \sigma_{p,\inn}^{t+1}
    } \in \bbR^{2(t+2) \times 1} . \label{eqn:Ufqin3}
\end{align}

Plugging \Cref{eqn:UU3,eqn:Ufqin3} in the right-hand side of \Cref{eqn:step3_TODO}, using the fact that the second term in \Cref{eqn:step3_TODO} has a Gaussian limit due to \Cref{prop:wto4}, denoted by $ \sfU_{p,\inn}^{t+2} $, independent of everything else, recalling from the induction hypothesis \Cref{eqn:concl1} in \Cref{itm:hyp2} that
\begin{align}
    V &\wto \matrix{\sfF_{p,\inn,2} & \sfF_{p,\inn,1}^0 & \cdots & \sfF_{p,\inn,1}^{t+1} & \sfP_\inn^1 & \cdots & \sfP_\inn^{t+1}} , \notag 
\end{align}
and applying \Cref{prop:wto3}, we have that the conditional law of $ p_\inn^{t+2} $ converges in Wasserstein to the law of 
\begin{align}
    & \matrix{
        \sfF_{p,\inn,2} & \sfF_{p,\inn,1}^0 & \cdots & \sfF_{p,\inn,1}^{t+1} & \sfP_\inn^1 & \cdots & \sfP_\inn^{t+1}
    }
    \matrix{
        \rho^{-1} & 0_{1\times(t+2)} & 0_{1\times(t+1)} \\
        0_{(t+2)\times1} & (\Sigma_{q,\inn}^{t+1})^{-1} & 0_{(t+2) \times (t+1)} \\
        0_{(t+1)\times1} & 0_{(t+1) \times (t+2)} & (\Sigma_{p,\inn}^{t+1})^{-1}
    } \matrix{
        0 \\ 0_{t+2} \\ \sigma_{p,\inn}^{t+1}
    } + \sfU_{p,\inn}^{t+2} \notag \\
    &= \sum_{s = 1}^{t+1} \sfP_\inn^s \paren{ (\Sigma_{p,\inn}^{t+1})^{-1} \sigma_{p,\inn}^{t+1} }_s + \sfU_{p,\inn}^{t+2} . \label{eqn:Pint+2}
\end{align}
The random variable above is jointly Gaussian with $ (\sfP_\inn^{s+1})_{0\le s\le t} $ and they are independent of $ \sfA $ and $ (\sfQ_\inn^s)_{0\le s\le t+1}, (\sfQ_\out^s)_{0\le s\le t+1} $. 
It can be verified in a way similar to \Cref{eqn:check_cov} that the above random variable has the desired covariance matrix $ \Sigma_{p,\inn}^{t+2} $ with $ (\sfP_\inn^{s+1})_{0\le s\le t} $ and therefore is precisely the random variable $ \sfP_\inn^{t+2} $. 
Furthermore, using the Schur complement formula, positive definiteness of $ \Sigma_{p,\inn}^{t+2} $ shown around \Cref{eqn:Sigmap>0}, the independence between the two terms in \Cref{eqn:Pint+2} and applying arguments similar to \Cref{eqn:var}, the variance of $ \sfU_{p,\inn}^{t+2} $ is strictly positive. 
We omit these details to avoid repetition. 
This proves the first statement of \Cref{eqn:hyp4_out} in \Cref{itm:hyp4}. 

The above reasoning establishes the convergence result: 
\begin{align}
    & \matrix{
        q_{\inn,2} & \bq_{\inn,1}^{t+1} & 
        \bf_{q,\inn}^{t+1} & 
        \bp_\inn^{t+2} & 
        f_{p,\inn,2} & \bf_{p,\inn,1}^{t+1} & 
        w_{p,\inn} & w_{q,\inn}
    } \notag \\
    &\qquad\wto \matrix{
        \sfQ_{\inn,2} & \bm\sfQ_{\inn,1}^{t+1} &
        \bm\sfF_{q,\inn}^{t+1} &
        \bm\sfP_\inn^{t+2} &
        \sfF_{p,\inn,2} & \bm\sfF_{p,\inn,1}^{t+1} &
        \sfW_{p,\inn}^\top & \sfW_{q,\inn}^\top
    } . \label{eqn:pin_conv} 
\end{align}

\paragraph{Induction step for $ p_\out^{t+2} $.}
Redefine $ U,V $ to be 
\begin{align}
    U &= \matrix{
        q_\out^0 & \cdots & q_\out^{t+1} & f_{q,\out,2} & f_{q,\out,1}^0 & \cdots & f_{q,\out,1}^t
    } \in \bbR^{n\times2(t+2)} , \label{eqn:U4} \\
    V &= \matrix{
        f_{p,\out}^0 & \cdots & f_{p,\out}^{t+1} & p_{\out,2} & p_{\out,1}^1 & \cdots & p_{\out,1}^{t+1}
    } \in \bbR^{n\times2(t+2)} . \label{eqn:V4}
\end{align}
The conditional law of $ p_\out^{t+2} $ given all past iterates and side information is the same as that of $ p_\out^{t+2} $ given the condition $ V = O U $. 
We then use \Cref{prop:haar_cond} to characterize the latter conditional distribution: for all sufficiently large $n$, 
\begin{align}
    \left. p_{\out}^{t+2} \mid \brace{ V = O U } \right. 
    = \left. O f_{q,\out}^{t+1} \mid \brace{ V = O U } \right. 
    \eqqlaw V (U^\top U)^{-1} U^\top f_{q,\out}^{t+1} + \Pi_{V^\perp} \wt{O} \Pi_{U^\perp}^\top f_{q,\out}^{t+1} , \label{eqn:step4_TODO} 
\end{align}
where $ \wt{O} \sim \haar(\bbO(n - 2(t+2))) $ is independent of everything else. 
Again, we compute the limits $ n^{-1} U^\top U $ and $ n^{-1} U^\top f_{q,\out}^{t+1} $. 
In particular, \Cref{eqn:UU4} below justifies, for all sufficiently large $n$, the invertibility of $ U^\top U $ which shows up on the right of \Cref{eqn:step4_TODO}. 

First consider $ n^{-1} U^\top U $. 
Note that $U$ defined in \Cref{eqn:U4} has one additional column $ q_\out^{t+1} $ compared to $U$ in \Cref{eqn:U2}. 
So to compute the limit of $ n^{-1} U^\top U $, we only need to examine the inner product between $ q_\out^{t+1} $ and other columns of $U$. 
By the convergence result \Cref{eqn:qout_concl} already shown, for any $ 0\le s\le t $, 
\begin{align}
&& 
    \frac{1}{n} \inprod{q_\out^{t+1}}{q_\out^s} &\to (\sigma_{q,\out}^t)_s , & 
    \frac{1}{n} \inprod{q_\out^{t+1}}{q_\out^{t+1}} &\to s_{q,\out}^{t+1} , &
& \notag \\
&&
    \frac{1}{n} \inprod{q_\out^{t+1}}{f_{q,\out,2}} &\to 0 , & 
    \frac{1}{n} \inprod{q_\out^{t+1}}{f_{q,\out,1}^s} &\to 0 , & 
& \notag 
\end{align}
where the last two limits can be treated similarly to \Cref{eqn:qoutt+1_TODO}. 
Combining these limits with \Cref{eqn:UU_qoutt+1} yields
\begin{align}
    \frac{1}{n} U^\top U &\to \matrix{
        \Sigma_{q,\out}^{t+1} & 0_{(t+2)\times1} & 0_{(t+2)\times(t+1)} \\
        0_{1\times(t+2)} & \sigma^2 & 0_{1\times(t+1)} \\
        0_{(t+1)\times(t+2)} & 0_{(t+1)\times1} & \Sigma_{p,\out}^{t+1}
    } \in \bbR^{2(t+2)\times2(t+2)} . \label{eqn:UU4} 
\end{align}
The matrix on the right is strictly positive definite by the positive definiteness of $ \Sigma_{q,\out}^{t+1} $ shown around \Cref{eqn:Sigma_q>0}, positivity of $ \sigma^2 $ noted in \Cref{eqn:sigma2} and the positive definiteness of $ \Sigma_{p,\out}^{t+1} $ from the induction hypothesis in \Cref{itm:hyp3}. 

Next consider $ n^{-1} U^\top f_{q,\out}^{t+1} $. 
For any $ 0\le s\le t+1 $ and $ 0\le r \le t $, 
\begin{align}
&& 
    \frac{1}{n} \inprod{f_{q,\out,1}^{t+1}}{q_\out^s} &\to 0 , & 
    \frac{1}{n} \inprod{f_{q,\out,1}^{t+1}}{f_{q,\out,2}} &\to 0 , & 
    \frac{1}{n} \inprod{f_{q,\out,1}^{t+1}}{f_{q,\out,1}^r} &\to (\sigma_{p,\out}^{t+1})_r , & 
& \label{eqn:fqout1}
\end{align}
where the first limit can be treated similarly to the first term in \Cref{eqn:qoutt+1_TODO}.
The second limit follows from $ \expt{ \sfF_{q,\out,1}^{t+1} \sfF_{q,\out,2} } = 0 $ shown in \Cref{eqn:pout12}; 
the third one directly follows from the convergence result \Cref{eqn:fqout_conv} shown previously. 
Moreover, it can be established similarly to the second term in \Cref{eqn:qoutt+1_TODO} that 
\begin{align}
    \frac{1}{n} \inprod{f_{q,\out,2}}{q_\out^s} &\to 0 . \label{eqn:fqout2}
\end{align}
Recall also the second equation in \Cref{eqn:qoutt+1_cov1}. 
This, together with \Cref{eqn:fqout1,eqn:fqout2}, implies
\begin{align}
    \frac{1}{n} U^\top f_{q,\out}^{t+1} &\to \matrix{
        0_{t+2} & 0_{t+2} \\ 
        0 & \sigma^2 \\ 
        \sigma_{p,\out}^{t+1} & 0_{t+2}
    } \in \bbR^{2(t+2)\times2} . \label{eqn:Ufqout} 
\end{align}

By \Cref{prop:wto4}, the second term in \Cref{eqn:step4_TODO} admits a Gaussian limit, denoted by $ \matrix{ \sfU_{p,\out}^{t+2} & 0 } \in \bbR^{1\times2} $. 
The second entry vanishes since every column of $ \Pi_{U^\perp} $ is orthogonal to the column $ f_{q,\out,2} $ in $ U $ and therefore $ \Pi_{U^\perp}^\top f_{q,\out,2} = 0_{n-2(t+2)} $. 
Combining this observation with \Cref{eqn:UU4,eqn:Ufqout}, recalling from the induction hypothesis \Cref{eqn:concl2} in \Cref{itm:hyp2} that 
\begin{align}
    V &\wto \matrix{\sfF_{p,\out}^0 & \cdots & \sfF_{p,\out}^{t+1} & \sfP_{\out,2} & \sfP_{\out,1}^1 & \cdots & \sfP_{\out,1}^{t+1}} , \notag 
\end{align}
and applying \Cref{prop:wto3} give us the Wasserstein limit of the right-hand side of \Cref{eqn:step4_TODO}: 
\begin{align}
    & \left. p_{\out}^{t+2} \mid \brace{V = OU} \right. \notag \\
    &\wto \matrix{\sfF_{p,\out}^0 & \cdots & \sfF_{p,\out}^{t+1} & \sfP_{\out,2} & \sfP_{\out,1}^1 & \cdots & \sfP_{\out,1}^{t+1}} \matrix{
        \Sigma_{q,\out}^{t+1} & 0_{(t+2)\times1} & 0_{(t+2)\times(t+1)} \\
        0_{1\times(t+2)} & \sigma^2 & 0_{1\times(t+1)} \\
        0_{(t+1)\times(t+2)} & 0_{(t+1)\times1} & \Sigma_{p,\out}^{t+1}
    }^{-1} \notag \\
    &\qquad \times \matrix{
        0_{t+2} & 0_{t+2} \\ 
        0 & \sigma^2 \\ 
        \sigma_{p,\out}^{t+1} & 0_{t+2}
    } + \matrix{ \sfU_{p,\out}^{t+2} & 0 } \notag \\
    &= \matrix{
        \displaystyle \sum_{s=0}^t \sfP_{\out,1}^{s+1} \paren{ (\Sigma_{p,\out}^{t+1})^{-1} \sigma_{p,\out}^{t+1} }_{s+1} + \sfU_{p,\out}^{t+2} & \sfP_{\out,2}
    } . \notag 
\end{align}
We claim that the first entry on the right-hand side is precisely $ \sfP_{\out,1}^{t+2} $. 
To justify this, we take notice of the following facts. 
The random variable in the first entry is jointly Gaussian with $ (\sfP_{\out,1}^{s+1})_{0\le s\le t} $, with covariance matrix $ \Sigma_{p,\out}^{t+2} $ (following calculations similar to \Cref{eqn:qout_cov}). 
These jointly Gaussian random variables are independent of $ \sfA $ and $ (\sfQ_\inn^s)_{0\le s\le t+1}, (\sfQ_\out^s)_{0\le s\le t+1} $. 
Furthermore, they are also independent of $ \sfP_{\out,2} $ since both $ (\sfP_{\out,1}^{s+1})_{0\le s\le t} $ and $ \sfU_{p,\out}^{t+2} $ are independent of $ \sfP_{\out,2} $. 
As before, the variance of $ \sfU_{p,\out}^{t+2} $ is strictly positive, where we use the positive definiteness of $ \Sigma_{p,\out}^{t+2} $ shown around \Cref{eqn:Sigmap>0}. 
This proves the second statement of \Cref{eqn:hyp4_out} in \Cref{itm:hyp4}. 

Therefore, we have shown 
\begin{align}
    & \matrix{
        \bq_\out^{t+1} & 
        f_{q,\out,2} & \bf_{q,\out,1}^{t} & 
        p_{\out,2} & \bp_{\out,1}^{t+2} & 
        \bf_{p,\out}^{t+1} & 
        w_{p,\out} & w_{q,\out}
    } \notag \\
    &\qquad\wto \matrix{
        \bm\sfQ_\out^{t+1} &
        \sfF_{q,\out,2} & \bm\sfF_{q,\out,1}^{t} &
        \sfP_{\out,2} & \bm\sfP_{\out,1}^{t+2} &
        \bm\sfF_{p,\out}^{t+1} &
        \sfW_{p,\out}^\top & \sfW_{q,\out}^\top
    } . 
    \label{eqn:pout_conv}
\end{align}

\paragraph{Induction step for $ f_{p,\inn}^{t+2}, f_{p,\out}^{t+2} $.}
By the update rules \Cref{eqn:OQ4}, the convergence results \Cref{eqn:pin_conv,eqn:pout_conv} shown in the last two steps, and \Cref{prop:wto2}, we immediately get 
\begin{align}
    & \matrix{
        q_{\inn,2} & \bq_{\inn,1}^{t+1} & 
        \bf_{q,\inn}^{t+1} & 
        \bp_\inn^{t+2} & 
        f_{p,\inn,2} & \bf_{p,\inn,1}^{t+2} & 
        w_{p,\inn} & w_{q,\inn}
    } \notag \\
    &\qquad\wto \matrix{
        \sfQ_{\inn,2} & \bm\sfQ_{\inn,1}^{t+1} &
        \bm\sfF_{q,\inn}^{t+1} &
        \bm\sfP_\inn^{t+2} &
        \sfF_{p,\inn,2} & \bm\sfF_{p,\inn,1}^{t+2} &
        \sfW_{p,\inn}^\top & \sfW_{q,\inn}^\top
    } , \notag \\
    & \matrix{
        \bq_\out^{t+1} & 
        f_{q,\out,2} & \bf_{q,\out,1}^{t+1} & 
        p_{\out,2} & \bp_{\out,1}^{t+2} & 
        \bf_{p,\out}^{t+2} & 
        w_{p,\out} & w_{q,\out}
    } \notag \\
    &\qquad\wto \matrix{
        \bm\sfQ_\out^{t+1} &
        \sfF_{q,\out,2} & \bm\sfF_{q,\out,1}^{t+1} &
        \sfP_{\out,2} & \bm\sfP_{\out,1}^{t+2} &
        \bm\sfF_{p,\out}^{t+2} &
        \sfW_{p,\out}^\top & \sfW_{q,\out}^\top
    } , \notag 
\end{align}
which finally completes the whole induction argument. 
\end{proof}

\subsection{Proof of \Cref{cor:degenerate}}
\label{app:pf_cor:degenerate}

We perturb the random matrices $ \Phi_t(X^\top X), \wt{\Phi}_t(X), \Psi_t(XX^\top), \wt{\Psi}_t(X) $ and the denoising functions $ f_{t+1}, g_{t+1} $ in the original GVAMP iteration \Cref{eqn:GVAMP} and show that the perturbed version satisfies \Cref{asmp:nondegenerate_abs,asmp:regularity,asmp:4matrices} so that \Cref{thm:SE_abs_VAMP} is applicable. 
We then show that these two iterations and also their associated state evolution parameters are close to each other and conclude the proof by sending the perturbation to zero. 

Fix a small constant $ \zeta > 0 $ and an arbitrary integer $ \Gamma \in \bbZ_{\ge0} $. 
Let $ M_t \iid \GOE(d) $, $ N_t \iid \GOE(n) $, $ \wt{M}_t \iid \Gin(n,d) $, $ \wt{N}_t \iid \Gin(n,d) $, $ s^t \iid \cN(0_d, I_d) $, $ q^t \iid \cN(0_n, I_n) $ (for all $ t\in\brace{0,1,\cdots,\Gamma+1} $) be independent of each other and of everything else.\footnote{Slightly abuse terminology in \Cref{sec:strong_free}, we denote by $ \Gin(n,d) $ the (rectangular) Ginibre ensemble which is the probability measure on $ n\times d $ real rectangular matrices $ M $ such that $ M_{i,j} \iid \cN(0,1/d) $ for all $ (i,j)\in[n]\times[d] $. } 
Define the augmented side information matrices $ \ul{\Theta} = \matrix{\Theta & s^0 & \cdots & s^{\Gamma+1}} \in \bbR^{d\times(d'+\Gamma+2)} $ and $ \ul{\Xi} = \matrix{\Xi & q^0 & \cdots & q^{\Gamma+1}} \in \bbR^{n\times(n'+\Gamma+2)} $. 
Let the perturbed initializers be given by 
\begin{align}
&&
    \wt{r}_\zeta^0 &= \wt{r}^0 + \zeta s^0, & 
    \wt{p}_\zeta^0 &= \wt{p}^0 + \zeta q^0 , & 
& \label{eqn:pert_init}
\end{align}
where $ \wt{r}^0, \wt{p}^0 $ are the initializers of \Cref{eqn:GVAMP}. 
Defining
\begin{align}
&&
    f_0^\zeta(b, \theta, \matrix{s_0 & \cdots & s_{\Gamma+1}}) &= f_0(b, \theta) + \zeta s_0 , &
    g_0^\zeta(z, e, \xi, \matrix{q_0 & \cdots & q_{\Gamma+1}}) &= g_0(z, e, \xi) + \zeta q_0 , & 
& \notag 
\end{align}
it is clear that the perturbed initializers satisfy \Cref{asmp:init} in which the Wasserstein limits of $ \wt{r}_\zeta^0, \wt{p}_\zeta^0 $ are given by 
\begin{align}
&&
    \wt{\sfR}_0^\zeta &= \wt{\sfR}_0 + \zeta \sfS_0 , & 
    \wt{\sfP}_0^\zeta &= \wt{\sfP}_0 + \zeta \sfQ_0 , &
& \label{eqn:art_SE_init}
\end{align}
where $ \sfS_0, \sfQ_0 $ are i.i.d.\ standard Gaussians independent of $ \sfB_*, \sfZ, \sfE, \sfTheta, \sfXi $. 
For each $ 0\le t\le \Gamma $, define the perturbed denoising functions $ f^\zeta_{t+1} \colon \bbR^{2+d'+\Gamma+2} \to \bbR $ and $ g^\zeta_{t+1} \colon \bbR^{3+n'+\Gamma+2} \to \bbR $ as
\begin{align}
\begin{split}
    f^\zeta_{t+1}(r_\zeta; b, \matrix{\theta & s_0 & \cdots & s_{\Gamma+1}}) &= f_{t+1}(r_\zeta; b, \theta) + \zeta s_{t+1} - \expt{f_{t+1}'(\sfR_{t+1}^\zeta; \sfB_*, \sfTheta)} r_\zeta , \\
    g^\zeta_{t+1}(p_\zeta; z, e, \matrix{\xi & q_0 & \cdots & q_{\Gamma+1}}) &= g_{t+1}(p_\zeta; z, e, \xi) + \zeta q_{t+1} - \expt{g_{t+1}'(\sfP_{t+1}^\zeta; \sfZ, \sfE, \sfXi)} p_\zeta , 
\end{split}
\label{eqn:pert_fg}
\end{align}
where $ r_\zeta,p_\zeta,b,z,e,s_0,\cdots,s_{\Gamma+1},q_0,\cdots,q_{\Gamma+1} \in \bbR $, $\theta\in\bbR^{d'}$, $\xi\in\bbR^{n'}$, and $ f_{t+1}, g_{t+1} $ are the denoising functions in \Cref{eqn:GVAMP}. 
The random variables $ \sfR_{t+1}^\zeta, \sfP_{t+1}^\zeta $ are given in \Cref{eqn:pert_SE} below. 
In particular, the negative correction terms in the definition \Cref{eqn:pert_fg} ensures that the divergence-free condition required in \Cref{asmp:div_free} is satisfied by $ f_{t+1}^\zeta, g_{t+1}^\zeta $ when the expectation of their derivatives with respect to $ r_\zeta, p_\zeta $ are evaluated using the laws of $ \sfR_{t+1}^\zeta, \sfP_{t+1}^\zeta $, respectively. 
Since $ f_{t+1}, g_{t+1} $ satisfy \Cref{asmp:remove}, $ f_{t+1}^\zeta, g_{t+1}^\zeta $ satisfy \Cref{asmp:regularity}. 
For any $ \beta_0, \cdots, \beta_t \in \bbR $, suppose 
\begin{align}
    \wt{\sfR}_{t+1}^\zeta - \wt{b}_{t+1}^\zeta \sfB_* &= \sum_{s = 0}^t \beta_s \paren{\wt{\sfR}_s^\zeta - \wt{b}_s^\zeta \sfB_*} . \label{eqn:contradiction}
\end{align}
Rearranging terms, this is equivalent to 
\begin{align}
    \zeta \sfS_{t+1} &= - f_{t+1}(\sfR_{t+1}^\zeta; \sfB_*, \sfTheta) + \expt{f_{t+1}'(\sfR_{t+1}^\zeta; \sfB_*, \sfTheta)} + \wt{b}_{t+1}^\zeta \sfB_* + \sum_{s = 0}^t \beta_s \paren{\wt{\sfR}_s^\zeta - \wt{b}_s^\zeta \sfB_*} . \label{eqn:contradiction2}
\end{align}
Since $ \zeta>0 $ and $ \sfS_{t+1} \sim \cN(0,1) $ is independent of all terms on the right of \Cref{eqn:contradiction2}, \Cref{eqn:contradiction} cannot hold almost surely. 
A similar argument applies to $ \wt{\sfP}_{t+1}^\zeta $. 
Therefore $ f_{t+1}^\zeta, g_{t+1}^\zeta $ satisfy \Cref{asmp:nondegenerate_abs}. 

Using the above quantities, we construct a perturbed GVAMP iteration. 
Recall from \Cref{asmp:design} that $ X = O \Lambda Q^\top $. 
Initialized with $ \wt{r}_\zeta^0, \wt{p}_\zeta^0 $ in \Cref{eqn:pert_init}, for every $ 0\le t\le \Gamma $, the perturbed GVAMP iterates are updated as follows: 
\begin{subequations}
\label{eqn:GVAMP_pert}    
\begin{align}
&&
    r_\zeta^{t+1} &= Q (\Phi_t(\Lambda^\top \Lambda) + \zeta M_t) Q^\top \wt{r}_\zeta^t + Q (\wt{\Phi}_t(\Lambda) + \zeta \wt{M}_t)^\top O^\top \wt{p}_\zeta^t , & 
    \wt{r}_\zeta^{t+1} &= f^\zeta_{t+1}(r_\zeta^{t+1}; \beta_*, \ul{\Theta}) , & 
& \label{eqn:GVAMP_pert1} \\ 
&&
    p_\zeta^{t+1} &= O (\Psi_t(\Lambda \Lambda^\top) + \zeta N_t) O^\top \wt{p}_\zeta^t + O (\wt{\Psi}_t(\Lambda) + \zeta \wt{N}_t) Q^\top \wt{r}_\zeta^t , & 
    \wt{p}_\zeta^{t+1} &= g^\zeta_{t+1}(p_\zeta^{t+1}; z, \eps, \ul{\Xi}) . &
& \label{eqn:GVAMP_pert2}
\end{align}
\end{subequations}

Let $ m_t \in \bbR^d, n_t\in\bbR^n $ be the vectors of eigenvalues of $ M_t, N_t $, respectively. 
Let $ \wt{m}_t \in \bbR^d , \wt{n}_t \in \bbR^n $ be the vectors of singular values of $ \wt{M}_t, \wt{N}_t $, respectively, appended with zeros if needed, as in \Cref{eqn:lambda_nd}. 
Then we identify 
\begin{align}
&&
    \phi^\zeta_{t} &= \Phi_t(\lambda_d \circ \lambda_d) + \zeta m_t , & 
    \wt{\phi}^\zeta_{t} &= \wt{\Phi}_t(\lambda_d) + \zeta \wt{m}_t , & 
    \psi^\zeta_{t} &= \Psi_t(\lambda_n \circ \lambda_n) + \zeta n_t , & 
    \wt{\psi}^\zeta_{t} &= \wt{\Psi}_t(\lambda_n) + \zeta \wt{n}_t . & 
& \label{eqn:pert_phipsi}
\end{align}
Using orthogonal invariance and the convergence of empirical spectral distributions of GOE and the Ginibre ensemble, \Cref{asmp:tr_free} and \Cref{prop:wto2}, we obtain \Cref{eqn:conv_phi_psi} in which the limiting random variables are given by
\begin{align}
&&
    \sfPhi^\zeta_{t} &= \Phi_t(\sfLambda_d^2) + \zeta \sfM_t , & 
    \wt{\sfPhi}^\zeta_{t} &= \wt{\Phi}_t(\sfLambda_d) + \zeta \wt{\sfM}_t , & 
    \sfPsi^\zeta_{t} &= \Psi_t(\sfLambda_n^2) + \zeta \sfN_t , & 
    \wt{\sfPsi}^\zeta_{t} &= \wt{\Psi}_t(\sfLambda_n) + \zeta \wt{\sfN}_t , & 
& \label{eqn:pert_PhiPsi}
\end{align}
where $ \sfM_t, \sfN_t $ follow the semicircle distribution, $ \wt{\sfM}_t, \wt{\sfN}_t $ follow the (square root pushforward of) Marchenko--Pastur distribution with aspect ration $ \delta $, moreover, they are independent of each other and of $ \sfLambda $ for all $t\ge0$. 
Since the semicircle distribution has mean zero, so do $ \sfPhi^\zeta_{t}, \sfPsi^\zeta_{t} $. 
Since $ \zeta>0 $ and the variances of $ \sfM_t, \wt{\sfM}_t, \sfN_t, \wt{\sfN}_t $ are all positive, it holds that $ \min\brace{\expt{(\sfPhi_t^\zeta)^2}, \var{\wt{\sfPhi}_t^\zeta}} > 0$, $\min\brace{\expt{(\sfPsi_t^\zeta)^2}, \var{\wt{\sfPsi}_t^\zeta}} > 0 $. 
Therefore, \Cref{asmp:4matrices} is satisfied by the vectors in \Cref{eqn:pert_phipsi} and the associated random variables in \Cref{eqn:pert_PhiPsi}. 

Now \Cref{thm:SE_abs_VAMP} is applicable to the perturbed GVAMP iteration \Cref{eqn:GVAMP_pert}. 
Its state evolution parameters 
\begin{align}
&&
    & a_{t+1}^\zeta ,&
    & b_{t+1}^\zeta ,&
    & \wt{a}^\zeta_t ,& 
    & \wt{b}^\zeta_t ,& 
    & \wt{\sfP}^\zeta_{t+1} ,& 
    & \wt{\sfR}^\zeta_{t+1} ,& 
    & \sfJ^\zeta_{t+1} ,& 
    & \sfK^\zeta_{t+1} ,& 
    & \sfR_{t+1}^\zeta ,& 
    & \sfP_{t+1}^\zeta ,&
    & \forall 0\le t\le \Gamma ,&
& \label{eqn:pert_SE}
\end{align}
can be defined through \Cref{eqn:SE_GVAMP} with the modification that $ \Phi_t(\sfLambda_d^2), \wt{\Phi}_t(\sfLambda_d), \Psi_t(\sfLambda_n^2), \wt{\Psi}_t(\sfLambda_n) $, $ f_{t+1}(\sfR_{t+1}; \sfB_*, \sfTheta) $, $ g_{t+1}(\sfP_{t+1}; \sfZ, \sfE, \sfXi) $ are replaced with the perturbed versions $ \sfPhi_t^\zeta, \wt{\sfPhi}_t^\zeta, \sfPsi_t^\zeta, \wt{\sfPsi}_t^\zeta $, $ f_{t+1}^\zeta(\sfR_{t+1}^\zeta; \sfB_*, \ul{\sfTheta}) $, $ g_{t+1}^\zeta(\sfP_{t+1}^\zeta; \sfZ, \sfE, \ul{\sfXi}) $ introduced in \Cref{eqn:pert_fg,eqn:pert_PhiPsi}. 

We will inductively show the following statements. 
For every $t\ge0$ and every $ 1\le r,s\le t $, 
\begin{align}
&&
    \limsup_{d\to\infty} d^{-1} \normtwo{r^{t}}^2 &< \infty , & 
    \limsup_{n\to\infty} n^{-1} \normtwo{p^{t}}^2 &< \infty , & 
& \label{eqn:finite} \\
&&
    \lim_{\zeta\to0} \limsup_{d\to\infty} d^{-1} \normtwo{r^{t} - r_\zeta^{t}}^2 &= 0 , & 
    \lim_{\zeta\to0} \limsup_{n\to\infty} n^{-1} \normtwo{p^{t} - p_\zeta^{t}}^2 &= 0 , & 
& \label{eqn:close} \\
&&
    \lim_{\zeta\to0} a_{t}^\zeta &= a_{t} , & 
    \lim_{\zeta\to0} b_{t}^\zeta &= b_{t} , & 
& \label{eqn:converge_ab} \\
&&
    \lim_{\zeta\to0} \expt{\sfJ_{r}^\zeta \sfJ_{s}^\zeta} &= \expt{\sfJ_{r} \sfJ_{s}} , & 
    \lim_{\zeta\to0} \expt{\sfK_{r}^\zeta \sfK_{s}^\zeta} &= \expt{\sfK_{r} \sfK_{s}} , & 
& \label{eqn:converge_JK} \\
&&
    \lim_{\zeta\to0} \expt{f_{t}'(\sfR_{t}^\zeta; \sfB_*, \sfTheta)} &= 0 , & 
    \lim_{\zeta\to0} \expt{g_{t}'(\sfP_{t}^\zeta; \sfZ, \sfE, \sfXi)} &= 0 , & 
& \label{eqn:vanish} \\
&&
    \limsup_{d\to\infty} d^{-1} \normtwo{\wt{r}^{t}}^2 &< \infty , & 
    \limsup_{n\to\infty} n^{-1} \normtwo{\wt{p}^{t}}^2 &< \infty , & 
& \label{eqn:finite_tilde} \\
&&
    \lim_{\zeta\to0} \limsup_{d\to\infty} d^{-1} \normtwo{\wt{r}^{t} - \wt{r}_\zeta^{t}}^2 &= 0 , & 
    \lim_{\zeta\to0} \limsup_{n\to\infty} n^{-1} \normtwo{\wt{p}^{t} - \wt{p}_\zeta^{t}}^2 &= 0 , & 
& \label{eqn:close_tilde} \\
&& 
    \lim_{\zeta\to0} \wt{a}_t^\zeta &= \wt{a}_t , & 
    \lim_{\zeta\to0} \wt{b}_t^\zeta &= \wt{b}_t , & 
& \label{eqn:converge_abtilde} 
\end{align}
where \Cref{eqn:finite,eqn:close,eqn:vanish,eqn:converge_ab,eqn:converge_JK} are vacuous for $t=0$. 
Furthermore, for every $t\ge0$, 
\begin{align}
\matrix{
    r^1 & \cdots & r^{t} & 
    \wt{r}^0 & \cdots & \wt{r}^{t} & 
    \beta_* & \Theta 
} 
&\xrightarrow{W_2} \matrix{
    \sfR_1 & \cdots & \sfR_{t} & 
    \wt{\sfR}_0 & \cdots & \wt{\sfR}_{t} & 
    \sfB_* & \sfTheta^\top
} , \label{eqn:converge_r} \\
\matrix{
    p^1 & \cdots & p^{t} & 
    \wt{p}^0 & \cdots & \wt{p}^{t} & 
    z & \eps & \Xi
} 
&\xrightarrow{W_2} \matrix{
    \sfP_1 & \cdots & \sfP_{t} & 
    \wt{\sfP}_0 & \cdots & \wt{\sfP}_{t} & 
    \sfZ & \sfE & \sfXi^\top
} , \label{eqn:converge_p}
\end{align}
where for $t=0$, these are just $ \matrix{\wt{r}^0 & \beta_* & \Theta} \xrightarrow{W_2} \matrix{\wt{\sfR}_0 & \sfB_* & \sfTheta^\top} $, $ \matrix{\wt{p}^0 & z & \eps & \Xi} \xrightarrow{W_2} \matrix{\wt{\sfP}_0 & \sfZ & \sfE & \sfXi^\top} $. 
This will imply \Cref{cor:degenerate}. 

\paragraph{\textsc{Base case.}}
For $ t = 0 $, the Wasserstein convergence in \Cref{eqn:converge_r} is given in \Cref{asmp:init}. 
It also implies $ d^{-1} \normtwo{\wt{r}^0}^2 \to \expt{\wt{\sfR}_0^2} < \infty $ almost surely. 
By the choice of the perturbed initializer \Cref{eqn:pert_init}, $ d^{-1} \normtwo{\wt{r}_\zeta^0 - \wt{r}^0}^2 = \zeta^2 d^{-1} \normtwo{s^0}^2 \to 0 $ almost surely as $ d\to\infty, \zeta\to0 $ in that order. 
According to the definition \Cref{eqn:ab_tilde} and the identification of $ \wt{\sfR}_0^\zeta $ in \Cref{eqn:art_SE_init}, $ \wt{b}_0^\zeta = \rho^{-1} \expt{(\wt{\sfR}_0 + \zeta \sfS_0) \sfB_*} = \rho^{-1} \expt{\wt{\sfR}_0 \sfB_*} = \wt{b}_0 $. 
The same arguments apply to $ \wt{p}^0, \wt{p}^0_\zeta, \wt{a}_0^\zeta $. 
This justifies \Cref{eqn:finite_tilde,eqn:close_tilde,eqn:converge_abtilde}. 

\paragraph{\textsc{Induction step.}}
Assuming that \Cref{eqn:finite,eqn:close,eqn:finite_tilde,eqn:close_tilde,eqn:vanish,eqn:converge_ab,eqn:converge_JK} hold up to time $t$, we verify their validity for time $t+1$. 

Recalling the update rule \Cref{eqn:GVAMP1} and using the triangle inequality, 
\begin{align}
    \normtwo{r^{t+1}} &\le \normtwo{\Phi_t(X^\top X)} \normtwo{\wt{r}^t} + \normtwo{\wt{\Phi}_t(X)} \normtwo{\wt{p}^t} . \notag 
\end{align}
By \Cref{asmp:remove_PhiPsi} and the induction hypothesis \Cref{eqn:finite_tilde}, we have 
\begin{align}
    \limsup_{d\to\infty} d^{-1} \normtwo{r^{t+1}}^2 < \infty . \label{eqn:finite_rt+1}
\end{align}
The same argument implies $ \limsup_{n\to\infty} n^{-1} \normtwo{p^{t+1}}^2 < \infty $, thereby establishing \Cref{eqn:finite} for $ t+1 $. 

Comparing the update rules for $ r^{t+1}, r^{t+1}_\zeta $ in \Cref{eqn:GVAMP_pert1,eqn:GVAMP1} and applying the triangle inequality, 
\begin{align}
    \normtwo{r^{t+1} - r^{t+1}_\zeta} &\le \normtwo{\Phi_t(\Lambda^\top \Lambda)} \normtwo{\wt{r}^t - \wt{r}^t_\zeta} + \zeta \normtwo{M_t} \normtwo{\wt{r}^t_\zeta} + \normtwo{\wt{\Phi}_t(\Lambda)} \normtwo{\wt{p}^t - \wt{p}^t_\zeta} + \zeta \normtwo{\wt{M}_t} \normtwo{\wt{p}^t_\zeta} . \notag 
\end{align}
Combining \Cref{asmp:remove_PhiPsi}, the induction hypotheses \Cref{eqn:finite_tilde,eqn:close_tilde}, and the facts $ \normtwo{M_t} \to 2 , \normtwo{\wt{M}_t} \to 1 + \sqrt{\delta} $, we obtain 
\begin{align}
    \lim_{\zeta\to0} \limsup_{d\to\infty} d^{-1} \normtwo{r^{t+1} - r^{t+1}_\zeta}^2 < \infty . \label{eqn:close_rt+1}
\end{align}
An analogous argument shows $ \lim_{\zeta\to0} \limsup_{n\to\infty} n^{-1} \normtwo{p^{t+1} - p^{t+1}_\zeta}^2 < \infty $, thereby justifying \Cref{eqn:close} for $t+1$. 

By the definition \Cref{eqn:ab} and the identification of $ \wt{\sfPhi}_t $ in \Cref{eqn:pert_PhiPsi}, 
\begin{align}
    b_{t+1}^\zeta &= \expt{ (\wt{\Phi}_t(\sfLambda_d) + \zeta \wt{\sfM}_t) \sfLambda_d } \wt{a}_t^\zeta
    = \expt{\wt{\Phi}_t(\sfLambda_d) \sfLambda_d} \wt{a}_t^\zeta
    \to b_{t+1} , \label{eqn:bzeta} 
\end{align}
where the convergence in the $ \zeta\to0 $ limit is by the induction hypothesis \Cref{eqn:converge_abtilde}. 
A similar argument shows $ a_{t+1}^\zeta \to a_{t+1} $ as $ \zeta\to0 $, thereby establishing \Cref{eqn:converge_ab} for $t+1$. 

For any $ 0\le s\le t $, recalling \Cref{eqn:K_cov} and using the triangle inequality, 
\begin{align}
\begin{split}
    & \abs{\expt{\sfK_{s+1}^\zeta \sfK_{t+1}^\zeta} - \expt{\sfK_{s+1} \sfK_{t+1}}} \\
    &\le \abs{\expt{\sfPhi_s^\zeta \sfPhi_t^\zeta} \expt{\wt{\sfR}_s^\zeta \wt{\sfR}_t^\zeta} - \expt{\Phi_s(\sfLambda_d^2) \Phi_t(\sfLambda_d^2)} \expt{\wt{\sfR}_s \wt{\sfR}_t}} 
    + \rho \abs{b_{s+1}^\zeta b_{t+1}^\zeta - b_{s+1} b_{t+1}} \\
    &+ \abs{\expt{\wt{\sfPhi}_s^\zeta \wt{\sfPhi}_t^\zeta} \paren{\expt{\wt{\sfP}_s^\zeta \wt{\sfP}_t^\zeta} - \wt{a}_s^\zeta \wt{a}_t^\zeta \sigma^2} - \expt{\wt{\Phi}_s(\sfLambda_d) \wt{\Phi}_t(\sfLambda_d)} \paren{\expt{\wt{\sfP}_s \wt{\sfP}_t} - \wt{a}_s \wt{a}_t \sigma^2}} \\
    &+ \rho \abs{\expt{\wt{\sfPhi}_s^\zeta \wt{\sfPhi}_t^\zeta \sfLambda_d^2} \wt{a}_s^\zeta \wt{a}_t^\zeta - \expt{\wt{\Phi}_s(\sfLambda_d) \wt{\Phi}_t(\sfLambda_d) \sfLambda_d^2} \wt{a}_s \wt{a}_t} \\
    &+ \rho \abs{\expt{\sfPhi_s^\zeta \wt{\sfPhi}_t^\zeta \sfLambda_d} \wt{a}_s^\zeta \wt{b}_t^\zeta - \expt{\Phi_s(\sfLambda_d^2) \wt{\Phi}_t(\sfLambda_d) \sfLambda_d} \wt{a}_s \wt{b}_t} 
    + \rho \abs{\expt{\sfPhi_t^\zeta \wt{\sfPhi}_s^\zeta \sfLambda_d} \wt{a}_t^\zeta \wt{b}_s^\zeta - \expt{\Phi_t(\sfLambda_d^2) \wt{\Phi}_s(\sfLambda_d) \sfLambda_d} \wt{a}_t \wt{b}_s} . 
\end{split}
\label{eqn:Kzeta} 
\end{align}
By \Cref{eqn:pert_PhiPsi}, we have
\begin{align}
    \lim_{\zeta\to0} \abs{\expt{\sfPhi_s^\zeta \sfPhi_t^\zeta} - \expt{\Phi_s(\sfLambda_d^2) \Phi_t(\sfLambda_d^2)}}
    &\le \lim_{\zeta\to0} \zeta^2 \expt{\sfM_t^2} = 0 , \notag 
\end{align}
where inequality becomes an equality if and only if $s=t$, and the last equality holds since the semicircle distribution has finite variance. 
Similarly, 
$ \expt{\wt{\sfPhi}_s^\zeta \wt{\sfPhi}_t^\zeta} \to \expt{\wt{\Phi}_s(\sfLambda_d) \wt{\Phi}_t(\sfLambda_d)} $, 
$ \expt{\wt{\sfPhi}_s^\zeta \wt{\sfPhi}_t^\zeta \sfLambda_d^2} \to \expt{\wt{\Phi}_s(\sfLambda_d) \wt{\Phi}_t(\sfLambda_d) \sfLambda_d^2} $, 
$ \expt{\sfPhi_s^\zeta \wt{\sfPhi}_t^\zeta \sfLambda_d} \to \expt{\Phi_s(\sfLambda_d^2) \wt{\Phi}_t(\sfLambda_d) \sfLambda_d} $, 
$ \expt{\sfPhi_t^\zeta \wt{\sfPhi}_s^\zeta \sfLambda_d} \to \expt{\Phi_t(\sfLambda_d^2) \wt{\Phi}_s(\sfLambda_d) \sfLambda_d} $, 
as $ \zeta\to0 $. 
Moreover, recalling \Cref{eqn:pert_fg} and using the triangle inequality, we have
\begin{align}
    & \abs{\expt{\wt{\sfR}_s^\zeta \wt{\sfR}_t^\zeta} - \expt{\wt{\sfR}_s \wt{\sfR}_t}} \label{eqn:rr} \\
    &\le \abs{\expt{f_s(b_s^\zeta \sfB_* + \sfK_s^\zeta; \sfB_*, \sfTheta) f_t(b_t^\zeta \sfB_* + \sfK_t^\zeta; \sfB_*, \sfTheta)} - \expt{f_s(b_s \sfB_* + \sfK_s; \sfB_*, \sfTheta) f_t(b_t \sfB_* + \sfK_t; \sfB_*, \sfTheta)}}
    \notag \\
    &\quad + \zeta^2 + \abs{\expt{f_s'(\sfR_s^\zeta; \sfB_*, \sfTheta)}} \cdot \abs{\expt{f_t'(\sfR_t^\zeta; \sfB_*, \sfTheta)}} \cdot \abs{\expt{\sfR_s^\zeta \sfR_t^\zeta}} \notag \\
    &\quad + \abs{\expt{f_s'(\sfR_s^\zeta; \sfB_*, \sfTheta)}} \cdot \abs{\expt{f_t(\sfR_t^\zeta; \sfB_*, \sfTheta) \sfR_s^\zeta}}
    + \abs{\expt{f_t'(\sfR_t^\zeta; \sfB_*, \sfTheta)}} \cdot \abs{\expt{f_s(\sfR_s^\zeta; \sfB_*, \sfTheta) \sfR_t^\zeta}} . \notag 
\end{align}
By the Lipschitz property of $ f_t $ from \Cref{asmp:remove}, 
\begin{align}
    & \abs{\expt{f_t(b_t^\zeta \sfB_* + \sfK_t^\zeta; \sfB_*, \sfTheta) - f_t(b_t \sfB_* + \sfK_t; \sfB_*, \sfTheta)}} \notag \\
    &\le \expt{\abs{f_t\paren{b_t^\zeta \sfB_* + \expt{(\sfK_t^\zeta)^2}^{1/2}\sfG; \sfB_*, \sfTheta} - f_t\paren{b_t \sfB_* + \expt{\sfK_t^2}^{1/2} \sfG; \sfB_*, \sfTheta}}} \notag \\
    &\le C \expt{\abs{(b_t^\zeta - b_t) \sfB_* + \paren{\expt{(\sfK_t^\zeta)^2}^{1/2} - \expt{\sfK_t^2}^{1/2}} \sfG}} \notag \\
    &\le C \paren{\abs{b_t^\zeta - b_t} \rho^{1/2} + \abs{\expt{(\sfK_t^\zeta)^2} - \expt{\sfK_t^2}}^{1/2}}
    \xrightarrow[\zeta\to0]{} 0 , \label{eqn:ff} 
\end{align}
where $ \sfG \sim \cN(0,1) $ is independent of $ \sfB_* $ and the last line follows from the induction hypotheses \Cref{eqn:converge_ab,eqn:converge_JK}. 
Combining \Cref{eqn:ff} and the hypotheses \Cref{eqn:converge_ab,eqn:converge_JK,eqn:vanish}, we have that \Cref{eqn:rr} converges to $0$ as $\zeta\to0$. 
Using all the above observations in \Cref{eqn:Kzeta}, we obtain 
\begin{align}
    \expt{\sfK_{s+1}^\zeta \sfK_{t+1}^\zeta} \to \expt{\sfK_{s+1} \sfK_{t+1}} . \label{eqn:Kzeta_t+1}
\end{align}
An analogous argument shows $ \expt{\sfJ_{s+1}^\zeta \sfJ_{t+1}^\zeta} \to \expt{\sfJ_{s+1} \sfJ_{t+1}} $ as $ \zeta\to0 $ for every $ 0\le s\le t $. 
Now \Cref{eqn:converge_JK} is established for $ t+1 $. 

Next, by \Cref{eqn:bzeta,eqn:Kzeta_t+1} (with $ s=t $) shown above, 
\begin{align}
    \expt{f_{t+1}'(\sfR_{t+1}^\zeta; \sfB_*, \sfTheta)} - \expt{f_{t+1}'(\sfR_{t+1}; \sfB_*, \sfTheta)}
    &= \expt{f_{t+1}'(b_{t+1}^\zeta \sfB_* + \sfK_{t+1}^\zeta; \sfB_*, \sfTheta)} - \expt{f_{t+1}'(b_{t+1} \sfB_* + \sfK_{t+1}; \sfB_*, \sfTheta)}
    \to 0 , \notag 
\end{align}
as $ \zeta\to0 $. 
Since $ \expt{f_{t+1}'(\sfR_{t+1}; \sfB_*, \sfTheta)} = 0 $ by \Cref{asmp:div_free}, we have 
\begin{align}
    \lim_{\zeta\to0} \expt{f_{t+1}'(\sfR_{t+1}^\zeta; \sfB_*, \sfTheta)} &= 0 . \label{eqn:vanish_t+1}
\end{align}
A similar argument allows us to prove the other limit in \Cref{eqn:vanish} for $ t+1 $. 

For $ \wt{r}^{t+1} $, by the Lipschitz property of $ f_{t+1} $ in \Cref{asmp:remove}, 
\begin{align}
    \sum_{i = 1}^d \paren{ f_{t+1}(r^{t+1}_i; \beta_{*,i}, \Theta_{i,:}) - f_{t+1}(0; \beta_{*,i}, \Theta_{i,:}) }^2
    &\le C^2 \normtwo{r^{t+1}}^2 . \notag 
\end{align}
Combining this with \Cref{eqn:finite_rt+1,eqn:close_rt+1} shown above, we obtain 
\begin{align}
    \limsup_{d\to\infty} \frac{1}{\sqrt{d}} \normtwo{\wt{r}^{t+1}} 
    &= \limsup_{d\to\infty} \frac{1}{\sqrt{d}} \paren{\normtwo{f_{t+1}(r^{t+1}; \beta_*, \Theta)} + \zeta \normtwo{s^{t+1}} + \abs{\expt{f_{t+1}'(\sfR_{t+1}^\zeta; \sfB_*, \sfTheta)}} \normtwo{r^{t+1}_\zeta}}
    < \infty . \label{eqn:finite_rtilde}
\end{align}
The same argument implies $ \limsup_{n\to\infty} n^{-1} \normtwo{\wt{p}^{t+1}}^2 < \infty $, thereby establishing \Cref{eqn:finite_tilde}. 

Recalling the design of the perturbed denoiser \Cref{eqn:pert_fg} and applying the triangle inequality, 
\begin{align}
    \normtwo{\wt{r}^{t+1} - \wt{r}^{t+1}_\zeta} &\le \normtwo{f_{t+1}(r^{t_+1}; \beta_*, \Theta) - f_{t+1}(r^{t+1}_\zeta; \beta_*, \Theta)} + \zeta \normtwo{s^{t+1}} + \abs{\expt{f_{t+1}'(\sfR_{t+1}^\zeta; \sfB_*; \sfTheta)}} \normtwo{r^{t+1}_\zeta} \notag \\
    &\le C \normtwo{r^{t+1} - r^{t+1}_\zeta} + \zeta \normtwo{s^{t+1}} + \abs{\expt{f_{t+1}'(\sfR_{t+1}^\zeta; \sfB_*; \sfTheta)}} \paren{\normtwo{r^{t+1}} + \normtwo{r^{t+1} - r^{t+1}_\zeta}} , \notag 
\end{align}
where the last line is by the Lipschitz property of $f_{t+1}$ in \Cref{asmp:remove}. 
Using \Cref{eqn:close_rt+1,eqn:finite_rt+1,eqn:vanish_t+1} shown above, $ d^{-1} \normtwo{s^{t+1}}^2 \to 1 $, we obtain 
\begin{align}
    \lim_{\zeta\to0} \limsup_{d\to\infty} d^{-1} \normtwo{\wt{r}^{t+1} - \wt{r}^{t+1}_\zeta}^2 = 0 . \label{eqn:close_rtilde}
\end{align}
A similar argument implies $ \lim_{\zeta\to0} \limsup_{n\to\infty} n^{-1} \normtwo{\wt{p}^{t+1} - \wt{p}^{t+1}_\zeta}^2 = 0 $, thereby verifying \Cref{eqn:close_tilde} for $t+1$. 

From the definition \Cref{eqn:ab_tilde}, 
\begin{align}
    \wt{b}_{t+1}^\zeta &= \rho^{-1} \expt{\wt{\sfR}_{t+1}^\zeta \sfB_*} = \rho^{-1} \expt{f_{t+1}(b_{t+1}^\zeta \sfB_* + \sfK_{t+1}^\zeta; \sfB_*, \sfTheta) \sfB_*} - \rho^{-1} \expt{f_{t+1}'(\sfR_{t+1}^\zeta; \sfB_*, \sfTheta)} \expt{\sfR_{t+1}^\zeta \sfB_*}
    \to \wt{b}_{t+1} , \label{eqn:wtbzeta} 
\end{align}
where the convergence follows from \Cref{eqn:vanish_t+1,eqn:bzeta,eqn:Kzeta_t+1}. 
Similarly, $ \wt{a}_{t+1}^\zeta \to \wt{a}_{t+1} $, establishing \Cref{eqn:converge_abtilde} for $t+1$. 

Next, consider \Cref{eqn:converge_r} for $t+1$. 
Denote by $ R^{t+1} \in \bbR^{d\times(2(t+2)+d')} $ and $ \bm{\sfR}_{t+1} \in \bbR^{1\times(2(t+2)+d')} $ the left-/right-hand sides of \Cref{eqn:converge_r} for $t+1$. 
Denote also their perturbed versions by
\begin{align}
    R^{t+1}_\zeta &\coloneqq \matrix{r_\zeta^1 & \cdots & r_\zeta^{t+1} & \wt{r}_\zeta^0 & \cdots & \wt{r}_\zeta^{t+1} & \beta_* & \Theta} \in \bbR^{d\times(2(t+2)+d')} , \notag \\
    \bm{\sfR}_{t+1}^\zeta &\coloneqq \matrix{\sfR^\zeta_1 & \cdots & \sfR^\zeta_{t+1} & \wt{\sfR}^\zeta_0 & \cdots & \wt{\sfR}^\zeta_{t+1} & \sfB_* & \sfTheta} \in \bbR^{1\times(2(t+2)+d')} . \notag
\end{align}
For any $ h \in \PL_{(2(t+2)+d') \to 1}^2 $, we apply the triangle inequality to bound
\begin{align}
    \limsup_{d\to\infty} \abs{\frac{1}{d} \sum_{i = 1}^d h(R_{i,:}^{t+1}) - h(\bm{\sfR}_{t+1})} 
    &\le \lim_{\zeta\to0} \limsup_{d\to\infty} \abs{\frac{1}{d} \sum_{i = 1}^d \brack{h(R^{t+1}_{i,:}) - h((R^{t+1}_\zeta)_{i,:})}} \label{eqn:R1} \\
    &\quad + \limsup_{d\to\infty} \abs{\frac{1}{d} \sum_{i = 1}^d h((R^{t+1}_\zeta)_{i,:}) - \expt{h(\bm{\sfR}_{t+1}^\zeta)}} \label{eqn:R2} \\
    &\quad + \lim_{\zeta\to0} \abs{\expt{h(\bm{\sfR}_{t+1}^\zeta)} - \expt{h(\bm{\sfR}_{t+1})}} . \label{eqn:R3}
\end{align}
Using the pseudo-Lipschitz property of $ h $, the triangle inequality and Cauchy--Schwarz, we upper bound the first term in \Cref{eqn:R1} by 
\begin{align}
    &\abs{\frac{1}{d} \sum_{i = 1}^d \brack{h(R^{t+1}_{i,:}) - h((R^{t+1}_\zeta)_{i,:})}}
    \le \frac{C}{d} \sum_{i = 1}^d \paren{1 + \normtwo{R^{t+1}_{i,:}} + \normtwo{(R^{t+1}_\zeta)_{i,:}}} \normtwo{R^{t+1}_{i,:} - (R^{t+1}_\zeta)_{i,:}} \notag \\
    &\le C \sqrt{\frac{3}{d} \sum_{i = 1}^d \paren{1 + 4 \normtwo{R^{t+1}_{i,:}}^2 + \normtwo{R^{t+1}_{i,:} - (R^{t+1}_\zeta)_{i,:}}^2}} \sqrt{\frac{1}{d} \sum_{i = 1}^d \normtwo{R^{t+1}_{i,:} - (R^{t+1}_\zeta)_{i,:}}^2} , \notag 
\end{align}
which vanishes as $ d\to\infty, \zeta\to0 $ in that order, using the induction hypotheses \Cref{eqn:finite,eqn:close,eqn:finite_tilde,eqn:close_tilde} along with their counterparts \Cref{eqn:finite_rt+1,eqn:close_rt+1,eqn:finite_rtilde,eqn:close_rtilde} for $t+1$ shown above. 
The second term in \Cref{eqn:R2} equals zero by the Wasserstein convergence result in \Cref{thm:SE_abs_VAMP} applied to $ R^{t+1}_\zeta $ for any fixed $ \zeta>0 $. 
The third term in \Cref{eqn:R3} is zero due to the convergence of perturbed state evolution parameters to the original ones guaranteed by the induction hypotheses \Cref{eqn:converge_ab,eqn:converge_JK,eqn:vanish,eqn:converge_abtilde} together with their counterparts \Cref{eqn:bzeta,eqn:Kzeta_t+1,eqn:vanish_t+1,eqn:wtbzeta} for $ t+1 $. 
Taking \Cref{eqn:R1,eqn:R2,eqn:R3} collectively, we conclude \Cref{eqn:converge_r} for $t+1$. 
The convergence result \Cref{eqn:converge_p} for $t+1$ follows from a similar argument which is omitted. 
This completes the induction and therefore the proof of \Cref{cor:degenerate}. 

\subsection{Proof of \Cref{thm:spec_GVAMP}}
\label{app:pf_thm:spec_GVAMP}

Fix a large constant $\Gamma \ge 1$. 
Consider an auxiliary GVAMP iteration: 
\begin{subequations}
\label{eqn:art_GVAMP}
\begin{align}
&&
    \ul{r}^{t+1} &= \ul{\Phi}_t(X^\top X) \wt{\ul{r}}^t + \wt{\ul{\Phi}}_t(X)^\top \wt{\ul{p}}^t , & 
    \wt{\ul{r}}^{t+1} &= \ul{f}_{t+1}(\ul{r}^{t+1}) , & 
& \label{eqn:art_GVAMP1} \\ 
&&
    \ul{p}^{t+1} &= \ul{\Psi}_t(XX^\top) \wt{\ul{p}}^t + \wt{\ul{\Psi}}_t(X) \wt{\ul{r}}^t , & 
    \wt{\ul{p}}^{t+1} &= \ul{g}_{t+1}(\ul{p}^{t+1}; y) . &
& \label{eqn:art_GVAMP2}
\end{align}
For $ 0\le t\le \Gamma-1 $, we choose 
\begin{align}
&&
    \ul{\Phi}_t(x) &= 0 , &
    \wt{\ul{\Phi}}_t(x) &= \ol{\kappa}_2^{-1} x , & 
    \ul{\Psi}_t(x) &= \ol{\kappa}_2^{-1} x - 1 , & 
    \wt{\ul{\Psi}}_t(x) &= 0 , & 
    \ul{f}_t(\ul{r}) &= 0 , & 
    \ul{g}_t(\ul{p}; y) &= \ol{g}(y) \ul{p} , & 
& \label{eqn:denoiser1} 
\end{align}
such that 
\begin{align}
\begin{split}
    \ul{r}^{t+1} &= \ol{\kappa}_2^{-1} X^\top \wt{\ul{p}}^t , \\
    \ul{p}^{t+1} &= \paren{ \ol{\kappa}_2^{-1} X X^\top - I_n } \wt{\ul{p}}^t , \qquad 
    \wt{\ul{p}}^{t+1} = \ol{G} \ul{p}^{t+1} , 
\end{split}
\label{eqn:art_GVAMP0}
\end{align}
where $ \ol{G} = \diag(\ol{g}(y)) $. 
For $ t \ge \Gamma $, we choose
\begin{align}
&&
    \ul{\Phi}_t &= \Phi_{t-\Gamma} , & 
    \wt{\ul{\Phi}}_t &= \wt{\Phi}_{t-\Gamma} , & 
    \ul{\Psi}_t &= \Psi_{t-\Gamma} , & 
    \wt{\ul{\Psi}}_t &= \wt{\Psi}_{t-\Gamma} , & 
& \label{eqn:denoiser2} 
\end{align}
and
\begin{align}
&&
    \ul{f}_t(\ul{r}) &= \begin{cases}
        f_0(c_r \, \ul{r}) , & t = \Gamma \\
        f_{t-\Gamma}(\ul{r}) , & t > \Gamma 
    \end{cases} , & 
    \ul{g}_t(\ul{p};y) &= \begin{cases}
        g_0(c_p \, \ul{p}; y) , & t = \Gamma \\
        g_{t-\Gamma}(\ul{p};y) , & t > \Gamma 
    \end{cases} , & 
& \label{eqn:denoiser2_fg}
\end{align}
where the constants $c_r, c_p$ are defined in \Cref{eqn:cr_cp}. 
We initialize \Cref{eqn:art_GVAMP0} with 
\begin{align}
    \wt{\ul{p}}^0 = \gamma^{-\Gamma} w_3 \, \ol{G} \paren{\frac{1}{\sqrt{\rho + w_2}} z + \ul{j}^0} , \label{eqn:wtp0}
\end{align}
where $ \ul{j}^0 \sim \cN\paren{0_n, \frac{w_1}{\rho + w_2} I_n} $ is independent of everything else. 
Note that $ \wt{r}^0 $ is not needed to initialize \Cref{eqn:art_GVAMP0}. 
\end{subequations}

\Cref{lem:SE_aux} below offers a state evolution result for the auxiliary GVAMP iteration and is a direct consequence of the general result of \Cref{cor:degenerate}. 

\begin{subequations}
\label{eqn:art_SE}    
Recall $\sfZ, \sfY$ from \Cref{eqn:BZEY}. 
Let 
\begin{align}
    \wt{\ul{\sfP}}_0 &= \gamma^{-\Gamma} w_3 \, \ol{g}(\sfY) \paren{\frac{1}{\sqrt{\rho + w_2}} \sfZ + \ul{\sfJ}_0} , \label{eqn:wtP0}
\end{align}
where $ \ul{\sfJ}_0 \sim \cN\paren{0, \frac{w_1}{\rho + w_2}} $ is independent of $\sfZ, \sfY$. 
Let $ (\ul{\sfJ}_{t+1})_{t\ge0}, (\ul{\sfK}_{t+1})_{t\ge0} $ be two centered Gaussian processes that are independent of each other and of $ (\sfB_*, \sfZ, \sfE, \sfLambda) $ with covariance structures given below. 

\begin{enumerate}
\item For $ 0\le r,s\le \Gamma-1 $ and $ 1\le t\le \Gamma $, 
\begin{align}
    \ul{a}_t &= \gamma^{-(\Gamma - t)} \frac{w_3}{\sqrt{\rho + w_2}} , \qquad \expt{\ul{\sfJ}_t^2} = \gamma^{-2(\Gamma - t)} \frac{w_3^2 w_1}{\rho + w_2} , \label{eqn:art_first1} \\
    \ul{b}_t &= \gamma^{-(\Gamma-t)} \frac{1}{\sqrt{\rho + w_2}} , \qquad \expt{\ul{\sfK}_t^2} = \gamma^{-2(\Gamma - t)} \frac{w_2}{\rho + w_2} , \label{eqn:art_first2} \\
    \expt{\ul{\sfK}_{r+1} \ul{\sfK}_{s+1}} 
    &= \frac{\delta}{\ol{\kappa}_2} \expt{\ol{g}(\sfY)^2} \expt{\ul{\sfJ}_r \ul{\sfJ}_s} + \delta \rho \brace{ \expt{\ol{g}(\sfY)^2 \ol{\sfZ}^2} + \frac{\ol{\kappa}_4}{\ol{\kappa}_2^2} \expt{\ol{g}(\sfY) \ol{\sfZ}^2}^2 } \ul{a}_r \ul{a}_s , \label{eqn:art_first3} \\
    \expt{\ul{\sfJ}_{r+1} \ul{\sfJ}_{s+1}} 
    &= \paren{\frac{\ol{\kappa}_4}{\ol{\kappa}_2^2} + \delta} \expt{\ol{g}(\sfY)^2} \expt{\ul{\sfJ}_r \ul{\sfJ}_s} \notag \\
    &\quad + \rho \brace{
        \paren{\frac{\ol{\kappa}_4}{\ol{\kappa}_2} + \delta \ol{\kappa}_2} \expt{\ol{g}(\sfY)^2 \ol{\sfZ}^2} 
        + \paren{\frac{\ol{\kappa}_6}{\ol{\kappa}_2^2} - \frac{\ol{\kappa}_4^2}{\ol{\kappa}_2^3} + \delta \frac{\ol{\kappa}_4}{\ol{\kappa}_2}} \expt{\ol{g}(\sfY) \ol{\sfZ}^2}^2
    } \ul{a}_r \ul{a}_s . \label{eqn:art_first4} 
\end{align}
We also define $ \ul{\sfR}_t = \ul{b}_t \sfB_* + \ul{\sfK}_t, \ul{\sfP}_t = \ul{a}_t \sfZ + \ul{\sfJ}_t $ for every $1\le t \le \Gamma$, and $ \wt{\ul{\sfP}}_t = \ol{g}(\sfY) \ul{\sfP}_t, \wt{\ul{a}}_t = \expt{\ol{g}(\sfY) \ol{\sfZ}^2} \ul{a}_t $ for every $ 1\le t\le\Gamma - 1 $. 
Note that due to the choice of denoisers in \Cref{eqn:denoiser1} for the first phase, we do not need to define $ \wt{\ul{\sfR}}_t, \wt{\ul{b}}_t $ for $ 0\le t\le \Gamma-1 $. 

\item For $ r,s \ge \Gamma $, the correlations $ \expt{\ul{\sfJ}_{r+1} \ul{\sfJ}_{s+1}}, \expt{\ul{\sfK}_{r+1} \ul{\sfK}_{s+1}} $ together with the coefficients $ \wt{\ul{a}}_t, \wt{\ul{b}}_t, \ul{a}_{t+1}, \ul{b}_{t+1} $ and the random variables $ \ul{\sfR}_{t+1}, \ul{\sfP}_{t+1}, \wt{\ul{\sfR}}_{t}, \wt{\ul{\sfP}}_{t} $ (for $t\ge\Gamma$) formally follow the same recursion \Cref{eqn:SE_GVAMP} where the denoisers in \Cref{eqn:denoiser2,eqn:denoiser2_fg} are used. 

\item For $ 0\le r\le \Gamma-1 < s $, 
\begin{align}
    \expt{\ul{\sfJ}_{r+1} \ul{\sfJ}_{s+1}} &= \ol{\kappa}_2^{-1} \expt{\ul{\Psi}_s(\sfLambda_n^2) \sfLambda_n^2} \expt{\wt{\ul{\sfP}}_r \wt{\ul{\sfP}}_s} - 2 \expt{\ul{\Psi}_s(\sfLambda_n^2) \sfLambda_n^2} \wt{\ul{a}}_r \wt{\ul{a}}_s \rho + \ol{\kappa}_2^{-1} \expt{\ul{\Psi}_s(\sfLambda_n^2) \sfLambda_n^4} \wt{\ul{a}}_r \wt{\ul{a}}_s \rho \notag \\
    &\quad + \expt{(\ol{\kappa}_2^{-1} \sfLambda_n^2 - 1) \wt{\ul{\Psi}}_s(\sfLambda_n) \sfLambda_n} \wt{\ul{a}}_r \wt{\ul{b}}_s \rho - \ul{a}_{r+1} \ul{a}_{s+1} \sigma^2 , \\
    \expt{\ul{\sfK}_{r+1} \ul{\sfK}_{s+1}} &= \ol{\kappa}_2^{-1} \expt{\wt{\ul{\Phi}}_s(\sfLambda_d) \sfLambda_d} \expt{\wt{\ul{\sfP}}_r \wt{\ul{\sfP}}_s} + \expt{(\ol{\kappa}_2^{-1} - 1) \wt{\ul{\Phi}}_s(\sfLambda_d) \sfLambda_d} \wt{\ul{a}}_r \wt{\ul{a}}_s \rho \notag \\
    &\quad + \ol{\kappa}_2^{-1} \expt{\ul{\Phi}_s(\sfLambda_d^2) \sfLambda_d^2} \wt{\ul{a}}_r \wt{\ul{b}}_s \rho - \ul{b}_{r+1} \ul{b}_{s+1} \rho . 
\end{align}
\end{enumerate}
\end{subequations}

\begin{lemma}
\label{lem:SE_aux}
Consider the auxiliary GVAMP iteration \Cref{eqn:art_GVAMP}. 
Then for any fixed $t\ge0$, 
\begin{align}
\matrix{
    \ul{r}^1 & \cdots & \ul{r}^{t+1} & 
    \wt{\ul{r}}^\Gamma & \cdots & \wt{\ul{r}}^{t+1} & 
    \beta_* 
} 
&\xrightarrow{W_2} \matrix{
    \ul{\sfR}_1 & \cdots & \ul{\sfR}_{t+1} & 
    \wt{\ul{\sfR}}_\Gamma & \cdots & \wt{\ul{\sfR}}_{t+1} & 
    \sfB_* 
} , \notag \\
\matrix{
    \ul{p}^1 & \cdots & \ul{p}^{t+1} & 
    \wt{\ul{p}}^0 & \cdots & \wt{\ul{p}}^{t+1} & 
    z & \eps
} 
&\xrightarrow{W_2} \matrix{
    \ul{\sfP}_1 & \cdots & \ul{\sfP}_{t+1} & 
    \wt{\ul{\sfP}}_0 & \cdots & \wt{\ul{\sfP}}_{t+1} & 
    \sfZ & \sfE
} , \notag
\end{align}
where the random variables on the right are defined in \Cref{eqn:art_SE}. 
\end{lemma}

\begin{proof}
The lemma follows directly from the general result of \Cref{cor:degenerate}. 
In particular, specializing the state evolution recursion \Cref{eqn:SE_GVAMP} to the first phase \Cref{eqn:denoiser1} and using the moment-cumulant relation in \Cref{sec:cumulant} yield \Cref{eqn:art_first1,eqn:art_first2,eqn:art_first3,eqn:art_first4} which are consistent with \Cref{eqn:stein2,eqn:omega_tau,eqn:nu_nu,eqn:tau_tau} obtained from a reduction to \cite[Theorem 1]{VKM}. 
We leave out these tedious calculations but only make the clarification that for $1\le t\le \Gamma$, the notation
\begin{align}
&&
    & \ul{p}^t , &
    & \wt{\ul{p}}^t , &
    & \ul{a}_t , &
    & \ul{\sfJ}_t , &
    & \ul{\sfP}_t , &
    & \ul{r}^t , &
    & \ul{b}_t , &
    & \ul{\sfK}_t , &
    & \ul{\sfR}_t &
& \notag 
\end{align}
in \Cref{lem:SE_aux} corresponds to the notation
\begin{align}
&&
    & u^t , &
    & g^t , &
    & \nu_t , &
    & \tau_t \sfN_t , &
    & \sfZ_t , &
    & v^t , &
    & \chi_t , &
    & \omega_t \sfM_t , &
    & \sfB_t &
& \notag 
\end{align}
in \Cref{prop:SE}. 
Finally, we explain how to obtain the expressions \Cref{eqn:art_first1,eqn:art_first2}. 
Note that the only difference between the initial conditions \Cref{eqn:wtp0,eqn:lin_init} is in the additional $ \gamma^{-\Gamma} w_3 $ factor in \Cref{eqn:wtp0}. 
Inspecting the proof of \Cref{lem:SE_stay}, this does not change the validity of the following properties of the state evolution parameters: for every $ 2\le t\le \Gamma $, 
\begin{align}
    \frac{\ul{a}_t}{\ul{a}_{t-1}} &= \frac{\ul{b}_{t}}{\ul{b}_{t-1}}
    = \frac{\expt{\ul{\sfJ}_{t}^2}^{1/2}}{\expt{\ul{\sfJ}_{t-1}^2}^{1/2}}
    = \frac{\expt{\ul{\sfK}_{t}^2}^{1/2}}{\expt{\ul{\sfK}_{t-1}^2}^{1/2}} = \gamma , \notag 
\end{align}
and for every $ 1\le t\le \Gamma $, 
\begin{align}
&&
    \frac{\ul{a}_t}{\ul{b}_t} &= w_3 , & 
    \frac{\expt{\ul{\sfK}_t^2}^{1/2}}{\ul{b}_t} &= \sqrt{w_2} . & 
& \label{eqn:w_relation}
\end{align}
By induction on $1\le t\le \Gamma$, we have 
\begin{align}
    \paren{\ul{a}_t, \ul{b}_t, \expt{\ul{\sfJ}_t^2}^{1/2}, \expt{\ul{\sfK}_t^2}^{1/2}} = \gamma^t \paren{\ul{a}_0, \ul{b}_0, \expt{\ul{\sfJ}_0^2}^{1/2}, \expt{\ul{\sfK}_0^2}^{1/2}} . \label{eqn:gamma_recur}
\end{align}
According to the initial condition \Cref{eqn:wtP0}, 
\begin{align}
&&
    \ul{a}_0 &= \gamma^{-\Gamma} w_3 \frac{1}{\sqrt{\rho + w_2}} , & 
    \expt{\ul{\sfJ}_0^2}^{1/2} &= \gamma^{-\Gamma} w_3 \sqrt{\frac{w_1}{\rho + w_2}} . &
& \notag 
\end{align}
Plugging this in \Cref{eqn:gamma_recur} and using \Cref{eqn:w_relation}, we conclude that \Cref{eqn:art_first1,eqn:art_first2} hold. 
\end{proof}

\Cref{thm:spec_GVAMP} immediately follows from combining \Cref{eqn:art_r,eqn:art_p} in the lemma below with \Cref{lem:SE_aux}. 

\begin{lemma}
\label{lem:art_close_to_true}
For every $ \tau\ge0 $, 
\begin{align}
&&
    \lim_{\Gamma \to \infty} \abs{\ul{a}_\Gamma - \frac{a_0}{c_p}} = \lim_{\Gamma \to \infty} \abs{\ul{a}_{\Gamma+\tau+1} - a_{\tau+1}} &= 0 , & 
    \lim_{\Gamma \to \infty} \abs{\ul{b}_\Gamma - \frac{b_0}{c_r}} = \lim_{\Gamma \to \infty} \abs{\ul{b}_{\Gamma+\tau+1} - a_{\tau+1}} &= 0 , & 
& \label{eqn:art_ab} \\
&& 
    \lim_{\Gamma\to\infty} \abs{\wt{\ul{a}}_{\Gamma+\tau} - \wt{a}_{\tau}} &= 0 , &
    \lim_{\Gamma\to\infty} \abs{\wt{\ul{b}}_{\Gamma+\tau} - \wt{b}_{\tau}} &= 0 . &
& \label{eqn:art_abtilde}
\end{align}
For every $ \varsigma,\tau\ge1 $, 
\begin{align}
    \lim_{\Gamma\to\infty} \abs{\expt{\ul{\sfJ}_\Gamma^2} - \frac{1}{c_p^2} \expt{\sfJ_0^2}} = 
    \lim_{\Gamma\to\infty} \abs{\expt{\ul{\sfJ}_\Gamma \ul{\sfJ}_{\Gamma+\tau}} - \frac{1}{c_p} \expt{\sfJ_0 \sfJ_{\tau}}} =
    \lim_{\Gamma\to\infty} \abs{\expt{\ul{\sfJ}_{\Gamma+\varsigma} \ul{\sfJ}_{\Gamma+\tau}} - \expt{\sfJ_{\varsigma} \sfJ_{\tau}}} &= 0 , 
    \label{eqn:art_J} \\
    \lim_{\Gamma\to\infty} \abs{\expt{\ul{\sfK}_\Gamma^2} - \frac{1}{c_r^2} \expt{\sfK_0^2}} =
    \lim_{\Gamma\to\infty} \abs{\expt{\ul{\sfK}_\Gamma \ul{\sfK}_{\Gamma+\tau}} - \frac{1}{c_r} \expt{\sfK_0 \sfK_{\tau}}} =
    \lim_{\Gamma\to\infty} \abs{\expt{\ul{\sfK}_{\Gamma+\varsigma} \ul{\sfK}_{\Gamma+\tau}} - \expt{\sfK_{\varsigma} \sfK_{\tau}}} &= 0 . 
    \label{eqn:art_K}
\end{align}
For any fixed $\tau\ge0$ and any $ \phi_r \in \PL^2_{2(\tau+1)+1\to1}, \phi_p \in \PL^2_{2(\tau+2)\to1} $, almost surely
\begin{align}
    \lim_{\Gamma\to\infty} \lim_{d\to\infty} \abs{\frac{1}{d} \sum_{i = 1}^d \paren{\phi_r(c_r \ul{r}^\Gamma_i, \cdots, \ul{r}^{\Gamma+\tau}_i, \wt{\ul{r}}^\Gamma_i, \cdots, \wt{\ul{r}}^{\Gamma+\tau}_i, \beta_{*,i}) - \phi_r(r^0_i, \cdots, r^{\tau}_i, \wt{r}^0_i, \cdots, \wt{r}^{\tau}_i, \beta_{*,i})}} &= 0 , \label{eqn:art_r} \\
    \lim_{\Gamma\to\infty} \lim_{n\to\infty} \abs{\frac{1}{n} \sum_{i = 1}^n \paren{\phi_p(c_p \ul{p}^\Gamma_i, \cdots, \ul{p}^{\Gamma+\tau}_i, \wt{\ul{p}}^\Gamma_i, \cdots, \wt{\ul{p}}^{\Gamma+\tau}_i, z_i, \eps_i) - \phi_p(p^0_i, \cdots, p^{\tau}_i, \wt{p}^0_i, \cdots, \wt{p}^{\tau}_i, z_i, \eps_i)}} &= 0 . \label{eqn:art_p}
\end{align}
\end{lemma}

\begin{proof}
The proof is by induction on time. 

\paragraph{\textsc{Base case.}}
We take time $ \Gamma $ in the auxiliary GVAMP iteration \Cref{eqn:art_GVAMP} as the base case. 
By \Cref{eqn:art_first1} and the definition \Cref{eqn:cr_cp} of $ c_p $, 
\begin{align}
    c_p \, \ul{a}_\Gamma &= \frac{\sqrt{\rho + w_2}}{w_1 w_3} \frac{w_3}{\sqrt{\rho + w_2}} = \frac{1}{w_1} \explain{\Cref{eqn:a0_b0}} a_0 . \label{eqn:a=}
\end{align}
It can be shown similarly that 
\begin{align}
    c_r \, \ul{b}_\Gamma &= b_0 . \label{eqn:b=}
\end{align}
This proves the first equalities in \Cref{eqn:art_ab}. 
Again by \Cref{eqn:art_first1}, 
\begin{align}
    c_p^2 \expt{\ul{\sfJ}_\Gamma^2} &= \frac{\rho + w_2}{(w_1 w_3)^2} \frac{w_3^2 w_1}{\rho + w_2} = \frac{1}{w_1} \explain{\Cref{eqn:a0_b0}} a_0 . \label{eqn:J=} 
\end{align}
A similar calculation shows 
\begin{align}
    c_r^2 \expt{\ul{\sfK}_\Gamma^2} &= b_0 , \label{eqn:K=}
\end{align}
thereby establishing the first equalities in \Cref{eqn:art_J,eqn:art_K}. 
By the design of $ \ul{g}_\Gamma $ in \Cref{eqn:denoiser2_fg} along with \Cref{eqn:a=,eqn:J=} shown above, 
\begin{align}
    \wt{\ul{a}}_\Gamma &= \sigma^{-2} \expt{\ul{g}_\Gamma\paren{c_p \paren{\ul{a}_\Gamma \sfZ + \expt{\ul{\sfJ}_\Gamma^2}^{1/2} \sfG}; \sfY} \sfZ}
    = \sigma^{-2} \expt{g_0\paren{a_0 \sfZ + \sqrt{a_0} \sfG; \sfY} \sfZ}
    \explain{\Cref{eqn:wta0_wtb0}} \wt{a}_0 , \notag 
\end{align}
where $ \sfG \sim \cN(0,1) $ independent of $ \sfZ, \sfY $. 
Similarly, $ \wt{\ul{b}}_\Gamma = \wt{b}_0 $. 
Therefore, \Cref{eqn:art_abtilde} holds for $ \tau = 0 $. 

Let $ \phi_r \in \PL_{3\to1}^2 $. 
Then we have
\begin{align}
    & \abs{ \frac{1}{d} \sum_{i = 1}^d \paren{\phi_r(c_r \ul{r}_i^\Gamma, \wt{\ul{r}}_i^\Gamma, \beta_{*,i}) - \phi_r(r_i^0, \wt{r}_i^0, \beta_{*,i})} } \notag \\
    &\le \frac{C}{d} \sum_{i = 1}^d \paren{1 + \normtwo{\matrix{c_r \ul{r}_i^\Gamma & \wt{\ul{r}}_i^\Gamma & \beta_{*,i}}} + \normtwo{\matrix{r_i^0 & \wt{r}_i^0 & \beta_{*,i}}}} \normtwo{\matrix{c_r \ul{r}_i^\Gamma & \wt{\ul{r}}_i^\Gamma} - \matrix{r_i^0 & \wt{r}_i^0}} \notag \\
    &\le C \sqrt{3} \brack{ 1 + d^{-1} \paren{\normtwo{c_r \ul{r}^\Gamma}^2 + \normtwo{\wt{\ul{r}}^\Gamma}^2 + \normtwo{r^0}^2 + \normtwo{\wt{r}^0}^2 + 2 \normtwo{\beta_*}^2} }^{1/2} 
    \brack{d^{-1} \paren{\normtwo{c_r\ul{r}^\Gamma - r^0}^2 + \normtwo{\wt{\ul{r}}^\Gamma - \wt{r}^0}^2}}^{1/2} , \label{eqn:phir} 
\end{align}
where the last line is by Cauchy--Schwarz. 
Inspecting the proof of \Cref{lem:align}, we find that despite of an additional multiplicative factor $ \gamma^{-\Gamma} w_3 $ in the initializer \Cref{eqn:wtp0}, the following result still holds: 
\begin{align}
    \lim_{\Gamma\to\infty} \lim_{d\to\infty} \frac{\abs{\inprod{\ul{r}^\Gamma}{v_1(D)}}}{\normtwo{\ul{r}^\Gamma}} &= 1 , \notag 
\end{align}
almost surely. 
Moreover, due to the rescaling of the initializer \Cref{eqn:wtp0}, it holds almost surely that for any fixed $\Gamma$, 
\begin{align}
    \lim_{d\to\infty} d^{-1} \normtwo{\ul{r}^\Gamma}^2 &= \ul{b}_\Gamma^2 \rho + \expt{\ul{\sfK}_\Gamma^2} 
    \explain{\Cref{eqn:art_first2}} \frac{1}{\rho + w_2} \rho + \frac{w_2}{\rho + w_2}
    = 1 . \label{eqn:ulr=1} 
\end{align}
Combining the above two results, we have that almost surely 
\begin{align}
    \lim_{\Gamma\to\infty} \lim_{d\to\infty} d^{-1} \normtwo{\ul{r}^\Gamma - v_1(D)}^2 &= 0 . \label{eqn:close_v1}
\end{align}
In view of the choice of the initializer \Cref{eqn:spec_GVAMP_init} for the original GVAMP, this is equivalent to 
\begin{align}
    \lim_{\Gamma\to\infty} \lim_{d\to\infty} d^{-1} \normtwo{c_r \, \ul{r}^\Gamma - r^0}^2 &= 0 ,
    \qquad \textnormal{a.s.} \label{eqn:cr_r0}
\end{align}
By the choice of $ \ul{f}_\Gamma $ in \Cref{eqn:denoiser2_fg} and its Lipschitz property guaranteed by \Cref{asmp:remove}, \Cref{eqn:cr_r0} in turn immediately implies
\begin{align}
    d^{-1/2} \normtwo{\wt{\ul{r}}^\Gamma - \wt{r}^0}
    &= d^{-1/2} \normtwo{f_0(c_r\,\ul{r}^\Gamma) - f_0(r^0)}
    \le C d^{-1/2} \normtwo{c_r\,\ul{r}^\Gamma - r^0} \to 0 , \label{eqn:close_wtr} 
\end{align}
almost surely in the limit $ d\to\infty $ followed by $ \Gamma\to\infty $. 
The state evolution result in \Cref{lem:SE_aux} implies for any fixed $\Gamma$, almost surely, 
\begin{align}
    \lim_{d\to\infty} d^{-1} \normtwo{\ul{\wt{r}}^\Gamma}^2 &= \expt{\ul{f}_\Gamma\paren{c_r \paren{\ul{b}_\Gamma \sfB_* + \expt{\ul{\sfK}^2_\Gamma}^{1/2} \sfG}}^2}
    = \expt{f_0\paren{b_0 \sfB_* + \sqrt{b_0} \sfG}^2} < \infty , \label{eqn:finite_wtr} 
\end{align}
where $ \sfG \sim \cN(0,1) $ is independent of $ \sfB_* $ and the last equality follows from \Cref{eqn:b=,eqn:K=}. 
Since the right-hand side above does not depend on $\Gamma$, finiteness holds uniformly over $\Gamma$. 
Combining this with \Cref{eqn:close_wtr}, we have $ \lim_{d\to\infty} d^{-1} \normtwo{\wt{r}^0}^2 < \infty $. 
Combining \Cref{eqn:ulr=1,eqn:cr_r0}, we also have $ \lim_{d\to\infty} d^{-1} \normtwo{r^0}^2 < \infty $. 
Using all the above results together with \Cref{asmp:signal} in \Cref{eqn:phir} shows \Cref{eqn:art_r} for $\tau=0$. 

Let $ \phi_p \in \PL_{4\to1}^2 $. 
Similarly to \Cref{eqn:phir}, we have
\begin{align}
    & \abs{\frac{1}{n} \sum_{i = 1}^n \paren{\phi_p(c_p \ul{p}^\Gamma_i, \wt{\ul{p}}^\Gamma_i, z_i, \eps_i) - \phi_p(p^0_i, \wt{p}^0_i, z_i, \eps_i)}} \notag \\
    &\le C \sqrt{3} \brack{1 + n^{-1} \paren{\normtwo{c_p \ul{p}^\Gamma}^2 + \normtwo{\wt{\ul{p}}^\Gamma}^2 + \normtwo{p^0}^2 + \normtwo{\wt{p}^0}^2 + 2 \normtwo{z}^2 + 2 \normtwo{\eps}^2}}^{1/2} 
    \brack{1 + n^{-1} \paren{\normtwo{c_p \ul{p}^\Gamma - p^0}^2 + \normtwo{\wt{\ul{p}}^\Gamma - \wt{p}^0}^2}}^{1/2} . \label{eqn:phip} 
\end{align}
We examine each term on the right individually. 
First consider $ n^{-1} \normtwo{c_p \ul{p}^\Gamma - p^0}^2 $. 
From the update rule \Cref{eqn:art_GVAMP0} in the first phase of the auxiliary GVAMP, we have
\begin{align}
    \ul{p}^\Gamma &= (\ol{\kappa}_2^{-1} XX^\top - I_n) \ol{G} \ul{p}^{\Gamma - 1}
    = X \ul{r}^\Gamma - \ol{G} \ul{p}^{\Gamma - 1} , \notag 
\end{align}
or equivalently, 
\begin{align}
    \ul{p}^\Gamma + \ol{G} \ul{p}^{\Gamma - 1} &= X \ul{r}^\Gamma . \label{eqn:update} 
\end{align}
We claim 
\begin{align}
    \lim_{\Gamma\to\infty} \lim_{n\to\infty} n^{-1} \normtwo{\ul{p}^\Gamma - \gamma \ul{p}^{\Gamma - 1}}^2 &= 0 . \label{eqn:close_pp} 
\end{align}
To show this, by the state evolution result in \Cref{lem:SE_aux}, almost surely, 
\begin{align}
    \lim_{n\to\infty} n^{-1} \normtwo{\ul{p}^\Gamma - \gamma \ul{p}^{\Gamma - 1}}^2
    &= \ul{a}_\Gamma^2 \sigma^2 + \expt{\ul{\sfJ}_\Gamma^2} + \gamma^2 \ul{a}_{\Gamma}^2 \sigma^2 + \gamma^2 \expt{\ul{\sfJ}_{\Gamma-1}^2}
    - 2 \gamma \paren{\ul{a}_\Gamma \ul{a}_{\Gamma - 1} \sigma^2 + \expt{\ul{\sfJ}_\Gamma \ul{\sfJ}_{\Gamma - 1}}} \notag \\
    &\explain{\Cref{eqn:art_first1}} 2 \gamma^2 \expt{\ul{\sfJ}_{\Gamma - 1}^2} - 2 \gamma \expt{\ul{\sfJ}_\Gamma \ul{\sfJ}_{\Gamma - 1}} . \label{eqn:pp_TODO} 
\end{align}
It was shown in \Cref{eqn:lim_nt} that 
\begin{align}
    \lim_{\Gamma\to\infty} \frac{\expt{\ul{\sfJ}_{\Gamma} \ul{\sfJ}_{\Gamma - 1}}}{\expt{\ul{\sfJ}_\Gamma^2}^{1/2} \expt{\ul{\sfJ}_{\Gamma - 1}^2}^{1/2}}
    &\explain{\Cref{eqn:art_first1}} \lim_{\Gamma\to\infty} \frac{\expt{\ul{\sfJ}_{\Gamma} \ul{\sfJ}_{\Gamma - 1}}}{\gamma \expt{\ul{\sfJ}_{\Gamma - 1}^2}}
    = 1 . \label{eqn:JJ} 
\end{align}
Note that \Cref{eqn:lim_nt} was established under the initialization scheme \Cref{eqn:lin_init}. 
Nevertheless, it is easily checked that \Cref{eqn:JJ} continues to hold under the initialization scheme \Cref{eqn:wtp0} which differs from \Cref{eqn:lin_init} by a multiplicative factor $ \gamma^{-\Gamma} w_3 $. 
Using \Cref{eqn:JJ} in \Cref{eqn:pp_TODO} allows us to conclude \Cref{eqn:close_pp}. 
Now combining \Cref{eqn:close_pp,eqn:update,eqn:close_v1} gives 
\begin{align}
    \lim_{\Gamma\to\infty} \lim_{n\to\infty} n^{-1} \normtwo{(I_n + \gamma^{-1} \ol{G}) \ul{p}^\Gamma - X v_1(D)}^2 &= 0 . \notag 
\end{align}
In view of the definition of $ p^0 $ in \Cref{eqn:spec_GVAMP_init}, this implies
\begin{align}
    \lim_{\Gamma\to\infty} \lim_{n\to\infty} n^{-1} \normtwo{c_p \, \ul{p}^\Gamma - p^0}^2 &= 0 . \label{eqn:close_p0} 
\end{align}
By the state evolution result in \Cref{lem:SE_aux},
\begin{align}
    \lim_{n\to\infty} n^{-1} \normtwo{\ul{p}^\Gamma}^2 &= \ul{a}_\Gamma^2 \sigma^2 + \expt{\ul{\sfJ}_\Gamma^2}
    \explain{\Cref{eqn:art_first1}} \frac{w_3^2 (\sigma^2 + w_1)}{\rho + w_2} < \infty . \label{eqn:finite_ulpGamma} 
\end{align}
Moreover, finiteness holds uniformly in $\Gamma$ and therefore also in the $\Gamma\to\infty$ limit. 
\Cref{eqn:close_p0,eqn:finite_ulpGamma} jointly imply $ \lim_{n\to\infty} n^{-1} \normtwo{p^0}^2 < \infty $. 
They also imply, by similar reasoning leading to \Cref{eqn:close_wtr,eqn:finite_wtr}, that
\begin{align}
&&
    \lim_{\Gamma\to\infty} \lim_{n\to\infty} n^{-1} \normtwo{\wt{\ul{p}}^\Gamma - \wt{p}^0}^2 &= 0 , & 
    \lim_{\Gamma\to\infty} \lim_{n\to\infty} n^{-1} \normtwo{\wt{\ul{p}}^\Gamma}^2 &< \infty . & 
& \label{eqn:finite_wtp}
\end{align}
As a consequence, we further have $ \lim_{n\to\infty} n^{-1} \normtwo{\wt{p}^0}^2 < \infty $. 
Furthermore, by \Cref{asmp:signal,asmp:design,asmp:noise}, 
\begin{align}
&&
    n^{-1} \normtwo{z}^2 \to \ol{\kappa}_2 \rho &< \infty , &
    n^{-1} \normtwo{\eps}^2 \to \expt{\sfE^2} &< \infty , & 
& \label{eqn:finite_ze}
\end{align}
almost surely as $n\to\infty$. 
Collecting \Cref{eqn:finite_wtp,eqn:finite_ulpGamma,eqn:close_p0,eqn:finite_ze}, returning to \Cref{eqn:phip}, we conclude \Cref{eqn:art_p} for $\tau = 0$. 

\paragraph{\textsc{Induction step.}}
Now fix any $t\ge1$. 
We put the following statements as induction hypotheses: \Cref{eqn:art_ab} for every $ 0\le\tau\le t-1 $, \Cref{eqn:art_abtilde} for every $ 0\le\tau\le t $, \Cref{eqn:art_J,eqn:art_K} for every $ 1\le\varsigma,\tau\le t $, \Cref{eqn:art_r,eqn:art_p} for every $ 0\le\tau\le t $. 
The remainder of the proof amounts to verifying the analogous statements for $ t+1 $. 

First consider $ \ul{a}_{\Gamma+t+1} $: 
\begin{align}
    \ol{\kappa}_2 \abs{\ul{a}_{\Gamma+t+1} - a_{t+1}}
    &\explain{\Cref{eqn:ab}} \abs{\expt{\ul{\Psi}_{\Gamma+t}(\sfLambda_n^2) \sfLambda_n^2} \wt{\ul{a}}_{\Gamma+t} - \expt{\Psi_t(\sfLambda_n^2) \sfLambda_n^2} \wt{a}_t
    + \expt{\wt{\ul{\Psi}}_{\Gamma+t}(\sfLambda_n) \sfLambda_n} \wt{\ul{b}}_{\Gamma+t} - \expt{\wt{\Psi}_t(\sfLambda_n) \sfLambda_n} \wt{b}_t} \notag \\
    &\le \abs{\expt{\Psi_t(\sfLambda_n^2) \sfLambda_n^2}} \cdot \abs{\wt{\ul{a}}_{\Gamma+t} - \wt{a}_t}
    + \abs{\expt{\wt{\Psi}_t(\sfLambda_n) \sfLambda_n}} \cdot \abs{\wt{\ul{b}}_{\Gamma+t} - \wt{a}_t}
    \to 0 , \label{eqn:art_at+1} 
\end{align}
as $\Gamma\to\infty$, where the last line uses 
the choice of $ \ul{\Psi}_{\Gamma+t}, \wt{\ul{\Psi}}_{\Gamma+t} $ in \Cref{eqn:denoiser2}, 
the triangle inequality,
and the induction hypotheses for $ \wt{\ul{a}}_{\Gamma+t}, \wt{\ul{b}}_{\Gamma+t} $ in \Cref{eqn:art_abtilde}. 
A similar argument proves $ \abs{\ul{b}_{\Gamma+t+1} - b_{t+1}} \to 0 $ as $\Gamma\to\infty$, thereby establishing \Cref{eqn:art_ab} for $t+1$. 

Then, for any $ 0\le \tau\le t $, we consider $ \expt{\ul{\sfJ}_{\Gamma+t+1} \ul{\sfJ}_{\Gamma+\tau+1}} $: 
\begin{align}
    & \abs{\expt{\ul{\sfJ}_{\Gamma+t+1} \ul{\sfJ}_{\Gamma+\tau+1}} - \expt{\sfJ_{t+1} \sfJ_{\tau+1}}} \label{eqn:art_Jt+1} \\
    &\overset{\Cref{eqn:J_cov}}{\le} \abs{\Psi_t(\sfLambda_n^2) \Psi_\tau(\sfLambda_n^2)} \paren{\abs{\expt{\wt{\ul{\sfP}}_{\Gamma+t} \wt{\ul{\sfP}}_{\Gamma+\tau}} - \expt{\wt{\sfP}_t \wt{\sfP}_\tau}} + \sigma^2 \abs{\wt{\ul{a}}_{\Gamma+t} \wt{\ul{a}}_{\Gamma+\tau} - \wt{a}_t \wt{a}_\tau}} \notag \\
    &\quad + \abs{\expt{\Psi_t(\sfLambda_n^2) \Psi_\tau(\sfLambda_n^2) \sfLambda_n^2}} \cdot \rho \cdot \abs{\wt{\ul{a}}_{\Gamma+t} \wt{\ul{a}}_{\Gamma+\tau} - \wt{a}_t \wt{a}_\tau} \notag \\
    &\quad + \abs{\expt{\wt{\Psi}_t(\sfLambda_n) \wt{\Psi}_\tau(\sfLambda_n)}} \cdot \abs{\expt{\wt{\ul{\sfR}}_{\Gamma+t} \wt{\ul{\sfR}}_{\Gamma+\tau}} - \expt{\wt{\sfR}_t \wt{\sfR}_\tau}}
    + \sigma^2 \abs{\ul{a}_{\Gamma+t+1} \ul{a}_{\Gamma+\tau+1} - a_{t+1} a_{\tau+1}} \notag \\
    &\quad + \rho \abs{\Psi_t(\sfLambda_n^2) \wt{\Psi}_\tau(\sfLambda) \sfLambda_n} \abs{\wt{\ul{a}}_{\Gamma+t} \wt{\ul{b}}_{\Gamma+\tau} - \wt{a}_t \wt{b}_\tau}
    + \rho \abs{\Psi_\tau(\sfLambda_n^2) \wt{\Psi}_t(\sfLambda) \sfLambda_n} \abs{\wt{\ul{a}}_{\Gamma+\tau} \wt{\ul{b}}_{\Gamma+t} - \wt{a}_\tau \wt{b}_t} . \notag 
\end{align}
Recalling the definition \Cref{eqn:denoiser2_fg} and using the induction hypotheses \Cref{eqn:art_ab,eqn:art_J}, 
\begin{align}
    \abs{\expt{\wt{\ul{\sfP}}_{\Gamma+t} \wt{\ul{\sfP}}_{\Gamma+\tau}} - \expt{\wt{\sfP}_t \wt{\sfP}_\tau}} 
    &= \abs{\expt{g_t(\ul{a}_{\Gamma+t} \sfZ + \ul{\sfJ}_{\Gamma+t}; \sfY) g_\tau(\ul{a}_{\Gamma+\tau} \sfZ + \ul{\sfJ}_{\Gamma+\tau}; \sfY)} - \expt{g_t(a_t \sfZ + \sfJ_t; \sfY) g_\tau(a_\tau \sfZ + \sfJ_\tau; \sfY)}}
    \to 0 , \notag
\end{align}
as $\Gamma\to\infty$. 
Similarly, $ \expt{\wt{\ul{\sfR}}_{\Gamma+t} \wt{\ul{\sfR}}_{\Gamma+\tau}} - \expt{\wt{\sfR}_t \wt{\sfR}_\tau} \xrightarrow[\Gamma\to\infty]{} 0 $. 
Using the above observations and the induction hypotheses \Cref{eqn:art_ab,eqn:art_abtilde,eqn:art_J}, we have that \Cref{eqn:art_Jt+1} converges to $0$ as $\Gamma\to\infty$. 
A similar argument can be used to show $ \abs{\expt{\ul{\sfJ}_{\Gamma} \ul{\sfJ}_{\Gamma+t+1}} - c_p^{-1} \expt{\sfJ_0 \sfJ_{t+1}}} \xrightarrow[\Gamma\to\infty]{} 0 $, thereby verifying \Cref{eqn:art_J} for $t+1$. 
The verification o \Cref{eqn:art_K} for $t+1$ is analogous and omitted. 

Next, consider $ \wt{\ul{a}}_{\Gamma+t+1} $: 
\begin{align}
    \wt{\ul{a}}_{\Gamma+t+1} &\explain{\Cref{eqn:ab_tilde}} \sigma^{-2} \expt{\ul{g}_{\Gamma+t+1}(\ul{\sfP}_{\Gamma+t+1}; \sfY) \sfZ} 
    = \sigma^{-2} \expt{g_{t+1}(\ul{a}_{\Gamma+t+1} \sfZ + \ul{\sfJ}_{\Gamma+t+1}; \sfY)} \notag \\
    &\xrightarrow[\Gamma\to\infty]{} \sigma^{-2} \expt{g_{t+1}(a_{t+1} \sfZ + \sfJ_{t+1}; \sfY)} 
    = \wt{a}_{t+1} , \notag 
\end{align}
where the convergence follows from \Cref{eqn:art_at+1,eqn:art_Jt+1} shown above. 
A similar argument shows $ \wt{\ul{b}}_{\Gamma+t+1} \to \wt{\ul{b}}_{t+1} $ as $\Gamma\to\infty$, thereby justifying \Cref{eqn:art_abtilde} for $t+1$. 

Finally, we aim to show \Cref{eqn:art_p} for $t+1$. 
Denote
\begin{align}
\begin{split}
    \ul{P}^{\Gamma+t+1} &\coloneqq \matrix{c_p \ul{p}^\Gamma & \ul{p}^{\Gamma+1} & \cdots & \ul{p}^{\Gamma+t+1} & \wt{\ul{p}}^\Gamma & \cdots & \wt{\ul{p}}^{\Gamma+t+1} & z & \eps} \in \bbR^{d\times 2(t+3)} , \\
    P^{t+1} &\coloneqq \matrix{p^0 & p^1 & \cdots & p^{t+1} & \wt{p}^0 & \cdots & \wt{p}^{t+1} & z & \eps} \in \bbR^{d\times 2(t+3)} . 
\end{split}
\label{eqn:PP}
\end{align}
Writing $ \phi \equiv \phi_p $ for simplicity, using the pseudo-Lipschitz property of $\phi$ and applying Cauchy--Schwarz, we bound the left-hand side of \Cref{eqn:art_p} (for $t+1$) as follows: 
\begin{align}
    \abs{\frac{1}{n} \sum_{i = 1}^n \phi(\ul{P}^{\Gamma+t+1}_{i,:}) - \frac{1}{n} \sum_{i = 1}^n \phi(P^{t+1}_{i,:})}
    &\le C \sqrt{\frac{3}{n} \sum_{i = 1}^n \paren{1 + \normtwo{\ul{P}^{\Gamma+t+1}_{i,:}}^2 + \normtwo{P^{t+1}_{i,:}}^2}} \sqrt{\frac{1}{n} \sum_{i = 1}^n \normtwo{\ul{P}^{\Gamma+t+1}_{i,:} - P^{t+1}_{i,:}}^2} . \label{eqn:finite_TODO} 
\end{align}
By definition \Cref{eqn:PP}, 
\begin{align}
    \frac{1}{n} \sum_{i = 1}^n \normtwo{\ul{P}^{\Gamma+t+1}_{i,:} - P^{t+1}_{i,:}}^2
    &= \frac{1}{n} \paren{\normtwo{c_p \ul{p}^\Gamma - p^0}^2 + \sum_{i = 1}^{t+1} \normtwo{\ul{p}^{\Gamma+i} - p^i}^2 + \sum_{i = 0}^{t+1} \normtwo{\wt{\ul{p}}^{\Gamma+i} - \wt{p}^{i}}^2} . \label{eqn:diff} 
\end{align}
The induction hypothesis \Cref{eqn:art_p} implies that all terms on the right above converge to $0$ as in the limit $n\to\infty$ followed by $\Gamma\to\infty$, except two terms $ n^{-1} \normtwo{\ul{p}^{\Gamma+t+1} - p^{t+1}}^2 $, $ n^{-1} \normtwo{\wt{\ul{p}}^{\Gamma+t+1} - \wt{p}^{t+1}}^2 $ which we separately examine below. 
According to the update rules \Cref{eqn:art_GVAMP2,eqn:GVAMP2}, and the choice of denoisers in \Cref{eqn:denoiser2,eqn:denoiser2_fg}, 
\begin{align}
    \frac{1}{\sqrt{n}} \normtwo{\ul{p}^{\Gamma+t+1} - p^{t+1}}
    &\le \frac{1}{\sqrt{n}} \normtwo{\Psi_{t}(XX^\top)} \normtwo{g_t(\ul{p}^{\Gamma+t}; y) - g_t(p^t; y)} + \frac{1}{\sqrt{n}} \normtwo{\wt{\Psi}_t(X)} \normtwo{f_t(\ul{r}^{\Gamma+t}) - f_t(r^t)} \notag \\
    &\le \frac{C}{\sqrt{n}} \normtwo{\Psi_{t}(XX^\top)} \normtwo{\ul{p}^{\Gamma+t} - p^t} + \frac{C}{\sqrt{n}} \normtwo{\wt{\Psi}_t(X)} \normtwo{\ul{r}^{\Gamma+t} - r^t} , \notag
\end{align}
where the triangle inequality and the Lipschitz property of $ f_t, g_t $ (see \Cref{asmp:remove}) have been used. 
Then \Cref{asmp:remove_PhiPsi} and the induction hypotheses \Cref{eqn:art_r,eqn:art_p} ensure that 
\begin{align}
    \lim_{\Gamma\to\infty} \lim_{n\to\infty} n^{-1} \normtwo{\ul{p}^{\Gamma+t+1} - p^{t+1}}^2 &= 0 . \label{eqn:close_pt+1}
\end{align}
By the choice of the denoisers \Cref{eqn:denoiser2_fg} and the Lipschitz property of $ f_{t+1}, g_{t+1} $, this immediately implies
\begin{align}
    \frac{1}{\sqrt{n}} \normtwo{\wt{\ul{p}}^{\Gamma+t+1} - \wt{p}^{t+1}} &\le \frac{C}{\sqrt{n}} \normtwo{\ul{p}^{\Gamma+t+1} - p^{t+1}} \to 0 , \label{eqn:close_wtpt+1} 
\end{align}
as $ n\to\infty, \Gamma\to\infty $ in that order. 
Consequently, \Cref{eqn:diff} converges to $0$ in the same sequential limit. 
We then show that $ \lim_{\Gamma\to\infty} \limsup_{n\to\infty} n^{-1} \sum_{i = 1}^n \normtwo{\ul{P}^{\Gamma+t+1}_{i,:}}^2 $ is finite. 
To this end, 
\begin{align}
    \frac{1}{n} \sum_{i = 1}^n \normtwo{\ul{P}^{\Gamma+t+1}_{i,:}}^2
    &= \frac{1}{n} \paren{ c_p^2 \normtwo{\ul{p}^\Gamma}^2 + \sum_{i = 1}^{t+1} \normtwo{\ul{p}^{\Gamma+i}}^2 + \sum_{i = 0}^{t+1} \normtwo{\wt{\ul{p}}^{\Gamma+i}}^2 + \normtwo{z}^2 + \normtwo{\eps}^2 } . \label{eqn:finite_P} 
\end{align}
By the state evolution result in \Cref{lem:SE_aux} for the auxiliary GVAMP iteration \Cref{eqn:art_GVAMP}, for any $\Gamma\ge0$ fixed and any $ 0\le i\le t+1 $, 
\begin{align}
&&
    \frac{1}{n} \normtwo{\ul{p}^{\Gamma+i}}^2 &\to \expt{\ul{\sfP}_{\Gamma+i}^2} < \infty , & 
    \frac{1}{n} \normtwo{\wt{\ul{p}}^{\Gamma+i}}^2 &\to \expt{\wt{\ul{\sfP}}_{\Gamma+i}^2} < \infty , & 
& \notag 
\end{align}
almost surely, as $n\to\infty$. 
By the induction hypotheses \Cref{eqn:art_ab,eqn:art_J}, 
\begin{align}
&&
    \lim_{\Gamma\to\infty} \lim_{n\to\infty} \frac{1}{n} \normtwo{\ul{p}^{\Gamma+i}}^2 &= \expt{\sfP_{i}^2} < \infty , & 
    \lim_{\Gamma\to\infty} \lim_{n\to\infty} \frac{1}{n} \normtwo{\wt{\ul{p}}^{\Gamma+i}}^2 &= \expt{\wt{\sfP}_{i}^2} < \infty . & 
& \label{eqn:finite_ulp} 
\end{align}
Combining this with the induction hypothesis \Cref{eqn:art_p} along with \Cref{eqn:close_wtpt+1,eqn:close_pt+1} shown above, and applying the triangle inequality, we also obtain the finiteness of 
\begin{align}
&& 
    \lim_{\Gamma\to\infty} \lim_{n\to\infty} \frac{1}{n} \normtwo{p^i}^2 & < \infty , & 
    \lim_{\Gamma\to\infty} \lim_{n\to\infty} \frac{1}{n} \normtwo{\wt{p}^i}^2 & < \infty , &
& \label{eqn:finite_p}
\end{align}
for every $ 0\le i\le t+1 $. 
Taking \Cref{eqn:finite_p,eqn:finite_ulp,eqn:finite_ze} collectively, we see that both $ n^{-1} \sum_{i = 1}^n \normtwo{\ul{P}_{i,:}^{\Gamma+t+1}}^2 $ and $ n^{-1} \sum_{i = 1}^n \normtwo{P_{i,:}^{t+1}}^2 $ are bounded in the double limit. 
Finally, combining this with the fact that \Cref{eqn:diff} converges to $0$ already shown, we have that the left-hand side of \Cref{eqn:finite_TODO} also converges to $0$. 
This therefore proves \Cref{eqn:art_r} for $t+1$. 
The proof for \Cref{eqn:art_p} is syntactically symmetric and omitted. 
This closes the induction and the desired statements in \Cref{lem:art_close_to_true} follows. 
\end{proof}

\section{Proof of \Cref{thm:SE_Bayes_GVAMP}}
\label{app:pf_thm:SE_Bayes_GVAMP}

We show that the canonical form of Bayes-GVAMP in \Cref{eqn:BGVAMP} is indeed a special case of the generic GVAMP in \Cref{eqn:GVAMP} by verifying the trace-free and divergence-free conditions in \Cref{asmp:tr_free,asmp:div_free}. 

\begin{lemma}
\label{lem:BGVAMP_free}
The functions $ \Phi_t, \Psi_t $ in \Cref{eqn:B_Phi,eqn:B_Psi} are trace-free in the sense of \Cref{asmp:tr_free} and the functions $ f_t, g_t $ are divergence-free in the sense of \Cref{asmp:div_free}. 
\end{lemma}

\begin{proof}
By the definitions of $ v_2^t, c_2^t $ in \Cref{eqn:B_v2_c2}, 
\begin{align}
    \expt{\Phi_t(\sfLambda_d^2)} &= \frac{1}{v_2^t} \expt{\frac{1}{\gamma_2^t + \tau_2^t \sfLambda_d^2}} - 1
    = \frac{v_2^t}{v_2^t} - 1 = 0 , \notag \\
    \expt{\Psi_t(\sfLambda_n^2)} &= \frac{1}{c_2^t} \expt{\frac{\sfLambda_n^2}{\gamma_2^t + \tau_2^t \sfLambda_n^2}} - 1 = \frac{c_2^t}{c_2^t} - 1 = 0 . \notag 
\end{align}
By the definitions of $ v_ 1^t, c_1^t $ in \Cref{eqn:B_v1_c1}, 
\begin{align}
    \expt{f_t'(\sfR_t)} &= \frac{1}{v_1^t} \expt{g_{x1}'(\sfR_t; \gamma_1^t)} - 1 = \frac{v_1^t}{v_1^t} - 1 = 0 , \notag \\
    \expt{g_t'(\sfP_t; \sfY)} &= \frac{1}{c_1^t} \expt{g_{z1}'(\sfP_t, \sfY; \tau_1^t)} - 1 = \frac{c_1^t}{c_1^t} - 1 = 0 . \notag 
\end{align}
This completes the proof. 
\end{proof}

Then \Cref{thm:SE_Bayes_GVAMP} follows immediately from the lemma below which states that the scalars $ \gamma_1^t, \tau_1^t, v_1^t, c_1^t, \gamma_2^t, \tau_2^t, v_2^t, c_2^t $ are in one-to-one correspondence with the state evolution parameters $ \wt{a}_t, \wt{b}_t, a_t, b_t, \expt{\sfJ_t^2}, \expt{\sfK_t^2} $. 

\begin{lemma}
\label{lem:SE_BGVAMP}
Let $ \gamma_1^t, \tau_1^t, \gamma_2^t, \tau_2^t $ be defined via \Cref{eqn:Bayes_GVAMP} and let $ \wt{a}_t, \wt{b}_t, a_t, b_t, \expt{\sfJ_t^2}, \expt{\sfK_t^2} $ be defined via \Cref{eqn:SE_GVAMP} with $ \Phi_t, \wt{\Phi}_t, \Psi_t, \wt{\Psi}_t, f_t, g_t $ given in \Cref{eqn:B_Phi,eqn:B_Psi,eqn:B_fg}. 
Then these quantities satisfy
\begin{align}
&&
    \wt{a}_t &= \tau_2^t , & 
    \wt{b}_t &= \gamma_2^t - \frac{1}{\rho} , & 
    a_{t+1} &= \expt{\sfJ_{t+1}^2} = \tau_1^{t+1} - \frac{1}{\sigma^2} , & 
    b_{t+1} &= \expt{\sfK_{t+1}^2} = \gamma_1^{t+1} . & 
& \label{eqn:SE_relation} 
\end{align}
\end{lemma}

\begin{proof}
The proof is by induction on $t$. 
In preparation for performing the induction, we first make a few observations that hold for every $t$. 

\paragraph{\textsc{Preparation.}}
By the definitions \Cref{eqn:B_Psi} of $ \Psi_t, \wt{\Psi}_t $, it holds 
\begin{align}
    \Psi_t(x^2) &= x \wt{\Psi}_t(x) - 1 . \label{eqn:Psi_wtPsi} 
\end{align}
The trace-free property of $ \Psi_t $ (see \Cref{lem:BGVAMP_free}) then implies 
\begin{align}
    \expt{ \sfLambda_n \wt{\Psi}_t(\sfLambda_n) } &= 1 . \notag 
\end{align}
We use these observations to simplify the expressions of a few expectations involving $ \Psi_t, \wt{\Psi}_t $: 
\begin{subequations}
\label{eqn:expt_Psi0}
\begin{align}
    \expt{ \Psi_t(\sfLambda_n^2)^2 } &= \expt{ \paren{ \sfLambda_n \wt{\Psi}_t(\sfLambda_n) - 1 }^2 } = \expt{ \sfLambda_n^2 \wt{\Psi}_t(\sfLambda_n)^2 } - 1 , \\
    \expt{ \Psi_t(\sfLambda_n^2)^2 \sfLambda_n^2 }
    &= \expt{ \paren{ \sfLambda_n \wt{\Psi}_t(\sfLambda_n) - 1 }^2 \sfLambda_n^2 }
    = \expt{ \sfLambda_n^4 \wt{\Psi}_t(\sfLambda_n)^2 } + \expt{ \sfLambda_n^2 } - 2 \expt{ \sfLambda_n^3 \wt{\Psi}_t(\sfLambda_n) } , \\
    \expt{ \Psi_t(\sfLambda_n^2) \wt{\Psi}_t(\sfLambda_n) \sfLambda_n } &= \expt{ \paren{ \sfLambda_n \wt{\Psi}_t(\sfLambda_n) - 1 } \wt{\Psi}_t(\sfLambda_n) \sfLambda_n }
    = \expt{ \sfLambda_n^2 \wt{\Psi}_t(\sfLambda_n)^2 } - 1 . 
\end{align}
\end{subequations}
Plugging in the formula of $ \wt{\Psi}_t $ in \Cref{eqn:B_Psi}, we write some terms on the right-hand side, together with $ \expt{\wt{\Psi}_t(\sfLambda_n^2)} $, more explicitly as
\begin{subequations}
\label{eqn:expt_Psi}
\begin{align}
&& 
    \expt{\wt{\Psi}_t(\sfLambda_n)^2} &= \frac{1}{(c_2^t)^2} \expt{ \frac{\sfLambda_n^2}{(\gamma_2^t + \tau_2^t \sfLambda_n^2)^2} } , & 
    \expt{ \sfLambda_n^2 \wt{\Psi}_t(\sfLambda_n)^2 }
    &= \frac{1}{(c_2^t)^2} \expt{ \frac{\sfLambda_n^4}{(\gamma_2^t + \tau_2^t \sfLambda_n^2)^2} } , & 
& \\
&&
    \expt{\sfLambda_n^4 \wt{\Psi}_t(\sfLambda_n)^2} &= \frac{1}{(c_2^t)^2} \expt{ \frac{\sfLambda_n^6}{(\gamma_2^t + \tau_2^t \sfLambda_n^2)^2} } , &
    \expt{\sfLambda_n^3 \wt{\Psi}_t(\sfLambda_n)} &= \frac{1}{c_2^t} \expt{ \frac{\sfLambda_n^4}{\gamma_2^t + \tau_2^t \sfLambda_n^2} } , &
&  
\end{align}
\end{subequations}
We then express the expectations on the right-hand side above in terms of $ c_2^t $ and $ \expt{ \frac{\sfLambda_n^4}{(\gamma_2^t + \tau_2^t \sfLambda_n^2)^2}} $: 
\begin{align}
    \expt{ \frac{\sfLambda_n^4}{\gamma_2^t + \tau_2^t \sfLambda_n^2} } 
    &= \frac{1}{\tau_2^t} \expt{\sfLambda_n^2 - \gamma_2^t \frac{\sfLambda_n^2}{\gamma_2^t + \tau_2^t \sfLambda_n^2}} 
    = \frac{\ol{\kappa}_2}{\tau_2^t} - \frac{\gamma_2^t c_2^t}{\tau_2^t} , \notag \\
    \expt{ \frac{\sfLambda_n^6}{(\gamma_2^t + \tau_2^t \sfLambda_n^2)^2} }
    &= \frac{1}{\tau_2^t} \expt{ \frac{\sfLambda_n^4}{\gamma_2^t + \tau_2^t \sfLambda_n^2} } - \frac{\gamma_2^t}{\tau_2^t} \expt{ \frac{\sfLambda_n^4}{(\gamma_2^t + \tau_2^t \sfLambda_n^2)^2} } 
    = \frac{\ol{\kappa}_2}{(\tau_2^t)^2} - \frac{\gamma_2^t c_2^t}{(\tau_2^t)^2} - \frac{\gamma_2^t}{\tau_2^t} \expt{ \frac{\sfLambda_n^4}{(\gamma_2^t + \tau_2^t \sfLambda_n^2)^2} } , \notag \\
    \expt{ \frac{\sfLambda_n^2}{(\gamma_2^t + \tau_2^t \sfLambda_n^2)^2} }
    &= \frac{1}{\gamma_2^t} \expt{ \frac{\sfLambda_n^2}{\gamma_2^t + \tau_2^t \sfLambda_n^2} - \tau_2^t \frac{\sfLambda_n^4}{(\gamma_2^t + \tau_2^t \sfLambda_n^2)^2} } 
    = \frac{c_2^t}{\gamma_2^t} - \frac{\tau_2^t}{\gamma_2^t} \expt{ \frac{\sfLambda_n^4}{(\gamma_2^t + \tau_2^t \sfLambda_n^2)^2} } . \label{eqn:22} 
\end{align}
Putting these back in \Cref{eqn:expt_Psi}, we arrive at 
\begin{subequations}
\label{eqn:expect_Psi_explicit}
\begin{align}
    \expt{\wt{\Psi}_t(\sfLambda_n)^2} &= \frac{1}{\gamma_2^t c_2^t} - \frac{\tau_2^t}{\gamma_2^t (c_2^t)^2} \expt{ \frac{\sfLambda_n^4}{(\gamma_2^t + \tau_2^t \sfLambda_n^2)^2} } , \\
    \expt{ \sfLambda_n^2 \wt{\Psi}_t(\sfLambda_n)^2 }
    &= \frac{1}{(c_2^t)^2} \expt{ \frac{\sfLambda_n^4}{(\gamma_2^t + \tau_2^t \sfLambda_n^2)^2} } , \\
    \expt{\sfLambda_n^4 \wt{\Psi}_t(\sfLambda_n)^2} 
    &= \frac{\ol{\kappa}_2}{(\tau_2^t c_2^t)^2} - \frac{\gamma_2^t}{(\tau_2^t)^2 c_2^t} - \frac{\gamma_2^t}{\tau_2^t (c_2^t)^2} \expt{ \frac{\sfLambda_n^4}{(\gamma_2^t + \tau_2^t \sfLambda_n^2)^2} } , \\
    \expt{\sfLambda_n^3 \wt{\Psi}_t(\sfLambda_n)} 
    &= \frac{\ol{\kappa}_2}{\tau_2^t c_2^t} - \frac{\gamma_2^t}{\tau_2^t} . 
\end{align}
\end{subequations}

We perform analogous manipulations to several quantities involving $ \Phi_t, \wt{\Phi}_t $. 
First note that by the definitions of both functions in \Cref{eqn:B_Phi} and the definition of $ \gamma_1^{t+1} $ in \Cref{eqn:B_gamma1_tau1}, the following identity holds
\begin{align}
    \gamma_2^t \Phi_t(x^2) + \tau_2^t x \wt{\Phi}_t(x)
    &= \frac{1}{v_2^t} - \gamma_2^t = \gamma_1^{t+1} . \label{eqn:Phi_wtPhi} 
\end{align}
By rearranging terms and using the trace-free property of $ \Phi_t $ (see \Cref{lem:BGVAMP_free}), this implies
\begin{align}
&& 
    \Phi_t(x^2) &= \frac{1}{\gamma_2^t} \paren{ \gamma_1^{t+1} - \tau_2^t x \wt{\Phi}_t(x) } , & 
    \expt{ \sfLambda_d \wt{\Phi}_t(\sfLambda_d) } &= \frac{\gamma_1^{t+1}}{\tau_2^t} . & 
& \notag 
\end{align}
We use the above observations to compute a few expectations involving these $ \Phi_t, \wt{\Phi}_t $: 
\begin{subequations}
\label{eqn:expect_Phi0}
\begin{align}
    \expt{ \Phi_t(\sfLambda_d^2)^2 } &= \expt{ \paren{ \frac{1}{\gamma_2^t} \paren{ \gamma_1^{t+1} - \tau_2^t \sfLambda_d \wt{\Phi}_t(\sfLambda_d) } }^2 } \notag \\
    &= \frac{(\gamma_1^{t+1})^2}{(\gamma_2^t)^2} + \frac{(\tau_2^t)^2}{(\gamma_2^t)^2} \expt{ \sfLambda_d^2 \wt{\Phi}_t(\sfLambda_d)^2 } - 2 \frac{\gamma_1^{t+1} \tau_2^t}{(\gamma_2^t)^2} \cdot \frac{\gamma_1^{t+1}}{\tau_2^t} 
    = \frac{(\tau_2^t)^2}{(\gamma_2^t)^2} \expt{ \wt{\Phi}_t(\sfLambda_d)^2 \sfLambda_d^2 } - \frac{(\gamma_1^{t+1})^2}{(\gamma_2^t)^2} , \\
    \expt{ \Phi_t(\sfLambda_d^2) \wt{\Phi}_t(\sfLambda_d) \sfLambda_d } &= \expt{ \frac{1}{\gamma_2^t} \paren{ \gamma_1^{t+1} - \tau_2^t \sfLambda_d \wt{\Phi}_t(\sfLambda_d) } \wt{\Phi}_t(\sfLambda_d) \sfLambda_d } \notag \\
    &= \frac{(\gamma_1^{t+1})^2}{\gamma_2^t \tau_2^t} - \frac{\tau_2^t}{\gamma_2^t} \expt{ \wt{\Phi}_t(\sfLambda_d)^2 \sfLambda_d^2 } . 
\end{align}
\end{subequations}
We then derive more explicit expressions for $ \expt{ \wt{\Phi}_t(\sfLambda_d)^2 \sfLambda_d^2 } $ that shows up on the right-hand side and $ \expt{\wt{\Phi}_t(\sfLambda_d)^2} $. 
By the definition \Cref{eqn:B_Phi} of $ \wt{\Phi}_t $, 
\begin{align}
&&
    \expt{ \wt{\Phi}_t(\sfLambda_d)^2 } &= \frac{1}{(v_2^t)^2} \expt{ \frac{\sfLambda_d^2}{(\gamma_2^t + \tau_2^t \sfLambda_d^2)^2} } , & 
    \expt{ \wt{\Phi}_t(\sfLambda_d)^2 \sfLambda_d^2 } &= \frac{1}{(v_2^t)^2} \expt{ \frac{\sfLambda_d^4}{(\gamma_2^t + \tau_2^t \sfLambda_d^2)^2} } , &
& \notag 
\end{align}
By elementary algebra, the expectations of the fractions on the right-hand side can be written in terms of $ v_2^t $ and $ \expt{\frac{1}{(\gamma_2^t + \tau_2^t \sfLambda_d^2)^2}} $. 
Specifically, we have
\begin{align}
    \expt{ \frac{\sfLambda_d^2}{(\gamma_2^t + \tau_2^t \sfLambda_d^2)^2} }
    &= \frac{1}{\tau_2^t} \expt{ \frac{1}{\gamma_2^t + \tau_2^t \sfLambda_d^2} - \gamma_2^t \frac{1}{(\gamma_2^t + \tau_2^t \sfLambda_d^2)^2} }
    = \frac{v_2^t}{\tau_2^t} - \frac{\gamma_2^t}{\tau_2^t} \expt{ \frac{1}{(\gamma_2^t + \tau_2^t \sfLambda_d^2)^2} } , \label{eqn:Lamd12} 
\end{align}
and 
\begin{align}
    (\gamma_2^t)^2 \expt{ \frac{1}{(\gamma_2^t + \tau_2^t \sfLambda_d)^2} } 
    &= \expt{ \paren{1 - \frac{\tau_2^t \sfLambda_d^2}{\gamma_2^t + \tau_2^t \sfLambda_d^2}}^2 }
    = 1 + (\tau_2^t)^2 \expt{ \frac{\sfLambda_d^4}{(\gamma_2^t + \tau_2^t \sfLambda_d^2)^2} } - 2 \expt{ \frac{\tau_2^t \sfLambda_d^2}{\gamma_2^t + \tau_2^t \sfLambda_d^2} } \notag \\
    &= 1 + (\tau_2^t)^2 \expt{ \frac{\sfLambda_d^4}{(\gamma_2^t + \tau_2^t \sfLambda_d^2)^2} } - 2 \paren{ 1 - \gamma_2^t \expt{ \frac{1}{\gamma_2^t + \tau_2^t \sfLambda_d^2} } } \notag \\
    &= (\tau_2^t)^2 \expt{ \frac{\sfLambda_d^4}{(\gamma_2^t + \tau_2^t \sfLambda_d^2)^2} } + 2 \gamma_2^t v_2^t - 1 , \notag 
\end{align}
which implies
\begin{align}
    \expt{ \frac{\sfLambda_d^4}{(\gamma_2^t + \tau_2^t \sfLambda_d^2)^2} } 
    &= \frac{(\gamma_2^t)^2}{(\tau_2^t)^2} \expt{ \frac{1}{(\gamma_2^t + \tau_2^t \sfLambda_d)^2} }  - \frac{2 \gamma_2^t v_2^t}{(\tau_2^t)^2} + \frac{1}{(\tau_2^t)^2} . \label{eqn:Lamd22} 
\end{align}
Dividing both terms in \Cref{eqn:Lamd12,eqn:Lamd22} by $ (v_2^t)^2 $ brings us to the following expressions of $ \expt{\wt{\Phi}_t(\sfLambda_d)^2}, \expt{\wt{\Phi}_t(\sfLambda_d)^2 \sfLambda_d^2} $: 
\begin{subequations}
\label{eqn:expect_Phi}
\begin{align}
    \expt{\wt{\Phi}_t(\sfLambda_d)^2} &= \frac{1}{\tau_2^t v_2^t} - \frac{\gamma_2^t}{\tau_2^t (v_2^t)^2} \expt{ \frac{1}{(\gamma_2^t + \tau_2^t \sfLambda_d^2)^2} } , \notag \\
    \expt{ \wt{\Phi}_t(\sfLambda_d)^2 \sfLambda_d^2 } &= \frac{(\gamma_2^t)^2}{(\tau_2^t v_2^t)^2} \expt{ \frac{1}{(\gamma_2^t + \tau_2^t \sfLambda_d^2)^2} }  - \frac{2 \gamma_2^t}{(\tau_2^t)^2 v_2^t} + \frac{1}{(\tau_2^t v_2^t)^2} . \notag 
\end{align}
\end{subequations}

We also make an observation regarding $ v_1^t, c_1^t $. 
By their definitions in \Cref{eqn:B_v1_c1}, the expressions of $ g_{x1}', g_{z1}' $ in \Cref{eqn:g'} and the tower property of conditional expectation, we have 
\begin{align}
    v_1^t &= \expt{\var{\sfB_* \mid \sfR_t}} = \expt{\expt{\sfB_*^2 \mid \sfR_t} - \expt{\sfB_* \mid \sfR_t}^2} = \rho - \expt{\sfB_* \expt{\sfB_* \mid \sfR_t}} , \label{eqn:v1_alt} \\
    c_1^t &= \expt{\var{\sfZ \mid \sfP_t, \sfY}} = \expt{\expt{\sfZ^2 \mid \sfP_t, \sfY} - \expt{\sfZ \mid \sfP_t, \sfY}^2} = \sigma^2 - \expt{\sfZ \expt{\sfZ \mid \sfP_t, \sfY}} , \label{eqn:c1_alt}
\end{align}
implying 
\begin{align}
&&
    \expt{\sfB_* \expt{\sfB_* \mid \sfR_t}} &= \rho - v_1^t , &
    \expt{\sfZ \expt{\sfZ \mid \sfP_t, \sfY}} &= \sigma^2 - c_1^t . & 
& \label{eqn:cond_expt}
\end{align}

We are now ready to present the induction proof. 
We first perform the induction step since the verification of the base case $t=0$ follows similar arguments in which the initial condition \Cref{eqn:RP0_BGVAMP} instead of the induction hypotheses is used. 

\paragraph{\textsc{Induction step.}}
Now assuming the validity of the conclusion \Cref{eqn:SE_relation} for time up to $t-1$, we verify its validity for time $t$. 

Let us start with $ \wt{a}_t $:
\begin{align}
    \sigma^2 \wt{a}_t &\explain{\Cref{eqn:ab_tilde}} \expt{\wt{\sfP}_t \sfZ}
    \explain{\Cref{eqn:wtP_wtR}} \expt{ g_t(\sfP_t; \sfY) \sfZ }
    \explain{\Cref{eqn:B_fg}} \frac{1}{c_1^t} \expt{ g_{z1}(\sfP_t, \sfY; \tau_1^t) \sfZ } - \expt{\sfP_t \sfZ}
    \explain{\Cref{eqn:gz1,eqn:PR}} \frac{1}{c_1^t} \expt{ \expt{\sfZ \mid \sfP_t, \sfY} \sfZ } - a_t \sigma^2 \notag \\
    &= \frac{1}{c_1^t} (\sigma^2 - c_1^t) - \paren{ \tau_1^t - \frac{1}{\sigma^2} } \sigma^2
    = \sigma^2 \paren{ \frac{1}{c_1^t} - \tau_1^t }
    \explain{\Cref{eqn:B_gamma2_tau2}} \sigma^2 \tau_2^t . \notag 
\end{align}
implying 
\begin{align}
    \wt{a}_t &= \tau_2^t , 
    \label{eqn:ind_wta}
\end{align}
as desired. 
In the above display, the first equality in the second line uses the observation \Cref{eqn:cond_expt} and the induction hypothesis $ a_t = \tau_1^t - \sigma^{-2} $. 

Similarly, 
\begin{align}
    \rho \wt{b}_t &\explain{\Cref{eqn:ab_tilde}} \expt{\wt{\sfR}_t \sfB_*}
    \explain{\Cref{eqn:wtP_wtR}} \expt{ f_t(\sfR_t) \sfB_* }
    \explain{\Cref{eqn:B_fg}} \frac{1}{v_1^t} \expt{g_{x1}(\sfR_t; \gamma_1^t) \sfB_*} - \expt{\sfR_t \sfB_*}
    \explain{\Cref{eqn:PR,eqn:gx1}} \frac{1}{v_1^t} \expt{ \expt{\sfB_* \mid \sfR_t} \sfB_* } - b_t \rho , \notag \\
    &= \frac{1}{v_1^t} (\rho - v_1^t) - \gamma_1^t \rho 
    = \rho \paren{ \frac{1}{v_1^t} - \gamma_1^t } - 1
    \explain{\Cref{eqn:B_gamma2_tau2}} \rho \gamma_2^t - 1 , \notag 
\end{align}
where the first equality in the second line uses \Cref{eqn:cond_expt} and the induction hypothesis $ b_t = \gamma_1^t $. 
This implies
\begin{align}
    \wt{b}_t &= \gamma_2^t - \frac{1}{\rho} . 
    \label{eqn:ind_wtb}
\end{align}

We then consider $ b_{t+1} $: 
\begin{align}
    b_{t+1} 
    &\explain{\Cref{eqn:ab}} \expt{ \wt{\Phi}_t(\sfLambda_d) \sfLambda_d } \wt{a}_t
    \explain{\Cref{eqn:B_Phi,eqn:ind_wta}} \frac{1}{v_2^t} \cdot \expt{ \frac{\tau_2^t \sfLambda_d^2}{\gamma_2^t + \tau_2^t \sfLambda_d^2} } \notag \\
    &= \frac{1}{v_2^t} \paren{ 1 - \gamma_2^t \cdot \expt{ \frac{1}{\gamma_2^t + \tau_2^t \sfLambda_d^2} } }
    \explain{\Cref{eqn:B_v2_c2}} \frac{1}{v_2^t} \paren{ 1 - \gamma_2^t \cdot v_2^t } 
    = \frac{1}{v_2^t} - \gamma_2^t
    \explain{\Cref{eqn:B_gamma1_tau1}} \gamma_1^{t+1} , \label{eqn:ind_b} 
\end{align}
as desired. 
Next, move to $ a_{t+1} $: 
\begin{align}
    a_{t+1} &\explain{\Cref{eqn:ab}} \frac{1}{\ol{\kappa}_2} \paren{ \expt{ \wt{\Psi}_t(\sfLambda_n) \sfLambda_n } \wt{b}_t + \expt{ \Psi_t(\sfLambda_n^2) \sfLambda_n^2 } \wt{a}_t } \notag \\
    &\explain{\Cref{eqn:ind_wta,eqn:ind_wtb,eqn:B_Psi,eqn:Psi_wtPsi}} \frac{1}{\ol{\kappa}_2} \paren{ \expt{ \frac{1}{c_2^t} \cdot \frac{\sfLambda_n^2}{\gamma_2^t + \tau_2^t \sfLambda_n^2} } \paren{ \gamma_2^t - \frac{1}{\rho} } + \expt{ \paren{ \frac{1}{c_2^t} \cdot \frac{\sfLambda_n^2}{\gamma_2^t + \tau_2^t \sfLambda_n^2} - 1 } \sfLambda_n^2 } \tau_2^t } \notag \\
    &\explain{\Cref{eqn:B_v2_c2}} \frac{1}{\ol{\kappa}_2} \brace{ \frac{1}{c_2^t} \cdot c_2^t \cdot \paren{\gamma_2^t - \frac{1}{\rho}} + \frac{\tau_2^t}{c_2^t} \cdot \frac{1}{\tau_2^t} \paren{ \expt{\sfLambda_n^2} - \gamma_2^t \expt{ \frac{\sfLambda_n^2}{\gamma_2^t + \tau_2^t \sfLambda_n^2} } } - \tau_2^t \expt{\sfLambda_n^2} } \notag \\
    &\explain{\Cref{eqn:B_v2_c2}} \frac{1}{\ol{\kappa}_2} \brace{ \gamma_2^t - \frac{1}{\rho} + \frac{1}{c_2^t} \paren{ \ol{\kappa}_2 - \gamma_2^t c_2^t } - \tau_2^t \ol{\kappa}_2 }
    = \frac{1}{\ol{\kappa}_2} \brace{ \frac{1}{c_2^t} \ol{\kappa}_2  - \tau_2^t \ol{\kappa}_2 - \frac{1}{\rho} }
    \explain{\Cref{eqn:sigma2}} \frac{1}{c_2^t} - \tau_2^t - \frac{1}{\sigma^2} 
    \explain{\Cref{eqn:B_gamma1_tau1}} \tau_1^{t+1} - \frac{1}{\sigma^2} . \label{eqn:ind_a} 
\end{align}

Before proceeding to $ \expt{\sfJ_{t+1}^2}, \expt{\sfK_{t+1}^2} $, using the induction hypotheses, we derive expressions of $ \expt{\wt{\sfP}_t^2}, \expt{\wt{\sfR}_t^2} $ in terms of $ \tau_2^t, \gamma_2^t $. 
By the tower property of conditional expectation, 
\begin{align}
    \expt{\wt{\sfP}_t^2} &= \expt{\paren{ \frac{1}{c_1^t} \expt{\sfZ \mid \sfP_t, \sfY} - \sfP_t }^2}
    = \frac{1}{(c_1^t)^2} \expt{ \sfZ \expt{\sfZ \mid \sfP_t, \sfY} } + \expt{\sfP_t^2} - \frac{2}{c_1^t} \expt{\sfP_t \sfZ} \notag \\
    &= \frac{1}{(c_1^t)^2} (\sigma^2 - c_1^t) + \brack{ \paren{\tau_1^t - \frac{1}{\sigma^2}}^2 \sigma^2 + \paren{\tau_1^t - \frac{1}{\sigma^2}} } - \frac{2}{c_1^t} \paren{ \tau_1^t - \frac{1}{\sigma^2} } \sigma^2 \notag \\
    &= \sigma^2 \paren{ \frac{1}{c_1^t} - \tau_1^t } \paren{ \frac{1}{c_1^t} - \tau_1^t + \frac{1}{\sigma^2} } 
    \explain{\Cref{eqn:B_gamma2_tau2}} \sigma^2 \tau_2^t \paren{ \tau_2^t + \frac{1}{\sigma^2} } , \label{eqn:wtP_tau} 
\end{align}
where the second line follows from the induction hypotheses $ a_t = \expt{\sfJ_t^2} = \tau_1^t - \sigma^{-2} $. 

Similarly, 
\begin{align}
    \expt{\wt{\sfR}_t^2} &= \expt{ \paren{\frac{1}{v_1^t} \expt{\sfB_* \mid \sfR_t} - \sfR_t}^2 }
    = \frac{1}{(v_1^t)^2} \expt{ \sfB_* \expt{\sfB_* \mid \sfR_t} } + \expt{\sfR_t^2} - \frac{2}{v_1^t} \expt{\sfB_* \sfR_t} \notag \\
    &= \frac{1}{(v_1^t)^2} (\rho - v_1^t) + \brack{ (\gamma_1^t)^2 \rho + \gamma_1^t } - \frac{2}{v_1^t} \cdot \gamma_1^t \rho 
    = \frac{\rho}{v_1^t} \paren{\frac{1}{v_1^t} - \gamma_1^t} - \paren{\frac{1}{v_1^t} - \gamma_1^t} - \rho \gamma_1^t \paren{\frac{1}{v_1^t} - \gamma_1^t} \notag \\
    &= \rho \gamma_2^t \paren{ \frac{1}{v_1^t} - \gamma_1^t - \frac{1}{\rho} }
    \explain{\Cref{eqn:B_gamma2_tau2}} \rho \gamma_2^t \paren{ \gamma_2^t - \frac{1}{\rho} } , \label{eqn:wtR_gamma} 
\end{align}
where the second line follows from the induction hypotheses $ b_t = \expt{\sfK_t^2} = \gamma_1^t $. 

Equipped with the expressions \Cref{eqn:wtP_tau,eqn:wtR_gamma}, we now verify $ \expt{\sfJ_{t+1}^2} = \tau_1^{t+1} - \sigma^{-2}, \expt{\sfK_{t+1}^2} = \gamma_1^{t+1} $. 
First consider $ \expt{\sfJ_{t+1}^2} $: 
\begin{align}
    \expt{\sfJ_{t+1}^2}
    &\explain{\Cref{eqn:J_cov}} \expt{ \Psi_t(\sfLambda_n^2)^2 } \paren{ \expt{\wt{\sfP}_t^2} - \wt{a}_t^2 \sigma^2 }
    + \expt{ \Psi_t(\sfLambda_n^2)^2 \sfLambda_n^2 } \wt{a}_t^2 \rho 
    + \expt{ \wt{\Psi}_t(\sfLambda_n)^2 } \expt{ \wt{\sfR}_t^2 } 
    - a_{t+1}^2 \sigma^2 \notag \\
    &\quad + 2 \expt{ \Psi_t(\sfLambda_n^2) \wt{\Psi}_t(\sfLambda_n) \sfLambda_n } \wt{a}_t \wt{b}_t \rho , \notag \\
    &= \paren{ \expt{ \sfLambda_n^2 \wt{\Psi}_t(\sfLambda_n)^2 } - 1 } \brack{ \sigma^2 \tau_2^t \paren{ \tau_2^t + \frac{1}{\sigma^2} } - (\tau_2^t)^2 \sigma^2 } \notag \\
    &\quad + \paren{ \expt{ \sfLambda_n^4 \wt{\Psi}_t(\sfLambda_n)^2 } + \expt{ \sfLambda_n^2 } - 2 \expt{ \sfLambda_n^3 \wt{\Psi}_t(\sfLambda_n) } } (\tau_2^t)^2 \rho \notag \\
    &\quad + \expt{ \wt{\Psi}_t(\sfLambda_n)^2 } \rho \gamma_2^t \paren{ \gamma_2^t - \frac{1}{\rho} }
    - \paren{\tau_1^{t+1} - \frac{1}{\sigma^2}}^2 \sigma^2 \notag \\
    &\quad + 2 \paren{ \expt{ \sfLambda_n^2 \wt{\Psi}_t(\sfLambda_n)^2 } - 1 } \tau_2^t \paren{\gamma_2^t - \frac{1}{\rho}} \rho \notag \\
    &\explain{\Cref{eqn:B_gamma1_tau1}} \expt{ \sfLambda_n^2 \wt{\Psi}_t(\sfLambda_n)^2 } \paren{ \tau_2^t + 2 \tau_2^t \paren{\gamma_2^t - \frac{1}{\rho}} \rho } \notag \\
    &\quad + \paren{ \expt{ \sfLambda_n^4 \wt{\Psi}_t(\sfLambda_n)^2 } + \expt{ \sfLambda_n^2 } - 2 \expt{ \sfLambda_n^3 \wt{\Psi}_t(\sfLambda_n) } } (\tau_2^t)^2 \rho
    + \expt{ \wt{\Psi}_t(\sfLambda_n)^2 } \rho \gamma_2^t \paren{ \gamma_2^t - \frac{1}{\rho} } \notag \\
    &\quad - \paren{ \tau_2^t + 2 \tau_2^t \paren{\gamma_2^t - \frac{1}{\rho}} \rho + \paren{\frac{1}{c_2^t} - \tau_2^t - \frac{1}{\sigma^2}}^2 \sigma^2 } . \label{eqn:others1} 
\end{align}
In the above display, the second equality uses the expressions \Cref{eqn:expt_Psi0} of several expectations involving $ \Psi_t, \wt{\Psi}_t $, the expressions \Cref{eqn:wtP_tau,eqn:wtR_gamma} of $ \expt{\wt{\sfP}_t^2}, \expt{\wt{\sfR}_t^2} $ and the expressions of $ \wt{a}_t, \wt{b}_t, a_{t+1} $ just proved in \Cref{eqn:ind_wta,eqn:ind_wtb,eqn:ind_a}. 

Let us first examine terms involving expectations, excluding the last line \Cref{eqn:others1}. 
\begin{align}
    & \expt{ \sfLambda_n^2 \wt{\Psi}_t(\sfLambda_n)^2 } \paren{ \tau_2^t + 2 \tau_2^t \paren{\gamma_2^t - \frac{1}{\rho}} \rho } \notag \\
    &\quad + \paren{ \expt{ \sfLambda_n^4 \wt{\Psi}_t(\sfLambda_n)^2 } + \ol{\kappa}_2 - 2 \expt{ \sfLambda_n^3 \wt{\Psi}_t(\sfLambda_n) } } (\tau_2^t)^2 \rho
    + \expt{ \wt{\Psi}_t(\sfLambda_n)^2 } \rho \gamma_2^t \paren{ \gamma_2^t - \frac{1}{\rho} } \notag \\
    &\explain{\Cref{eqn:expect_Psi_explicit}} \paren{ \frac{1}{(c_2^t)^2} \expt{ \frac{\sfLambda_n^4}{(\gamma_2^t + \tau_2^t \sfLambda_n^2)^2} } } \paren{ \tau_2^t + 2 \tau_2^t \paren{\gamma_2^t - \frac{1}{\rho}} \rho } \notag \\
    &\quad + \paren{ \frac{\ol{\kappa}_2}{(\tau_2^t c_2^t)^2} - \frac{\gamma_2^t}{(\tau_2^t)^2 c_2^t} - \frac{\gamma_2^t}{\tau_2^t (c_2^t)^2} \expt{ \frac{\sfLambda_n^4}{(\gamma_2^t + \tau_2^t \sfLambda_n^2)^2} } + \ol{\kappa}_2 - 2 \paren{ \frac{\ol{\kappa}_2}{\tau_2^t c_2^t} - \frac{\gamma_2^t}{\tau_2^t} } } (\tau_2^t)^2 \rho \label{eqn:others0} \\
    &\quad + \paren{ \frac{1}{\gamma_2^t c_2^t} - \frac{\tau_2^t}{\gamma_2^t (c_2^t)^2} \expt{ \frac{\sfLambda_n^4}{(\gamma_2^t + \tau_2^t \sfLambda_n^2)^2} } } \rho \gamma_2^t \paren{ \gamma_2^t - \frac{1}{\rho} } . \label{eqn:others2} 
\end{align}
Elementary algebra shows that the coefficients in front of $ \expt{ \frac{\sfLambda_n^4}{(\gamma_2^t + \tau_2^t \sfLambda_n^2)^2} } $ add up to zero: 
\begin{align}
    & \frac{1}{(c_2^t)^2} \paren{ \tau_2^t + 2 \tau_2^t \paren{\gamma_2^t - \frac{1}{\rho}} \rho }
    - \frac{\gamma_2^t}{\tau_2^t (c_2^t)^2} (\tau_2^t)^2 \rho
    - \frac{\tau_2^t}{\gamma_2^t (c_2^t)^2} \rho \gamma_2^t \paren{ \gamma_2^t - \frac{1}{\rho} }
    = 0 . \notag 
\end{align}
Therefore all terms involving that expectation are cancelled. 
Moreover, the remaining terms in \Cref{eqn:others0,eqn:others2,eqn:others1} that do not involve $ \expt{ \frac{\sfLambda_n^4}{(\gamma_2^t + \tau_2^t \sfLambda_n^2)^2} } $ add up to 
\begin{align}
    & \paren{ \frac{\ol{\kappa}_2}{(\tau_2^t c_2^t)^2} - \frac{\gamma_2^t}{(\tau_2^t)^2 c_2^t} + \ol{\kappa}_2 - 2 \paren{ \frac{\ol{\kappa}_2}{\tau_2^t c_2^t} - \frac{\gamma_2^t}{\tau_2^t} } } (\tau_2^t)^2 \rho 
    + \frac{1}{\gamma_2^t c_2^t} \cdot \rho \gamma_2^t \paren{ \gamma_2^t - \frac{1}{\rho} } \notag \\
    & - \paren{ \tau_2^t + 2 \tau_2^t \paren{\gamma_2^t - \frac{1}{\rho}} \rho + \paren{\frac{1}{c_2^t} - \tau_2^t - \frac{1}{\sigma^2}}^2 \sigma^2 } \notag \\
    &= \frac{1}{c_2^t} - \tau_2^t - \frac{1}{\sigma^2} 
    \explain{\Cref{eqn:B_gamma1_tau1}} \tau_1^{t+1} - \frac{1}{\sigma^2} . \notag 
\end{align}
So we have shown 
\begin{align}
    \expt{\sfJ_{t+1}^2} = \tau_1^{t+1} - \sigma^{-2} . 
    \label{eqn:ind_J}
\end{align}

Finally, we turn to $ \expt{\sfK_{t+1}^2} $: 
\begin{align}
    \expt{ \sfK_{t+1}^2 }
    &\explain{\Cref{eqn:K_cov}} \expt{ \Phi_t(\sfLambda_d^2)^2 } \expt{\wt{\sfR}_t^2} 
    - b_{t+1}^2 \rho 
    + \expt{ \wt{\Phi}_t(\sfLambda_d)^2 } \paren{ \expt{\wt{\sfP}_t^2} - \wt{a}_t^2 \sigma^2 }
    + \expt{ \wt{\Phi}_t(\sfLambda_d)^2 \sfLambda_d^2 } \wt{a}_t^2 \rho \notag \\
    &\quad + 2 \expt{ \Phi_t(\sfLambda_d^2) \wt{\Phi}_t(\sfLambda_d) \sfLambda_d } \wt{a}_t \wt{b}_t \rho \notag \\
    &= \paren{ \frac{(\tau_2^t)^2}{(\gamma_2^t)^2} \expt{ \wt{\Phi}_t(\sfLambda_d)^2 \sfLambda_d^2 } - \frac{(\gamma_1^{t+1})^2}{(\gamma_2^t)^2} } \rho \gamma_2^t \paren{\gamma_2^t - \frac{1}{\rho}} 
    - (\gamma_1^{t+1})^2 \rho \notag \\
    &\quad + \expt{ \wt{\Phi}_t(\sfLambda_d)^2 } \paren{ \sigma^2 \tau_2^t \paren{ \tau_2^t + \frac{1}{\sigma^2} } - (\tau_2^t)^2 \sigma^2 } 
    + \expt{ \wt{\Phi}_t(\sfLambda_d)^2 \sfLambda_d^2 } (\tau_2^t)^2 \rho \notag \\
    &\quad + 2 \paren{ \frac{(\gamma_1^{t+1})^2}{\gamma_2^t \tau_2^t} - \frac{\tau_2^t}{\gamma_2^t} \expt{ \wt{\Phi}_t(\sfLambda_d)^2 \sfLambda_d^2 } } \tau_2^t \paren{ \gamma_2^t - \frac{1}{\rho} } \rho \notag \\
    &\explain{\Cref{eqn:B_gamma1_tau1}} \expt{ \wt{\Phi}_t(\sfLambda_d)^2 \sfLambda_d^2 } \paren{ \frac{(\tau_2^t)^2}{(\gamma_2^t)^2} \rho \gamma_2^t \paren{\gamma_2^t - \frac{1}{\rho}} + (\tau_2^t)^2 \rho
    - 2 \frac{\tau_2^t}{\gamma_2^t} \tau_2^t \paren{ \gamma_2^t - \frac{1}{\rho} } \rho } 
    + \expt{\wt{\Phi}_t(\sfLambda_d)^2} \tau_2^t \label{eqn:Phi_line1} \\
    &\quad - \paren{ \frac{(1/v_2^t - \gamma_2^t)^2}{(\gamma_2^t)^2} \rho \gamma_2^t \paren{\gamma_2^t - \frac{1}{\rho}}
    + (1/v_2^t - \gamma_2^t)^2 \rho - 2 \frac{(1/v_2^t - \gamma_2^t)^2}{\gamma_2^t \tau_2^t} \tau_2^t \paren{\gamma_2^t - \frac{1}{\rho}} \rho } , \label{eqn:Phi_line2}
\end{align}
where the second equality uses the expressions of $ \expt{\wt{\sfP}_t^2}, \expt{\wt{\sfR}_t^2} $ in \Cref{eqn:wtP_tau,eqn:wtR_gamma}, the expressions of $ \wt{a}_t, \wt{b}_t, b_{t+1} $ derived in \Cref{eqn:ind_wta,eqn:ind_wtb,eqn:ind_b} and the expressions of $ \expt{\Phi_t(\sfLambda_d^2)^2}, \expt{\Phi_t(\sfLambda_d^2) \wt{\Phi}_t(\sfLambda_d) \sfLambda_d} $ in \Cref{eqn:expect_Phi0}. 

We use \Cref{eqn:expect_Phi} to expand the terms in the first line \Cref{eqn:Phi_line1}: 
\begin{align}
    & \expt{ \wt{\Phi}_t(\sfLambda_d)^2 \sfLambda_d^2 } \paren{ \frac{(\tau_2^t)^2}{(\gamma_2^t)^2} \rho \gamma_2^t \paren{\gamma_2^t - \frac{1}{\rho}} + (\tau_2^t)^2 \rho
    - 2 \frac{\tau_2^t}{\gamma_2^t} \tau_2^t \paren{ \gamma_2^t - \frac{1}{\rho} } \rho } 
    + \expt{\wt{\Phi}_t(\sfLambda_d)^2} \tau_2^t \notag \\
    &= \paren{ \frac{(\gamma_2^t)^2}{(\tau_2^t v_2^t)^2} \expt{ \frac{1}{(\gamma_2^t + \tau_2^t \sfLambda_d)^2} }  - \frac{2 \gamma_2^t}{(\tau_2^t)^2 v_2^t} + \frac{1}{(\tau_2^t v_2^t)^2} } \paren{ \frac{(\tau_2^t)^2}{(\gamma_2^t)^2} \rho \gamma_2^t \paren{\gamma_2^t - \frac{1}{\rho}} + (\tau_2^t)^2 \rho
    - 2 \frac{\tau_2^t}{\gamma_2^t} \tau_2^t \paren{ \gamma_2^t - \frac{1}{\rho} } \rho } \label{eqn:Phi_other1} \\
    &\quad + \paren{ \frac{1}{\tau_2^t v_2^t} - \frac{\gamma_2^t}{\tau_2^t (v_2^t)^2} \expt{ \frac{1}{(\gamma_2^t + \tau_2^t \sfLambda_d^2)^2} } } \tau_2^t . \label{eqn:Phi_other2}
\end{align}
It can be verified that all terms involving $ \expt{ \frac{1}{(\gamma_2^t + \tau_2^t \sfLambda_d^2)^2} } $ are cancelled. 
Furthermore, the remaining terms in \Cref{eqn:Phi_other1,eqn:Phi_other2,eqn:Phi_line2} that do not involve $ \expt{ \frac{1}{(\gamma_2^t + \tau_2^t \sfLambda_d^2)^2} } $ collectively give $ \frac{1}{v_2^t} - \gamma_2^t = \gamma_1^{t+1} $, where the equality is by \Cref{eqn:B_gamma1_tau1}. 
So we have shown that 
\begin{align}
    \expt{\sfK_{t+1}^2} = \gamma_1^{t+1} 
    \label{eqn:ind_K}
\end{align}
which completes the induction step. 

To conclude the proof of the result, it remains to check the validity of \Cref{eqn:SE_relation} in the base case $t=0$.

\paragraph{\textsc{Base case.}}
Using \Cref{thm:spec_GVAMP} for $t=0$ and the initial condition \Cref{eqn:RP0_BGVAMP}, 
\begin{align}
&&
    \matrix{r^0 & \beta_*} &\xrightarrow{W_2} \matrix{\sfR_0 & \sfB_*} , & 
    \matrix{p^0 & z & \eps} &\xrightarrow{W_2} \matrix{\sfP_0 & \sfZ & \sfE} . & 
& \label{eqn:RP0_BGVAMP_pf}
\end{align}
Now to show 
\begin{align}
&&
    \wt{a}_0 &= \tau_2^0 , &
    \wt{b}_0 &= \gamma_2^0 - \rho^{-1} , &
& \label{eqn:ind_wtab0}
\end{align}
we repeat the same computations that lead to \Cref{eqn:ind_wta,eqn:ind_wtb} but use 
\begin{align}
&&
    \expt{\sfP_0 \sfZ} &= \tau_1^0 \sigma^2 - 1, & 
    \expt{\sfR_0 \sfB_*} &= \gamma_1^0 \rho & 
& \label{eqn:P0Z_R0B}
\end{align}
guaranteed by \Cref{eqn:RP0_BGVAMP_pf} instead of the induction hypotheses. 
Using this in \Cref{eqn:ind_b,eqn:ind_a} for $t=0$ then yields 
\begin{align}
&&
    b_1 &= \gamma_1^1, &
    a_1 &= \tau_1^1 - \sigma^{-2}. &
& \label{eqn:ind_ab1}
\end{align}
Moreover, using \Cref{eqn:P0Z_R0B} in place of the induction hypotheses in \Cref{eqn:wtP_tau,eqn:wtR_gamma} for $t=0$, we have 
\begin{align}
&&
    \expt{\wt{\sfP}_0^2} &= \sigma^2 \tau_2^0 (\tau_2^0 + \sigma^{-2}) , &
    \expt{\wt{\sfR}_0^2} &= \rho \gamma_2^0 (\gamma_2^0 - \rho^{-1}) . & 
& \label{eqn:wtP0_wtR0}
\end{align}
Using \Cref{eqn:ind_wtab0,eqn:ind_ab1,eqn:wtP0_wtR0} shown above together with other observations in the \textsc{Preparation} step (which do not require the induction hypotheses) in the calculations leading to \Cref{eqn:ind_J,eqn:ind_K}, we obtain the desired $ \expt{\sfJ_1^2} = \tau_1^1 - \sigma^{-2} , \expt{\sfK_1^2} = \gamma_1^1 $. 
This completes the verification for the $t=0$ case of \Cref{eqn:SE_relation} and thereby concludes the whole proof. 
\end{proof}

\section{Proof of \Cref{thm:replica}}
\label{app:pf_thm:replica}

Before diving in the proof, we first note from the definition \Cref{eqn:Qz} of $ Q_z $ and the definition \Cref{eqn:sigma2} of $ \sigma^2 $ that 
\begin{align}
    Q_z &= \frac{\rho \expt{\sfLambda_d^2}}{\delta} = \rho \expt{ \sfLambda_n^2 } = \rho \ol{m}_2 = \rho \ol{\kappa}_2 = \sigma^2 . \label{eqn:Qz_sigma2} 
\end{align}
Therefore $ \wh{Q}_z = \sigma^{-2} $. 

To prove the result, we will argue that under the change of variables in \Cref{eqn:change_var}, the fixed point equations for the variables in \Cref{eqn:SE_FP} coincide with the replica saddle point equations \Cref{eqn:replica}. 
The former equations can be obtained by removing the time indices from the recursion for $ (v_1^t, c_1^t, \gamma_2^t, \tau_2^t, v_2^t, c_2^t, \gamma_1^{t+1}, \tau_1^{t+1}) $ in \Cref{eqn:Bayes_GVAMP}: 
\begin{subequations}
\label{eqn:B_FP}
\begin{align}
&&
    v_1 &= \expt{g_{x1}'(\sfR; \gamma_1)} , & 
    c_1 &= \expt{g_{z1}'(\sfP, \sfY; \tau_1)} , & 
& \label{eqn:BFP_v1_c1} \\
&&
    \gamma_1 + \gamma_2 &= \frac{1}{v_1} = \frac{1}{v_2} , & 
    \tau_1 + \tau_2 &= \frac{1}{c_1} = \frac{1}{c_2} , & 
& \label{eqn:BFP_gamma_tau} \\
&&
    v_2 &= \expt{\frac{1}{\tau_2 \sfLambda_d^2 + \gamma_2}} , & 
    c_2 &= \expt{\frac{\sfLambda_n^2}{\tau_2 \sfLambda_n^2 + \gamma_2}} , & 
& \label{eqn:BFP_v2_c2}
\end{align}
\end{subequations}
We will check the equivalence between \Cref{eqn:replica,eqn:B_FP} equation by equation. 

Given the observation \Cref{eqn:Qz_sigma2}, the equivalences between \Cref{eqn:replica3} and the first equation in \Cref{eqn:BFP_gamma_tau}, between \Cref{eqn:replica4} and the second equation in \Cref{eqn:BFP_gamma_tau}, between \Cref{eqn:replica5} and the first equation in \Cref{eqn:BFP_v2_c2}, between \Cref{eqn:replica6} and the second equation in \Cref{eqn:BFP_v2_c2} are apparent. 
The rest of the proof is devoted to checking the equivalences between \Cref{eqn:replica1} and the first equation in \Cref{eqn:BFP_v1_c1}, between \Cref{eqn:replica2} and the second equation in \Cref{eqn:BFP_v1_c1}. 

First consider \Cref{eqn:replica1}. 
Recalling the definition of $ \cZ_0(b,a) $ in \Cref{eqn:Z0_Zout}, we compute the partial derivative of the log partition function evaluated at $ (\sqrt{\wh{q}_x} \xi, \wh{q}_x) $: 
\begin{align}
    \partial_b \log \cZ_0\paren{ \sqrt{\wh{q}_x} \xi, \wh{q}_x }
    &= \frac{1}{\cZ_0\paren{ \sqrt{\wh{q}_x} \xi, \wh{q}_x }} 
    \int x P_{\sfB_*}(x) \exp\paren{ -\frac{1}{2} \wh{q}_x x^2 + \sqrt{\wh{q}_x} \xi x } \diff x . \notag 
\end{align}
This allows us to evaluate the expectation on the right-hand side of \Cref{eqn:replica1}: 
\begin{align}
    & \expt{ \cZ_0\paren{\sqrt{\wh{q}_x}\,\sfXi, \wh{q}_x} \paren{ \partial_b \log \cZ_0\paren{\sqrt{\wh{q}_x}\,\sfXi, \wh{q}_x} }^2 } \notag \\
    &= \expt{ \frac{1}{\cZ_0\paren{ \sqrt{\wh{q}_x} \sfXi, \wh{q}_x }} 
    \paren{ \int x P_{\sfB_*}(x) \exp\paren{ -\frac{1}{2} \wh{q}_x x^2 + \sqrt{\wh{q}_x} \sfXi x } \diff x }^2 } \notag \\
    &= \int \frac{\exp\paren{-\frac{1}{2} \xi^2}}{\sqrt{2\pi}} \frac{\paren{ \int x P_{\sfB_*}(x) \exp\paren{ -\frac{1}{2} \wh{q}_x x^2 + \sqrt{\wh{q}_x} \xi x } \diff x }^2}{\int P_{\sfB_*}(x) \exp\paren{ -\frac{1}{2} \wh{q}_x x^2 + \sqrt{\wh{q}_x} \xi x } \diff x} 
     \diff \xi \notag \\
     &= \frac{1}{\sqrt{\wh{q}_x}} \int \frac{\paren{ \int x P_{\sfB_*}(x) \cdot \frac{1}{\sqrt{2\pi/\wh{q}_x}} \exp\paren{ -\frac{1}{2} \wh{q}_x \paren{ x - \frac{\xi}{\sqrt{\wh{q}_x}} }^2 } \diff x }^2}{\int P_{\sfB_*}(x) \cdot \frac{1}{\sqrt{2\pi/\wh{q}_x}} \exp\paren{ -\frac{1}{2} \wh{q}_x \paren{ x - \frac{\xi}{\sqrt{\wh{q}_x}} }^2 } \diff x} 
     \diff \xi \notag \\
     &= \expt{ \expt{ \sfB_* \mid \sfR }^2 } , 
     \qquad \frac{\sfR}{\wh{q}_x} = \sfB_* + \frac{1}{\sqrt{\wh{q}_x}} \ol{\sfK} , \ (\sfB_*, \ol{\sfK}) \sim P_{\sfB_*} \ot \cN(0,1) , \notag 
\end{align}
where we made the change of variable $ \sfR = \sqrt{\wh{q}_x} \sfXi $ in the last step. 

On the other hand, from \Cref{eqn:v1_alt}, we see that the first equation of \Cref{eqn:BFP_v1_c1} is equivalent to 
\begin{align}
    v_1 &= \rho - \expt{\sfB_* \expt{\sfB_* \mid \sfR}} . \label{eqn:BFP_v1_alt}
\end{align}
In view of the change of variable $ v_1 = \rho - q_x $ and the tower property of conditional expectation, we conclude the equivalence between \Cref{eqn:replica1,eqn:BFP_v1_alt}. 

Next, consider \Cref{eqn:replica2}. 
Recalling the definition of $ \cZ_\out(y; \xi, v) $, we compute the partial derivative of the log partition function: 
\begin{align}
    \partial_\xi \log \cZ_\out(y; \xi, v)
    &= \frac{1}{\cZ_\out(y; \xi, v)} \expt{ \partial_\xi Q\paren{ y \mid \sqrt{v} \sfOmega + \xi } } \notag \\
    &= \frac{1}{\cZ_\out(y; \xi, v)} \frac{1}{\sqrt{v}} \expt{ \partial_z Q\paren{ y \mid \sqrt{v} \sfOmega + \xi } } 
    = \frac{1}{\cZ_\out(y; \xi, v)} \frac{1}{\sqrt{v}} \expt{ \sfOmega Q\paren{ y \mid \sqrt{v} \sfOmega + \xi } } . \notag 
\end{align}
With the shorthand $ \fra = \frac{1}{\sqrt{\wh{Q}_z + \wh{q}_z}} $ and $ \frb = \sqrt{\frac{\wh{q}_z}{\wh{Q}_z (\wh{Q}_z + \wh{q}_z)}} $, we can write the expectation on the right-hand side of \Cref{eqn:replica2} as
\begin{align}
    & \expt{ \int \cZ_\out\paren{ y; \frb\sfXi, \fra^2 } \paren{ \partial_\xi \log \cZ_\out\paren{ y; \frb\sfXi, \fra^2 } }^2 \diff y } \notag \\
    &= \frac{1}{\fra^2} \expt{ \int \frac{\expt{ \sfOmega Q\paren{y \mid \fra \sfOmega + \frb \sfXi} }^2}{\expt{ Q\paren{y \mid \fra \sfOmega + \frb \sfXi} }} \diff y } \notag \\
    &= \frac{1}{\fra^2} \iint p_{\sfXi}(\xi) \frac{\paren{ \int \omega p_{\sfOmega}(\omega) p_{\lr{\sfY \mid \sfOmega, \sfXi}}\paren{ y \mid \omega, \xi } }^2 \diff \omega}{\int p_{\sfOmega}(\omega) p_{\lr{\sfY \mid \sfOmega, \sfXi}}\paren{y \mid \omega, \xi} \diff \omega} \diff y \diff \xi \notag \\
    &= \frac{1}{\fra^2} \iint p_{\sfXi}(\xi) p_{\lr{\sfY \mid \sfXi}}\paren{y \mid \xi} \paren{\int \omega \frac{p_{\lr{\sfOmega, \sfY \mid \sfXi}}\paren{\omega, y \mid \xi}}{p_{\lr{\sfY \mid \sfXi}}\paren{y \mid \xi}} \diff \omega }^2 \diff y \diff \xi \notag \\
    &= \frac{1}{\fra^2} \iint p_{\sfXi, \sfY}(\xi, y) \paren{ \int \omega p_{\lr{\sfOmega \mid \sfY, \sfXi}}\paren{\omega \mid y, \xi} \diff \omega }^2 \diff y \diff \xi \notag \\
    &= \frac{1}{\fra^2} \expt{ \expt{\sfOmega \mid \sfXi, \sfY}^2 } ,
    \qquad \textnormal{where } (\sfXi, \sfOmega) \sim \cN(0,1) \ot \cN(0,1) , \ \sfY \sim Q\paren{\cdot \mid \fra \sfOmega + \frb \sfXi} \notag \\
    &= \frac{1}{\fra^2} \expt{ \expt{ (\sfM - \frb \sfXi) \mid \sfXi, \sfY }^2 } , \qquad \textnormal{where } \sfM = \fra \sfOmega + \frb \sfXi \notag \\
    &= \frac{1}{\fra^2} \paren{ \expt{ \expt{\sfM \mid \sfXi, \sfY}^2 } + \frb^2 - 2 \frb \expt{ \sfM \sfXi } } \notag \\
    &= \frac{1}{\fra^2} \paren{ \expt{ \expt{\sfM \mid \sfXi, \sfY}^2 } - \frb^2 } , \notag 
\end{align}
where the penultimate line is due to the tower property of conditional expectation. 

To connect this to the second equation in \Cref{eqn:BFP_v1_c1}, we make the following observation. 
By \Cref{thm:SE_Bayes_GVAMP}, at fixed point, $ \sfP $ has distribution 
\begin{align}
&&
    \sfP &= \paren{\tau_1 - \frac{1}{\sigma^2}} \sfZ + \sqrt{\tau_1 - \frac{1}{\sigma^2}} \ol{\sfK} , & 
    & \textnormal{where } (\sfZ, \ol{\sfK}) \sim \cN(0,\sigma^2) \ot \cN(0,1) . &
& \notag 
\end{align}
In particular $ \sfP \sim \cN(0, (\tau_1 - \sigma^{-2}) \tau_1 \sigma^2) $. 
From this, one can deduce the reverse channel
\begin{align}
&&
    \sfZ &= \frac{1}{\sqrt{\tau_1}} \sfOmega + \frac{1}{\tau_1} \sfP , & 
    & \textnormal{where } (\sfOmega, \sfP) \sim \cN(0,1) \ot \cN(0, (\tau_1 - \sigma^{-2}) \tau_1 \sigma^2) . & 
& \label{eqn:sfZ_alt} 
\end{align}

On the other hand, recalling the shorthand $ \fra,\frb $, using the observation $ \wh{Q}_z = \sigma^{-2} $ and the change of variables $ \wh{q}_z = \tau_1 - \sigma^{-2} $, we have
\begin{align}
&&
    \sfM &= \frac{1}{\sqrt{\wh{Q}_z + \wh{q}_z}} \sfOmega + \sqrt{\frac{\wh{q}_z}{\wh{Q}_z (\wh{Q}_z + \wh{q}_z)}} \sfXi 
    = \frac{1}{\sqrt{\tau_1}} \sfOmega + \sqrt{\sigma^2 - \frac{1}{\tau_1}} \sfXi . & 
& \label{eqn:sfM} 
\end{align}
Contrasting \Cref{eqn:sfM} with \Cref{eqn:sfZ_alt}, we find that 
\begin{align}
    \paren{\sfP, \sfZ} &\eqqlaw \paren{\sqrt{(\tau_1 - \sigma^2) \tau_1 \sigma^2} \sfXi, \sfM} . \notag 
\end{align}
Therefore, 
\begin{align}
    \frac{1}{\fra^2} \paren{ \expt{ \expt{\sfM \mid \sfXi, \sfY}^2 } - \frb^2 } 
    &= \tau_1 \paren{\expt{\expt{\sfM \mid \sqrt{(\tau_1 - \sigma^2) \tau_1 \sigma^2} \sfXi, \sfY}^2} - \paren{\sigma^2 - \frac{1}{\tau_1}}} \notag \\
    &= \tau_1 \paren{\expt{\expt{\sfZ \mid \sfP, \sfY}^2} - \paren{\sigma^2 - \frac{1}{\tau_1}}} . \label{eqn:EM_EZ}
\end{align}

With the relation \Cref{eqn:EM_EZ} above and the change of variables $ q_z = \sigma^2 - c_1 , \wh{q}_z = \tau_1 - \sigma^{-2} $, we can use the tower property of conditional expectation to write \Cref{eqn:replica2} as
\begin{align}
    \sigma^2 - c_1 &= \frac{1}{\tau_1} \brace{\frac{\tau_1 - \sigma^{-2}}{\sigma^{-2}} + \tau_1 \paren{\expt{\sfZ \expt{\sfZ \mid \sfP, \sfY}} - \paren{\sigma^2 - \frac{1}{\tau_1}}}} . \notag 
\end{align}
Upon rearrangement, this is precisely 
\begin{align}
    \sigma^2 - c_1 &= \expt{\sfZ \expt{\sfZ \mid \sfP, \sfY}} , \notag 
\end{align}
which, in view of \Cref{eqn:c1_alt}, is in turn equivalent to the second equation in \Cref{eqn:BFP_v1_c1}. 

This completes the proof of the equivalence between \Cref{eqn:B_FP,eqn:replica} and therefore of \Cref{thm:replica}. 

\section{Details of examples and numerical experiments}
\label{sec:details}

\subsection{Replica saddle point equations}
\label{sec:replica_eg}

For fixed $ \gamma_1\ge0, \tau_1\ge\sigma^{-2} $, let
\begin{subequations}
\label{eqn:replica_law}
\begin{align}
&&
    (\sfB_*, \ol{\sfJ}) &\sim P_{\sfB_*} \ot \cN(0,1) , & 
    \sfR &= \gamma_1 \sfB_* + \sqrt{\gamma_1} \, \ol{\sfJ} , & 
    &&
& \label{eqn:law_R} \\
&&
    (\sfZ, \ol{\sfK}) &\sim \cN(0,\sigma^2) \ot \cN(0,1) , & 
    \sfP &= (\tau_1 - \sigma^{-2}) \sfZ + \sqrt{\tau_1 - \sigma^{-2}} \, \ol{\sfK} , & 
    \sfY &\sim Q\paren{\cdot \mid \sfZ} . & 
& \label{eqn:law_P}
\end{align}
\end{subequations}
With respect to the above laws, define functions
\begin{align}
&&
    \cE_x(\gamma_1) &= \expt{\sfB_* \expt{\sfB_* \mid \sfR}} , & 
    \cE_z(\tau_1) &= \expt{\sfZ \expt{\sfZ \mid \sfP, \sfY}} . & 
& \label{eqn:def_Ex_Ez}
\end{align}
Recalling the construction of $ \sfLambda $ in \Cref{eqn:constr_Lambda} and the definition of $\sfLambda_d$ in \Cref{eqn:sfLambda_nd}, define the function 
\begin{align}
    \cS(\tau_2, \gamma_2) &\coloneqq \expt{(\tau_2 \sfLambda_d^2 + \gamma_2)^{-1}}
    = \begin{cases}
        \delta \expt{\paren{\tau_2 \wt{\sfLambda}^2 + \gamma_2}^{-1}} + (1-\delta) \gamma_2^{-1} , & \delta \le 1 \\
        \expt{\paren{\tau_2 \delta \wt{\sfLambda}^2 + \gamma_2}^{-1}} , & \delta > 1 
    \end{cases} . \label{eqn:def_S}
\end{align}

\begin{proposition}
The fixed point $(v_1,v_2,\gamma_1)$ of \Cref{eqn:Bayes_GVAMP} satisfies $v_1=v_2=v$ and $ (v,\gamma_1) \in [0,\rho] \times [0,\infty) $ solves the following system of equations
\begin{align}
&&
    v &= \rho - \cE_x(\gamma_1) , & 
    \gamma_1 &= \frac{\delta}{v} \brack{1 - \tau_1(v,\gamma_1) \paren{\sigma^2 - \cE_z(\tau_1(v,\gamma_1))} } , & 
& \label{eqn:replica_vgamma1} 
\end{align}
where $ \tau_1(v,\gamma_1) \ge \sigma^{-2} $ is the unique solution to 
\begin{align}
    v &= \cS\paren{\frac{\tau_1 \gamma_1 v}{\delta - \gamma_1 v}, \frac{1}{v} - \gamma_1} . \label{eqn:replica_tau1}
\end{align}

Furthermore, for the special case of phase retrieval $ q(z,\eps) = \abs{z} $ with Gaussian prior $ P_{\sfB_*} = \cN(0,\rho) $, $v\in[0,\rho]$ can be expressed more explicitly as the solution to $ v = \cF(v) $, for a function $ \cF $ defined as
\begin{align}
    \cF(v) &= \frac{\delta \rho v}{\rho - v} \paren{1 - \tau_1(v) \paren{\sigma^2 - \expt{\sfY^2 \tanh(\sfP \sfY)^2}}} , \label{eqn:def_F} 
\end{align}
where 
\begin{align}
&&
    (\sfP, \ol{\sfK}) &\sim \cN(0,(\tau_1(v)-\sigma^{-2})\tau_1(v)\sigma^{2}) \ot \cN(0,1) , & 
    \sfY &= \abs{\frac{1}{\tau_1(v)} \sfP + \frac{1}{\sqrt{\tau_1(v)}} \ol{\sfK}} , & 
& \label{eqn:Ez_law}
\end{align}
and $ \tau_1(v) \ge \sigma^{-2} $ (with slight abuse of notation for the solution $ \tau_1(v,\gamma_1) $ to \Cref{eqn:replica_tau1}) is the unique solution to
\begin{align}
    v &= \cS\paren{\frac{(\rho - v) \tau_1}{\rho \delta - (\rho - v)}, \frac{1}{\rho}} . \notag 
\end{align}
\end{proposition}

\begin{proof}
The replica saddle point equations \Cref{eqn:B_FP} read: 
\begin{subequations}
\label{eqn:replica_FP}
\begin{align}
&&
    v &= \rho - \expt{\sfB_* \expt{\sfB_* \mid \sfR}} , & 
    c &= \sigma^2 - \expt{\sfZ \expt{\sfZ \mid \sfP, \sfY}} , & 
& \label{eqn:replica_FP1} \\
&&
    \frac{1}{v} &= \gamma_1 + \gamma_2 , & 
    \frac{1}{c} &= \tau_1 + \tau_2 , & 
& \label{eqn:replica_FP2} \\
&&
    v &= \expt{\frac{1}{\tau_2 \sfLambda_d^2 + \gamma_2}} , & 
    c &= \expt{\frac{\sfLambda_n^2}{\tau_2 \sfLambda_n^2 + \gamma_2}} , & 
& \label{eqn:replica_FP3}
\end{align}
\end{subequations}
where we use the notation $ v \coloneqq v_1 = v_2, c \coloneqq c_1 = c_2 $, 
and the random variables in the first line have joint laws specified in \Cref{eqn:replica_law}. 
\Cref{eqn:replica_FP} is intended to be solved in the domain 
\begin{align}
    (\gamma_1, \tau_1, \gamma_2, \tau_2, v, c) &\in [0,\infty) \times [\sigma^{-2}, \infty) \times [\rho^{-1}, \infty) \times [0,\infty) \times [0,\rho] \times [0,\sigma^2] . \notag 
\end{align}
Note that $c$ can be expressed by other variables as
\begin{align}
    c &= \expt{\frac{\sfLambda_n^2}{\tau_2 \sfLambda_n^2 + \gamma_2}} = \frac{1}{\delta} \expt{\frac{\sfLambda_d^2}{\tau_2 \sfLambda_d^2 + \gamma_2}}
    = \frac{1}{\delta \tau_2} (1 - \gamma_2 v) . \label{eqn:replica_c}
\end{align}
Using $ \gamma_2 = 1/v - \gamma_1, \tau_2 = 1/c - \tau_1 $ from \Cref{eqn:replica_FP2} to eliminate $ \gamma_2, \tau_2 $ in \Cref{eqn:replica_c}, upon rearrangements, we obtain
\begin{align}
    c &= \frac{1}{\tau_1} \paren{1 - \frac{\gamma_1 v}{\delta}} . \label{eqn:replica_cc}
\end{align}
Substituting the right-hand side of \Cref{eqn:replica_cc} in the second equation of \Cref{eqn:replica_FP2} for $c$, we have 
\begin{align}
    \tau_2 &= \frac{1}{c} - \tau_1 = \frac{\tau_1 \gamma_1 v}{\delta - \gamma_1 v} . \label{eqn:tmp_tau2}
\end{align}
A similar substitution in the second equation of \Cref{eqn:replica_FP1} gives, upon rearrangements,  
\begin{align}
    \gamma_1 &= \frac{\delta}{v} \brack{1 - \tau_1 \paren{\sigma^2 - \cE_z(\tau_1)} } , \label{eqn:tmp_gamma1}
\end{align}
where $ \cE_z(\tau_1) $ is defined in \Cref{eqn:def_Ex_Ez}. 
Further recalling $\cE_x(\gamma_1)$ from \Cref{eqn:def_Ex_Ez} and $\cS(\tau_2,\gamma_2)$ from \Cref{eqn:def_S}, we collect \Cref{eqn:tmp_tau2,eqn:tmp_gamma1} to arrive at the desired system of equations \Cref{eqn:replica_vgamma1,eqn:replica_tau1}. 

For the phase retrieval model $q(z,\eps) = \abs{z}$ with Gaussian prior $ P_{\sfB_*} = \cN(0,\rho) $,
we have 
\begin{align}
    \expt{\sfB_* \mid \sfR = r} &= \frac{\rho r}{\gamma_1 \rho + 1} , \label{eqn:gx1_pr} \\
    \expt{\sfZ \mid \sfP = p, \sfY = y} 
    &= \frac{\int \varphi_{\sigma}(z) \varphi_{\sqrt{\tau_1 - \sigma^{-2}}}(p - (\tau_1 - \sigma^{-2}) z) \delta(y - \abs{z}) z \diff z}{\int \varphi_{\sigma}(z') \varphi_{\sqrt{\tau_1 - \sigma^{-2}}}(p - (\tau_1 - \sigma^{-2}) z') \delta(y - \abs{z'}) \diff z'}
    = y \tanh(py) , \label{eqn:gz1_pr}
\end{align}
and therefore, by the tower property of conditional expectation, 
\begin{align}
&&
    \cE_x(\gamma_1) &= \frac{\expt{\sfR^2}}{(\gamma_1\rho + 1)^2}
    = \frac{\rho^2 \gamma_1}{\gamma_1 \rho + 1} , & 
    \cE_z(\tau_1) &= \expt{\sfY^2 \tanh(\sfP \sfY)^2} , & 
& \label{eqn:Exz_Gauss} 
\end{align}
where the second expectation is taken over the joint law in \Cref{eqn:Ez_law} which is equivalent to \Cref{eqn:law_P} specialized to the phase retrieval model. 
Using the expression of $ \cE_x(\gamma_1) $ in \Cref{eqn:Exz_Gauss}, we can solve for $ \gamma_1 $ in terms of $v$ from the first equation in \Cref{eqn:replica_vgamma1} which yields 
\begin{align}
    \gamma_1 &= \frac{\rho - v}{\rho v} . \label{eqn:solve_gamma1}
\end{align}
This then leads to the more succinct fixed point equation $ v = \cF(v) $ where $ \cF $ is defined as \Cref{eqn:def_F}. 
\end{proof}

\subsection{Synthetic data}
\label{sec:details_synthetic}

\paragraph{Phase retrieval.}
\label{sec:details_phase_retrieval}

Under the setting of \Cref{sec:phase_real}, it is easy to verify that the function $ \ol{g} $ in \Cref{eqn:olg_opt} is given by
\begin{align}
\ol{g}(y) &= y^2 - 1 , \label{eqn:g_phase_retrieval} 
\end{align}
where $ \gamma $ is defined as per \Cref{eqn:T}. 
The corresponding preprocessing function (see \Cref{asmp:preprocess}) is therefore $ \cT(y) = 1 - \frac{\gamma}{y^2 - 1 + \gamma} $, whereas the preprocessing function from \cite{Luo_Alghamdi_Lu} specializes to $ \cT(y) = 1 - \frac{1}{y^2} $. 
The weak recovery condition \Cref{eqn:thr_mmt} becomes
\begin{align}
\delta &> 
\frac{\ol{m}_2'^2}{\ol{m}_4'} \frac{3}{2} . \notag 
\end{align}
The vector denoisers $ g_{x1}, g_{z1} $ (see \Cref{eqn:gx1,eqn:gz1}) used in Bayes-GVAMP admit explicit expressions which have been derived in \Cref{eqn:gx1_pr,eqn:gz1_pr}. 

\paragraph{Poisson regression.}
\label{sec:details_pois_regr}

In this setting, elementary calculations show that $ \ol{g} $ in \Cref{eqn:olg_opt} evaluates to 
\begin{align}
    \ol{g}(y) &= \frac{2}{3} \paren{\frac{1}{2} + y} - 1 . \notag 
\end{align}
The preprocessing functions from \Cref{asmp:preprocess} and \cite{Luo_Alghamdi_Lu} can then be evaluated accordingly. 
The criticality condition \Cref{eqn:thr_mmt} becomes
\begin{align}
    \delta &> \frac{\ol{m}_2'^2}{\ol{m}_4'} \frac{7}{4} . \notag     
\end{align}

For Bayes-GVAMP, though the denoiser $ g_{x1} $ is still given by the simple Gaussian MMSE estimator \Cref{eqn:gx1_pr}, $ g_{z1} $ in \Cref{eqn:gz1} and its second moment (required by \Cref{eqn:g'}) are more tricky to evaluate: 
\begin{align}
&&
    g_{z1}(p,y;\tau_1) &= \expt{\sfZ \mid \sfP = p, \sfY = y} , & 
    \expt{g_{z1}(\sfP,\sfY;\tau_1)^2} &= \expt{\expt{\sfZ \mid \sfP, \sfY}^2} , & 
& \label{eqn:gz1_12}
\end{align}
where 
\begin{align}
&&
    (\sfZ,\ol{\sfK}) &\sim \cN(0,\sigma^2) \ot \cN(0,1) , & 
    \sfP &= (\tau_1 - \sigma^{-2}) \sfZ + \sqrt{\tau_1 - \sigma^{-2}} \, \ol{\sfK} , & 
    \sfY \sim \pois(\sfZ^2) . & 
& \notag 
\end{align}
Under the above law, \Cref{eqn:gz1_12} can be written equivalently as 
\begin{align}
&& 
    g_{z1}(p,y;\tau_1) &= \frac{\expt{\sfG(p;\tau_1)^{2y+1}}}{\expt{\sfG(p;\tau_1)^{2y}}} , & 
    \expt{g_{z1}(\sfP,\sfY;\tau_1)^2} &= \expt{\paren{\frac{\expt{\sfG(\sfP;\tau_1)^{2\sfY+1}}}{\expt{\sfG(\sfP;\tau_1)^{2\sfY}}}}^2} , & 
& \label{eqn:gz1_m12} 
\end{align}
where the random variable $ \sfG(p;\tau_1) $ is distributed according to
\begin{align}
    \sfG(p;\tau_1) &\sim \cN\paren{\frac{p}{\tau_1 + 2}, \frac{1}{\tau_1 + 2}} . \notag 
\end{align}
The difficulty in evaluating \Cref{eqn:gz1_m12} is two-fold: 
the ratio in $ g_{z1} $ involves moments of a non-centered Gaussian whose values can quickly reach machine precision for moderate values of $y\in\bbZ_{\ge0}$; 
the second moment of $ g_{z1} $ involves nested expectations over $ 3 $ variables $ \sfZ, \ol{\sfK}, \sfY $. 
To overcome the difficulty, we make use of the following observation. 
\begin{proposition}
\label{prop:stein_recur}
Fix $\mu,a,b\in\bbR$ and $\sigma>0$. 
Let $ \sfG \sim \cN(\mu,\sigma^2) $ and denote for any $k\ge0$, 
\begin{align}
    R_{k+1}(\mu,\sigma^2) &= \frac{\expt{(a+b\sfG)^{k+1}}}{\expt{(a+b\sfG)^k}} . \notag 
\end{align}
Then we have $ R_1(\mu,\sigma^2) = a+b\mu $; and for $ k\in\bbZ_{\ge1} $, 
\begin{align}
    R_{k+1}(\mu,\sigma^2) &= a + b \paren{\mu + \sigma^2 k b \frac{1}{R_k(\mu,\sigma^2)}} . \notag 
\end{align}
\end{proposition}
The proof follows immediately from Stein's lemma applied to the numerator of $ R_{k+1} $; we omit the details. 

\Cref{prop:stein_recur} allows us to evaluate $ g_{z1} $ recursively (over $y$) without evaluating its numerator and denominator individually. 
For the second moment of $ g_{z1} $, we further combine \Cref{prop:stein_recur} with the $2$-dimensional Gauss--Hermite quadrature \cite[Section 4.6.1]{Numerical_Recipes} and employ the following approximation scheme (with suitably large $N,Y$): 
\begin{align}
    \expt{g_{z1}(\sfP,\sfY;\tau_1)^2} &\approx \sum_{z = 1}^N \sum_{k = 1}^N w_N(x_z) w_N(x_k) \sum_{y = 0}^Y \frac{(\sigma x_z)^{2y} e^{-(\sigma x_z)^2}}{y!} R_{2y+1}\paren{\frac{(\tau_1 - \sigma^{-2}) \sigma x_z + \sqrt{\tau_1 - \sigma^{-2}} \, x_k}{\tau_1+2}, \frac{1}{\tau_1+2}}^2 , \notag 
\end{align}
where the weights $ w_N(\cdot) $ are given by 
\begin{align}
    w_N(x) &= \frac{2^{N-1} N!}{N^2 H_{N-1}(x)^2} , \notag 
\end{align}
with $H_N$ denoting the degree-$N$ physicists' Hermite polynomial and $(x_i/\sqrt{2})_{1\le i\le N}$ its roots. 
In the numerics, we find it efficient and accurate to take the degree of the quadrature to be $N=50$ and the upper bound $Y$ on $y$ to be $10$ standard deviations larger than the mean of the Poisson. 
This approximation scheme is deterministic and avoids stochastic errors of Monte Carlo method which we find non-negligible even with reasonably large sample sizes. 
Similar approximation scheme involving the recursion in \Cref{prop:stein_recur} and quadrature is used to evaluate analogous quantities in GAMP and its state evolution \cite{Mondelli_Venkataramanan}. 

\subsection{Real data}
\label{sec:details_real}

\paragraph{Coded diffraction patterns.}
\label{sec:details_CDP}

Coded diffraction patterns over the reals used in \Cref{sec:CDP} are constructed similarly to their complex analogues in \cite[Section 2.2]{candes_CDP}. 
For $ \delta \in \bbZ_{\ge1} $, $ X \in \bbR^{d\delta \times d} $ is defined as 
\begin{align}
    X &= \matrix{
    C_d^\top D_1 & \cdots & C_d^\top D_\delta
    }^\top , \notag 
\end{align}
where $ C_d \in \bbO(d) $ is the discrete cosine transform (DCT) matrix and each $ D_i \in \bbR^{d\times d} $ (for $ 1\le i\le \delta $) is an independent diagonal matrix with diagonal elements i.i.d.\ according to the law of a certain random variable $ \sfD $. 
For binary and ternary CDP, $ \law(\sfD) $ is taken to be 
\begin{align}
&&
    \law(\sfD) &= \bern(1/2) & 
    & \textnormal{and} & 
    \law(\sfD) &= \frac{1}{2} \textnormal{Rademacher} + \frac{1}{2} \delta_0 & 
& \notag 
\end{align}
respectively. 
Recalling \Cref{eqn:sfLambda_nd}, by orthogonality of $C_d$, we have 
\begin{align}
    \sfLambda_d^2 &\eqqlaw \sum_{i = 1}^\delta \sfD_i^2 , \notag 
\end{align}
where $(\sfD_i)_{1\le i\le \delta}$ are i.i.d.\ copies of $\sfD$. 
However, note that $ X $ is neither left- nor right-orthogonally invariant in law. 

\paragraph{GTEx.}
\label{sec:details_GTEx}

The design matrices used in \Cref{sec:GTEx} are obtained by sequentially applying the following standard preprocessing steps \cite{PLINK} to the GTEx datasets `Skin - Sun Exposed (Lower leg)' and `Muscle - Skeletal' \cite{GTEx} respectively. 
\begin{enumerate}
    \item Remove rows consisting of all zeros;
    \item Build a matrix by sequentially including each row from above only if it has an overlap at most $0.3$ with all existing rows;
    \item Standardize all rows collected from the previous step. 
\end{enumerate}
The resulting matrices are entirely deterministic and, in particular, do not satisfy \Cref{asmp:design}. 
The only randomness involved in the experiments in \Cref{sec:GTEx} lies in the generation of $ \beta_* \sim \cN(0_d,I_d) $. 

\section{Auxiliary lemmas}
\label{app:aux}

\Cref{prop:wto1,prop:wto2,prop:wto3,prop:wto5} are a few facts regarding Wasserstein convergence of empirical measures collected from \cite[Appendices E and F.1]{fan2020approximate} and \cite[Appendix G.1]{Li_Fan_Sen_Wu}. 

\begin{proposition}
\label{prop:wto1}
Fix $ p\ge1, k\ge1, \ell\ge0 $. 
Let $ E \in \bbR^{n\times \ell} $ be a deterministic matrix satisfying $ E \xrightarrow{W_p} \sfE^\top $. 
If $ V \in \bbR^{n\times k} $ have i.i.d.\ rows according to the law of $ \sfV \in \bbR^k $ with finite mixed moments up to order $p$, then $ \matrix{ V & E } \xrightarrow{W_p} \matrix{ \sfV^\top & \sfE^\top } $ where $ \sfV $ is independent of $ \sfE $. 
\end{proposition}

\begin{proposition}
\label{prop:wto2}
Fix $ p,q\ge1 $ and $ k,\ell\ge1 $. 
Suppose that $ V \in \bbR^{n\times k} $ satisfies $ V \xrightarrow{W_{p+q}} \sfV^\top $. 
Let $ g \in \PG^q_{k\to\ell} $ (see \Cref{def:poly_grow}) be continuous. 
Then $ g(V) \xrightarrow{W_p} g(\sfV)^\top $, where $ g $ is applied row-wise to $V$. 
\end{proposition}

\begin{proposition}
\label{prop:wto5}
Fix $p\ge1,k\ge1$. 
Suppose that $V\in\bbR^{n\times k}$ satisfies $ V \xrightarrow{W_p} \sfV^\top $. 
Let $ f \in \PG^p_{k\to1} $ be continuous everywhere except on a set of probability $0$ under the law of $ \sfV $. 
Then almost surely, 
\begin{align}
    \lim_{n\to\infty} \frac{1}{n} \sum_{i = 1}^n f(V_{i,:}) &= \expt{f(\sfV)} . \notag 
\end{align}
\end{proposition}

\begin{proposition}
\label{prop:wto3}
Fix $ p,k,\ell\ge1 $. 
Suppose $ V \in \bbR^{n\times k} $, $ W \in \bbR^{n\times\ell} $, and $ M_n, M \in \bbR^{k\times\ell} $ satisfy $ V \xrightarrow{W_p} \sfV^\top $, $ W \xrightarrow{W_p} 0_{1\times\ell} $, and $ M_n\to M $ entrywise, almost surely. 
Then $ VM_n + W \xrightarrow{W_p} \sfV^\top M $. 
\end{proposition}

\Cref{prop:wto4} below is a generalization of \cite[Proposition F.2]{fan2020approximate} from $m=1$ to $m\ge1$ fixed. 

\begin{proposition}
\label{prop:wto4}
Fix $p\ge1$, $k,\ell\ge0$ and $ m\ge1 $. 
Let $ O \sim \haar(\bbO(n - \ell)) $. 
Let $ E \in \bbR^{n\times k} $, $ V \in \bbR^{(n - \ell) \times m} $, $ \Pi \in \bbR^{n\times(n - \ell)} $ be deterministic matrices satisfying the following: $ E \xrightarrow{W_p} \sfE^\top $, $ n^{-1} V^\top V \to \Sigma $ entry-wise, and $\Pi $ has orthonormal columns. 
Then $ \matrix{ \Pi O V & E } \xrightarrow{W_p} \matrix{ \sfZ^\top & \sfE^\top } $ where $ \sfZ \sim \cN(0_m,\Sigma) $ is independent of $ \sfE $. 

Let $j \ge 0$ be fixed, and $d$ a dimension that grows jointly with $n$. 
Fix a constant $j\ge0$. 
Let $ F \in \bbR^{d\times j} $ be deterministic satisfying $ F \xrightarrow{W_p} \sfF^\top $. 
Let $ U \in \bbR^{d\times m} $ be the first $d$ rows of $ \Pi O V $ if $d\le n$, or $ \Pi O V $ appended by $d-n$ i.i.d.\ rows distributed as $ \cN(0_m,\Sigma) $ if $d>n$. 
Then $ \matrix{U & F} \xrightarrow{W_p} \matrix{\sfU^\top & \sfF^\top} $ where $ \sfU \sim \cN(0_m,\Sigma) $ is independent of $\sfF$. 
\end{proposition}

\begin{proof}
Denote $ n' = n-\ell $ and suppose that it is sufficiently large. 
Note that left-multiplying $V$ by $O$ uniformly randomizes the left singular space of $V$ while retaining its singular values. 
Writing $ V $ in terms of its SVD  we have $V = L \matrix{\Lambda \\ 0_{(n'-m)\times m}} R^\top $. We can express the law of $ OV $ as 
\begin{align}
    OV &\eqqlaw O \matrix{\Lambda \\ 0_{(n'-m)\times m}} R^\top
    \eqqlaw O \matrix{R & 0_{m\times(n'-m)} \\ 0_{(n'-m)\times m} & I_{n'-m}} \matrix{\Lambda \\ 0_{(n'-m)\times m}} R^\top 
    = O' R \Lambda R^\top 
    = O' (V^\top V)^{1/2} , \notag 
\end{align}
where $ O' $ contains the first $m$ columns of $O$. 
We use the well-known representation of the law of a partial Haar matrix (see, e.g., \cite[Section I-D]{Thrampoulidis_Hassibi}): $ O' \eqqlaw H (H^\top H)^{-1/2} $, where $ H \in \bbR^{n'\times m} $ consists of i.i.d.\ standard Gaussian entries. 
Since $ \Pi $ has orthonormal columns, $ \Pi^\top G \eqqlaw H $ where $ G \in \bbR^{n\times m} $ consists of i.i.d.\ standard Gaussian entries. 
Using both observations above, we have 
\begin{align}
    \Pi OV &\eqqlaw \Pi \Pi^\top G (G^\top \Pi \Pi^\top G)^{-1/2} (V^\top V)^{1/2} \notag \\
    &= G (G^\top \Pi \Pi^\top G)^{-1/2} (V^\top V)^{1/2} - P G (G^\top \Pi \Pi^\top G)^{-1/2} (V^\top V)^{1/2} , \label{eqn:PiOV} 
\end{align}
where $ P = I_n - \Pi \Pi^\top $. 
By \Cref{prop:wto1}, 
\begin{align}
    \matrix{G & E} &\xrightarrow{W_p} \matrix{\sfG^\top & \sfE^\top} , \label{eqn:GE}
\end{align}
where $ \sfG \sim \cN(0_m,I_m) $ is independent of $ \sfE $. 
Since $ \Pi^\top G $ has i.i.d.\ $ \cN(0,1) $ entries, we have $ n^{-1} G^\top \Pi \Pi^\top G \to I_m $ entry-wise, almost surely. 
Therefore, recalling the assumption $ n^{-1} V^\top V \to \Sigma $, we have 
\begin{align}
    & (G^\top \Pi \Pi^\top G)^{-1/2} (V^\top V)^{1/2} \to \Sigma^{1/2} \label{eqn:OV1}
\end{align}
entry-wise, almost surely. 

Next, we show that 
\begin{align}
    PG &\xrightarrow{W_p} 0_m^\top . \label{eqn:PG0_TODO}
\end{align}
The law of $ P G $ can be written as $ P G \eqqlaw \sum_{j = 1}^\ell \pi_j g_j^\top $ where $ \pi_1, \cdots, \pi_m \in \bbR^n $ form an orthonormal basis of the column span of $P$ and $ g_1, \cdots, g_m $ are i.i.d.\ $ \cN(0_m, I_m) $ vectors. 
Write each $ \pi_j $ as $ \matrix{\pi_{1,j} & \cdots & \pi_{n,j}}^\top $. 
For any fixed $p\ge1$, 
\begin{align}
    \frac{1}{n} \sum_{i = 1}^n \normtwo{\sum_{j = 1}^\ell \pi_{i,j} g_j^\top}^p &\le \frac{1}{n} \sum_{i = 1}^n \paren{\sum_{j = 1}^\ell \normtwo{\pi_{i,j} g_j^\top}}^p 
    \le \frac{1}{n} \sum_{i = 1}^n \brack{\ell^{1-1/p} \paren{\sum_{j = 1}^\ell \normtwo{\pi_{i,j} g_j^\top}^p}^{1/p}}^p \notag \\
    &= \frac{\ell^{p-1}}{n} \sum_{j = 1}^\ell \sum_{i = 1}^n \abs{\pi_{i,j}}^p \normtwo{g_j}^p 
    = \ell^{p-1} \sum_{j = 1}^\ell \normtwo{g_j}^p \frac{\norm{p}{\pi_j}^p}{n} , \notag 
\end{align}
where we have used subadditivity of $\ell_2$ norm and the inequality $ \norm{p}{x} \le \norm{r}{x} \le d^{1/r - 1/p} \norm{p}{x} $ for any vector $ x \in \bbR^d $ and $ 0 \le r \le p $. 
If $ p\ge2 $, then $ n^{-1} \norm{p}{\pi_j}^p \le n^{-1} \normtwo{\pi_j}^p = n^{-1} $. 
If $ 1\le p<2 $, then $ n^{-1} \norm{p}{\pi_j}^p \le n^{-1} \paren{n^{1/p-1/2} \normtwo{\pi_j}}^p = n^{-p/2} $. 
In both cases, $ n^{-1} \norm{p}{\pi_j}^p \to 0 $. 
Therefore, almost surely, 
\begin{align}
    \frac{1}{n} \sum_{i = 1}^n \normtwo{(PG)_{i,:}}^p &\to 0 . \label{eqn:PG0} 
\end{align}
Furthermore, for any bounded Lipschitz $ f \colon \bbR^m \to \bbR $, by \Cref{eqn:PG0}, almost surely, 
\begin{align}
    \abs{\frac{1}{n} \sum_{i = 1}^n \brack{f((PG)_{i,:}) - f(0_m^\top)}} &\le \frac{1}{n} \sum_{i = 1}^n C \normtwo{(PG)_{i,:} - 0_m^\top} \to 0 . \label{eqn:PG0_lip} 
\end{align}
Taking \Cref{eqn:PG0,eqn:PG0_lip} collectively, we have almost surely, 
\begin{align}
    \frac{1}{n} \sum_{i = 1}^n f((PG)_{i,:}) &\to f(0_m^\top) , \notag 
\end{align}
for both $ f(v) = \normtwo{v}^p $ and any bounded Lipschitz $ f $. 
Therefore \Cref{eqn:PG0_TODO} holds. 
Combining \Cref{eqn:PG0_TODO,eqn:OV1,eqn:PiOV,eqn:GE} and applying \Cref{prop:wto3} yield the first convergence result for $ \matrix{\Pi O V & E} $: 
\begin{align}
    \matrix{\Pi O V & E}
    &\eqqlaw \matrix{G & E} \matrix{(G^\top \Pi \Pi^\top G)^{-1/2} (V^\top V)^{1/2} & 0_{m\times k} \\ 0_{k\times m} & I_k} - \matrix{P G (G^\top \Pi \Pi^\top G)^{-1/2} (V^\top V)^{1/2} & 0_{n\times k}}
    \xrightarrow{W_p} \matrix{\sfZ^\top & \sfE^\top} , \notag 
\end{align}
where $ \sfZ \sim \cN(0_m, \Sigma) $ is independent of $ \sfE $. 

Next consider $U$. 
Construct $ \wt{G}, W_1, W_2 \in \bbR^{d\times m} $ as follows: 
\begin{align}
&&
    \wt{G} &= \begin{cases}
        G_{1:d,:} , & d\le n \\
        \matrix{G \\ G'} , & d > n
    \end{cases} , & 
    W_1 &= \begin{cases}
        \brack{P G (G^\top \Pi \Pi^\top G)^{-1/2} (V^\top V)^{1/2}}_{1:d,:} , & d\le n \\
        \matrix{P G (G^\top \Pi \Pi^\top G)^{-1/2} (V^\top V)^{1/2} \\ 0_{(d-n)\times m}} , & d>n
    \end{cases} , & 
& \notag \\
&&
    &&
    W_2 &= \begin{cases}
        0_{d\times m} , & d\le n \\
        \matrix{0_{n\times m} \\ G'\brack{(G^\top \Pi \Pi^\top G)^{-1/2} (V^\top V)^{1/2} - \Sigma}} , & d>n
    \end{cases} , & 
& \notag 
\end{align}
where $ G' \in \bbR^{(d-n)\times m} $ (if $ d>n $) contains i.i.d.\ $\cN(0,1)$ entries, independent of everything else. 
By construction, it is easy to verify that
\begin{align}
    U &= \wt{G} (G^\top \Pi \Pi^\top G)^{-1/2} (V^\top V)^{1/2} - W_1 - W_2 . \notag 
\end{align}
Again by \Cref{prop:wto1}, 
\begin{align}
    \matrix{\wt{G} & F} &\xrightarrow{W_p} \matrix{\wt{\sfG}^\top & \sfF^\top} , \label{eqn:GF} 
\end{align}
where $ \wt{\sfG} \sim \cN(0_m,I_m) $ is independent of $ \sfF $. 
The same argument as before shows $ W_1 \xrightarrow{W_p} 0_m^\top $. 
Also, $ W_2 \xrightarrow{W_p} 0_m^\top $. 
To verify this, we only need to consider the case $ d>n $. 
Almost surely, 
\begin{align}
    \frac{1}{d} \sum_{i = 1}^{d-n} \normtwo{G'_{i,:} \brack{(G^\top \Pi \Pi^\top G)^{-1/2} (V^\top V)^{1/2} - \Sigma}}^p
    &\le \normtwo{(G^\top \Pi \Pi^\top G)^{-1/2} (V^\top V)^{1/2} - \Sigma} \frac{1}{d} \sum_{i = 1}^{d-n} \normtwo{G'_{i,:}}^p 
    \to 0 . \notag 
\end{align}
Using this, one can then verify \Cref{eqn:wass_f} for any bounded Lipschitz $f$, thereby justifying the $W_p$-convergence of $ W_2 $; the argument is similar to \Cref{eqn:PG0_lip}. 
Combining these with \Cref{eqn:GF} and applying \Cref{prop:wto3} then give the convergence result for $ \matrix{U & F} $. 
\end{proof}

The following result identifies the distribution of a Haar matrix conditioned on it solving a linear system. 

\begin{proposition}[{\cite[Lemma 4]{Rangan_Schniter_Fletcher}}]
\label{prop:haar_cond}
Fix $k\ge1$. 
Let $ X, Y\in\bbR^{n\times k} $ be deterministic matrices with rank $k$ such that $ QY = X $ for some $ Q\in\bbO(n) $. 
If $ O\in\haar(\bbO(n)) $, then 
\begin{align}
O \mid \brace{ OY = X }
&\eqqlaw X (X^\top X)^{-1} Y^\top + \Pi_{X^\perp} \wt{O} \Pi_{Y^\perp}^\top = X (Y^\top Y)^{-1} Y^\top + \Pi_{X^\perp} \wt{O} \Pi_{Y^\perp}^\top , \notag 
\end{align}
where $ \wt{O} \sim \haar(\bbO(n-k)) $ is independent of everything else and $ \Pi_{X^\perp}, \Pi_{Y^\perp} \in \bbR^{n\times(n-k)} $ are (any) matrices with orthonormal columns spanning $ \spn\brace{\brace{X_{:,i}}_{i = 1}^k}^\perp, \spn\brace{\brace{Y_{:,i}}_{i = 1}^k}^\perp $, respectively.  
\end{proposition}

\end{document}